\documentclass{jams-l}

\usepackage{amssymb}
\usepackage{amsmath}

\usepackage{graphicx}

\usepackage[cmtip,all]{xy}

\usepackage{stmaryrd}
\usepackage{comment}
\usepackage{tikz}
\usepackage{caption}
\usepackage{subcaption}
\usepackage{float}
\usepackage{wrapfig}
\usepackage{bm}
\usepackage{dashbox}
\usepackage{multirow}
\usepackage{xcolor}
\usepackage{hyperref}

\usetikzlibrary{quotes,angles}

\newtheorem{theorem}{Theorem}[section]
\newtheorem{lemma}[theorem]{Lemma}
\newtheorem{proposition}[theorem]{Proposition}
\newtheorem{proposition-definition}[theorem]{Proposition-Definition}
\newtheorem{corollary}[theorem]{Corollary}

\theoremstyle{definition}
\newtheorem{definition}[theorem]{Definition}

\theoremstyle{remark}
\newtheorem{remark}[theorem]{Remark}

\numberwithin{equation}{section}

\DeclareMathOperator{\id}{id}
\DeclareMathOperator{\arcsinh}{Arcsinh}

\newcommand\scalemath[2]{\scalebox{#1}{\mbox{\ensuremath{\displaystyle #2}}}}

\begin{document}

\title[Pseudo-Laplacian on a cuspidal end: Dirichlet conditions]{Pseudo-Laplacian on a cuspidal end with a flat unitary line bundle: Dirichlet boundary conditions}


\author{Mathieu Dutour}
\address{University of Alberta, Department of Mathematical and Statistical Sciences, Central Academic Building, Edmonton, Alberta, Canada T6G2G1}
\curraddr{}
\email{dutour@ualberta.ca}
\thanks{}



\date{}

\dedicatory{}

\begin{abstract}

    A cuspidal end is a type of metric singularity, described as a product $S^1 \times \left] a, +\infty \right[$ with the Poincaré metric. The underlying set can also be seen as $\mathbb{R} \times \left] a, +\infty \right[$ subject to the action of the translation $T : \left( x,y \right) \longrightarrow \left( x+1, y \right)$. On it, one may consider a holomorphic line bundle $L$, coming from a unitary character of the group generated by $T$. The complex modulus induces a flat metric on $L$, and a pseudo-Laplacian $\Delta_{L,0}$ can be associated to the Chern connection, with Dirichlet boundary conditions. The aim of this paper is to find the asymptotic behavior of the zeta-regularized determinant $\det \left( \Delta_{L,0} + \mu \right)$, as $\mu > 0$ goes to infinity for any $a$, and also as $a$ goes to infinity for $\mu = 0$. \vspace{10pt}

\end{abstract}

\maketitle

\tableofcontents


\section{Introduction}

    \subsection{Description of the situation}
    
        This paper is devoted to the spectral study of special types of metric singularities on Riemann surfaces, called \textit{cusps}, with flat unitary holomorphic line bundles on them. Such a situation naturally arises when considering a modular curve defined by a Fuchsian group of the first kind $\Gamma$, and a vector bundle induced by a unitary representation of $\Gamma$. The computations of determinants made here can then be used to obtain a \textit{Deligne--Riemann--Roch isometry} extending \cite{MR4167014}, where Freixas i Montplet and von Pippich deal with the case of the trivial line bundle. The following introduction is detailed, so as to facilitate the reading of the more technical parts of this paper.

        \subsubsection{Metric singularities}
        \label{SubSubSecMetricSingularities}
    
            The underlying set of a cusp is defined as the quotient of $\mathbb{R} \times \left] a, + \infty \right[$ by the action of the translation $\left( x,y \right) \mapsto \left( x+1, y \right)$, or alternatively as the product $S^1 \times \left] a, + \infty \right[$, endowed with the Poincaré metric
        
            \[ \begin{array}{lll}
                    \mathrm{d}s^2_{\mathrm{hyp}} & = & \frac{\mathrm{d}x^2 + \mathrm{d}y^2}{y^2}
                \end{array} . \]
        
            \noindent Using the coordinate $z = \exp \left( 2 i \pi \left( x+iy \right) \right)$, a cusp can also be seen as a punctured disk of radius $ \varepsilon =\exp \left( - 2 \pi a \right)$, whose center corresponds to the singularity. In this description, the Poincaré metric becomes
        
            \[ \begin{array}{lll}
                    \mathrm{d}s^2_{\mathrm{hyp}} & = & \frac{\left\vert \mathrm{d} z \right \vert^2}{\left( \left \vert z \right \vert \log \left \vert z \right\vert \right)^2}
                \end{array} . \]
        
            \noindent This metric cannot be extended into a smooth metric at the center of the disk, which is the meaning of the term ``singularity'' here. We also need to consider a flat unitary line bundle over a cuspidal end. Such an object is induced by a unitary character $\chi : \mathbb{Z} \longrightarrow \mathbb{C}^{\ast}$, which provides an action of $\mathbb{Z}$ onto the trivial $\mathbb{C}$-bundle of rank $1$ over $\mathbb{R} \times \left] a, + \infty \right[$ defined by
        
            \[ \begin{array}{lll}
                    k \cdot \left( \left( x,y \right), \lambda \right) & = & \left( \left( x+k, y \right), \chi \left( k \right)  \lambda \right) 
                \end{array} \]
        
            \noindent for $k \in \mathbb{Z}$, as well as $\left( x,y \right) \in \mathbb{R} \times \left] a, +\infty \right[$ and $\lambda \in \mathbb{C}$ . Under this group action, the quotient $\mathbb{Z} \backslash \left( \left( \mathbb{R} \times \left] a, + \infty \right[ \right) \times \mathbb{C} \right)$ is a flat unitary line bundle $L$ over the cuspidal end, which is entirely determined by the complex number of modulus $1$
        
            \[ \begin{array}{lll}
                \chi \left( 1 \right) & = & e^{2 i \pi \alpha}
            \end{array} , \]
        
            \noindent with $\alpha$ being a real number well-defined modulo $1$. To simplify, we identify $\alpha$ with its representative in $\left[ 0, 1 \right[$. We can extend $L$ to a holomorphic line bundle over the cusp, \textit{i.e.} over the center of the disk in the coordinate $z$, using Deligne's canonical extension (see \cite{MR0417174, MR2393625}). The complex modulus on $\mathbb{C}$, being compatible with the action of $\mathbb{Z}$, induces a metric on $L$, called the \textit{canonical flat metric}, which cannot, in general, be extended smoothly over the cusp.

        \subsubsection{Pseudo-Laplacian}
        
            In this paper, we consider the \textit{pseudo-Laplacian} with Dirichlet boundary condition, studied by Colin de Verdière in \cite{MR688031, MR699488}. The value of the representative $\alpha \in \left[ 0,1 \right[$ splits the discussion into two parts.
            
            \medskip
            
            \hspace{10pt} $\bullet$ First, consider the case of a (metrically) non-trivial line bundle $L$, which corresponds to having $\alpha > 0$. The Chern Laplacian, acting on compactly supported smooth sections of $L$, is a symmetric operator, and its \textit{Friedrichs extension} is a self-adjoint operator, called the \textit{Chern Laplacian with Dirichlet boundary condition}. This operator does not have an essential spectrum. For the purpose of this paper, in order to be consistent with the case of a trivial bundle, this Laplacian is renamed the \textit{pseudo-Laplacian with Dirichlet boundary condition}, and is denoted by $\Delta_{L,0}$.
            
            \medskip
            
            \hspace{10pt} $\bullet$ Should $L$ be (metrically) trivial, which corresponds to having $\alpha = 0$, the Chern Laplacian has an essential spectrum, which must be removed from consideration before we can compute a determinant. This is achieved by considering the orthogonal decomposition
            
            \[ \begin{array}{lll}
                    \scriptstyle L^2 \left( S^1 \times \left] a, +\infty \right[, \frac{\mathrm{d}x^2 + \mathrm{d}y^2}{y^2} \right) & \scriptstyle = & \scriptstyle L^2 \left( S^1 \times \left] a, +\infty \right[, \frac{\mathrm{d}x^2 + \mathrm{d}y^2}{y^2} \right)_0 \; \; \oplus \; \; L^2 \left( \left] a, + \infty \right[, \frac{\mathrm{d}y^2}{y^2} \right)
                \end{array} , \]
            
            \noindent where the subscript $0$ on the right-hand side means ``with vanishing constant Fourier coefficient''. The pseudo-Laplacian with Dirichlet boundary condition $\Delta_{L,0}$ is the operator induced by the Chern Laplacian with Dirichlet boundary condition and this decomposition. Its determinant can be seen as the \textit{relative determinant}, a notion introduced by Müller in \cite{MR1617554}, of the Chern Laplacian with Dirichlet boundary condition $\Delta_L$ and of the Laplacian $-y^2 \mathrm{d}^2/\mathrm{d}y^2$ on $\left] a, + \infty \right[$ also with Dirichlet boundary condition.

    \subsection{Statement of the main result}
    
        This paper is devoted to two results related to the \textit{zeta-regularized determinant} of the pseudo-Laplacian with Dirichlet boundary condition. In the course of proving these formulas, we have to adapt in subsection \ref{SubSecWeylLaw} some computations from \cite{MR699488}, and find a slightly different argument, to take the presence of a line bundle into account.
        
        \medskip
        
        \hspace{10pt} $\bullet$ Our first result, in theorems \ref{ThmAlphaNeq0Mu} and \ref{ThmAlpha0Mu}, is a $\mu$-aymptotic expansion
            \begin{equation}
            \label{AsymptoticStudyMu}
                \begin{array}{lll}
                    \log \det \left( \Delta_{L,0} + \mu \right) & = & \fbox{$\mu$-divergent part} + \fbox{$\mu$-constant term} + o \left( 1 \right)
                \end{array}
            \end{equation}
            
            \noindent for the logarithm of the determinant of the pseudo-Laplacian (with Dirichlet boundary condition), as $\mu$ goes to infinity through strictly positive real values. This type of evaluation can be used to compute the constant in Mayer-Vietoris type formulas with parameter, in a way similar to \cite[Sec. 3.19 \& 4.8]{MR1165865} and \cite{MR1890995}.

        \medskip
        
        \hspace{10pt} $\bullet$ Our second result, in theorems \ref{ThmAlphaNeq0A} and \ref{ThmAlpha0A}, is an $a$-asymptotic expansion
            \begin{equation}
            \label{AsymptoticStudyA}
                \begin{array}{lll}
                    \log \det \Delta_{L,0} & = & \fbox{$a$-divergent part} + \fbox{$a$-constant term} + o \left( 1 \right)
                \end{array}
            \end{equation}
            
            \noindent for the logarithm of the determinant of the pseudo-Laplacian (with Dirichlet boundary condition), as the height $a$ of the cuspidal end goes to infinity, \textit{i.e.} as the cusp shrinks, without parameter $\mu$. This computation generalizes the case of the trivial line bundle, studied in \cite[Sec. 6]{MR4167014}.

    \subsection{Presentation of the paper}
    
        The technical nature of this paper makes it important to have an overview of the methods we use. This is achieved by splitting the reasoning into three parts: the first two, devoted to understanding the spectrum of the pseudo-Laplacian, serve as preparation for the third and most intricate part, where we obtain the asymptotic expansions \eqref{AsymptoticStudyMu} and \eqref{AsymptoticStudyA}.
    
        \subsubsection{Step $1$: preliminary work on the pseudo-Laplacian}
        
            Let us first go through the paragraphs comprising the first main step of this paper.
            
            \medskip
            
            \hspace{10pt} $\bullet$ In subsections \ref{SubSecCuspidalEnds} and \ref{SubSecFlatUnitaryLineBundles}, the definitions of cuspidal ends and flat unitary holomorphic line bundles on them are given.
            
            \medskip
            
            \hspace{10pt} $\bullet$ Subsection \ref{SubSecPseudoLaplacian} is devoted to the precise definition of the pseudo-Laplacian with Dirichlet boundary condition, including its domain in terms of Sobolev spaces, using the \textit{Friedrichs extension process}. This last notion is explained in appendix \ref{AppSelfAdjointOp}.
            
            \medskip
            
            \hspace{10pt} $\bullet$ The last part of this first step is subsection \ref{SubSecWeylLaw}, in which a Weyl law
            
                \begin{equation}
                \label{WeylLaw}
                    \begin{array}{lll}
                        N \left( \Delta_{L,0}, \lambda \right) & \leqslant & C \lambda
                    \end{array}
                \end{equation}
                
                \noindent is proved for any $\lambda > 0$, with $C > 0$ being a constant, mainly following \cite{MR699488}. To make one of the arguments used by Colin de Verdière more detailed, a Poincaré inequality is proved in lemma \ref{LemmaPoincaréInequality}, which results from the Banach--Alaoglu theorem. Unlike more standard versions of Weyl laws, the left-hand side of \eqref{WeylLaw} involves the \textit{spectral counting function}, defined using the Inf-Sup principle (see theorem \ref{ThmInfSup}), which exists even for self-adjoint positive definite operators with non-discrete spectra. This inequality proves that the pseudo-Laplacian has no essential spectrum, and that its spectral zeta function is holomorphic on the half-plane $\Re s > 1$.

        \subsubsection{Step $2$: localizing the eigenvalues}
        
            In the second step, comprised of subsections \ref{SubSecSpectralProblem} and \ref{SubSecLocalizationEigenvalues}, we study the eigenvalues of the pseudo-Laplacian with Dirichlet boundary condition, which amounts to solving the spectral problem
            
            \begin{equation}
            \tag{$S_0$}
            \label{SystPsi}
                \left \{ \begin{array}{llll}
                    - y^2 \left( \frac{\partial^2}{\partial x^2} + \frac{\partial^2}{\partial y^2} \right) \psi \left( x, y \right) & = & \multicolumn{2}{l}{\lambda \psi \left( x, y \right)} \\[0.7em]
                
                    \psi \left( x+1, y \right) & = & \multicolumn{2}{l}{e^{2 i \pi \alpha} \psi \left( x, y \right)} \\[0.7em]
                
                    \int_{S^1 \times \left] a, +\infty \right[} \; \left \vert \psi \right \vert^2 & < & + \infty & \text{(integrability condition)} \\[0.7em]
                
                    \psi \left( x, a \right) & = & 0 & \text{(Dirichlet boundary condition)} \\[0.7em]
                
                    \int_{S^1} \; \psi \left( x, y \right) \; \mathrm{d}x & = & 0 &  \text{for almost all } y > a \text{ if } \alpha = 0
                \end{array} \right.
            \end{equation}
            
            \noindent for smooth functions $\psi$ on $\mathbb{R} \times \left] a, +\infty \right[$. Using the change of function
            
            \[ \begin{array}{lll}
                    \varphi \left( x, y \right) & = & e^{- 2 i \pi \alpha x} \psi \left( x, y \right)
                \end{array} , \]
            
            \noindent we get a smooth and $1$-periodic in the first variable function, which we write as
            
            \[ \begin{array}{lll}
                    \varphi \left( x, y \right) & = & \sum\limits_{k \in \mathbb{Z}} a_k \left( y \right) e^{2 i \pi k x}
                \end{array} . \]
            
            \noindent Hence \eqref{SystPsi} gives us a differential equation for each $a_k$, which can be solved for every $k \in \mathbb{Z}$ and gives, up to multiplication by a constant depending only on $k$,
            
            \[ \begin{array}{lll}
                    a_k \left( y \right) & = & \sqrt{y} K_{s-1/2} \left( 2 \pi \left \vert k + \alpha \right \vert y \right)
                \end{array} , \]
            
            \noindent where $K$ denotes a \textit{modified Bessel function of the second kind}, for which the reader is referred to appendix \ref{AppModifiedBessel}, and $s$ is determined by $\lambda = s \left( 1-s \right)$. With the boundary condition $\psi \left( x,a \right) = 0$, the only possible eigenvalues are characterized by
            
            \[ \begin{array}{lllll}
                    a_k \left( a \right) & = & K_{s-1/2} \left( 2 \pi \left \vert k + \alpha \right \vert a \right) & = & 0
                \end{array} . \]
            
            \noindent In order to understand where the eigenvalues of $\Delta_{L,0}$ are located, we need more information on the zeros of the holomorphic function
            
            \[ \begin{array}{lll}
                    s & \longmapsto & K_{s-1/2} \left( 2 \pi \left \vert k + \alpha \right \vert a \right)
                \end{array} . \]
            
            \noindent This is the purpose of proposition \ref{PropZerosModifiedBessel}, which states that the function above only vanishes on a discrete subset of the line $\Re s = 1/2$, with $s=1/2$ being excluded. Such a result is proved by adapting Saharian's argument from \cite[Appendix A]{MR2439244}.

        \subsubsection{Step $3$: asymptotic studies}
        
            Using steps $1$ and $2$, we recover in subsection~\ref{SubSecIntegralRepresentationSpectralZeta} the spectral zeta function $\zeta_{L,\mu}$ of the pseudo-Laplacian with Dirichlet boundary conditions by using the argument principle with
            
            \[ \begin{array}{lll}
                    s & \longmapsto & K_{s-1/2} \left( 2 \pi \left \vert k + \alpha \right \vert a \right)
                \end{array} \]
            
            \noindent The holomorphy of the function above, as well as the simplicity of its zeros are paramount for this method to work. We get the following integral representation
            
            \[ \begin{array}{lll}
                    \zeta_{L,\mu} \left( s \right) & = & \frac{1}{2 i \pi} \sum\limits_k \int_{i \gamma_{\vartheta}} \left( \frac{1}{4} - t^2 + \mu \right)^{-s} \frac{\partial}{\partial t} \log K_t \left( 2 \pi \left \vert k + \alpha \right \vert a \right) \mathrm{d}t
                \end{array} , \]
            
            \noindent on the half-plane $\Re s > 1$, where the contour $\gamma_{\vartheta}$ surrounds the half-line of positive real numbers. The sum ranges over all integers $k \in \mathbb{Z}$, with $k=0$ being excluded by the ``vanishing constant Fourier coefficient'' condition if we have $\alpha = 0$. To make the computation possible, we want to let $\vartheta$ go to $\pi/2$. Avoiding convergence problems requires care, and we show in proposition \ref{PropSpectralZetaIntegralRepresentation2} that we have
            
            \begin{equation}
            \label{SpectralZetaIntegralRepresentation}
                \begin{array}{lll}
                    \zeta_{L, \mu} \left( s \right) & = & \frac{\sin \left( \pi s \right)}{\pi} \; \sum\limits_k \; \int_{\sqrt{\frac{1}{4}+\mu}}^{+\infty} \; \left( t^2 - \left( \frac{1}{4} + \mu \right) \right)^{-s} \; f_{\mu,k} \left( t \right) \; \text{d}t
                \end{array} ,
            \end{equation}
            
            \noindent with the function $f_{\mu,k}$ being given in definition \ref{Deffmuk} by
            
            \begin{equation}
            \label{fmuk}
                \scalemath{0.94}{\begin{array}{lll}
                    f_{\mu,k} \left( t \right) & = & \frac{\partial}{\partial t} \log K_t \left( 2\pi \left \vert k + \alpha \right \vert a \right) - \frac{2t}{\sqrt{4 \mu + 1}} \frac{\partial}{\partial t}_{\left \vert t = \sqrt{\frac{1}{4} + \mu} \right.} \log K_t \left( 2 \pi \left \vert k + \alpha \right \vert a \right).
                \end{array}}
            \end{equation}
            
            \noindent This last integral representation holds on the strip $1 < \Re s < 2$. In the course of subsections \ref{SubSecSplittingIntervalOfIntegration}, \ref{SubSecStudyLmuk}, and \ref{SubSecStudyMmuk}, the spectral zeta function undergoes several decompositions. A summary of these splittings is presented in subsection \ref{SubSecDiagramSplittings}. Starting with the integral representation \eqref{SpectralZetaIntegralRepresentation}, we set, in definition \ref{DefImuk},

            \[ \begin{array}{lll}
                    I_{\mu, k} \left( s \right) & = & \frac{\sin \left( \pi s \right)}{\pi} \int_{\sqrt{\frac{1}{4}+\mu}}^{+\infty} \; \left( t^2 - \left( \frac{1}{4} + \mu \right) \right)^{-s} \; f_{\mu,k} \left( t \right) \; \text{d}t
                \end{array} , \]
            
            \noindent on the strip $1 < \Re s < 2$. By analogy with a technical trick due to Freixas i Montplet and used in \cite[Sec. 6.1]{MR4167014}, the following decomposition of the interval of integration
            
            \[ \begin{array}{lll}
                    \left] \sqrt{\frac{1}{4} + \mu}, \, + \infty \right[ & = & \left] \sqrt{\frac{1}{4} + \mu}, \; 2 \left \vert k \right \vert^{\delta} \sqrt{\frac{1}{4} + \mu} \right[ \sqcup \left[ 2 \left \vert k \right \vert^{\delta} \sqrt{\frac{1}{4} + \mu}, \; +\infty \right[
                \end{array} \]
            
            \noindent is suggested, for non-zero integers $k$, for some parameter $\delta > 0$. In section \ref{SecDetPseudoLaplacian}, we find several inequalities which $\delta$ must satisfy. A ``small enough'' parameter $\delta > 0$ is taken in the end. When $\alpha$ does not vanish, we also consider the case $k=0$, where the interval can be split at any point. The integral $I_{\mu,k} \left( s \right)$ is then written as $I_{\mu,k} \left( s \right) = L_{\mu,k} \left( s \right) + M_{\mu,k} \left( s \right)$, with
            
            \begin{equation}
            \label{LmukMmuk}
                \begin{array}{lll}
                    L_{\mu,k} \left( s \right) & = & \frac{\sin \left( \pi s \right)}{\pi} \int_{\sqrt{\frac{1}{4}+\mu}}^{2 \left \vert k \right \vert^{\delta} \sqrt{\frac{1}{4} + \mu}} \; \left( t^2 - \left( \frac{1}{4} + \mu \right) \right)^{-s} \; f_{\mu,k} \left( t \right) \; \text{d}t, \\[2em]
                    
                    M_{\mu,k} \left( s \right) & = & \frac{\sin \left( \pi s \right)}{\pi} \int_{2 \left \vert k \right \vert^{\delta} \sqrt{\frac{1}{4} + \mu}}^{+ \infty} \; \left( t^2 - \left( \frac{1}{4} + \mu \right) \right)^{-s} \; f_{\mu,k} \left( t \right) \; \text{d}t.
                \end{array}
            \end{equation}
            
            \noindent The study of series with general term $L_{\mu,k} \left( s \right)$ is the focus of subsection \ref{SubSecStudyLmuk}, and the same is done with $M_{\mu,k} \left( s \right)$ in subsection \ref{SubSecStudyMmuk}.
            
            \medskip
            
            \hspace{10pt} $\bullet$ We begin studying $L_{\mu,k} \left( s \right)$ with proposition \ref{PropGlobalStudy}, which allows us to perform an integration by parts, resulting in the splitting $L_{\mu,k} \left( s \right) = A_{\mu,k} \left( s \right) + B_{\mu,k} \left( s \right)$, with
            
                \begin{equation}
                \label{IntroAB}
                    \begin{array}{lll}
                        A_{\mu,k} \left( s \right) & = & \frac{\sin \left( \pi s \right)}{\pi} \left( 4 \mu + 1 \right)^{-s} \left( \left \vert k \right \vert^{2 \delta} - \frac{1}{4} \right)^{-s} F_{\mu,k} \left( 2 \left \vert k \right \vert^{\delta} \sqrt{\frac{1}{4}+\mu} \right), \\[2em]
                    
                        B_{\mu,k} \left( s \right) & = & 2s \, \frac{\sin \left( \pi s \right)}{\pi} \, \int_{\sqrt{\frac{1}{4}+\mu}}^{2 \left \vert k \right \vert^{\delta} \sqrt{\frac{1}{4}+\mu}} \; t \left( t^2 - \left( \frac{1}{4}+\mu \right) \right)^{-s-1} \; F_{\mu,k} \left( t \right) \; \text{d}t.
                    \end{array}
                \end{equation}
            
                \noindent The function $F_{\mu,k}$ is defined for every $k \in \mathbb{Z}$ as a primitive of $f_{\mu,k}$ by
            
                \begin{equation}
                \label{Fmuk}
                    \begin{array}{lll}
                        F_{\mu,k} \left( t \right) & = & \log K_t \left( 2 \pi \left \vert k + \alpha \right \vert a \right) - \log K_{\sqrt{\frac{1}{4} + \mu}} \left( 2 \pi \left \vert k + \alpha \right \vert a \right) \\[1em]
                    
                        && \qquad \qquad \qquad - \frac{t^2 - \left( 1/4 + \mu \right)}{\sqrt{4 \mu + 1}} \frac{\partial}{\partial t}_{\left \vert t = \sqrt{\frac{1}{4} + \mu} \right.} \log K_t \left( 2 \pi \left \vert k + \alpha \right \vert a \right).
                    \end{array}
                \end{equation}

                \noindent The simpler of the two terms from \eqref{IntroAB} is $B_{\mu,k} \left( s \right)$. In proposition \ref{PropBmuk}, it is proved that the series with general term $B_{\mu,k} \left( s \right)$ has a holomorphic continuation to an open region of the complex plane containing $0$, and that its derivative at $s=0$ vanishes. For $A_{\mu,k} \left( s \right)$, using \eqref{IntroAB} and \eqref{Fmuk}, we see in proposition \ref{PropAmukAsymptA} that the function
                
                \[ \begin{array}{lll}
                        s & \longmapsto & \sum\limits_k A_{\mu,k} \left( s \right)
                    \end{array} \]
                
                \noindent has a holomorphic continuation near $0$, whose derivative at $s=0$ for $\mu=0$ satisfies
                
                \begin{equation}
                \label{AsymptAmukA}
                    \begin{array}{lll}
                        \frac{\partial}{\partial s}_{\left \vert s=0 \right.} \sum\limits_k A_{0,k} \left( s \right) & = & O \left( \frac{1}{a^2} \right)
                    \end{array}
                \end{equation}
                
                \noindent as $a$ goes to infinity. Note that the left-hand side of \eqref{AsymptAmukA} refers to ``the derivative of the continuation of''. This central result does not give the $\mu$-asymptotic behavior needed for \eqref{AsymptoticStudyMu}, as it relies on proposition \ref{PropGlobalStudy}, which uses the parameter asymptotics of the modified Bessel functions of the second kind, from proposition \ref{PropAsymptoticsModifiedBesselParameter}. It is seen in this last result that the remainder $\gamma_k$ would behave poorly with respect to $\mu$. To avoid that problem, we must, after having written
                
                \[ \begin{array}{lll}
                        \scriptstyle A_{\mu,k} \left( s \right) & \scriptstyle = & \scriptstyle \frac{\sin \left( \pi s \right)}{\pi} \left( 4 \mu+1 \right)^{-s} \left( \left \vert k \right \vert^{2 \delta} - \frac{1}{4} \right)^{-s} \Big[ \log K_{2 \left \vert k \right \vert^{\delta} \sqrt{\frac{1}{4} + \mu}} \left( 2 \pi \left \vert k + \alpha \right \vert a \right) - \log K_{\sqrt{\frac{1}{4} + \mu}} \left( 2 \pi \left \vert k + \alpha \right \vert a \right) \\[1em]
                        
                        && \qquad \qquad \qquad \qquad \qquad \quad \scriptstyle - \sqrt{4 \mu + 1} \left( \left \vert k \right \vert^{2 \delta} - \frac{1}{4} \right) \frac{\partial}{\partial t}_{\left \vert t = \sqrt{1/4 + \mu} \right.} \log K_t \left( 2 \pi \left \vert k + \alpha \right \vert a \right) \Big],
                    \end{array} \]
                
                \noindent use propositions \ref{PropLogK1}, \ref{PropLogK2}, and \ref{PropLogK3}, which give order-asymptotic expansions of the Bessel functions. The series with general term $A_{\mu,k} \left( s \right)$ is split into eleven parts.
            
                \medskip
            
                \hspace{30pt} $\blacktriangleright$ Let us begin with part $11$. The series with general term
                
                    \begin{equation}  
                    \label{PartialDerivativeTerm}
                        \scalemath{0.98}{\begin{array}{c}
                            - \frac{\sin \left( \pi s \right)}{\pi} \left( 4 \mu+1 \right)^{-s+1/2} \left( \left \vert k \right \vert^{2 \delta} - \frac{1}{4} \right)^{-s+1} \frac{\partial}{\partial t}_{\left \vert t = \sqrt{1/4 + \mu} \right.} \log K_t \left( 2 \pi \left \vert k + \alpha \right \vert a \right)
                        \end{array}}
                    \end{equation}
                    
                    \noindent is too complicated to be studied directly. We can use proposition \ref{PropGlobalStudy} to prove that the series has a holomorphic continuation near $0$, if we prove a similar result for all the other ten parts. This does not yield the $\mu$-asymptotic expension however. Fortunately, it is not the $\mu$-expansion of every term which matters to get \eqref{AsymptoticStudyMu}, but their sum. In paragraph \ref{SubSubSecStudyRmuk}, we end up studying the series with general term
                    
                    \[ \begin{array}{c}
                            \frac{1}{1-s} \cdot \frac{\sin \left( \pi s \right)}{\pi} \left( 4 \mu+1 \right)^{-s+1/2} \left( \left \vert k \right \vert^{2 \delta} - \frac{1}{4} \right)^{-s+1} \frac{\partial}{\partial t}_{\left \vert t = \sqrt{1/4 + \mu} \right.} \log K_t \left( 2 \pi \left \vert k + \alpha \right \vert a \right)
                        \end{array} , \]
                    
                    \noindent and we are faced with the same problem for the $\mu$-expansion. Since the sum of these two terms is what matters, we should consider the series with general term
                    
                    \[ \begin{array}{c}
                            \frac{s}{1-s} \cdot \frac{\sin \left( \pi s \right)}{\pi} \left( 4 \mu+1 \right)^{-s+1/2} \left( \left \vert k \right \vert^{2 \delta} - \frac{1}{4} \right)^{-s+1} \frac{\partial}{\partial t}_{\left \vert t = \sqrt{1/4 + \mu} \right.} \log K_t \left( 2 \pi \left \vert k + \alpha \right \vert a \right)
                        \end{array} . \]
                    
                    \noindent The extra factor $s$ would ideally make this derivative vanish entirely. For this to happen, however, the (continuation of) the series with general term \eqref{PartialDerivativeTerm} would need to vanish at $s=0$. This does not happen, and we need to remove some explicit terms from it, which are found in propositions \ref{PropA7} and \ref{PropA8}, and correspond to moments when the factor $\sin \left( \pi s \right)$ has to be used to cancel a simple pole.

                \medskip
                    
                \hspace{30pt} $\blacktriangleright$ Parts $1$ and $2$ of paragraph \ref{SubSubSecStudyAmuk} deal with the series involving the remainder term $\widetilde{\eta_2}$. They are comprised of propositions \ref{PropA1}, \ref{PropA1-0}, \ref{PropA2}, and \ref{PropA2-0}, and can be taken care of using the estimates detailed in corollary \ref{CorAsymptoticExpansionModifiedBessel}.

                \medskip
                
                \hspace{30pt} $\blacktriangleright$ Part $3$, in proposition \ref{PropA3}, proves, using a Taylor expansion and the known behavior of the Riemann zeta function, that the series with general term
                
                    \begin{equation}
                    \label{ThirdPart}
                        \begin{array}{c}
                            - \frac{\sin \left( \pi s \right)}{\pi} \left( 4 \mu+1 \right)^{-s} \left( \left \vert k \right \vert^{2 \delta} - \frac{1}{4} \right)^{-s} \sqrt{\left( 2 \pi \left \vert k + \alpha \right \vert a \right)^2 + \left( 4 \mu + 1 \right) \left \vert k \right \vert^{2 \delta}}
                        \end{array}
                    \end{equation}
                    
                    \noindent has a holomorphic continuation to an open neighborhood of $0$, and provides an expression for its derivative at $s=0$. This involves a derivative which cannot be computed asymptotically in $\mu$, but cancels one found in proposition \ref{PropMTilde2}.

                \medskip
                
                \hspace{30pt} $\blacktriangleright$ Parts $4$, $5$, and $6$, which are mainly comprised of propositions \ref{PropA4},~\ref{PropA5}, and \ref{PropA6}, are also taken care of using Taylor expansions, also with derivatives left uncomputed as they are canceled by derivatives found in propositions \ref{PropMTilde3} and \ref{PropMTilde4}.

                \medskip
                
                \hspace{30pt} $\blacktriangleright$ We now come to part $7$, which is the first in a series of much more complicated ones. In proposition \ref{PropA7}, we must handle the series with general term
                
                    \begin{equation}
                    \label{SeventhPart}
                        \begin{array}{c}
                            \frac{\sin \left( \pi s \right)}{\pi} \left( 4 \mu+1 \right)^{-s} \left( \left \vert k \right \vert^{2 \delta} - \frac{1}{4} \right)^{-s} \sqrt{\left( 2 \pi \left \vert k + \alpha \right \vert a \right)^2 + \frac{1}{4} + \mu}
                        \end{array} .
                    \end{equation}
                    
                    \noindent Proving that this series has a holomorphic continuation near $0$ could be done using a Taylor expansion, but unlike what we have done before, there are no cancellation with other terms found later on. Thus, the asymptotic expansion as $\mu$ goes to infinity of the derivative at $s=0$ must be found without any uncomputed term. A Taylor expansion cannot give us that. Let us see why on the simpler example of the \textit{Hurwitz zeta function}, which presents the same difficulty, defined for $\mu \geqslant 0$ by
                    
                    \[ \begin{array}{lll}
                            \zeta_H \left( s, 1+\mu \right) & = & \sum\limits_{k = 1}^{+ \infty} \left( k + \mu \right)^{-s}
                        \end{array} . \]
                    
                    \noindent Suppose we want to prove that $\zeta_H$ has a holomorphic continuation near $0$, and find an asymptotic expansion of its derivative at $s=0$ as $\mu$ goes to infinity. We have
                    
                    \[ \scalemath{0.95}{\begin{array}{lllll}
                            \zeta_H \left( s, 1+\mu \right) & = & \zeta \left( s \right) - \frac{\mu}{k} s \zeta \left( s+1 \right) + s \left( s+1 \right) \sum\limits_{k=1}^{+ \infty} k^{-s} \int_0^{\mu/k} \left( 1+x \right)^{-s-2} \left( \frac{\mu}{k} - x \right) \mathrm{d}x
                        \end{array}} \]
                    
                    \noindent using a Taylor expansion in $1/k$. The Hurwitz zeta function therefore has a holomorphic continuation near $0$, but increasing the convergence in $k$ has made a divergence in $\mu$ appear, and we cannot compute the derivative at $s=0$ asymptotically in this manner. One way to solve that problem is to find an integral representation for $\zeta_H$, but that cannot work for our more complicated examples, since we had to simplify an already existing representation. Another solution, called the \textit{Ramanujan summation}, is presented by Candelpherger in \cite{MR3677185}. This method, which is close to the Euler-Maclaurin and Abel-Plana formulas, is presented in appendix \ref{AppRamanujan}, and the asymptotic study of special values of $\zeta_H$ is made there as an example. We study the series with general term \eqref{SeventhPart} in proposition \ref{PropA7}. While doing that, we find one of the terms which must be removed from the partial derivative in part $11$.

                \medskip
                
                \hspace{30pt} $\blacktriangleright$ The remaining parts $8$, $9$, and $10$ of \ref{SubSubSecStudyAmuk} can be dealt with using the Ramanujan summation as well. Only the relevant results are given in this paper, as writting all the details would take significantly more space.

                \medskip
                
                This concludes \ref{SubSubSecStudyAmuk} and the study of the series with general term $A_{\mu,k} \left( s \right)$.

        \medskip
        
        \hspace{10pt} $\bullet$ We can now comment subsection \ref{SubSecStudyMmuk}, whose purpose is to study $M_{\mu,k} \left( s \right)$, as defined in \eqref{LmukMmuk}. We begin by writting 
        
            \[ \begin{array}{lll}
                    M_{\mu,k} \left( s \right) = \widetilde{M}_{\mu,k} \left( s \right) + R_{\mu,k} \left( s \right)
                \end{array} \]
            
            \noindent according to the definition of $f_{\mu,k}$ given by \eqref{fmuk}. More precisely, we set
        
            \[ \begin{array}{lll}
                    \widetilde{M}_{\mu,k} \left( s \right) & = & \frac{\sin \left( \pi s \right)}{\pi} \int_{2 \left \vert k \right \vert^{\delta} \sqrt{\frac{1}{4} + \mu}}^{+ \infty} \left( t^2 - \left( \frac{1}{4} + \mu \right) \right)^{-s} \frac{\partial}{\partial t} \log K_t \left( 2\pi \left \vert k + \alpha \right \vert a \right) \mathrm{d}t, \\[1.5em]
                    
                    R_{\mu,k} \left( s \right) & = & - \frac{\sin \left( \pi s \right)}{\pi} \cdot \frac{2}{\sqrt{4 \mu + 1}} \int_{2 \left \vert k \right \vert^{\delta} \sqrt{\frac{1}{4} + \mu}}^{+ \infty} t \left( t^2 - \left( \frac{1}{4} + \mu \right) \right)^{-s} \\[0.8em]
                    
                    && \qquad \qquad \qquad \qquad \qquad \qquad \quad \cdot \frac{\partial}{\partial t}_{\left \vert t = \sqrt{\frac{1}{4} + \mu} \right.} \log K_t \left( 2 \pi \left \vert k + \alpha \right \vert a \right) \mathrm{d}t.
                \end{array} \]
            
            \noindent First, we deal with $R_{\mu,k} \left( s \right)$ in \ref{SubSubSecStudyRmuk}, since the integral appearing in this term can be computed. This is done in lemma \ref{LemRmuk}, and we end up with an expression close to the term treated in proposition \ref{PropA11}, as we have already noted when we described the eleventh part of \ref{SubSubSecStudyAmuk}. Up to removing some explicit terms, we can cancel the derivative at $s=0$ of the series with general term $R_{\mu,k} \left( s \right)$, and these terms which have been removed are studied in proposition \ref{PropR}. To complete the study of this term, we need to find the $a$-asymptotic behavior for $\mu = 0$. To that effect, we use the asymptotics of the \textit{exponential integral} function $\mathbb{E}_1$ the explicit expression
            
            \[ \begin{array}{lll}
                    \frac{\partial}{\partial t}_{\left \vert t = 1/2 \right.} \log K_t \left( 2 \pi \left \vert k + \alpha \right \vert a \right) & = & \mathbb{E}_1 \left( 4 \pi \left \vert k + \alpha \right \vert a \right) e^{4 \pi \left \vert k + \alpha \right \vert a}
                \end{array} . \]
            
            \noindent Only the study of $\widetilde{M}_{\mu,k} \left( s \right)$ remains, and we now compute the logarithmic derivative of the Bessel function to get, in lemma \ref{LemLogDerivativeModifiedBessel},
            
            \[ \begin{array}{lll}
                    \scriptstyle \frac{\partial}{\partial t} \log K_t \left( 2 \pi \left \vert k + \alpha \right \vert a \right) & \scriptstyle = & \scriptstyle \arcsinh \left( \frac{t}{2 \pi \left \vert k + \alpha \right \vert a} \right) - \frac{1}{2} \cdot \frac{t}{t^2 + 4 \pi^2 \left( k + \alpha \right)^2 a^2} \\[1em]
                        
                    && \qquad \scriptstyle - \frac{\partial}{\partial t} \left( \frac{1}{t} U_1 \left( p \left( \frac{2 \pi \left \vert k + \alpha \right \vert a}{t} \right) \right) \right) + \frac{\partial}{\partial t} \left( \frac{1}{t^2} \widetilde{\eta_2} \left( t, \frac{1}{t} \cdot 2 \pi \left \vert k + \alpha \right \vert a \right) \right).
                \end{array} \]
            
            \noindent We split the remaining work into four parts induced by the decomposition above, and recall that we need to prove, for each part, the existence of a contination, find the $\mu$-asymptotic behavior for all $a > 0$, and the $a$-asymptotic behavior for $\mu = 0$ of the derivative at $s=0$.

            \medskip
            
            \hspace{30pt} $\blacktriangleright$ In part $1$, and more precisely in proposition \ref{PropMTilde1}, we take care of the remainder, which is the series with general term
            
                \[ \begin{array}{c}
                        \frac{\sin \left( \pi s \right)}{\pi} \int_{2 \left \vert k \right \vert^{\delta} \sqrt{\frac{1}{4} + \mu}}^{+ \infty} \left( t^2 - \left( \frac{1}{4} + \mu \right) \right)^{-s} \frac{\partial}{\partial t} \left( \frac{1}{t^2} \widetilde{\eta_2} \left( t, \frac{1}{t} \cdot 2 \pi \left \vert k + \alpha \right \vert a \right) \right) \mathrm{d}t
                    \end{array} . \]
                
                \noindent We use an integration by parts and the upper bounds on $\widetilde{\eta_2}$ found in \ref{CorAsymptoticExpansionModifiedBessel}.

            \medskip
            
            \hspace{30pt} $\blacktriangleright$ Part $2$, in proposition \ref{PropMTilde2}, studies the series with general term
                
                \[ \begin{array}{c}
                        \frac{\sin \left( \pi s \right)}{\pi} \int_{2 \left \vert k \right \vert^{\delta} \sqrt{\frac{1}{4} + \mu}}^{+ \infty} \left( t^2 - \left( \frac{1}{4} + \mu \right) \right)^{-s} \arcsinh \left( \frac{t}{2 \pi \left \vert k + \alpha \right \vert a} \right) \mathrm{d}t
                    \end{array} . \]
                    
                \noindent The guiding principle in this study is to make step-by-step simplifications, and, as much as possible, to make a factor $s$ appear. First, we use the binomial formula, which is recalled as proposition \ref{PropBinomialFormula}, on the complex power, yielding
                
                \[ \begin{array}{lll}
                        \left( t^2 - \left( \frac{1}{4} + \mu \right) \right)^{-s} & = & \sum\limits_{j = 0}^{+ \infty} \frac{\left( s \right)_j}{j!} \left( \frac{1}{4} + \mu \right)^j t^{-2 \left( s + j \right)}
                    \end{array} . \]
                
                \noindent The sum and the integral can be interchanged, and we can perform an integration by parts, to replace $\arcsinh$ by a fraction, in order to have an integral similar to the one from corollary \ref{CorIntegralFormula}. Since the Pochhammer symbol vanishes at $s=0$ for all $j \geqslant 1$, we can now prove that the sum over $j \geqslant 2$ has a continuation around $0$, and that its derivative at $s=0$ vanishes. We are are thus reduced to dealing with the terms corresponding to $j=0$ and $j=1$, separately. Each of these is computed using hypergeometric functions. What follows is a cumbersome calculation, which involves various formulas related to hypergeometric functions, all of which are presented in appendix \ref{AppHypergeometric}. The $\mu$ and $a$-asymptotic studies are then obtained almost simultaneously. In all of this, we find two derivatives which cannot be computed as $\mu$ goes to infinity. Fortunately, they cancel the derivatives left aside in propositions \ref{PropA7} and \ref{PropA8}.

            \medskip
            
            \hspace{30pt} $\blacktriangleright$ Parts $3$ and $4$ are dealt with in a similar fashion. Let us only mention that in part $2$, and more precisely in proposition \ref{PropMTilde3}, we find a derivative which cannot be computed as $\mu$ goes to infinity, and which is used to cancel the one left aside in proposition \ref{PropA9}. This is done by using \textit{Euler's integral formula}, recalled in proposition \ref{PropEulerIntegralFormula}, on hypergeometric functions.

    \subsection{Summary of the splittings}
    \label{SubSecDiagramSplittings}
    
        The following diagram sums up the splittings performed on the spectral zeta function, and points to the relevant results in the paper for the various parts.
        
        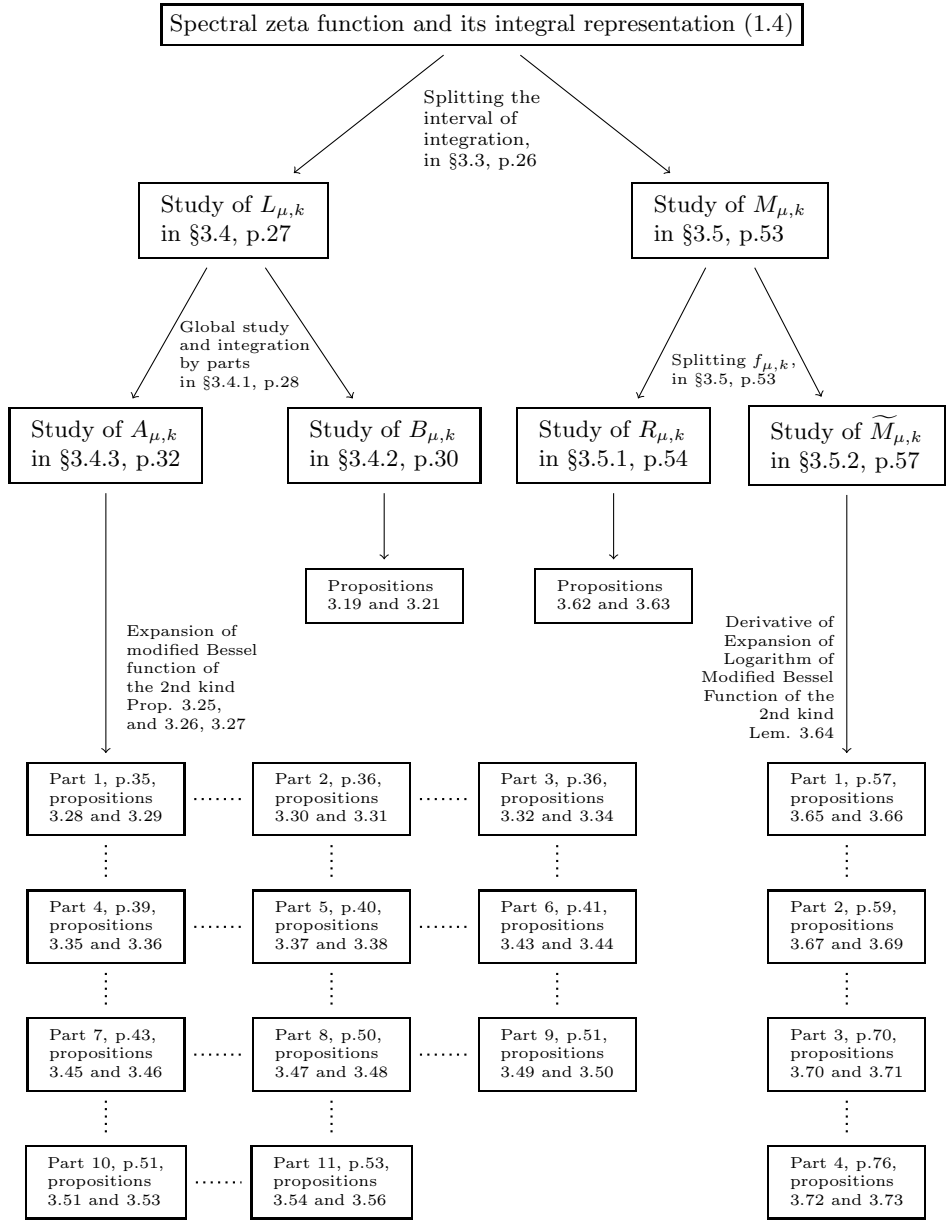
\begin{figure}[H]
    
            \centering
    
            \begin{tikzpicture}
        
                \node (A) at (-1,0) {\fbox{\small Spectral zeta function and its integral representation \eqref{SpectralZetaIntegralRepresentation}}};
                    
                \node (B) at (-4.3,-2.6) {\fbox{\small \begin{tabular}{l}Study of $L_{\mu,k}$ \\ in \S \ref{SubSecStudyLmuk}, p.\pageref{SubSecStudyLmuk} \end{tabular}}};
                \node (C) at (2.3,-2.6) {\fbox{\small \begin{tabular}{l}Study of $M_{\mu,k}$ \\ in \S \ref{SubSecStudyMmuk}, p.\pageref{SubSecStudyMmuk} \end{tabular}}};
                    
                \node (ABC) at (-1,-1.4) {\scriptsize \begin{tabular}{l} Splitting the \\ interval of \\ integration, \\ in \S \ref{SubSecSplittingIntervalOfIntegration}, p.\pageref{SubSecSplittingIntervalOfIntegration} \end{tabular}};
                    
                \node (D) at (-6,-5.6) {\fbox{\small \begin{tabular}{l} Study of $A_{\mu,k}$ \\ in \S \ref{SubSubSecStudyAmuk}, p.\pageref{SubSubSecStudyAmuk} \end{tabular}}};
                \node (E) at (-2.3,-5.6) {\fbox{\small \begin{tabular}{l} Study of $B_{\mu,k}$ \\ in \S \ref{SubSubSecStudyBmuk}, p.\pageref{SubSubSecStudyBmuk} \end{tabular}}};
                    
                \node (BDE) at (-4.15,-4.4) {\tiny \begin{tabular}{l} Global study \\ and integration \\ by parts \\ in \S \ref{SubSubSecGlobalStudy}, p.\pageref{SubSubSecGlobalStudy} \end{tabular}};
                    
                \node (F) at (0.75,-5.6) {\fbox{\small \begin{tabular}{l} Study of $R_{\mu,k}$ \\ in \S \ref{SubSubSecStudyRmuk}, p.\pageref{SubSubSecStudyRmuk} \end{tabular}}};
                \node (G) at (3.85,-5.6) {\fbox{\small \begin{tabular}{l} Study of $\widetilde{M}_{\mu,k}$ \\ in \S \ref{SubSubSecStudyMTildemuk}, p.\pageref{SubSubSecStudyMTildemuk} \end{tabular}}};
                    
                \node (AGF) at (2.35,-4.6) {\tiny \begin{tabular}{l} Splitting $f_{\mu,k}$, \\ in \S \ref{SubSecStudyMmuk}, p.\pageref{SubSecStudyMmuk} \end{tabular}};
                    
                \node (H) at (-2.3, -7.6) {\fbox{\tiny \begin{tabular}{l} Propositions \\\ref{PropBmuk} and \ref{PropBmu0} \end{tabular}}};
                    
                \node (I) at (0.75,-7.6) {\fbox{\tiny \begin{tabular}{l} Propositions \\\ref{PropR} and \ref{PropR-0} \end{tabular}}};
                    
                \node (J) at (-6,-10.3) {\fbox{\tiny \begin{tabular}{l} Part 1, p.\pageref{AFirstPart}, \\ propositions \\ \ref{PropA1} and \ref{PropA1-0} \end{tabular}}};
                    
                \node (K) at (-6,-12) {\fbox{\tiny \begin{tabular}{l} Part 4, p.\pageref{AFourthPart}, \\ propositions \\ \ref{PropA4} and \ref{PropA4-0} \end{tabular}}};
                    
                \node (L) at (-6,-13.7) {\fbox{\tiny \begin{tabular}{l} Part 7, p.\pageref{ASeventhPart}, \\ propositions \\ \ref{PropA7} and \ref{PropA7-0} \end{tabular}}};
                    
                \node (M) at (-6,-15.4) {\fbox{\tiny \begin{tabular}{l} Part 10, p.\pageref{ATenthPart}, \\ propositions \\ \ref{PropA10} and \ref{PropA10-0} \end{tabular}}};

                \node (N) at (-3,-10.3) {\fbox{\tiny \begin{tabular}{l} Part 2,  p.\pageref{ASecondPart}, \\ propositions \\ \ref{PropA2} and \ref{PropA2-0} \end{tabular}}};
                    
                \node (O) at (-3,-12) {\fbox{\tiny \begin{tabular}{l} Part 5, p.\pageref{AFifthPart}, \\ propositions \\ \ref{PropA5} and \ref{PropA5-0} \end{tabular}}};
                    
                \node (P) at (-3,-13.7) {\fbox{\tiny \begin{tabular}{l} Part 8, p.\pageref{AEighthPart}, \\ propositions \\ \ref{PropA8} and \ref{PropA8-0} \end{tabular}}};
                    
                \node (Q) at (-3,-15.4) {\fbox{\tiny \begin{tabular}{l} Part 11, p.\pageref{AEleventhPart}, \\ propositions \\ \ref{PropA11} and \ref{PropA11-0} \end{tabular}}};

                \node (R) at (0,-10.3) {\fbox{\tiny \begin{tabular}{l} Part 3, p.\pageref{AThirdPart}, \\ propositions \\ \ref{PropA3} and \ref{PropA3-0} \end{tabular}}};
                    
                \node (S) at (0,-12) {\fbox{\tiny \begin{tabular}{l} Part 6, p.\pageref{ASixthPart}, \\ propositions \\ \ref{PropA6} and \ref{PropA6-0} \end{tabular}}};
                    
                \node (T) at (0,-13.7) {\fbox{\tiny \begin{tabular}{l} Part 9, p.\pageref{ANinthPart}, \\ propositions \\ \ref{PropA9} and \ref{PropA9-0} \end{tabular}}};

                \node (U) at (3.85,-10.3) {\fbox{\tiny \begin{tabular}{l} Part 1,  p.\pageref{MTildeFirstPart}, \\ propositions \\ \ref{PropMTilde1} and \ref{PropMTilde1-0} \end{tabular}}};
                    
                \node (V) at (3.85,-12) {\fbox{\tiny \begin{tabular}{l} Part 2,  p.\pageref{MTildeSecondPart}, \\ propositions \\ \ref{PropMTilde2} and \ref{PropMTilde2-0} \end{tabular}}};
                    
                \node (W) at (3.85,-13.7) {\fbox{\tiny \begin{tabular}{l} Part 3,  p.\pageref{MTildeThirdPart}, \\ propositions \\ \ref{PropMTilde3} and \ref{PropMTilde3-0} \end{tabular}}};
                    
                \node (X) at (3.85,-15.4) {\fbox{\tiny \begin{tabular}{l} Part 4,  p.\pageref{MTildeFourthPart}, \\ propositions \\ \ref{PropMTilde4} and \ref{PropMTilde4-0} \end{tabular}}};

                \node (DJ) at (-4.85, -8.7) {\tiny \begin{tabular}{l} Expansion of \\ modified Bessel \\ function of \\ the $2$nd kind \\ Prop. \ref{PropLogK1}, \\ and \ref{PropLogK2}, \ref{PropLogK3} \end{tabular}};
                    
                \node (GU) at (2.8, -8.7) {\tiny \begin{tabular}{r} Derivative of \\ Expansion of \\ Logarithm of \\ Modified Bessel \\ Function of the \\ $2$nd kind \\ Lem. \ref{LemLogDerivativeModifiedBessel} \end{tabular}};

                \draw[->] (A)--(B);
                \draw[->] (A)--(C);
                    
                \draw[->] (B)--(D);
                \draw[->] (B)--(E);
                    
                \draw[->] (C)--(F);
                \draw[->] (C)--(G);
                    
                \draw[->] (E)--(H);
                \draw[->] (F)--(I);
                    
                \draw[->] (D)--(J);
                    
                \draw[->] (G)--(U);
                    
                \draw[dotted, thick] (J)--(N);
                \draw[dotted, thick] (N)--(R);
                \draw[dotted, thick] (J)--(K);
                \draw[dotted, thick] (N)--(O);
                \draw[dotted, thick] (R)--(S);
                \draw[dotted, thick] (K)--(O);
                \draw[dotted, thick] (O)--(S);
                \draw[dotted, thick] (K)--(L);
                \draw[dotted, thick] (O)--(P);
                \draw[dotted, thick] (S)--(T);
                \draw[dotted, thick] (L)--(P);
                \draw[dotted, thick] (P)--(T);
                \draw[dotted, thick] (L)--(M);
                \draw[dotted, thick] (P)--(Q);
                \draw[dotted, thick] (M)--(Q);

                \draw[dotted, thick] (U)--(V);
                \draw[dotted, thick] (V)--(W);
                \draw[dotted, thick] (W)--(X);

            \end{tikzpicture}
        
            \caption{The splittings of the spectral zeta function}
            \label{FigSplittings}
        
        \end{figure}
            
        \noindent In this diagram, plain arrows represent splittings, and dotted lines are meant to show which parts result from a given decomposition. Hence, parts $1$ to $11$ on the left-hand side are not linked to one another, and all result from the study of the terms $A_{\mu,k}$. In turn, they are unrelated to parts $1$ to $4$ on the right-hand side.

    \subsection{Acknowledgements}
    
        The work presented here was first done during my PhD at the \textit{Institut de Mathématiques de Jussieu-Paris Rive Gauche}, and funding from the DIM RDM IDF program of the \textit{région Île-de-France} should be acknowledged for that period. I would especially like to thank Gerard Freixas i Montplet, who introduced me to Arakelov geometry and the problem of obtaining Deligne--Riemann--Roch isometries, which led to this paper. As a PhD advisor, his comments, feedback, and more generally all the mathematical discussions we had proved invaluable. I would also like to thank Colin Guillarmou and Kai Köhler, who acted as referees for the doctoral thesis from which this work is extracted. Moreover, I want to thank José Burgos Gil and Xiaonan Ma for the discussions we had, and ideas they gave for future work. Finally, I would like to thank Siarhei Finski for the fruitful mathematical discussions we had, as well as Manish Patnaik, not only for our discussions, but also for his feedback.

\section{Description of the spectral problem}
\label{SecDescriptionSpectralProblem}

    \subsection{Cuspidal ends}
    \label{SubSecCuspidalEnds}
    
        This paper is focused on the study of metric singularities known as \textit{cuspidal ends}. The first task is to review what these are. Let $a > 0$.
        
        \begin{proposition}
                
            The translation $T$, defined by
            
            \[ \begin{array}{lllll}
                    T & : & \mathbb{R} \times \left] a, + \infty \right[ & \longrightarrow & \mathbb{R} \times \left] a, + \infty \right[ \\[0.5em]
                    && \left( x, y \right) & \longmapsto & \left( x+1, y \right)
                \end{array} , \]
            
            \noindent is a bijection of $\mathbb{R} \times \left] a, + \infty \right[$ of infinite order. The subgroup it generates is canonically identified to $\mathbb{Z}$ by $T \mapsto 1$, and acts on $\mathbb{R} \times \left] a, + \infty \right[$.
            
        \end{proposition}
        
        \begin{proof}
        
            This result is direct, the action of $T$ on $\mathbb{R} \times \left] a, + \infty \right[$ being the natural one.
        
        \end{proof}
        
        \begin{definition}
        
            The \textit{cuspidal end} of height $a$ is defined as the product $S^1 \times \left] a, + \infty \right[$, which is the quotient $\mathbb{Z} \backslash \left( \mathbb{R} \times \left] a, + \infty \right[ \right)$, endowed with the \textit{Poincaré metric}
            
            \[ \begin{array}{lll}
                    \mathrm{d}s^2_{\mathrm{hyp}} & = & \frac{\mathrm{d}x^2 + \mathrm{d}y^2}{y^2}
                \end{array} . \]
            
            \noindent It is implicitly assumed here that $S^1$ is parametrized by $x \in \left[ 0, 1 \right[$.
        
        \end{definition}
        
        \begin{proposition}
        
            The cuspidal end of height $a > 0$ is isometric to the punctured disk $D^{\times} \left( 0, \varepsilon \right)$ of radius $\varepsilon = \exp \left( - 2 \pi a \right)$ with the Poincar\'{e} metric
            
            \[ \begin{array}{lll}
                    \mathrm{d}s^2_{\mathrm{hyp}} & = & \frac{\left\vert \mathrm{d} z \right \vert^2}{\left( \left \vert z \right \vert \log \left \vert z \right\vert \right)^2}
                \end{array} . \]
        
        \end{proposition}
        
        \begin{proof}
        
            Let $\varepsilon = \exp \left( - 2 \pi a \right)$. The map
            
            \[ \begin{array}{lllll}
                    \varphi & : & \mathbb{R} \times \left] a, +\infty \right[ & \longrightarrow & D^{\times} \left( 0, \varepsilon \right) \\[0.5em]
                    
                    && \left( x, y \right) & \longmapsto & e^{2 i \pi \left( x + iy \right)}
                \end{array} \]
            
            \noindent is invariant by the action of $\mathbb{Z}$, and thus induces a map $S^1 \times \left] a, + \infty \right[ \longrightarrow D^{\times} \left( 0, \varepsilon \right)$, which is bijective. Using this map, and denoting by $z$ the coordinate on the punctured disk, we obtain the Poincaré metric on $D^{\times} \left( 0, \varepsilon \right)$.
        
        \end{proof}
        
        \begin{remark}
        
            The Poincaré metric is \textit{singular} at $z=0$, meaning it cannot be extended into a smooth metric on the full disk $D \left( 0, \varepsilon \right)$. This is because we have
            
            \[ \begin{array}{lll}
                    \lim\limits_{z \rightarrow 0} \frac{1}{\left( \left \vert z \right \vert \log \left \vert z \right\vert \right)^2} & = & + \infty
                \end{array} . \]
            
            \noindent This can be seen as a loss of control on the metric as we approach the cusp.
        
        \end{remark}

    \subsection{Flat unitary line bundles}
    \label{SubSecFlatUnitaryLineBundles}
    
        The second part of the setting we consider here is that of flat unitary line bundles on cuspidal ends.
        
        \begin{definition}
        
            Let $\chi : \mathbb{Z} \longrightarrow \mathbb{C}^{\ast}$ be a unitary character. The group $\mathbb{Z}$ acts on the trivial line bundle $\mathbb{C}$ over $\mathbb{R} \times \left] a, +\infty \right[$ by
            
            \[ \begin{array}{lllll}
                    k \cdot \left( \left( x, y \right), \lambda \right) & = & \left( \left( x+k, y \right), \chi \left( k \right) \lambda \right) & = & \left( \left( x+k, y \right), \chi \left( 1 \right)^k \lambda \right)
                \end{array} . \]
        
        \end{definition}
        
        \begin{remark}
        
            Since $\mathbb{Z}$ is generated by $1$, a unitary character of $\mathbb{Z}$ is entirely determined by its value at $1$, which takes the form
            
            \[ \begin{array}{lll}
                    \chi \left( 1 \right) & = & e^{2 i \pi \alpha}
                \end{array} , \]
            
            \noindent where $\alpha$ is a real number well-defined modulo $1$.
        
        \end{remark}
        
        \begin{remark}
        
            The usual complex modulus on $\mathbb{C}$, which is its canonical Hermitian metric, is compatible with the action of $\mathbb{Z}$.
        
        \end{remark}
        
        \begin{proposition-definition}
        
            The quotient $L = \mathbb{Z} \backslash \left( \left( \mathbb{R} \times \left] a, +\infty \right[ \right) \times \mathbb{C} \right)$ is a holomorphic line bundle on the cuspidal end. Furthermore, the Hermitian metric on $L$ induced by the modulus on $\mathbb{C}$ is a flat metric, called the canonical flat metric.
        
        \end{proposition-definition}
        
        \begin{proof}
        
            This is a classical result.
        
        \end{proof}
        
        \begin{remark}
        \label{RmkIdentificationSectionFunction}
        
            Using the definition of the action of $\mathbb{Z}$ on $\left( \mathbb{R} \times \left] a, +\infty \right[ \right) \times \mathbb{C}$, one notes that smooth sections of $L$ over the cuspidal end can be identified to smooth functions $f : \mathbb{R} \times \left] a, + \infty \right[ \longrightarrow \mathbb{C}$ such that we have $f \left( x+1, y \right) = e^{2 i \pi \alpha} f \left( x,y \right)$.
        
        \end{remark}

    \subsection{Pseudo-Laplacian}
    \label{SubSecPseudoLaplacian}
    
        The main operator we will be concerned with is a type of Laplacian, similar to the one used by Colin de Verdière in \cite{MR688031, MR699488}. Let us see how it is defined. The Laplacian acting on smooth sections of $L$ associated to the Chern connection (see \cite[Prop 3.12]{MR2451566}) is denoted by $\Delta_L$ and called the Chern Laplacian.
        
        \begin{remark}
        
            For any smooth section $f$ of $L$ over the cuspidal end, remark \ref{RmkIdentificationSectionFunction} provides an identification between the section $\Delta_L f$ of $L$ and the function
            
            \[ \begin{array}{lllll}
                    -y^2 \left( \frac{\partial^2}{\partial x^2} + \frac{\partial^2}{\partial y^2} \right) f & : & \mathbb{R} \times \left] a,+\infty \right[ & \longrightarrow & \mathbb{C}
                \end{array} . \]
            
            \noindent This is due to the flatness of the Chern connection.
        
        \end{remark}
        
        \begin{definition}
        
            The Sobolev-type space $L^2_{2, \text{Dir}} \left( S^1 \times \left] a, +\infty \right[, L \right)$ is defined as
            
            \[ \begin{array}{lll}
                    \multicolumn{3}{l}{L^2_{2, \text{Dir}} \left( S^1 \times \left] a, +\infty \right[, L \right)} \\[0.7em]
                    
                    \qquad & = & \left \{ u \in L^2 \left( S^1 \times \left] a, +\infty \right[, L \right), \; \Delta_L u \in L^2 \left( S^1 \times \left] a, +\infty \right[, L \right), \; \gamma u = 0 \right \},
                \end{array} \]
            
            \noindent where the Laplacian $\Delta_L$ is considered in the distributional sense on the right-hand side, and where $\gamma$ is the boundary trace operator, meaning the restriction to the boundary. It is called the $L^2_2$-Sobolev space with \textit{Dirichlet boundary conditions}.
        
        \end{definition}
        
        \begin{remark}
        
            Using more common notations, the Sobolev space defined above could be seen as an intersection $H^2 \cap H^1_0$. However, it is more important to have simpler notations here, as well as in works related to this paper, where boundary conditions will be more complicated.
        
        \end{remark}
        
        \begin{proposition}
        
            The Chern Laplacian $\Delta_L$, acting on smooth compactly supported sections of $L$ over $S^1 \times \left] a, +\infty \right[$, is a symmetric positive operator. Its Friedrichs extension is a positive $L^2$ self-adjoint operator
            
            \[ \begin{array}{lllll}
                    \Delta_L & : & L^2_{2, \mathrm{Dir}} \left( S^1 \times \left] a, +\infty \right[, L \right) & \longrightarrow & L^2 \left( S^1 \times \left] a, +\infty \right[, L \right)
                \end{array} , \]
            
            \noindent called the Chern Laplacian with Dirichlet boundary conditions.
        
        \end{proposition}
        
        \begin{proof}
        
            Consider the Chern Laplacian
            
            \[ \begin{array}{lllll}
                    \Delta_L & : & \mathcal{C}^{\infty}_0 \left( S^1 \times \left] a, + \infty \right[, L \right) & \longrightarrow & \mathcal{C}^{\infty}_0 \left( S^1 \times \left] a, + \infty \right[, L \right)
                \end{array} \]
            
            \noindent on smooth compactly supported sections of $L$. It is a positive symmetric operator. The closure $\overline{Q_{\Delta_L}}$ of the associated quadratic form is defined on the completion of the domain of $\Delta_L$ for the norm given in definition \ref{DefQuadFormNorm}. Hence, we have
            
            \[ \begin{array}{lll}
                    \mathrm{Dom} \, \overline{Q_{\Delta_L}} & = & L^2_{1,\mathrm{Dir}} \left( S^1 \times \left] a, + \infty \right[, L \right)
                \end{array} , \]
            
            \noindent this Sobolev space being defined in terms of the Chern connection. We now need to find the domain of the adjoint of $\Delta_L$. We have
            
            \[ \begin{array}{lllll}
                    \mathrm{Dom} \, \Delta_L^{\ast} & = & \left \{ u \in L^2 \left( S^1 \times \left] a, + \infty \right[, L \right), \; \Delta_L u \in L^2 \left( S^1 \times \left] a, + \infty \right[, L \right) \right \} \\[1em]
                    
                    & = & L^2_2 \left( S^1 \times \left] a, + \infty \right[, L \right)
                \end{array} \]
        
            \noindent Using remark \ref{RmkDomainFriedrichs}, the domain of the Friedrichs extension of $\Delta_L$ is given by the intersection of these last two domains. Still denoting the extension $\Delta_L$, we get
            
            \[ \begin{array}{lll}
                    \mathrm{Dom} \, \Delta_L & = & L^2_{1,\mathrm{Dir}} \left( S^1 \times \left] a, + \infty \right[, L \right) \cap L^2_2 \left( S^1 \times \left] a, + \infty \right[, L \right) \\[1em]
                    
                    & = & L^2_{2, \mathrm{Dir}} \left( S^1 \times \left] a, +\infty \right[, L \right).
                \end{array} \]
            
            \noindent This concludes the proof of the proposition.
        
        \end{proof}
        
        \begin{remark}
        
            In general, Sobolev spaces can be defined in more than one way, for instance in terms of the Chern connection, or of fractional powers of the Chern Laplacian. They coincide in the case we study here. For comprehensive comparisons of such spaces in delicate situations, see \cite{MR2343536}.
        
        \end{remark}
        
        As will be made clear later, the operator $\Delta_L$ and its eigenvalues behave well when the character $\chi$ is non-trivial. When $\chi$ is trivial, we modify it to get a ``pseudo-Laplacian'', similar to the one used by Colin de Verdière in \cite{MR688031, MR699488}. Let us temporarily assume, up until definition \ref{DefPseudoLaplacian}, that $\chi$ is trivial. Since $L$ is then metrically trivial, we omit it from the notation. Let $f$ be a smooth function on $S^1 \times \left] a, +\infty \right[$, seen as a smooth function on $\mathbb{R} \times \left] a, + \infty \right[$ such that we have $f \left( x+1, y \right) = f \left( x,y \right)$. Its Fourier decomposition is

        \[ \begin{array}{lll}
                f \left( x, y \right) & = & a_0 \left( y \right) + \sum\limits_{k \neq 0} \; a_k \left( y \right) e^{2 i k \pi x}.
            \end{array} \]
        
        \noindent Computing the Laplacian of $f$ then yields
        
        \[ \begin{array}{lll}
                \Delta f & = & \left( - y^2 \frac{\mathrm{d}^2}{\mathrm{d}y^2} \right) a_0 \left( y \right) + \left( - y^2 \left( \frac{\partial^2}{\partial x^2} + \frac{\partial^2}{\partial y^2} \right) \right) \left( \sum\limits_{k \neq 0} \; a_k \left( y \right) e^{2 i k \pi x} \right)
            \end{array} . \]
        
        \noindent The constant term in this expansion $a_0 \left( y \right)$ is of particular interest, and we set
        
        \[ \begin{array}{lllll}
                p & : & L^2 \left( S^1 \times \left] a, + \infty \right[ \right) & \longrightarrow & L^2 \left( S^1 \times \left] a, + \infty \right[ \right) \\[0.5em]
                
                && f & \longmapsto & a_0 \left( y \right)
            \end{array} . \]
        
        \noindent This map is surjective.
        
        \begin{definition}
        
            The kernel of $p$, \textit{i.e.} the space of $L^2$ functions with vanishing constant Fourier coefficient, is denoted by $L^2 \left( S^1 \times \left] a, +\infty \right[ \right)_0$.
        
        \end{definition}
        
        \begin{proposition}
        
            We have the following orthogonal decomposition
            
            \[ \begin{array}{lll}
                    L^2 \left( S^1 \times \left] a, +\infty \right[ \right) & = & L^2 \left( S^1 \times \left] a, +\infty \right[ \right)_0 \oplus L^2 \left( \left] a, +\infty \right[ \right)
                \end{array} . \]
            
            \noindent Furthermore, the Chern Laplacian with Dirichlet boundary condition splits
            
            \[ \begin{array}{lll}
                    \Delta & = & \Delta \oplus \left( - y^2 \frac{\mathrm{d}^2}{\mathrm{d}y^2} \right)
                \end{array} , \]
            
            \noindent where the Laplacian on the right-hand side acts on
            
            \[ \begin{array}{lll}
                    L^2_{2, \mathrm{Dir}} \left( S^1 \times \left] a, +\infty \right[ \right)_0 & = & L^2_{2, \mathrm{Dir}} \left( S^1 \times \left] a, +\infty \right[ \right) \cap L^2 \left( S^1 \times \left] a, +\infty \right[ \right)_0
                \end{array} . \]
        
        \end{proposition}
        
        For this last definition, let us go back to the more general setting of a possibly non trivial flat unitary line bundle $L$.
        
        \begin{definition}
        \label{DefPseudoLaplacian}
        
            The \textit{pseudo-Laplacian with Dirichlet boundary condition} $\Delta_{L,0}$ is defined to be the Chern Laplacian with Dirichlet boundary condition:
            
            \begin{itemize}
                \item acting on $L^2_{2, \mathrm{Dir}} \left( S^1 \times \left] a, +\infty \right[ \right)$, if the character $\chi$ is non-trivial;
                
                \item acting on $L^2_{2, \mathrm{Dir}} \left( S^1 \times \left] a, +\infty \right[ \right)_0$ if $\chi$, and thus $L$, is trivial.
            \end{itemize}
        
        \end{definition}
        
        \begin{remark}
        
            When $\chi$ is trivial, there is an added condition of vanishing constant Fourier coefficient. This will be important in subsection \ref{SubSecSpectralProblem} to make sense of the \textit{spectral zeta function} of the pseudo-Laplacian.
        
        \end{remark}

    \subsection{Weyl type law and the spectral zeta function}
    \label{SubSecWeylLaw}
    
        Before we can define the spectral zeta function of the pseudo-Laplacian with Dirichlet boundary conditions, we need some information on the distribution and the multiplicity of the eigenvalues. To that effect, we obtain the following Weyl type law.
        
        \begin{theorem}[Weyl type law]
        \label{ThmWeylLaw}
        
            There exists a constant $C > 0$ such that we have, for any strictly positive real number $\lambda$,
            
            \[ \begin{array}{lll}
                    N \left( \Delta_{L,0}, \lambda \right) & \leqslant & C \lambda
                \end{array} , \]
            
            \noindent where the spectral counting function is presented in definition \ref{DefSpectralCountingFunction}.
        
        \end{theorem}
        
        \begin{remark}
        
            This estimate is obtained using similar computations to those performed by Colin de Verdière in \cite[Thm 6]{MR699488}. However, to avoid having to prove separately that $\Delta_L$ has a compact resolvent, we count the real numbers appearing in the sequence yielded by the Inf-Sup theorem (see theorem \ref{ThmInfSup}), instead of simply the eigenvalues. Furthermore, the argument from Colin de Verdière is slightly modified, so as not to deal with a version of proposition \ref{PropSpectralCountingDirectSum} for infinite sums.
        
        \end{remark}
        
        \noindent Before moving to the proof of the Weyl type law, let us see its consequences.
        
        \begin{corollary}
        \label{CorPseudoLaplacianHasNoEssentialSpectrum}
        
            The pseudo-Laplacian with Dirichlet boundary condition $\Delta_{L,0}$ has no essential spectrum. All its eigenvalues are thus isolated and have finite multiplicity. Denoting them by $\left( \lambda_j \right)$ in ascending order with multiplicity, we have
            
            \[ \begin{array}{lll}
                    N \left( \Delta_{L,0}, \lambda \right) & = & \# \left \{ j, \; \lambda_j \leqslant \lambda \right \}
                \end{array} . \]
        
        \end{corollary}
        
        \begin{proof}
        
            It was made clear in remark \ref{RmkSpectralCountingEssentialSpectrum} that the existence of an essential spectrum is equivalent to the spectral counting function being infinite from a certain point forward. This cannot be by virtue of theorem \ref{ThmWeylLaw}.
        
        \end{proof}
        
        \begin{proposition-definition}
        
            For any real number $\mu \geqslant 0$, the function
            
            \[ \begin{array}{lllll}
                    \zeta_{L, \mu} & : & s & \longmapsto & \sum\limits_j \; \left( \lambda_j + \mu \right)^{-s}
                \end{array} \]
            
            \noindent is well-defined and holomorphic on the half-plane $\Re s > 1$. It is called the spectral zeta function associated to the pseudo-Laplacian with Dirichlet boundary condition.
        
        \end{proposition-definition}
        
        \begin{proof}
            
            We can rephrase theorem \ref{ThmWeylLaw} (\textit{i.e.} the Weyl type law) as
            
            \[ \begin{array}{lllll}
                    \beta_j & = & \sum\limits_{k=1}^j m_k & \leqslant & C \lambda_r
                \end{array} \]
            
            \noindent for any integer $r \in \llbracket \beta_{j-1}+1, \beta_j \rrbracket$, with $\lambda_r$ being constant when $r$ is chosen in this interval of integers. This last inequality in turn yields
            
            \[ \begin{array}{lllll}
                    \frac{1}{\lambda_r} & \leqslant & \frac{C}{\beta_j} & \leqslant & \frac{C}{r}
                \end{array}\]
            
            \noindent for any $r \in \llbracket \beta_{j-1}+1, \beta_j \rrbracket$. We then have, on the half-plane $\Re s > 1$,
            
            \[ \begin{array}{lllll}
                    \left \vert \sum\limits_{j=1}^{+ \infty} \; \frac{1}{\left( \lambda_j + \mu \right)^s} \right \vert & \leqslant & \sum\limits_{j=1}^{+ \infty} \; \frac{1}{\lambda_j^{\Re s}} & = & \sum\limits_{j=1}^{+ \infty} \; \sum\limits_{r = \beta_{j-1}}^{\beta_j} \; \frac{1}{\lambda_r^{\Re s}}
                \end{array} , \]
            
            \noindent with the convention $\beta_0 = 0$. Using corollary \ref{CorPseudoLaplacianHasNoEssentialSpectrum}, we get
            
            \[ \begin{array}{lllll}
                    \sum\limits_{j=1}^{+ \infty} \; \sum\limits_{r = \beta_{j-1}}^{\beta_j} \; \frac{1}{\lambda_r^{\Re s}} & \leqslant & \sum\limits_{j=1}^{+ \infty} \; \sum\limits_{r = \beta_{j-1}}^{\beta_j} \; \frac{C}{r^{\Re s}} & = & C \sum\limits_{j=1}^{+ \infty} \; \frac{1}{j^{\Re s}}
                \end{array} . \]
            
            \noindent This series converges absolutely on the half-plane $\Re s > 1$, thus proving the result.
        
        \end{proof}
        
        Let us now prove the Weyl type law stated in theorem \ref{ThmWeylLaw}.
        
        \begin{remark}
        
            The spectral counting function $N \left( \Delta_{L,0}, \lambda \right)$ is denoted by $N_a \left( \lambda \right)$ when the line bundle $L$ is trivial. The dependence in $a$ is then made more explicit.
        
        \end{remark}

        \noindent The following key lemma allows us to get rid of the line bundle $L$.
        
        \begin{lemma}
        \label{LemmaSpectralCountingLineBundle}
        
            Assume the character $\chi$ has finite order $n$. Then, for any $\lambda > 0$, we have the inequality $N \left( \Delta_{L,0}, \lambda \right) \leqslant N_{a/n} \left( \lambda \right)$.
        
        \end{lemma}
        
        \begin{proof}
        
            First, we note that the result is automatic if the line bundle $L$ is trivial. Let us therefore assume this is not the case, \textit{i.e.} that $\alpha$ is non-zero. The pseudo-Laplacian $\Delta_{L,0}$ is then the Chern Laplacian $\Delta_L$, as specified in definition \ref{DefPseudoLaplacian}. Since $\chi$ having finite order is equivalent to the rationality of $\alpha$, we have $n \alpha \in \mathbb{Z}$. In this proof, we identify any section of $L$ to a function $\psi : \mathbb{R} \times \left] a, + \infty \right[ \longrightarrow \mathbb{C}$ which satisfies $\psi \left( x+1, y \right) = e^{2 i \pi \alpha} \psi \left( x, y \right)$. Such a function is $n$-periodic in the first variable, and the map
            
            \[ \begin{array}{cccllll}
                    f & : & \mathcal{C}^{\infty} \left( S^1 \times \left] a, +\infty \right[, L \right) & \longrightarrow & \multicolumn{3}{l}{\mathcal{C}^{\infty} \left( S^1 \times \left] \frac{a}{n}, +\infty \right[ \right)} \\[0.5em]
                    
                    && \psi & \longmapsto & \left( x, y \right) & \longmapsto & \psi \left( nx, ny \right)
                \end{array} \]
            
            \noindent is an injection from the space of smooth sections of $L$ into the space of smooth functions (\textit{i.e.} smooth sections of the trivial line bundle), which commutes with the Chern Laplacians. For any smooth section $\psi$ of $L$, and any $y > a/n$, we have
            
            \[ \begin{array}{lllllll}
                    \scriptstyle \int_0^1 \; f \left( \psi \right) \left( x, y \right) \; \mathrm{d}x & \scriptstyle = & \scriptstyle \int_0^1 \; \psi \left( nx, ny \right) \; \mathrm{d}x & \scriptstyle = & \scriptstyle \frac{1}{n} \sum\limits_{j=0}^{n-1} \left( e^{2i\pi \alpha} \right)^j \int_0^1 \; \psi \left( x, ny \right) \; \mathrm{d}x & \scriptstyle = & \scriptstyle 0,
                \end{array} \]
            
            \noindent Thus $f$ takes values in the space of smooth sections with vanishing constant Fourier coefficient. Now, we have
            
            \[ \begin{array}{lll}
                    \scriptstyle \int_{a/n}^{+ \infty} \; \int_0^1 \; \left( \Delta f \left( \psi \right) \right) \left( x, y \right) \overline{f \left( \psi \right) \left( x, y \right)} \; \frac{\mathrm{d}x \, \mathrm{d}y}{y^2} & \scriptstyle = & \scriptstyle \int_{a/n}^{+ \infty} \; \int_0^1 \; \left( \Delta_L \psi \right) \left( nx, ny \right) \overline{\psi \left( nx, ny \right)} \; \frac{\mathrm{d}x \, \mathrm{d}y}{y^2} \\[0.8em]
                    
                    & = & \scriptstyle n \int_{a}^{+ \infty} \; \int_0^1 \; \left( \Delta_L \psi \right) \left( x, y \right) \overline{\psi \left( x, y \right)} \; \frac{\mathrm{d}x \, \mathrm{d}y}{y^2}.
                \end{array} \]
            
            \noindent Similarly, we have
            
            \[ \begin{array}{lll}
                    \int_{a/n}^{+ \infty} \; \int_0^1 \; f \left( \psi \right) \left( x, y \right) \overline{f \left( \psi \right) \left( x, y \right)} \; \frac{\mathrm{d}x \, \mathrm{d}y}{y^2} & = & n \int_{a}^{+ \infty} \; \int_0^1 \; \psi \left( x, y \right) \overline{\psi \left( x, y \right)} \; \frac{\mathrm{d}x \, \mathrm{d}y}{y^2}.
                \end{array} \]
            
            \noindent Taking the quotient of these two integrals, assuming $\psi$ is non-zero, we get
            
            \begin{equation}
            \label{EqRayleighQuotient}
                \begin{array}{lll}
                    \frac{\overline{Q_{\Delta}} \left( f \left( \psi \right), f \left( \psi \right) \right)}{\left< f \left( \psi \right), f \left( \psi \right) \right>} & = & \frac{\overline{Q_{\Delta_L}} \left( \psi, \psi \right)}{\left< \psi, \psi \right>}.
                \end{array}
            \end{equation}
            
            \noindent The map $f$ then extends into
            
            \[ \begin{array}{lllll}
                    f & : & L^2_{1,\mathrm{Dir}} \left( \mathbb{R} \times \left] a, + \infty \right[, L \right) & \longrightarrow & L^2_{1,\mathrm{Dir}} \left( \mathbb{R} \times \left] a, + \infty \right[ \right)_0
                \end{array} , \]
            
            \noindent Furthermore, equality \eqref{EqRayleighQuotient} remains valid for this continuation. This allows us to compare the spectral quantities $\mu_n$ and $\mu_n \left( \Delta_L \right)$, as we have
            
            \[ \begin{array}{lll}
                    \multicolumn{3}{l}{\inf\limits_{f \left( \psi_1 \right), \dots, f \left( \psi_k \right)} \sup \left \{ \frac{\overline{Q_{\Delta}} \left( f \left( \psi \right), f \left( \psi \right) \right)}{\left< f \left( \psi \right), f \left( \psi \right) \right>}, \; f \left( \psi \right) \in \mathrm{span} \left( f \left( \psi_1 \right), \dots, f \left( \psi_k \right) \right), \; f \left( \psi \right) \neq 0 \right \}} \\[1em]
                    
                    \qquad \qquad \qquad \qquad & \leqslant & \inf\limits_{\psi_1, \dots, \psi_k} \sup \left \{ \frac{\overline{Q_{\Delta_L}} \left( \psi, \psi \right)}{\left< \psi, \psi \right>}, \; \psi \in \mathrm{span} \left( \psi_1, \dots, \psi_k \right), \; \psi \neq 0 \right \}.
                \end{array} \]
            
            \noindent Replacing the lower bound on the left-hand side by one with respect to $k$ elements in the $L_1^2$ space, instead of just the image of $f$, we get
            
            \[ \begin{array}{lll}
                    \mu_k \left( a/n \right) & \leqslant & \mu_k \left( L \right)
                \end{array} . \]
            
            \noindent We can now use the definition of the spectral counting functions, and get the result.
        
        \end{proof}
        
        The proof of theorem \ref{ThmWeylLaw} is thus reduced to the case of the trivial line bundle. We will now adapt some arguments used by Colin de Verdière in \cite[Sec. 4]{MR699488}. This will require pseudo-Laplacians with Dirichlet and Neumann conditions on the cuspidal end, or on ``steps'' within it. Let us first see what all of this means, using a language close to \cite[Sec. I.5]{MR1102675}.
        
        \begin{definition}
        
            Let $\Lambda$ be a union of open intervals contained in $\left] a, + \infty \right[$. The \textit{pseudo-Laplacian with Dirichlet boundary conditions} $\Delta^D_{\Lambda}$ is defined as the Friedrichs extension of the Laplacian associated to the quadratic form
            
            \[ \begin{array}{lllll}
                    Q_D & : & \mathcal{C}_0^{\infty} \left( S^1 \times \Lambda \right)_0 \times \mathcal{C}_0^{\infty} \left( S^1 \times \Lambda \right)_0 & \longrightarrow & \mathbb{C} \\[0.5em]
                    
                    && \left( u,v \right) & \longmapsto & \int_{S^1 \times \Lambda} \; \nabla u \wedge \overline{\nabla v}
                \end{array} \]
                
            \noindent on smooth compactly supported functions with vanishing constant Fourier coefficient. Here $\nabla$ is the gradient, \textit{i.e.} the Chern connection for the trivial line bundle.
        
        \end{definition}
        
        \begin{remark}
        
            The domain of this Friedrichs extension and of the associated quadratic form are given by
            
            \[ \begin{array}{lllclll}
                    \mathrm{Dom} \, \Delta^D_{\Lambda} & = & L_{2,\mathrm{Dir}}^2 \left( S^1 \times \Lambda \right)_0 & \quad \text{and} \quad & \mathrm{Dom} \, \overline{Q_D} & = & L_{1,\mathrm{Dir}}^2 \left( S^1 \times \Lambda \right)_0.
                \end{array} \]
        
        \end{remark}
        
        \begin{definition}
        
            Let $\Lambda$ be a union of open intervals contained in $\left] a, + \infty \right[$. The \textit{pseudo-Laplacian with Neumann boundary conditions} $\Delta^N_{\Lambda}$ is defined as the Laplacian associated to the closed quadratic form
            
            \[ \begin{array}{lllll}
                    Q_N & : & L_1^2 \left( S^1 \times \Lambda \right)_0 \times L_1^2 \left( S^1 \times \Lambda \right)_0 & \longrightarrow & \mathbb{C} \\[0.5em]
                    
                    && \left( u,v \right) & \longmapsto & \int_{S^1 \times \Lambda} \; \nabla u \wedge \overline{\nabla v}
                \end{array} \]
                
            \noindent on the Sobolev space $L_1^2$ with vanishing constant Fourier coefficient. Here again $\nabla$ is the gradient, \textit{i.e.} the Chern connection for the trivial line bundle.
        
        \end{definition}
        
        \begin{remark}
        
            The quadratic forms associated to the Laplacians with Dirichlet or Neumann boundary conditions being the same, comparing them in the sense of the order $\preccurlyeq$ from definition \ref{DefOrderSelfAdjoint} is just a matter of inclusion of domains.
        
        \end{remark}
        
        \begin{lemma}
        \label{LemDNBracketing}
        
            Let $\Lambda$ be an open interval in $\left] a, + \infty \right[$. We have $\Delta^N_{\Lambda} \preccurlyeq \Delta^D_{\Lambda}$. Consequently, we have $N \left( \Delta^D_{\Lambda}, \lambda \right) \leqslant N \left( \Delta^N_{\Lambda}, \lambda \right)$ for any $\lambda > 0$.
        
        \end{lemma}
        
        \begin{proof}
        
            This proposition stems directly from the inclusion of Sobolev spaces
            
            \[ \begin{array}{lll}
                    L_{1,\mathrm{Dir}}^2 \left( S^1 \times \Lambda \right)_0 & \subset & L_1^2 \left( S^1 \times \Lambda \right)_0
                \end{array} . \]
            
            \noindent The second part of the result is a consequence of proposition \ref{PropSpectralCountingOrderRelation}.
        
        \end{proof}
        
        \begin{lemma}
        \label{LemNeumannBreakingInterval}
        
            Let $\Lambda_1$ and $\Lambda_2$ be two open intervals included in $\left] a, +\infty \right[$. We have
            
            \[ \begin{array}{lll}
                    \Delta^N_{\Lambda_1} \oplus \Delta^N_{\Lambda_2} & \preccurlyeq & \Delta^N_{\Lambda}
                \end{array} , \]
            
            \noindent where $\Lambda$ is given by $\Lambda = \mathring{\overline{\Lambda_1 \cup \Lambda_2}}$, and the operator on the left-hand side acts on the direct sum of the relevant domains. Consequently, for any $\lambda > 0$, we have
            
            \[ \begin{array}{lll}
                    N \left( \Delta^N_{\Lambda}, \lambda \right) & \leqslant & N \left( \Delta^N_{\Lambda_1}, \lambda \right) + N \left( \Delta^N_{\Lambda_2}, \lambda \right)
                \end{array} . \]
        
        \end{lemma}
        
        \begin{proof}
        
            The comparison between the two operators is a consequence of the inclusion
            
            \[ \begin{array}{lll}
                    L_1^2 \left( S^1 \times \Lambda \right)_0 & \hookrightarrow & L_1^2 \left( S^1 \times \Lambda_1 \right)_0 \oplus L_1^2 \left( S^1 \times \Lambda_2 \right)_0
                \end{array} \]
            
            \noindent given by the restriction to each $S^1 \times \Lambda_j$. Propositions \ref{PropSpectralCountingOrderRelation} and \ref{PropSpectralCountingDirectSum} give the rest.
        
        \end{proof}
        
        \begin{remark}
        
            There is a similar comparison for Dirichlet boundary conditions
            
            \[ \begin{array}{lll}
                    \Delta^D_{\Lambda} & \preccurlyeq & \Delta^D_{\Lambda_1} \oplus \Delta^D_{\Lambda_2}
                \end{array} , \]
            
            \noindent with $\Lambda_1$ and $\Lambda_2$ two pairwise disjoint union of open intervals, and $\Lambda = \mathring{\overline{\Lambda_1 \cup \Lambda_2}}$. This comes from the fact that $L_1^2$ functions on $S^1 \times \Lambda_1$ and $S^1 \times \Lambda_2$ with Dirichlet boundary conditions can be glued into an $L_1^2$ function on $S^1 \times \Lambda$.
        
        \end{remark}
        
        We now need some asymptotic control as $a$ goes to infinity over the first eigenvalue of the pseudo-Laplacian with Neumann boundary condition on $S^1 \times \left] a, +\infty \right[$. This is done by first obtaining an \textit{ad-hoc} version of the \textit{Poincaré inequality}, for functions in the Sobolev space
        
        \[ \begin{array}{lll}
                L_1^2 \left( S^1 \times \left] a, + \infty \right[, \mathrm{d}x^2 + \mathrm{d}y^2 \right)_0
            \end{array} \]
        
        \noindent for the Euclidean metric, with vanishing constant Fourier coefficient.
        
        \begin{lemma}
        \label{LemmaPoincaréInequality}
        
            For any function $\psi \in L_1^2 \left( S^1 \times \left] a, + \infty \right[, \mathrm{d}x^2 + \mathrm{d}y^2 \right)_0$, we have
            
             \[ \begin{array}{lll}
                    \int_{S^1} \int_{a}^{+ \infty} \; \left \Vert \nabla \psi \right \Vert^2 \; \mathrm{d}y \; \mathrm{d}x & \geqslant & K \int_{S^1} \int_a^{+ \infty} \; \left \vert \psi \right \vert^2 \; \mathrm{d}y \; \mathrm{d}x
                \end{array} , \]
            
            \noindent where $\nabla$ stands for the usual gradient.
        
        \end{lemma}
        
        \begin{proof}
        
            Let us assume, by contradiction, that there exists a sequence $\left( \psi_n \right)$ of functions in the Sobolev space above, such that we have, for every integer $n > 0$,
            
            \[ \begin{array}{lll}
                    \int_{S^1} \int_{a}^{+ \infty} \; \left \Vert \nabla \psi_n \right \Vert^2 \; \mathrm{d}y \; \mathrm{d}x & < & 2^{-n} \int_{S^1} \int_a^{+ \infty} \; \left \vert \psi_n \right \vert^2 \; \mathrm{d}y \; \mathrm{d}x.
                \end{array} \]
            
            \noindent Up to multiplying all these functions by some constants, we assume that we have
            
            \begin{equation}
            \label{EqNorm1Poincaré}
                \begin{array}{lll}
                    \int_{S^1} \int_a^{+ \infty} \; \left \vert \psi_n \right \vert^2 \; \mathrm{d}y \; \mathrm{d}x & = & 1.
                \end{array}
            \end{equation}
            
            \noindent In particular, we have
            
            \begin{equation}
            \label{EqConvGradientInPoincaré}
                \begin{array}{lll}
                    \int_{S^1} \int_{a}^{+ \infty} \; \left \Vert \nabla \psi_n \right \Vert^2 \; \mathrm{d}y \; \mathrm{d}x & \underset{n \rightarrow + \infty}{\longrightarrow} & 0.
                \end{array}
            \end{equation}
            
            \noindent Using the Banach--Alaoglu theorem (see for instance \cite[Thm 1.3.17]{MR1157815}, \cite[Sec. V.3]{MR1070713}), and up to taking a subsequence, the sequence $\left( \psi_n \right)$ converges weakly to an element
            
            \[ \begin{array}{lll}
                    \psi & \in & L_1^2 \left( S^1 \times \left] a, + \infty \right[, \mathrm{d}x^2 + \mathrm{d}y^2 \right)_0
                \end{array} . \]
            
            \noindent Note that we can identify weak and weak$^{\ast}$ convergences here, because we are dealing with a Hilbert space. Using \eqref{EqConvGradientInPoincaré}, we get $\nabla \psi = 0$, so $\psi$ is constant. This can be proved by writing the Fourier decomposition of $\psi$ and noting that a distribution with vanishing derivative in dimension $1$ is constant. Because $\psi$ has vanishing constant Fourier coefficient, it vanishes identically. This is absurd by \eqref{EqNorm1Poincaré}.
        
        \end{proof}
        
        \begin{lemma}
        \label{LemmaSpectralCountingSchrinkingCusp}
        
            There exists a constant $K > 0$ such that we have 
            
            \[ \begin{array}{lll}
                    \mu_1 \left( \Delta^N_{\left] a, +\infty \right[} \right) & \geqslant & a^2 K
                \end{array} . \]
            
            \noindent As a consequence, for any fixed real number $\lambda > 0$, we have
            
            \[ \begin{array}{llll}
                   N \left( \Delta^N_{\left] a, +\infty \right[}, \lambda \right) & = & 0 & \text{for every $a$ large enough}
                \end{array} . \]
        
        \end{lemma}
        
        \begin{proof}
        
            To simplify, denote by $D_a$ the domain of $\Delta^N_{\left] a, +\infty \right[}$. We have
            
            \[ \begin{array}{lll}
                    \mu_1 \left( \Delta^N_{\left] a, +\infty \right[} \right) & = & \inf\limits_{\psi \in D_a} \; \frac{\int_{S^1} \int_{a_{\phantom{a}}}^{+ \infty} \; \left \Vert \nabla \psi \right \Vert^2 \; \mathrm{d}y \; \mathrm{d}x}{\int_{S^1} \int_a^{+ \infty^{\phantom{a}}} \; \left \vert \psi \right \vert^2 \; \frac{\mathrm{d}y \; \mathrm{d}x}{y^2}}
                \end{array} . \]
            
            \noindent For any smooth function $\psi$ on $S^1 \times \left] a, + \infty \right[$ whose constant Fourier coefficient is zero, and which vanishes for $y$ large enough, we have, using lemma \ref{LemmaPoincaréInequality},
            
            \[ \begin{array}{lllll}
                    \frac{\int_{S^1} \int_{a_{\phantom{a}}}^{+ \infty} \; \left \Vert \nabla \psi \right \Vert^2 \; \mathrm{d}y \; \mathrm{d}x}{\int_{S^1} \int_a^{+ \infty^{\phantom{a}}} \; \left \vert \psi \right \vert^2 \; \frac{\mathrm{d}y \; \mathrm{d}x}{y^2}} & \geqslant & \alpha^2 \frac{\int_{S^1} \int_{a_{\phantom{a}}}^{+ \infty} \; \left \Vert \nabla \psi \right \Vert^2 \; \mathrm{d}y \; \mathrm{d}x}{\int_{S^1} \int_a^{+ \infty^{\phantom{a}}} \; \left \vert \psi \right \vert^2 \; \mathrm{d}y \; \mathrm{d}x} & \geqslant & a^2 K.
                \end{array} \]
            
            \noindent The space of such functions being dense in $D_a$, we get the first part of the lemma, and the rest follows immediately, as we then have 
            
            \[ \begin{array}{lll}
                    \mu_1 \left( \Delta^N_{\left] a, +\infty \right[} \right) & \underset{a \rightarrow + \infty}{\longrightarrow} & + \infty
                \end{array} . \]
        
        \end{proof}
        
        We need a few more considerations before we can move on to the Weyl type law.
        
        \begin{definition}
        
            Let $\Lambda$ be an open interval included in $\left] a, +\infty \right[$. The \textit{Euclidean pseudo-Laplacian with Neumann boundary conditions} $H^N_{\Lambda}$ is defined as the Laplacian associated to the closed quadratic form
            
            \[ \begin{array}{lllll}
                    Q_N^{\mathrm{eucl}} & : & \left( L_1^2 \left( S^1 \times \Lambda, \mathrm{d}x^2 + \mathrm{d}y^2 \right)_0 \right)^2 & \longrightarrow & \mathbb{C} \\[0.5em]
                    
                    && \left( u,v \right) & \longmapsto & \int_{S^1 \times \Lambda} \; \nabla u \wedge \overline{\nabla v}
                \end{array} \]
                
            \noindent on the Sobolev space $L_1^2$ with vanishing constant Fourier coefficient, for the euclidean metric. Here $\nabla$ is the usual euclidean gradient.
        
        \end{definition}
        
        \begin{lemma}
        \label{LemmaNeumannLaplacianPoincaréEuclidean}
        
            For any real number $\lambda > 0$, and any $0 < a < b$, we have
            
            \[ \begin{array}{lll}
                    N \left( \Delta^N_{\left] a,b \right[}, \lambda \right) & \leqslant &  N \left( H^N_{\left] a,b \right[}, \frac{\lambda}{a^2} \right)
                \end{array} \]
        
        \end{lemma}
        
        \begin{proof}
        
            The argument used here is in \cite[Lem. 4.2]{MR699488}, and relies on the inequality
            
            \[ \begin{array}{lll}
                    \frac{\int_{S^1} \int_{a_{\phantom{a}}}^{b} \; \left \Vert \nabla \psi \right \Vert^2 \; \mathrm{d}y \; \mathrm{d}x}{\int_{S^1} \int_a^{b^{\phantom{a}}} \; \left \vert \psi \right \vert^2 \; \frac{\mathrm{d}y \; \mathrm{d}x}{y^2}} & \geqslant & a^2 \frac{\int_{S^1} \int_{a_{\phantom{a}}}^{b} \; \left \Vert \nabla \psi \right \Vert^2 \; \mathrm{d}y \; \mathrm{d}x}{\int_{S^1} \int_a^{b^{\phantom{a}}} \; \left \vert \psi \right \vert^2 \; \mathrm{d}y \; \mathrm{d}x}
                \end{array} . \]
        
        \end{proof}
        
        \begin{lemma}
        \label{LemmaSpectralCountingNeumannEuclidean}
        
            For any real number $\lambda > 0$, and any $0 < a < b$, we have
            
            \[ \begin{array}{llll}
                    N \left( H^N_{\left] a,b \right[}, \lambda \right) & \leqslant &  \frac{b-a}{4 \pi} \lambda + \frac{1}{\pi} \sqrt{\lambda} & \text{if } \lambda \geqslant 4 \pi^2 \\[1em]
                    
                    & = & 0 & \text{otherwise}
                \end{array} . \]
        
        \end{lemma}
        
        \begin{proof}
        
            The proof of this result amounts to an explicit computation of the eigenvalues, and can be found in \cite[Lem 4.1]{MR699488}.
        
        \end{proof}
        
        This was the last ingredient needed to prove the main theorem of this section.
        
        \begin{proof}[Proof of theorem \ref{ThmWeylLaw}]
        
            Consider a real number $\lambda > 0$. Using lemma \ref{LemmaSpectralCountingLineBundle}, we have
            
            \[ \begin{array}{lll}
                    N \left( \Delta_{L,0}, \lambda \right) & \leqslant & N_{a/n} \left( \lambda \right)
                \end{array} . \]
                
            \noindent Extending by $0$ gives an injection

            \[ \begin{array}{lll}
                    L_{1, \mathrm{Dir}}^2 \left( S^1 \times \left] a, +\infty \right[ \right)_0 & \hookrightarrow & L_{1, \mathrm{Dir}}^2 \left( S^1 \times \left] \frac{a}{n}, +\infty \right[ \right)_0,
                \end{array} \]
            
            \noindent which yields $N_{a/n} \left( \lambda \right) \leqslant N_a \left( \lambda \right)$, using proposition \ref{PropSpectralCountingOrderRelation}. We further have
            
            \[ \begin{array}{lllll}
                    N_a \left( \lambda \right) & = & N \left( \Delta^D_{\left] a, + \infty \right[}, \lambda \right) & \leqslant & N \left( \Delta^N_{\left] a, + \infty \right[}, \lambda \right)
                \end{array} , \]
            
            \noindent since operators for the Dirichlet boundary condition and the Neumann one can be compared using lemma \ref{LemDNBracketing}. The idea presented in \cite[Thm 6]{MR699488} by Colin de Verdière can be used to break apart the interval $\left] a, +\infty \right[$, which yields, for a fixed $\delta > 0$,
            
            \[ \begin{array}{lll}
                    N \left( \Delta^N_{\left] a, + \infty \right[}, \lambda \right) & \leqslant & \sum\limits_{k=0}^{\ell - 1} \; N \left( \Delta^N_{\left] a+k \delta, a+\left( k+1 \right) \delta \right[}, \lambda \right) + N \left( \Delta^N_{\left] a + \ell, + \infty \right[}, \lambda \right)
                \end{array} \]
            
            \noindent after having applied lemma \ref{LemNeumannBreakingInterval} inductively. We now note that lemma \ref{LemmaSpectralCountingSchrinkingCusp} gives us an integer $\ell_a$, depending on $a$ in an unknown manner, such that we have
            
            \[ \begin{array}{lll}
                     N \left( \Delta^N_{\left] a + \ell, + \infty \right[}, \lambda \right) & = & 0
                \end{array} \]
            
            \noindent for any integer $\ell \geqslant \ell_a$. For such integers, we have
            
            \[ \begin{array}{llllll}
                    N \left( \Delta^N_{\left] a, + \infty \right[}, \lambda \right) & \leqslant & \sum\limits_{k=0}^{\ell - 1} \; N \left( \Delta^N_{\left] a+k \delta, a+\left( k+1 \right) \delta \right[}, \lambda \right) \\[1em]
                    
                    & \leqslant & \sum\limits_{k=0}^{+ \infty} \; N \left( \Delta^N_{\left] a+k \delta, a+\left( k+1 \right) \delta \right[}, \lambda \right) \\[1em]
                    
                    & \leqslant & \sum\limits_{k=0}^{+ \infty} \; N \left( \Delta^N_{\left] a+k \delta, a+\left( k+1 \right) \delta \right[}, \frac{\lambda}{\left( a + k \delta \right)^2} \right) & \text{by lemma \ref{LemmaNeumannLaplacianPoincaréEuclidean}}.
                \end{array} \]
            
            \noindent The inequality formed by the first and last term above does not depend on $\ell$, and thus removes the problem of not understanding how $\ell_a$ varies with respect to $a$. The added terms do not play much of a role, since lemma \ref{LemmaSpectralCountingNeumannEuclidean} says we have
            
            \[ \begin{array}{lll}
                    N \left( \Delta^N_{\left] a+k \delta, a+\left( k+1 \right) \delta \right[}, \frac{\lambda}{\left( a + k \delta \right)^2} \right) & = & 0
                \end{array} \]
            
            \noindent if the inequality $\lambda < 4 \pi^2 \left( a + k \delta \right)^2$ is satisfied. Let us set
            
            \[ \begin{array}{lll}
                    k_a & = & \left \lfloor \frac{1}{\delta} \left[ \frac{1}{2\pi} \sqrt{\lambda} - a \right] \right \rfloor
                \end{array} , \]
            
            \noindent where $\left \lfloor x \right \rfloor$ denotes the largest integer smaller than $x$. We get
            
            \[ \begin{array}{lll}
                    N \left( \Delta^N_{\left] a, + \infty \right[}, \lambda \right) & \leqslant & \sum\limits_{k=0}^{k_a} \; N \left( \Delta^N_{\left] a+k \delta, a+\left( k+1 \right) \delta \right[}, \frac{\lambda}{\left( a + k \delta \right)^2} \right)
                \end{array} . \]
            
            \noindent We will now evaluate this term, using the method described in \cite[Thm 6]{MR699488}. Each term summed on the right-hand side above being evaluated in lemma \ref{LemmaSpectralCountingNeumannEuclidean}, we have
            
            \[ \begin{array}{lll}
                    N \left( \Delta^N_{\left] a, + \infty \right[}, \lambda \right) & \leqslant & \frac{\delta \lambda}{4 \pi} \sum\limits_{k=0}^{k_a} \; \frac{1}{\left( a + k \delta \right)^2} + \frac{1}{\pi} \sqrt{\lambda} \sum\limits_{k=0}^{k_a} \; \frac{1}{a + k \delta}
                \end{array} . \]
            
            \noindent Comparing series and integrals, we get
            
            \[ \begin{array}{lllll}
                    \frac{\delta \lambda}{4 \pi} \sum\limits_{k=0}^{k_a} \; \frac{1}{\left( a + k \delta \right)^2} & = & \multicolumn{3}{l}{\frac{\delta \lambda}{4 \pi a^2} + \frac{\lambda}{4 \pi} \sum\limits_{k=1}^{k_a} \; \int_{a+\left( k - 1 \right) \delta}^{a + k \delta} \; \frac{1}{\left( a + k \delta \right)^2} \; \mathrm{d}y} \\[1em]
                    
                    & \leqslant & \frac{\delta \lambda}{4 \pi a^2} + \frac{\lambda}{4 \pi} \sum\limits_{k=1}^{+ \infty} \; \int_{a+\left( k - 1 \right) \delta}^{a + k \delta} \; \frac{1}{y^2} \; \mathrm{d}y & = & \frac{\delta \lambda}{4 \pi a^2} + \frac{\lambda}{4 \pi a}
                \end{array} \]
            
            \noindent for the first term, and, similarly
            
            \[ \begin{array}{lll}
                    \frac{1}{\pi} \sqrt{\lambda} \sum\limits_{k=0}^{k_a} \; \frac{1}{a + k \delta} & = & \frac{1}{\pi} \sqrt{\lambda} \left[ \frac{1}{a} + \frac{1}{\delta} \sum\limits_{k=1}^{k_a} \; \int_{a + \left( k-1 \right) \delta}^{a + k \delta} \; \frac{1}{a + k \delta} \mathrm{d}y \right] \\[1em]
                    
                    & \leqslant & \frac{1}{\pi} \sqrt{\lambda} \left[ \frac{1}{a} + \frac{1}{\delta} \log \left( 1 + \frac{k_a \delta}{a} \right) \right] \\[1em]
                    
                    & \leqslant & \frac{1}{\pi a} \sqrt{\lambda} + \frac{1}{2 \pi \delta} \sqrt{\lambda} \log \lambda.
                \end{array} \]
            
            \noindent Combining these results, we get
            
            \begin{equation}
            \label{IneqWeylLaw}
                \begin{array}{lll}
                    N \left( \Delta^N_{\left] a, + \infty \right[}, \lambda \right) & \leqslant & \frac{1}{\pi a} \sqrt{\lambda} + \frac{1}{4 \pi a} \left[ 1 + \frac{\delta}{a} \right] \lambda + \frac{1}{2 \pi \delta} \sqrt{\lambda} \log \lambda
                \end{array} .
            \end{equation}
            
            \noindent All these estimates put together yield
            
            \[ \begin{array}{lll}
                    N \left( \Delta_{L,0}, \lambda \right) & \leqslant & \frac{1}{4 \pi a} \left[ 5 + \frac{\delta}{a} \right] \lambda + \frac{1}{2 \pi \delta} \sqrt{\lambda} \log \lambda
                \end{array} \]
            
            \noindent for any $\lambda > 1$ by using the inequality $\sqrt{\lambda} < \lambda$. The theorem follows on $\left[ 1, + \infty \right[$ by using the fact that the function $\lambda \mapsto \lambda^{-1/2} \log \lambda$ is bounded on $\left[ 1, + \infty \right[$, and then on $\mathbb{R}_+^{\ast}$ by further adjusting the constant.
            
        \end{proof}

    \subsection{Spectral problem}
    \label{SubSecSpectralProblem}
    
        Everything has now been set up to study the eigenvalues of the pseudo-Laplacian, and compute asymptotics related to the determinant of the pseudo-Laplacian. We consider the representative of $\alpha$ modulo $1$ in $\left[ 0,1 \right[$, and denote it the same way. The spectral problem we want to solve is:
        
        \[ \left \{ \begin{array}{llll}
                - y^2 \left( \frac{\partial^2}{\partial x^2} + \frac{\partial^2}{\partial y^2} \right) \psi \left( x, y \right) & = & \multicolumn{2}{l}{\lambda \psi \left( x, y \right)} \\[0.7em]
                
                \psi \left( x+1, y \right) & = & \multicolumn{2}{l}{e^{2 i \pi \alpha} \psi \left( x, y \right)} \\[0.7em]
                
                \int_{S^1 \times \left] a, +\infty \right[} \; \left \vert \psi \right \vert^2 & < & + \infty & \text{(Integrability condition)} \\[0.7em]
                
                \psi \left( x, a \right) & = & 0 & \text{(Dirichlet boundary condition)} \\[0.7em]
                
                \int_{S^1} \; \psi \left( x, y \right) \; \mathrm{d}x & = & 0 &  \text{for almost all } y > a \text{ if } \alpha = 0
            \end{array} \right. \]
        
        The last condition above is equivalent to a vanishing constant Fourier coefficient, only to be considered when $\chi$ is trivial, or equivalently when we have $\alpha = 0$. Setting
        
        \[ \begin{array}{lll}
                \varphi \left( x, y \right) & = & e^{- 2 i \pi \alpha x} \psi \left( x, y \right)
            \end{array} , \]
        
        \noindent the spectral problem written above becomes
        
        \[ \left \{ \begin{array}{llll}
                - y^2 \left( \frac{\partial^2}{\partial x^2} + \frac{\partial^2}{\partial y^2} \right) \varphi \left( x, y \right) & = & \multicolumn{2}{l}{\left( \lambda - 4 \pi^2 \alpha^2 y^2 \right) \varphi \left( x, y \right) + 4 i \pi \alpha y^2 \frac{\partial \varphi}{\partial x}} \\[0.7em]
                
                \varphi \left( x+1, y \right) & = & \multicolumn{2}{l}{\varphi \left( x, y \right)} \\[0.7em]
                
                \int_{S^1 \times \left] a, +\infty \right[} \; \left \vert \varphi \right \vert^2 & < & \multicolumn{2}{l}{+ \infty} \\[0.7em]
                
                \varphi \left( x, a \right) & = & \multicolumn{2}{l}{0} \\[0.7em]
                
                \int_{S^1} \; \varphi \left( x, y \right) \; \mathrm{d}x & = & 0 &  \text{for almost all } y > a \text{ if } \alpha = 0
            \end{array} \right. \]
            
        The Laplacian being an elliptic operator, solutions to either problems are smooth. Furthermore, the second formulation of the spectral problem is easier to work with, as solutions are periodic in the first variable, and can thus be written as a sum of their Fourier series. We write such a function as
        
        \[ \begin{array}{lll}
                \varphi \left( x, y \right) & = & \sum\limits_{k \in \mathbb{Z}} a_k \left( y \right) e^{2 i \pi k x}
            \end{array} . \]
        
        The unicity of Fourier coefficients then implies that the partial differential equation on $\varphi$ is equivalent to the ordinary differential equations
        
        \[ \begin{array}{lll}
                \left[ y^2 \frac{\mathrm{d}^2}{\mathrm{d}y^2} + \lambda - 4 \pi^2 y^2 \left( k + \alpha \right)^2 \right] a_k \left( y \right) & = & 0
            \end{array} \]
        
        \noindent for every relative integer $k$, with the exception of $k = 0$ if $\alpha$ vanishes. The solution to the problem above is given, up to multiplicative constant, by
        
        \[ \begin{array}{lll}
                a_k \left( y \right) & = & \sqrt{y} K_{s-1/2} \left( 2 \pi \left \vert k + \alpha \right \vert y \right)
            \end{array} , \]
        
        \noindent where the possible values of $\lambda = s \left( 1 - s \right)$ are determined by the Dirichlet boundary condition $K_{s-1/2} \left( 2 \pi \left \vert k + \alpha \right \vert a \right) = 0$, and $K$ is a modified Bessel function of the second kind. More information on those can be found in appendix \ref{AppModifiedBessel}.

    \subsection{Localization of the eigenvalues}
    \label{SubSecLocalizationEigenvalues}
    
        As we have seen in the last paragraph, the spectral problem we consider can always be solved, if we leave aside the Dirichlet boundary condition. We will now get more information as to when a solution satisfying the boundary condition exists. The goal is to know when the function
        
        \[ \begin{array}{lll}
                s & \longmapsto & K_{s-1/2} \left( 2 \pi \left \vert k + \alpha \right \vert a \right)
            \end{array} \]
        
        \noindent vanishes. This is obtained by adaptating an argument from \cite[Appendix A]{MR2439244}, which is developped there by Saharian for Legendre functions.

        \begin{proposition}
        \label{PropZerosModifiedBessel}
        
            For any $k \in \mathbb{Z}$, except $k = 0$ if $\alpha$ vanishes, the function
            
            \[ \begin{array}{lll}
                    s & \longmapsto & K_{s-1/2} \left( 2 \pi \left \vert k + \alpha \right \vert a \right)
                \end{array} \]
            
            \noindent is holomorphic. Its zeros, which are all simple, are given by a discrete set
            
            \[ \begin{array}{lll}
                    \left \{ \frac{1}{2} + i r_{k,j} \right \} & \subset & \left \{ \frac{1}{2} + ir, r \in \mathbb{R}^{\ast} \right \}
                \end{array} . \]
            
            \noindent The eigenvalues corresponding to the spectral problem are given by $\lambda_{k,j} = 1/4 + r_{k,j}^2$.
        
        \end{proposition}
        
        \begin{proof}
        
            This proposition is a direct consequence of proposition \ref{PropAppZerosModifiedBessel}.
        
        \end{proof}
        
        \begin{remark}
        
            It should be noted, when comparing with section \ref{SubSecWeylLaw}, that the eigenvalues of $\Delta_{L,0}$, which were denoted by $\lambda_j$ in ascending order, have been reindexed as~$\lambda_{k,j}$ in this last proposition. The definition of the spectral zeta function should be adapted to reflect that, but the function does not change on the half-plane $\Re s > 1$.
        
        \end{remark}

\section{Determinant of the pseudo-Laplacian with Dirichlet boundary condition}
\label{SecDetPseudoLaplacian}

    Following the classical theory of zeta-regularized determinants, the idea is to set
    
    \[ \begin{array}{lll}
            \log \det \left( \Delta_{L,0} + \mu \right) & = & - \zeta_{L, \mu} ' \left( 0 \right)
        \end{array} . \]
    
    \noindent The spectral zeta function $\zeta_{L,\mu}$ being \textit{a priori} only defined and holomorphic on the half-plane $\Re s > 1$, one must show that it has a holomorphic continuation to some region containing the origin. Doing so is one of the purposes of this section, and the other two are obtaining:
    
    \begin{enumerate}
        \item an asymptotic expansion of $\log \det \left( \Delta_{L,0} + \mu \right)$ as $\mu$ goes to infinity;
        
        \item an asymptotic expansion of $\log \det \Delta_{L,0}$ as $a$ goes to infinity.
    \end{enumerate}
    
    \noindent This requires technical computations, based on the ideas developped by Freixas i Montplet and von Pippich in \cite{MR4167014} for the case of the trivial line bundle. The reader is referred to the introduction for an overview of the methods, and some comments.

    \subsection{Integral representation of the spectral zeta function}
    \label{SubSecIntegralRepresentationSpectralZeta}
    
        The first step in studying the zeta function $\zeta_{L,0}$ is to write it as some integral. The main ingredient for that is the \textit{argument principle}, for which the reader is referred to \cite[Sec. 3.4]{MR1976398}.
        
        \begin{definition}
        
            For any $\vartheta \in \left] 0, \pi/2 \right[$, the contour $\gamma_{\vartheta}$ is defined by
            
            \[ \begin{array}{lll}
                    \gamma_{\vartheta} & = & \left \{ r e^{i \vartheta}, r \geqslant 0 \right \} \cup \left \{ r e^{- i \vartheta}, r \geqslant 0 \right \}
                \end{array} . \]
        
        \end{definition}
        
        \begin{remark}
        
            This contour and its rotated counterpart, which are represented below with their orientation, will serve as contour integrations in the argument principle.
            
            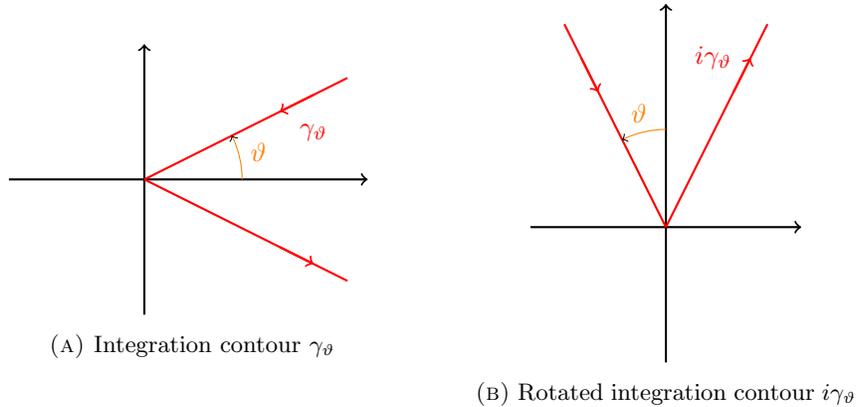
\begin{figure}[H]
    
                \begin{center}
        
                    \begin{subfigure}{0.4\textwidth}
                
                        \centering
                
                        \begin{tikzpicture}[scale=0.9]
        
                            \draw[thick] (-2,0) -- (0,0);
                            \draw[thick] (0,-2) -- (0,0);
                            \draw[thick,->] (0,0) coordinate (2) node {} -- (3.3,0) coordinate (1) node {};
                            \draw[thick, ->] (0,0) -- (0,2);
            
                            \draw[thick, red] (0,0) -- (3,3/2) coordinate (3) node {};
                            \draw[thick, red] (0,0) -- (3,-3/2) coordinate (4) node {};
            
                            \draw pic["$\color{orange} \vartheta$", draw=orange, ->, angle eccentricity=1.2, angle radius=1.3cm] {angle=1--2--3};
            
                            \node (P) at (2.5,0.7) {$\color{red} \gamma_{\vartheta}$};
            
                            \draw[thick, red, ->] (2.5,1.25) -- (2,1);
                            \draw[thick, red, ->] (2,-1) -- (2.5,-1.25);
                    
                        \end{tikzpicture}
                    
                        \caption{Integration contour $\gamma_{\vartheta}$}
                        \label{FigIntegrationContour}
        
                        \end{subfigure}
                        \hspace{1cm}
                        \begin{subfigure}{0.4\textwidth}
            
                        \centering
        
                        \begin{tikzpicture}[scale=0.9]
        
                            \draw[thick] (-2,0) -- (0,0);
                            \draw[thick] (0,-2) -- (0,0);
                            \draw[thick,->] (0,0) coordinate (2) node {} -- (0,3.3) coordinate (1) node {};
                            \draw[thick, ->] (0,0) -- (2,0);
            
                            \draw[thick, red] (0,0) -- (3/2,3) coordinate (4) node {};
                            \draw[thick, red] (0,0) -- (-3/2,3) coordinate (3) node {};
            
                            \draw pic["$\color{orange} \vartheta$", draw=orange, ->, angle eccentricity=1.2, angle radius=1.3cm] {angle=1--2--3};
            
                            \node (P) at (0.7,2.5) {$\color{red} i \gamma_{\vartheta}$};
            
                            \draw[thick, red, ->] (1,2) -- (1.25,2.5);
                            \draw[thick, red, ->] (-1.25,2.5) -- (-1,2);
                    
                        \end{tikzpicture}
                    
                        \caption{Rotated integration contour $i \gamma_{\vartheta}$}
                        \label{FigRotatedIntegrationContour}
        
                    \end{subfigure}
            
                \end{center}
        
                \caption{Integration contours}
                \label{FigIntegrationContours}
        
            \end{figure}
        
        \end{remark}
        
        \begin{proposition}
        \label{PropSpectralZetaIntegralRepresentation}
        
            On the half-plane $\Re s > 1$, we have
            
            \[ \begin{array}{lll}
                    \zeta_{L,\mu} \left( s \right) & = & \left \{ \begin{array}{ll} \frac{1}{2 i \pi} \sum\limits_{k \in \mathbb{Z}} \int_{i \gamma_{\vartheta}} \left( \frac{1}{4} - t^2 + \mu \right)^{-s} \frac{\partial}{\partial t} \log K_t \left( 2 \pi \left \vert k + \alpha \right \vert a \right) \mathrm{d}t & \hspace{-0.2em} \text{if } \alpha \neq 0 \\[1em] \frac{1}{2 i \pi} \sum\limits_{k \neq 0} \int_{i \gamma_{\vartheta}} \left( \frac{1}{4} - t^2 + \mu \right)^{-s} \frac{\partial}{\partial t} \log K_t \left( 2 \pi \left \vert k + \alpha \right \vert a \right) \mathrm{d}t & \hspace{-0.2em} \text{if } \alpha = 0 \end{array} \right.
                \end{array} . \]
        
        \end{proposition}
        
        \begin{proof}
        
            Let us assume that we have $\alpha \neq 0$, since the other case is similar. For every integer $k \in \mathbb{Z}$, the zeros of $\nu \mapsto K_{i \nu} \left( 2 \pi \left \vert k + \alpha \right \vert a \right)$ are simple, and denoted by $r_{k,j}$. The argument principle then states that we have, for any $k \in \mathbb{Z}$,
            
            \[ \begin{array}{lll}
                    \sum\limits_{j \geqslant 1} \left( \frac{1}{4} + r_{k,j}^2 + \mu \right)^{-s} & = & \frac{1}{2 i \pi} \int_{\gamma_{\vartheta}} \; \left( \frac{1}{4} + r^2 + \mu \right)^{-s} \frac{\partial}{\partial r} \log K_{ir} \left( 2 \pi \left \vert k+\alpha \right \vert a \right) \mathrm{d}r
                \end{array} . \]
            
            \noindent Using the Weyl type law, in order to make sense of the various series, we get
            
            \[ \begin{array}{lll}
                    \zeta_{L, \mu} \left( s \right) & = & \frac{1}{2 i \pi} \sum\limits_{k \in \mathbb{Z}} \; \int_{\gamma_{\vartheta}} \; \left( \frac{1}{4} + r^2 + \mu \right)^{-s} \frac{\partial}{\partial r} \log K_{ir} \left( 2 \pi \left \vert k+\alpha \right \vert a \right) \mathrm{d}r.
                \end{array} \]
            
            \noindent The change of variable $t = ir$ then gives the required formula.
        
        \end{proof}
        
        \begin{remark}
        
            The function $\zeta_{L,\mu}$ does not depend on the angle $\vartheta \in \left] 0, \pi/2 \right[$ chosen to define $\gamma_{\vartheta}$, though it cannot be $\pi/2$, as the complex power has to be well-defined.
        
        \end{remark}

    \subsection{\texorpdfstring{Letting $\vartheta$ go to $\pi/2$}{Letting vartheta go to pi/2}}
    \label{SubSecLettingThetaPiOver2}
    
        Let us now see how the angle $\vartheta$ can approach $\pi/2$.
        
        \begin{definition}
        \label{Deffmuk}
        
            For any integer $k$, we define the complex function $f_{\mu,k}$ on $\mathbb{C}$ by
            
            \[ \begin{array}{lll}
                    f_{\mu,k} \left( t \right) & = & \frac{\partial}{\partial t} \log K_t \left( 2\pi \left \vert k + \alpha \right \vert a \right) - \frac{2t}{\sqrt{4 \mu + 1}} \frac{\partial}{\partial t}_{\left \vert t = \sqrt{\frac{1}{4} + \mu} \right.} \log K_t \left( 2 \pi \left \vert k + \alpha \right \vert a \right).
                \end{array} \]
            
        \end{definition}
        
        The introduction of this function will be justified shortly, and is completely similar to what is done in \cite[Sec. 6.1]{MR4167014}. For now, let us note that we have
        
        \[ \begin{array}{lll}
                \int_{i \gamma_{\vartheta}} \left( \frac{1}{4} + \mu - t^2 \right)^{-s} t \; \mathrm{d}t & = & 0
            \end{array} \]
        
        \noindent on the half-plane $\Re s > 1$. We thus get
        
        \[ \begin{array}{lll}
                \int_{i \gamma_{\vartheta}} \left( \frac{1}{4} + \mu - t^2 \right)^{-s} \frac{\partial}{\partial t} \log K_t \left( 2 \pi \left \vert k + \alpha \right \vert a \right) \mathrm{d}t & = & \int_{i \gamma_{\vartheta}} \left( \frac{1}{4} + \mu - t^2 \right)^{-s} f_{\mu,k} \left( t \right) \mathrm{d}t
            \end{array} . \]
        
        \noindent The issue in letting $\vartheta$ go to $\pi/2$ lies with the complex power above, which needs to be well-defined. Let us see where $1/4 + \mu -t^2$ lands when $t$ goes through $i \gamma_{\vartheta}$.
        
        \begin{figure}[H]
        
            \centering
            
            \scalebox{0.8}{%
            \begin{minipage}{\linewidth}
            \centering
            \begin{tikzpicture}[scale=0.8]
        
                \draw[thick] (-3,0) -- (1,0) coordinate (1) node {};
                \draw[thick, ->] (1,0) -- (3.4,0) coordinate (3) node {};
                \draw[thick, ->] (0,-2.4) -- (0,2.4);
            
                \draw[thick, red] (1,0) -- (3,2) coordinate (2) node {};
                \draw[thick, red] (1,0) -- (3,-2);
            
                \draw pic["$\scriptstyle \color{orange} 2 \vartheta$", draw=orange, ->, angle eccentricity=1.2, angle radius=1.3cm] {angle=3--1--2};
            
                \draw[thick, red, ->] (2.5,1.5) -- (2.4,1.4);
                \draw[thick, red, ->] (2.4,-1.4) -- (2.5,-1.5);
            
                \node (P) at (4.5,-1) {$\scriptstyle \color{red} \frac{1}{4} + \mu - t^2, \; \; \; t \, \in \, i \gamma_{\vartheta}$};
            
                \node (Q) at (0.6,-0.5) {$\scriptstyle \color{red} \frac{1}{4} + \mu$};
        
            \end{tikzpicture}
            
            \end{minipage}
            }
            
            \caption{Variation on the contour $i \gamma_{\vartheta}$}
            \label{FigVariationT}

        \end{figure}
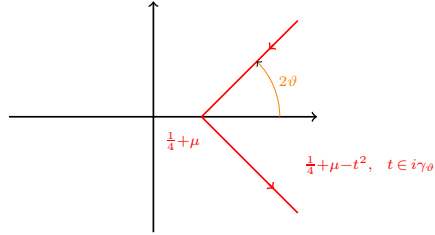
        
        \noindent As the figure above makes clear, the argument $1/4+\mu-t^2$ can collapse onto the half-line $\left] - \infty, 0 \right[$ when $\vartheta$ goes to $\pi/2$. To correct that, we split $i \gamma_{\vartheta}$ into four parts.
        
        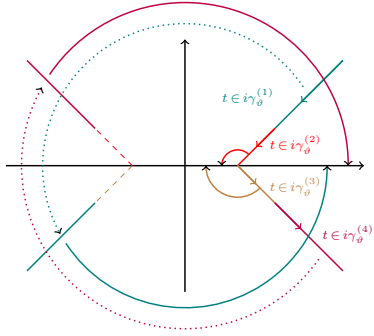
\begin{figure}[H]
    
            \centering
            \scalebox{0.7}{%
            \begin{minipage}{\linewidth}
            \centering
            \begin{tikzpicture}
        
                \draw[thick] (-3.4,0) coordinate (1) node {} -- (1,0) coordinate (2) node {};
                \draw[thick, ->] (1,0) -- (3.4,0) coordinate (3) node {};
                \draw[thick, ->] (0,-2.4) -- (0,2.4);
            
                \draw[thick, red] (1,0) -- (1+0.707106781186548,0.707106781186548) coordinate (4) node {};
                \draw[thick, teal] (1+0.707106781186548,0.707106781186548) -- (3,2) coordinate (5) node {};
                \draw[thick, brown] (1,0) -- (1+0.707106781186548,-0.707106781186548) coordinate (6) node {};
                \draw[thick, purple] (1+0.707106781186548,-0.707106781186548) -- (3,-2) coordinate (7) node {};
            
                \node (A) at (2.1,0.4) {$\scriptstyle \color{red} t \, \in \, i \gamma_{\vartheta}^{\left( 2 \right)}$};
                \node (B) at (1.2,1.3) {$\scriptstyle \color{teal} t \, \in \, i \gamma_{\vartheta}^{\left( 1 \right)}$};
                \node (C) at (2.1,-0.4) {$\scriptstyle \color{brown} t \, \in \, i \gamma_{\vartheta}^{\left( 3 \right)}$};
                \node (D) at (3.1,-1.3) {$\scriptstyle \color{purple} t \, \in \, i \gamma_{\vartheta}^{\left( 4 \right)}$};
            
                \draw[thick, red, ->] (1+0.707106781186548,0.707106781186548) -- (1+0.707106781186548/2,0.707106781186548/2);
                \draw[thick, teal, ->] (3,2) -- (1+0.707106781186548+1/2,0.707106781186548+1/2);
                \draw[thick, brown, ->] (1,0) -- (1+0.707106781186548/2,-0.707106781186548/2);
                \draw[thick, purple, ->] (1+0.707106781186548,-0.707106781186548) -- (1+0.707106781186548+1/2,-0.707106781186548-1/2);
            
                \draw[dashed, red] (-1,0) -- (-1-0.707106781186548,0.707106781186548) coordinate (8) node {};
                \draw[thick, purple] (-1-0.707106781186548,0.707106781186548) -- (-3,2) coordinate (9) node {};
                \draw[dashed, brown] (-1,0) -- (-1-0.707106781186548,-0.707106781186548) coordinate (10) node {};
                \draw[thick, teal] (-1-0.707106781186548,-0.707106781186548) -- (-3,-2) coordinate (11) node {};
            
                \node (P) at (0,-1) {};
                \node (Q) at (0,1) {};
                \node (O) at (0,0) {};
            
                \node (12) at (-1-0.707106781186548-1/2,-0.707106781186548-1/2) {};
            
                \node (13) at (-1-0.707106781186548-1/2,0.707106781186548+1/2) {};
            
                \node (14) at (-1-0.707106781186548-1/2+1/6,0.707106781186548+1/2+1/6) {};
            
                \node (15) at (-1-0.707106781186548-1/2+1/6,-0.707106781186548-1/2-1/8) {};
            
                \draw[thick] pic["", draw=red, ->, angle eccentricity=1.2, angle radius=0.3cm] {angle=4--2--O};
            
                \draw[thick] pic["", draw=brown, <-, angle eccentricity=1.2, angle radius=0.6cm] {angle=O--2--6};
            
                \draw[thick, dotted] pic["", draw=teal, ->, angle eccentricity=1, angle radius=2.7cm] {angle=5--O--12};
            
                \draw[thick, dotted] pic["", draw=purple, <-, angle eccentricity=1, angle radius=3.1cm] {angle=13--O--7};
            
                \draw[thick] pic["", draw=purple, <-, angle eccentricity=1, angle radius=3.1cm] {angle=3--O--14};
            
                \draw[thick] pic["", draw=teal, ->, angle eccentricity=1, angle radius=2.7cm] {angle=15--O--3};
        
            \end{tikzpicture}
            \end{minipage}
            }
            
            \caption{Modification of $i \gamma_{\vartheta}$ and limit as $\vartheta$ goes to $\frac{\pi}{2}$}
            \label{FigContourFourParts}
    
        \end{figure}
        
        \begin{definition}
        
            The four paths of integration $\gamma_{\vartheta}^{\left( 1 \right)}$, \dots, $\gamma_{\vartheta}^{\left( 4 \right)}$ are defined as follows:
            
            \[ \begin{array}{llllll}
                    \color{teal} \gamma_{\vartheta}^{\left( 1 \right)} & = & \left \{ r e^{i \vartheta} \in \gamma_{\vartheta}, \; r \, \geqslant \, \sqrt{\frac{1}{4} + \mu} \right \}, & \color{brown} \gamma_{\vartheta}^{\left( 3 \right)} & = & \left \{ r e^{-i \vartheta} \in \gamma_{\vartheta}, \; r \, < \, \sqrt{\frac{1}{4} + \mu} \right \}, \\[1em]

                    \color{red} \gamma_{\vartheta}^{\left( 2 \right)} & = & \left \{ r e^{i \vartheta} \in \gamma_{\vartheta}, \; r \, < \, \sqrt{\frac{1}{4} + \mu} \right \}, & \color{purple} \gamma_{\vartheta}^{\left( 4 \right)} & = & \left \{ r e^{- i \vartheta} \in \gamma_{\vartheta}, \; r \, \geqslant \, \sqrt{\frac{1}{4} + \mu} \right \}.

                \end{array} \]
        
        \end{definition}
        
        \begin{remark}
        
            Since we have partitioned $\gamma_{\vartheta}$ into four parts, the integral over $i \gamma_{\vartheta}$ can be written as a sum of four integrals.
        
        \end{remark}
        
        Going back to figure \ref{FigContourFourParts}, we note that parts $\color{red} \gamma_{\vartheta}^{\left( 2 \right)}$ and $\color{brown} \gamma_{\vartheta}^{\left( 3 \right)}$ are easier to deal with, as letting $\vartheta$ go to $\pi/2$ is not an issue. We have
        
        \[ \begin{array}{lllll}
                \multicolumn{5}{l}{\int_{i \gamma_{\vartheta}^{\left( 2 \right)}} \; \left( \frac{1}{4} + \mu - t^2 \right)^{-s} \; f_{\mu,k} \left( t \right) \text{d}t + \int_{i \gamma_{\vartheta}^{\left( 3 \right)}} \; \left( \frac{1}{4} + \mu - t^2 \right)^{-s} \; f_{\mu,k} \left( t \right) \text{d}t} \\[1em]
                
                \qquad \qquad \qquad \qquad \qquad \qquad & \underset{\vartheta \rightarrow \frac{\pi}{2}^{-}}{\longrightarrow} & \int_{- \sqrt{\frac{1}{4}+\mu}}^{\sqrt{\frac{1}{4}+\mu}} \; \left( \frac{1}{4} + \mu - t^2 \right)^{-s} \; f_{\mu,k} \left( t \right) \text{d}t & = & 0
            \end{array} \]
        
        \noindent where the last equality stems from the oddness of $f_{\mu,k}$. Thus $\color{teal} \gamma_{\vartheta}^{\left( 1 \right)}$ and $\color{purple} \gamma_{\vartheta}^{\left( 4 \right)}$ are the most interesting parts of the integration contour $\gamma_{\vartheta}$, and we have
        
        \[ \begin{array}{lll}
                \left( \frac{1}{4} + \mu - t^2 \right)^{-s} & = & \left \{ \begin{array}{lcl} e^{-is\pi} \left( t^2 - \left( \frac{1}{4} + \mu \right) \right)^{-s} & \text{on} & i\color{teal} \gamma_{\vartheta}^{\left( 1 \right)} \\[1em] e^{is\pi} \left( t^2 - \left( \frac{1}{4} + \mu \right) \right)^{-s} & \text{on} & i\color{purple} \gamma_{\vartheta}^{\left( 4 \right)} \end{array} \right.
            \end{array} . \]
        
        \noindent This last manipulation is represented on figure \ref{FigContourFourParts} by dotted arcs. It should also be noted that the difference of sign in the exponential above has to do with the choice of branch for the logarithm. We have
        
        \[ \begin{array}{lllll}
                \multicolumn{3}{l}{\int_{i \gamma_{\vartheta}^{\left( 1 \right)}} \; \left( \frac{1}{4} + \mu - t^2 \right)^{-s} \; f_{\mu,k} \left( t \right) \; \text{d}t} & = & e^{-is\pi} \int_{i \gamma_{\vartheta}^{\left( 1 \right)}} \; \left( t^2 - \left( \frac{1}{4} + \mu \right) \right)^{-s} \; f_{\mu,k} \left( t \right) \; \text{d}t \\[1em]
                
                \qquad \qquad \qquad \qquad & \underset{\vartheta \rightarrow \frac{\pi}{2}^{-}}{\longrightarrow} & \multicolumn{3}{l}{- e^{-is\pi} \int_{\sqrt{\frac{1}{4}+\mu}}^{+\infty} \; \left( t^2 - \left( \frac{1}{4} + \mu \right) \right)^{-s} \; f_{\mu,k} \left( t \right) \; \text{d}t,}
            \end{array} \]
        
        \noindent for the first part of $\gamma_{\vartheta}$, and
        
        \[ \begin{array}{lllll}
                \multicolumn{3}{l}{\int_{i \gamma_{\vartheta}^{\left( 4 \right)}} \; \left( \frac{1}{4} + \mu - t^2 \right)^{-s} \; f_{\mu,k} \left( t \right) \; \text{d}t} & = & e^{is\pi} \int_{i \gamma_{\vartheta}^{\left( 4 \right)}} \; \left( t^2 - \left( \frac{1}{4} + \mu \right) \right)^{-s} \; f_{\mu,k} \left( t \right) \; \text{d}t \\[1em]
                
                \qquad \qquad \qquad \qquad & \underset{\vartheta \rightarrow \frac{\pi}{2}^{-}}{\longrightarrow} & \multicolumn{3}{l}{e^{is\pi} \int_{\sqrt{\frac{1}{4}+\mu}}^{+\infty} \; \left( t^2 - \left( \frac{1}{4} + \mu \right) \right)^{-s} \; f_{\mu,k} \left( t \right) \; \text{d}t}
            \end{array} \]
        
        \noindent for the last one. Putting these results together, and using the fact that the integral
        
        \[ \begin{array}{c}
                \int_{i \gamma_{\vartheta}} \; \left( \frac{1}{4} + \mu - t^2 \right)^{-s} \; f_{\mu,k} \left( t \right) \; \text{d}t
            \end{array} \]
        
        \noindent is constant in $\vartheta$, we get the equality
        
        \[ \begin{array}{lll}
                \int_{i \gamma_{\vartheta}} \left( \frac{1}{4} + \mu - t^2 \right)^{-s} f_{\mu,k} \left( t \right) \text{d}t & = & 2 i \sin \left( \pi s \right) \int_{\sqrt{\frac{1}{4}+\mu}}^{+\infty} \left( t^2 - \left( \frac{1}{4} + \mu \right) \right)^{-s} f_{\mu,k} \left( t \right) \text{d}t
            \end{array} . \]
        
        So far, we have neglected to say for which complex numbers $s$ these equalities hold. We will do so now, and thus see why introducing the function $f_{\mu,k}$ was important. We have
        
        \[ \begin{array}{lll}
                \multicolumn{3}{l}{\left( t^2 - \left( \frac{1}{4} + \mu \right) \right)^{-s} f_{\mu,k} \left( t \right)} \\[1em]
                
                & = & \frac{1}{t - \sqrt{\frac{1}{4} + \mu}} \left[ \frac{\partial}{\partial t} \log K_t \left( 2 \pi \left \vert k + \alpha \right \vert a \right) - \frac{2t}{\sqrt{4 \mu + 1}}\frac{\partial}{\partial t}_{\vert t \, = \, \sqrt{\frac{1}{4}+\mu}} \log K_t \left( 2 \pi \left \vert k + \alpha \right \vert a \right) \right] \\[1em]
                
                && \qquad \qquad \qquad \qquad \qquad \qquad \qquad \qquad \qquad \qquad \cdot \frac{1}{\left( t^2 - \left( \frac{1}{4} + \mu \right) \right)^{s-1}} \cdot \frac{1}{t + \sqrt{\frac{1}{4}+\mu}} .
            \end{array} \]
        
        \noindent Noting that the first factor on the right-hand side above is a difference quotient, we see that $t \mapsto \left( t^2 - \left( 1/4 + \mu \right) \right)^{-s} f_{\mu,k} \left( t \right)$ is integrable at $\sqrt{1/4 + \mu}$ if and only if we have $\Re s < 2$. The integrability condition at $+ \infty$ is still $\Re s > 1$. We can summarize the discussion of this paragraph with the following proposition.
        
        \begin{proposition}
        \label{PropSpectralZetaIntegralRepresentation2}
        
            On the strip $1 \, < \, \Re s \, < \, 2$, the spectral function $\zeta_{L, \mu}$ is given by
            
            \[ \begin{array}{lll}
                    \zeta_{L, \mu} \left( s \right) & = & \left \{ \begin{array}{ll} \frac{\sin \left( \pi s \right)}{\pi} \; \sum\limits_{k \in \mathbb{Z}} \; \int_{\sqrt{\frac{1}{4}+\mu}}^{+\infty} \; \left( t^2 - \left( \frac{1}{4} + \mu \right) \right)^{-s} \; f_{\mu,k} \left( t \right) \; \text{d}t & \text{if } \alpha \neq 0 \\[1em] \frac{\sin \left( \pi s \right)}{\pi} \; \sum\limits_{k \in \mathbb{Z} \setminus \left \{ 0 \right \}} \; \int_{\sqrt{\frac{1}{4}+\mu}}^{+\infty} \; \left( t^2 - \left( \frac{1}{4} + \mu \right) \right)^{-s} \; f_{\mu,k} \left( t \right) \; \text{d}t & \text{if } \alpha = 0 \end{array} \right.
                \end{array} . \]
        
        \end{proposition}
        
        \begin{proof}
        
            This is a consequence of the discussion made right before the statement of this proposition.
        
        \end{proof}

    \subsection{Splitting the interval of integration}
    \label{SubSecSplittingIntervalOfIntegration}
    
        In proposition \ref{PropSpectralZetaIntegralRepresentation2}, we obtained an expression of the spectral zeta function $\zeta_{L, \mu}$ on a strip as a sum of integrals. In order to prove the existence of a holomorphic continuation, we will need to study these integrals.
        
        \begin{definition}
        \label{DefImuk}
        
            For any real number $\mu \geqslant 0$ and integer $k \in \mathbb{Z}$, the integral $I_{\mu,k}$  is defined on the strip $1 < \Re s < 2$ by
            
            \[ \begin{array}{lll}
                    I_{\mu, k} \left( s \right) & = & \displaystyle \frac{\sin \left( \pi s \right)}{\pi} \int_{\sqrt{\frac{1}{4}+\mu}}^{+\infty} \; \left( t^2 - \left( \frac{1}{4} + \mu \right) \right)^{-s} \; f_{\mu,k} \left( t \right) \; \text{d}t
                \end{array} , \]
            
            \noindent with the exception of $k = 0$ should we have $\alpha = 0$.
        
        \end{definition}
        
        The idea to study these terms is to apply the \textit{binomial formula}
        
        \[ \begin{array}{lll}
                \left( t^2 - \left( \frac{1}{4} + \mu \right) \right)^{-s} & = & \sum\limits_{j = 0}^{+ \infty} \; \frac{\left( s \right)_j}{j!} \left( \frac{1}{4} + \mu \right)^j \cdot \frac{1}{t^{2s+2j}}
            \end{array} , \]
            
        \noindent which one can obtain by applying proposition \ref{PropBinomialFormula}. This result holds on the interval of integration. However, in order to interchange the various sums and integrals, one needs ``some space'' between $\sqrt{1/4 + \mu}$ and $t$. This leads us to perform the following splitting, much in the fashion of \cite{MR4167014},
        
        \begin{equation}
        \label{EqSplittingkneq0}
            \begin{array}{lll}
                \left] \sqrt{\frac{1}{4} + \mu}, \, + \infty \right[ & = & \left] \sqrt{\frac{1}{4} + \mu}, \; 2 \left \vert k \right \vert^{\delta} \sqrt{\frac{1}{4} + \mu} \right[ \sqcup \left[ 2 \left \vert k \right \vert^{\delta} \sqrt{\frac{1}{4} + \mu}, \; +\infty \right[
            \end{array}
        \end{equation}
        
        \noindent for every integer $k \neq 0$. We have considered a real number $\delta > 0$ here, which will be adjusted throughout the rest of this section. Its sole purpose is to facilitate the convergence of series. In order to deal with the case $k = 0$, when $\alpha$ does not vanish, in a similar way, we write
        
        \begin{equation}
        \label{EqSplittingk0}
            \begin{array}{lll}
                \left] \sqrt{\frac{1}{4} + \mu}, \, + \infty \right[ & = & \left] \sqrt{\frac{1}{4} + \mu}, \; 2 \sqrt{\frac{1}{4} + \mu} \right[ \sqcup \left[ 2 \sqrt{\frac{1}{4} + \mu}, \; +\infty \right[
            \end{array} ,
        \end{equation}
        
        \noindent though the location of the splitting point here does not matter. We can now split every integral $I_{\mu,k}$ accordingly.
        
        \begin{definition}
        \label{DefLmuk}
        
            For every $k \in \mathbb{Z} \setminus \left \{ 0 \right \}$ and every $s \in \mathbb{C}$ with $1 < \Re s < 2$, we set
            
            \[ \begin{array}{lll}
                    L_{\mu,k} \left( s \right) & = & \displaystyle \frac{\sin \left( \pi s \right)}{\pi} \int_{\sqrt{\frac{1}{4}+\mu}}^{2 \left \vert k \right \vert^{\delta} \sqrt{\frac{1}{4} + \mu}} \; \left( t^2 - \left( \frac{1}{4} + \mu \right) \right)^{-s} \; f_{\mu,k} \left( t \right) \; \text{d}t
                \end{array} . \]
            
            \noindent Should $\alpha$ not vanish, we also set, on the same strip,
            
            \[ \begin{array}{lll}
                    L_{\mu,0} \left( s \right) & = & \displaystyle \frac{\sin \left( \pi s \right)}{\pi} \int_{\sqrt{\frac{1}{4}+\mu}}^{2 \sqrt{\frac{1}{4} + \mu}} \; \left( t^2 - \left( \frac{1}{4} + \mu \right) \right)^{-s} \; f_{\mu,0} \left( t \right) \; \text{d}t
                \end{array} . \]
        
        \end{definition}
        
        \begin{definition}
        \label{DefMmuk}
        
            For every $k \in \mathbb{Z} \setminus \left \{ 0 \right \}$ and every $s \in \mathbb{C}$ with $1 < \Re s < 2$, we set
            
            \[ \begin{array}{lll}
                    M_{\mu,k} \left( s \right) & = & \displaystyle \frac{\sin \left( \pi s \right)}{\pi} \int_{2 \left \vert k \right \vert^{\delta} \sqrt{\frac{1}{4} + \mu}}^{+ \infty} \; \left( t^2 - \left( \frac{1}{4} + \mu \right) \right)^{-s} \; f_{\mu,k} \left( t \right) \; \text{d}t
                \end{array} . \]
            
            \noindent Should $\alpha$ not vanish, we also set, on the same strip,
            
            \[ \begin{array}{lll}
                    M_{\mu,0} \left( s \right) & = & \displaystyle \frac{\sin \left( \pi s \right)}{\pi} \int_{2 \sqrt{\frac{1}{4} + \mu}}^{+ \infty} \; \left( t^2 - \left( \frac{1}{4} + \mu \right) \right)^{-s} \; f_{\mu,0} \left( t \right) \; \text{d}t
                \end{array} . \]
        
        \end{definition}
        
        \begin{remark}
        
            The last two definitions give $I_{\mu,k} \left( s \right) = L_{\mu,k} \left( s \right) + M_{\mu,k} \left( s \right)$, whenever these integrals make sense.
        
        \end{remark}

    \subsection{\texorpdfstring{Study of the integrals $L_{\mu,k}$}{Study of the integrals Lmu,k}}
    \label{SubSecStudyLmuk}
        
        We begin the study of the spectral zeta function by that of the integrals $L_{\mu,k}$ from definition \ref{DefLmuk}, as in \cite[Sec. 6.2]{MR4167014}. The difference is that we must keep the parameters $\mu$ and~$\alpha$.
        
        \subsubsection{Global study}
        \label{SubSubSecGlobalStudy}
        
            The first step is a global study of the integrals $L_{\mu,k}$, which will lead to splitting them into two different parts, which we will then study separately.
        
            \begin{definition}
            \label{DefFmuk}
        
                For any $\mu \geqslant 0$, we define the complex function $F_{\mu,k}$ by
            
                \[ \begin{array}{lll}
                        F_{\mu,k} \left( z \right) & = & \log K_z \left( 2 \pi \left \vert k + \alpha \right \vert a \right) - \log K_{\sqrt{\frac{1}{4} + \mu}} \left( 2 \pi \left \vert k + \alpha \right \vert a \right) \\[1em]
                    
                        && \qquad \qquad \qquad \qquad - \frac{z^2 - \left( 1/4 + \mu \right)}{\sqrt{4 \mu + 1}} \frac{\partial}{\partial z}_{\left \vert z = \sqrt{\frac{1}{4} + \mu} \right.} \log K_z \left( 2 \pi \left \vert k + \alpha \right \vert a \right)
                    \end{array} \]
                
                \noindent for $z$ in the angular sector $\left \vert z \right \vert < \pi / 4$, on which $K_z \left( 2 \pi \left \vert k + \alpha \right \vert a \right)$ does not vanish.
        
            \end{definition}
        
            The next result, is the same as \cite[Cor. 6.4]{MR4167014}, though it alone will not be sufficient. It is nevertheless a crucial step.
        
            \begin{proposition}
            \label{PropGlobalStudy}
        
                For any relative integer $k \neq 0$, we can write
            
                \[ \begin{array}{llll}
                        F_{\mu,k} \left( t \right) & = & \left( t^2 - \left( \frac{1}{4} + \mu \right) \right) \widetilde{F}_{\mu,k} \left( t \right) & \text{for } t \in \left[ \sqrt{\frac{1}{4} + \mu}, \; 2 \left \vert k \right \vert^{\delta} \sqrt{\frac{1}{4} + \mu} \right]
                    \end{array} \]
            
                \noindent where the function $\widetilde{F}_{\mu,k}$ is analytic in $t$ and satisfies a bound of the type
            
                \[ \begin{array}{ccc}
                        \left \vert \widetilde{F}_{\mu,k} \right \vert & \leqslant & C_{\mu} \cdot \frac{1}{\left \vert k \right \vert^{2-4\delta} a^2}
                    \end{array} , \]
            
                \noindent uniformly on the same interval, with a constant $C_{\mu} \, > \, 0$ depending only on $\mu$.
        
            \end{proposition}
        
            \begin{proof}
        
                We first note that the function $F_{\mu,k}$ has been defined so as to have
            
                \[ \begin{array}{lllll}
                        F_{\mu,k} \left( \pm \sqrt{\frac{1}{4}+\mu} \right) & = & F_{\mu,k}' \left( \pm \sqrt{\frac{1}{4}+\mu} \right) & = & 0
                    \end{array} . \]
            
                \noindent Since $F_{\mu,k}$ is holomorphic in $t$, it is of the form $F_{\mu,k} \left( t \right)= h_{\mu,k} \left( t^2 \right)$, where $h_{\mu,k}$ is holomorphic and such that we have
            
                \[ \begin{array}{lllll}
                        h_{\mu,k} \left( \frac{1}{4}+\mu \right) & = & h_{\mu,k}' \left( \frac{1}{4}+\mu \right) & = & 0
                    \end{array} . \]
            
                \noindent The Taylor-Lagrange theorem then allows us to write
            
                \[ \begin{array}{lllll}
                        F_{\mu,k} \left( t \right) & = & h_{\mu,k} \left( t^2 \right) & = & \frac{1}{2} \, \left( t^2 - \left( \frac{1}{4}+\mu \right) \right)^2 \; h_{\mu,k}'' \left( \xi_{\mu,t}^2 \right)
                    \end{array} , \]
            
                \noindent where $\xi_{\mu,t}$ is a real number with $\sqrt{1/4 + \mu} \leqslant \xi_{\mu,t} \leqslant 2 \left \vert k \right \vert^{\delta} \sqrt{1/4 + \mu}$. Note that we do not know how $\xi_{\mu,t}$ depends on $\mu$, $t$, or $k$. By differentiating $F_{\mu,k}$, we get
            
                \[ \begin{array}{lllclll}
                        F_{\mu,k}' \left( t \right) & = & 2t h_{\mu,k}' \left( t^2 \right), & & F_{\mu,k}'' \left( t \right) & = & 2 h_{\mu,k}' \left( t^2 \right) + 4t^2 h_{\mu,k}'' \left( t^2 \right)
                    \end{array} , \]
            
                \noindent and these two equalities can be combined to yield
            
                \[ \begin{array}{lll}
                        h_{\mu,k}'' \left( t^2 \right) & = & \frac{1}{4t^2} \, F_{\mu,k}'' \left( t \right) - \frac{1}{4t^3} \, F_{\mu,k}' \left( t \right) \\[0.5em]
                    
                        & = & \frac{1}{4t^2} \, \frac{\partial^2}{\partial t^2} \log K_t \left( 2 \pi \left \vert k + \alpha \right \vert a \right) - \frac{1}{4t^3} \, \frac{\partial}{\partial t} \log K_t \left( 2 \pi \left \vert k + \alpha \right \vert a \right)
                    \end{array} \]
            
                \noindent Therefore, we have
            
                \[ \scalemath{0.98}{\begin{array}{lll}
                        F_{\mu,k} \left( t \right) & = & \left( t^2 - \left( \frac{1}{4}+\mu \right) \right)^2 \left[ \frac{1}{4\xi_{\mu,t}^2} \, \frac{\partial^2}{\partial t^2}_{\vert t \, = \, \xi_{\mu,t}} \log K_t \left( 2 \pi \left \vert k + \alpha \right \vert a \right) \right. \\[1em]
                    
                        && \qquad \qquad \qquad \qquad \qquad \qquad - \left. \frac{1}{4\xi_{\mu,t}^3} \, \frac{\partial}{\partial t}_{\vert t \, = \, \xi_{\mu,t}} \log K_t \left( 2 \pi \left \vert k + \alpha \right \vert a \right) \right].
                    \end{array}} \]
            
                \noindent For any real number $\xi$ such that we have
            
                \[ \begin{array}{lll}
                        \xi & \in & \left[ \sqrt{\frac{1}{4}+\mu}, \; 2 \left \vert k \right \vert^{\delta} \sqrt{\frac{1}{4}+\mu} \right]
                    \end{array} , \]
            
                \noindent we denote by $D_{\xi}$ the disk centered at $\xi$ of radius $1/4$. The Bessel function $K_{\nu} \left( z \right)$ being entire in $\nu$ for any positive real number $z$, the Cauchy formula gives
            
                \[ \begin{array}{lll}
                        \frac{\partial}{\partial t}_{\vert t \, = \, \xi} \log K_t \left( 2 \pi \left \vert k + \alpha \right \vert a \right) & = & \frac{1}{2i\pi K_{\xi} \left( 2 \pi \left \vert k + \alpha \right \vert a \right)} \; \int_{\partial D_{\xi}} \; \frac{K_{\nu} \left( 2 \pi \left \vert k + \alpha \right \vert a \right)}{\left( \nu - \xi \right)^2} \; \text{d}\nu
                    \end{array} . \]
            
                \noindent Using proposition \ref{PropAsymptoticsModifiedBesselParameter}, we get
            
                \[ \scalemath{0.96}{\begin{array}{lll}
                        \multicolumn{3}{l}{\frac{\partial}{\partial t}_{\vert t \, = \, \xi} \log K_t \left( 2 \pi \left \vert k + \alpha \right \vert a \right)} \\[0.5em]
                    
                        \qquad & = & \frac{1}{2i\pi} \; \frac{K_{1/2} \left( 2 \pi \left \vert k + \alpha \right \vert a \right)}{K_{\xi} \left( 2 \pi \left \vert k + \alpha \right \vert a \right)} \; \int_{\partial D_{\xi}} \; \frac{1}{\left( \nu - \xi \right)^2} \left( 1 + \frac{A_1 \left( \nu \right)}{2 \pi \left \vert k + \alpha \right \vert a} + \gamma_2 \left( \nu, 2 \pi \left \vert k + \alpha \right \vert a \right) \right) \text{d}\nu \\[0.5em]
                    
                        & = & \frac{K_{1/2} \left( 2 \pi \left \vert k + \alpha \right \vert a \right)}{K_{\xi} \left( 2 \pi \left \vert k + \alpha \right \vert a \right)} \; \left( \frac{\xi}{2 \pi \left \vert k + \alpha \right \vert a} + O \left( \frac{1}{\left \vert k \right \vert^{2-4 \delta} a^2} \right) \right).
                    \end{array}} \]
            
                \noindent For that last point, we have used the following estimate for the remainder $\gamma_2$
            
                \[ \begin{array}{lllll}
                        \left \vert \gamma_2 \left( \nu, 2 \pi \left \vert k + \alpha \right \vert a \right) \right \vert & \leqslant & 2 \left \vert \frac{A_2 \left( \nu \right)}{\left( 2 \pi \left \vert k + \alpha \right \vert a \right)^2} \right \vert \exp \left( \frac{\left( \frac{1}{4}+\mu \right) \left \vert k \right \vert^{2 \delta}}{2 \pi \left \vert k + \alpha \right \vert a} \right)
                    \end{array} , \]
            
                \noindent which is uniformly bounded in $k$, but not in $\mu$. Similarly, we have
            
                \[ \scalemath{0.96}{\begin{array}{lll}
                        \frac{\partial^2}{\partial t^2}_{\vert t \, = \, \xi} \log K_t \left( 2 \pi \left \vert k + \alpha \right \vert a \right) & = & \frac{K_{1/2} \left( 2 \pi \left \vert k + \alpha \right \vert a \right)}{K_{\xi} \left( 2 \pi \left \vert k + \alpha \right \vert a \right)} \left( 1 + O \left( \frac{1}{\left \vert k \right \vert^{2-4 \delta}} \right) \right)
                    \end{array}} , \]
            
                \noindent which means that we have, still on the same interval,
            
                \[ \scalemath{0.96}{\begin{array}{lll}
                        \widetilde{F}_{\mu,k} \left( t \right) & = & \frac{K_{1/2} \left( 2 \pi \left \vert k + \alpha \right \vert a \right)}{K_{\xi} \left( 2 \pi \left \vert k + \alpha \right \vert a \right)} \cdot O \left( \frac{1}{\left \vert k \right \vert^{2-4 \delta} a^2} \right)
                    \end{array}} , \]
            
                \noindent with an implicit constant depending only on $\mu$. The asymptotics of the modified Bessel function of the second kind show that the factor on the right-hand side is bounded on the interval we consider, uniformly in $k$. This concludes the proof.
        
            \end{proof}
        
            \begin{remark}
            \label{RmkGlobalk0}
        
                Note that a similar result holds when $k = 0$, should $\alpha$ not vanish, with the interval being replaced by the appropriate one from \eqref{EqSplittingk0}.
        
            \end{remark}
        
            We will now further break apart the integrals $L_{\mu,k}$ in the following proposition.
        
            \begin{proposition}
            \label{PropSplittingLmuk}
        
                For any $k \in \mathbb{Z} \setminus \left \{ 0 \right \}$, and any real number $\mu \geqslant 0$, we have
            
                \[ \scalemath{0.96}{\begin{array}{ccl}
                        L_{\mu,k} \left( s \right) & = & \frac{\sin \left( \pi s \right)}{\pi} \; \left( \left( 4 \mu + 1 \right)^{-s} \, \left( \left \vert k \right \vert^{2 \delta} - \frac{1}{4} \right)^{-s} \; F_{\mu,k} \left( 2 \left \vert k \right \vert^{\delta} \sqrt{\frac{1}{4}+\mu} \right) \right. \\[1em]
                    
                        && \qquad \qquad \qquad \left. + 2s \int_{\sqrt{\frac{1}{4}+\mu}}^{2 \left \vert k \right \vert^{\delta} \sqrt{\frac{1}{4}+\mu}} \; t \left( t^2 - \left( \frac{1}{4}+\mu \right) \right)^{-s-1} \; F_{\mu,k} \left( t \right) \; \text{d}t \right),
                    \end{array}} \]
            
                \noindent on the strip $1 \, < \, \Re s \, < \, 2$. If $\alpha$ does not vanish, we also have, on the same strip,
            
                \[ \scalemath{0.96}{\begin{array}{ccl}
                        L_{\mu,0} \left( s \right) & = & \frac{\sin \left( \pi s \right)}{\pi} \; \left( \frac{1}{3^s} \, \left( \frac{1}{4}+\mu \right)^{-s} \; F_{\mu,0} \left( 2 \sqrt{\frac{1}{4}+\mu} \right) \right. \\[1em]
                    
                        && \qquad \qquad \qquad \left. + 2s \int_{\sqrt{\frac{1}{4}+\mu}}^{2 \sqrt{\frac{1}{4}+\mu}} \; t \left( t^2 - \left( \frac{1}{4}+\mu \right) \right)^{-s-1} \; F_{\mu,0} \left( t \right) \; \text{d}t \right).
                    \end{array}} \]
        
            \end{proposition}
        
            \begin{proof}
        
                Only the case $k \neq 0$ will be dealt with here, the other one being perfectly similar. We will further assume throughout this proof that $s$ satisfies $1 \, < \, \Re s \, < \, 2$. Performing an integration by parts on $L_{\mu,k}$, we get
            
                \[ \scalemath{0.9}{\begin{array}{lll}
                        L_{\mu,k} \left( s \right) & = & \frac{\sin \left( \pi s \right)}{\pi} \; \left(  \left[ \left( t^2 - \left( \frac{1}{4} + \mu \right) \right)^{-s} F_{\mu,k} \left( t \right) \right]_{\sqrt{\frac{1}{4}+\mu}}^{2 \left \vert k \right \vert^{\delta} \sqrt{\frac{1}{4}+\mu}} \right. \\[1.5em]
                    
                        && \qquad \qquad \qquad \left. + 2s \int_{\sqrt{\frac{1}{4}+\mu}}^{2 \left \vert k \right \vert^{\delta} \sqrt{\frac{1}{4}+\mu}} \; t \left( t^2 - \left( \frac{1}{4} + \mu \right) \right)^{-s-1} \; F_{\mu,k} \left( t \right) \; \text{d}t \right).
                    \end{array}} \]
            
                \noindent The first term on the right-hand side above can be explicitely computed, as we have
            
                \[ \scalemath{0.95}{\begin{array}{lllll}
                        \left( t^2 - \left( \frac{1}{4} + \mu \right) \right)^{-s} F_{\mu,k} \left( t \right) & = & \left( t^2 - \left( \frac{1}{4} + \mu \right) \right)^{-s+2} R_{\mu,k} \left( t \right) & \underset{t \rightarrow \sqrt{\frac{1}{4}+\mu}^{+}}{\longrightarrow} & 0,
                    \end{array}} \]
            
                \noindent using proposition \ref{PropGlobalStudy}. This completes the proof.
        
            \end{proof}
        
            This integration by parts, made possible by proposition \ref{PropGlobalStudy}, allows us to break each $L_{\mu,k}$ into two parts, to be studied separately. Let us properly define them.
        
            \begin{definition}
            \label{DefAmuk}
        
                For every $k \in \mathbb{Z} \setminus \left \{ 0 \right \}$ and every $s \in \mathbb{C}$ with $1 < \Re s < 2$, we set
            
                \[ \begin{array}{lll}
                        A_{\mu,k} \left( s \right) & = & \displaystyle \frac{\sin \left( \pi s \right)}{\pi} \left( 4 \mu + 1 \right)^{-s} \left( \left \vert k \right \vert^{2 \delta} - \frac{1}{4} \right)^{-s} F_{\mu,k} \left( 2 \left \vert k \right \vert^{\delta} \sqrt{\frac{1}{4}+\mu} \right)
                    \end{array} . \]
            
                \noindent Should $\alpha$ not vanish, we also set, on the same strip,
            
                \[ \begin{array}{lll}
                        A_{\mu,0} \left( s \right) & = & \displaystyle \frac{\sin \left( \pi s \right)}{\pi} \; \frac{1}{3^s} \, \left( \frac{1}{4}+\mu \right)^{-s} \; F_{\mu,0} \left( 2 \sqrt{\frac{1}{4}+\mu} \right)
                    \end{array} . \]
        
            \end{definition}
        
            \begin{definition}
            \label{DefBmuk}
        
                For every $k \in \mathbb{Z} \setminus \left \{ 0 \right \}$ and every $s \in \mathbb{C}$ with $1 < \Re s < 2$, we set
                
                \[ \begin{array}{lll}
                        B_{\mu,k} \left( s \right) & = & \displaystyle 2s \, \frac{\sin \left( \pi s \right)}{\pi} \, \int_{\sqrt{\frac{1}{4}+\mu}}^{2 \left \vert k \right \vert^{\delta} \sqrt{\frac{1}{4}+\mu}} \; t \left( t^2 - \left( \frac{1}{4}+\mu \right) \right)^{-s-1} \; F_{\mu,k} \left( t \right) \; \text{d}t
                    \end{array} . \]
            
                \noindent Should $\alpha$ not vanish, we also set, on the same strip,
            
                \[ \begin{array}{lll}
                        B_{\mu,0} \left( s \right) & = & \displaystyle 2s \, \frac{\sin \left( \pi s \right)}{\pi} \int_{\sqrt{\frac{1}{4}+\mu}}^{2 \sqrt{\frac{1}{4}+\mu}} \; t \left( t^2 - \left( \frac{1}{4}+\mu \right) \right)^{-s-1} \; F_{\mu,0} \left( t \right) \; \text{d}t
                    \end{array} . \]
        
            \end{definition}

        \subsubsection{Study of the terms $B_{\mu,k}$}
        \label{SubSubSecStudyBmuk}
        
            Using the splitting $L_{\mu,k} \left( s \right) = A_{\mu,k} \left( s \right) + B_{\mu,k} \left( s \right)$ obtained above, we will begin by analyzing the behavior of $B_{\mu,k}$, as it is by far the simplest of the two parts.
            
            \begin{proposition}
            \label{PropBmuk}
        
                For any real number $\mu \geqslant 0$, the function
            
                \[ \begin{array}{lll}
                        s & \longmapsto & \sum\limits_{\left \vert k \right \vert \geqslant 1} \; B_{\mu,k} \left( s \right),
                    \end{array} \]
            
                \noindent is holomorphic on the strip $4-1/\left( 2 \delta \right) < \Re s < 2$, and we have
            
                \[ \begin{array}{lll}
                        \frac{\partial}{\partial s}_{\vert s \, = \, 0} \; \; \sum\limits_{\left \vert k \right \vert \geqslant 1} \; B_{\mu,k} \left( s \right) & = & 0.
                    \end{array} \]
        
            \end{proposition}
            
            \begin{proof}
            
                For any non-zero integer $k$, and any real number
            
                \[ \begin{array}{ccc}
                        t & \in & \left] \sqrt{\frac{1}{4}+\mu}, \; 2 \left \vert k \right \vert^{\delta} \sqrt{\frac{1}{4}+\mu} \right[
                    \end{array} , \]
                
                \noindent we can use proposition \ref{PropGlobalStudy} to get the following estimate
            
                \[ \begin{array}{lll}
                        \left \vert \frac{t}{\left( t^2 - \left( \frac{1}{4}+\mu \right) \right)^{s+1}} \; F_{\mu,k} \left( t \right) \right \vert &\leqslant & C_{\mu} \cdot \frac{t}{\left( t^2 - \left( \frac{1}{4}+\mu \right) \right)^{\Re s-1}} \cdot \frac{1}{\left \vert k \right \vert^{2-4 \delta}} \cdot \frac{1}{a^2}.
                    \end{array} \]
            
                \noindent We now note that the right hand side of this inequality can be bounded uniformly in $s$ on any strip $\alpha < \Re s < \beta < 2$ with $\alpha$ and $\beta$ being fixed, possibly negative, real numbers, using the following inequalities
            
                \[ \begin{array}{lll}
                        \frac{1}{\left( t^2 - \left( \frac{1}{4}+\mu \right) \right)^{\Re s-1}} & \leqslant & \left \{ \begin{array}{lcl} \frac{1}{\left( t^2 - \left( \frac{1}{4}+\mu \right) \right)^{\alpha-1}} & \text{if} & t^2 - \left( \frac{1}{4}+\mu \right) > 1 \\[1em]
                    
                        \frac{1}{\left( t^2 - \left( \frac{1}{4}+\mu \right) \right)^{\beta-1}} & \text{if} & t^2 - \left( \frac{1}{4}+\mu \right) < 1 \end{array} \right.
                    \end{array} . \]
            
                \noindent For any such $\alpha$ and $\beta$, the dominated convergence theorem proves that the function
            
                \[ \begin{array}{lll}
                        s & \longmapsto & \int_{\sqrt{\frac{1}{4}+\mu}}^{2 \left \vert k \right \vert^{\delta} \sqrt{\frac{1}{4}+\mu}} \; t \left( t^2 - \left( \frac{1}{4}+\mu \right) \right)^{-s-1} \; F_{\mu,k} \left( t \right) \; \text{d}t
                    \end{array} \]
            
                \noindent is holomorphic on the strip $\alpha \, < \, \Re s \, < \, \beta$, which means that, due to the randomness of $\alpha$ and $\beta$, it is holomorphic on the half-plane $\Re s \, < \, 2$, where we further have
            
                \[ \begin{array}{lll}
                        \multicolumn{3}{l}{\left \vert \int_{\sqrt{\frac{1}{4}+\mu}}^{2 \left \vert k \right \vert^{\delta} \sqrt{\frac{1}{4}+\mu}} \; t \left( t^2 - \left( \frac{1}{4}+\mu \right) \right)^{-s-1} \; F_{\mu,k} \left( t \right) \; \text{d}t \right \vert} \\[1em]
                    
                        \qquad \qquad \qquad \qquad \qquad \quad & \leqslant & \frac{C_{\mu}}{\left \vert k \right \vert^{2-4\delta} a^2} \; \int_{\sqrt{\frac{1}{4}+\mu}}^{2 \left \vert k \right \vert^{\delta} \sqrt{\frac{1}{4}+\mu}} \; t \left( t^2 - \left( \frac{1}{4}+\mu \right) \right)^{- \Re s +1} \; \text{d}t \\[1em]
                    
                        & \leqslant & \frac{C_{\mu}}{\left \vert k \right \vert^{2-4\delta} a^2} \; \left[ \frac{1}{2 - \Re s} \left( t^2 - \left( \frac{1}{4}+\mu \right) \right)^{-\Re s + 2} \right]_{\sqrt{\frac{1}{4}+\mu}}^{2 \left \vert k \right \vert^{\delta} \sqrt{\frac{1}{4}+\mu}} \\[1em]
                    
                        & \leqslant & \frac{C_{\mu}}{2 - \Re s} \; \left( \frac{1}{4} + \mu \right)^{- \Re s + 2} \; \cdot \frac{1}{\left \vert k \right \vert^{2-4\delta} a^2} \cdot \frac{1}{\left( 4 \left \vert k \right \vert^{2 \delta} - 1 \right)^{\Re s - 2}}.
                    \end{array} \]
            
                \noindent The dominated convergence theorem then proves that the function
            
                \[ \begin{array}{ccc}
                        \displaystyle s & \longmapsto & \displaystyle \sum\limits_{\left \vert k \right \vert \geqslant 1} \; \int_{\sqrt{\frac{1}{4}+\mu}}^{2 \left \vert k \right \vert^{\delta} \sqrt{\frac{1}{4}+\mu}} \; t \left( t^2 - \left( \frac{1}{4}+\mu \right) \right)^{-s-1} \; F_{\mu,k} \left( t \right) \; \text{d}t
                    \end{array} \]
            
                \noindent is well-defined and holomorphic on the strip $4-1/\left( 2 \delta \right) < \Re s < 2$, which contains~$0$ if we have $0 < \delta < 1/8$, which we may assume. Hence the function
            
                \[ \begin{array}{lll}
                        \displaystyle s & \longmapsto & \displaystyle \sum\limits_{\left \vert k \right \vert \geqslant 1} \; B_{\mu,k} \left( s \right)
                    \end{array} \]
            
                \noindent is holomorphic around $0$, and we have
            
                \[ \begin{array}{lll}
                        \displaystyle \frac{\partial}{\partial s}_{\vert s \, = \, 0} \; \; \sum\limits_{\left \vert k \right \vert \geqslant 1} \; B_{\mu,k} \left( s \right) & = & \displaystyle 0,
                    \end{array} \]
            
                \noindent because the term $B_{\mu,k}$ involves the product of a function which we have shown was holomorphic around $0$ with the factor $s \sin \left( \pi s \right)$.
            
            \end{proof}
            
            \begin{remark}
            
                The proposition above only considers a sum over non-zero integers, regardless of whether or not $\alpha$ vanishes, so as to give a more uniform result. However, we still need to account for the case $k=0$ when $\alpha$ is non-zero.
            
            \end{remark}
            
            \begin{proposition}
            \label{PropBmu0}
        
                Assume we have $\alpha \neq 0$. For any $\mu \geqslant 0$, the function
            
                \[ \begin{array}{ccc}
                        \displaystyle s & \longmapsto & \displaystyle B_{\mu,0} \left( s \right)
                    \end{array} , \]
            
                \noindent is holomorphic on the half-plane $\Re s < 2$, and  we have
            
                \[ \begin{array}{ccc}
                        \displaystyle \frac{\partial}{\partial s}_{\vert s \, = \, 0} \; \; B_{\mu,0} \left( s \right) & = & \displaystyle 0.
                    \end{array} \]
        
            \end{proposition}
        
            \begin{proof}
        
                This is a simpler version of the argument used in proposition \ref{PropBmuk}.
        
            \end{proof}

        \subsubsection{Study of the terms $A_{\mu,k}$}
        \label{SubSubSecStudyAmuk}
            
            It must be noted outright that understanding the behavior of the series involving the terms $A_{\mu,k}$ introduced in definition \ref{DefAmuk} is significantly more complicated, and will involve a lot of computations.
            
            \begin{proposition}
            \label{PropComputA}
        
                Let $\mu \, \geqslant \, 0$. For every integer $k$, the function $s \longmapsto A_{\mu,k} \left( s \right)$ is holomorphic on $\mathbb{C}$, and its derivative satisfies
            
                \[ \scalemath{0.97}{\begin{array}{ccl}
                        \frac{\partial}{\partial s} \; A_{\mu,k} & = & \cos \left( \pi s \right) \left( 4 \mu + 1 \right)^{-s} \, \left( \left \vert k \right \vert^{2 \delta} - \frac{1}{4} \right)^{-s} \; F_{\mu,k} \left( 2 \left \vert k \right \vert^{\delta} \sqrt{\frac{1}{4}+\mu} \right) \\[1em]
                    
                        && \qquad \qquad - \left[ \frac{\sin \left( \pi s \right)}{\pi} \log \left( \left( 4 \mu + 1 \right) \left( \left \vert k \right \vert^{2 \delta} - \frac{1}{4} \right) \right) \; \left( 4 \mu + 1 \right)^{-s} \right. \\[1em]
                    
                        && \qquad \qquad \qquad \qquad \qquad \left. \cdot \left( \left \vert k \right \vert^{2 \delta} - \frac{1}{4} \right)^{-s} \; F_{\mu,k} \left( 2 \left \vert k \right \vert^{\delta} \sqrt{\frac{1}{4}+\mu} \right) \right],
                    \end{array}} \]
            
                \noindent whenever we have $k \, \neq \, 0$. Assuming $\alpha$ is different from zero, we further have
            
                \[ \scalemath{0.97}{\begin{array}{ccl}
                        \frac{\partial}{\partial s} \; A_{\mu,0} \left( s \right) & = & \cos \left( \pi s \right) \frac{1}{3^s} \, \left( \frac{1}{4}+\mu \right)^{-s} \; F_{\mu,0} \left( 2 \sqrt{\frac{1}{4}+\mu} \right) \\[1em]
                    
                        && \qquad - \frac{\sin \left( \pi s \right)}{\pi} \; \log \left( 3 \left( \frac{1}{4} + \mu \right) \right) \frac{1}{3^s} \, \left( \frac{1}{4}+\mu \right)^{-s} \; F_{\mu,0} \left( 2 \sqrt{\frac{1}{4}+\mu} \right).
                    \end{array}} \]
                
            \end{proposition}
            
            \begin{proof}
            
                The proof of this result directly stems from definition \ref{DefAmuk}.
            
            \end{proof}
            
            Recall that the aim of this paper is to get asymptotic expansions as $\mu$ goes to infinity for all $a > 0$, and as $a$ goes to infinity for $\mu = 0$. The next proposition deals with the second of these goals, as far as the terms $A_{\mu,k}$ are concerned.
            
            \begin{proposition}
            \label{PropAmukAsymptA}
            
                Let $\mu \geqslant 0$. The function
                
                \[ \begin{array}{ccc}
                        \displaystyle s & \longmapsto & \displaystyle \sum\limits_{\left \vert k \right \vert \geqslant 1} \; A_{\mu,k} \left( s \right)
                    \end{array} \]
            
                \noindent induces a holomorphic function on the half-plane $\Re s > 2 - 1/\left( 4 \delta \right)$ which contains $0$ if we have $\delta \, < \, 1/8$. On this half-plane, we can further differentiate term by term, and the derivative at $0$ for $\mu = 0$ satisfies, as $a$ goes to infinity,
            
                \[ \begin{array}{ccc}
                        \displaystyle \frac{\partial}{\partial s}_{\vert s \, = \, 0} \; \; \; \sum\limits_{\left \vert k \right \vert \geqslant 1} \; A_{0,k} \left( s \right) & = & \displaystyle O \left( \frac{1}{a^2} \right).
                    \end{array} \]
            
            \end{proposition}
            
            \begin{proof}
            
                For any $k \in \mathbb{Z} \setminus \left \{ 0 \right \}$ and any $\mu \geqslant 0$, proposition \ref{PropGlobalStudy} yields
                
                \[ \scalemath{0.97}{\begin{array}{lll}
                        \left \vert \left( \left \vert k \right \vert^{2 \delta} - \frac{1}{4} \right)^{-s} F_{\mu,k} \left( 2 \left \vert k \right \vert^{\delta} \sqrt{\frac{1}{4}+\mu} \right) \right \vert & \leqslant & \frac{C_{\mu}}{16 a^2 \left \vert k \right \vert^{2-4 \delta}} \left( \frac{1}{4} + \mu \right) \left( \left \vert k \right \vert^{2 \delta} - \frac{1}{4} \right)^{2 - \Re s}.
                    \end{array}} \]
            
                \noindent By the dominated convergence theorem, the sum of $A_{\mu,k}$ over non-zero integers is a holomorphic function on the half-plane $\Re s > 2 - 1 / \left( 4 \delta \right)$, and we can differentiate term by term. Evaluating the derivative of $A_{\mu,k}$ at $s \, = \, 0$ yields
            
                \[ \begin{array}{lll}
                        \frac{\partial}{\partial s}_{\vert s \, = \, 0} \; A_{\mu,k} \left( s \right) & = & F_{\mu,k} \left( 2 \left \vert k \right \vert^{\delta} \sqrt{\frac{1}{4}+\mu} \right)
                    \end{array} , \]
            
                \noindent and we can set $\mu \, = \, 0$, to get the estimate
            
                \[ \begin{array}{lll}
                        \displaystyle \left \vert \frac{\partial}{\partial s}_{\vert s \, = \, 0} \; A_{0,k} \left( s \right) \right \vert & \leqslant & \displaystyle \frac{1}{4} \; C_{0} \; \frac{1}{16 a^2} \cdot \frac{1}{\left \vert k \right \vert^{2 - 4 \delta}} \left( \left \vert k \right \vert^{2 \delta} - \frac{1}{4} \right)^2
                    \end{array} . \]
            
                \noindent This allows us to bound the derivative at $0$ of the series with general term $A_{0,k}$
            
                \[ \begin{array}{lllll}
                        \left \vert \frac{\partial}{\partial s}_{\vert s = 0} \; \sum\limits_{\left \vert k \right \vert \geqslant 1} \; A_{0,k} \left( s \right) \right \vert & \leqslant &  \frac{C_0}{4} \frac{1}{16 a^2} \; \sum\limits_{\left \vert k \right \vert \geqslant 1} \; \frac{1}{\left \vert k \right \vert^{2 - 4 \delta}} \left( \left \vert k \right \vert^{2 \delta} - \frac{1}{4} \right)^2
                    \end{array} . \]
            
                \noindent Since the series on the right hand side is absolutely convergent, we get
            
                \[ \begin{array}{lll}
                        \frac{\partial}{\partial s}_{\vert s = 0} \; \sum\limits_{\left \vert k \right \vert \geqslant 1} \; A_{0,k} \left( s \right) & = & O \left( \frac{1}{a^2} \right).
                    \end{array} \]
            
            \end{proof}
            
            As we did for $B_{\mu,k}$, we need to deal with the case $k = 0$ to complete the picture.
            
            \begin{proposition}
            \label{PropAmukAsymptA-0}
        
                For any $\mu \geqslant 0$, the derivative of the function $A_{0,0}$ satisfies
            
                \[ \begin{array}{lll}
                        \displaystyle \frac{\partial}{\partial s}_{\vert s = 0} \; A_{0,0} \left( s \right) & = & \displaystyle O \left( \frac{1}{a^2} \right)
                    \end{array} . \]
        
            \end{proposition}
            
            \begin{proof}
        
                The derivative of $A_{0,0}$ at $s=0$ being given by $F_{0,0} \left( 1 \right)$, we can prove the proposition by using remark \ref{RmkGlobalk0}.
        
        \end{proof}
        
        Having studied the $a$-asymptotics for $\mu = 0$, we turn our attention to the $\mu$-asymptotic behavior for all $a > 0$. We cannot proceed as in proposition \ref{PropAmukAsymptA}, since the upper-bound we use is not explicit in $\mu$. This complicates the study, and we need to split $A_{\mu,k}$, according to the asymptotic expansion given in corollary \ref{CorAsymptoticExpansionModifiedBessel}. From definitions \ref{DefAmuk} and \ref{DefFmuk}, we see that we need information on
        
        \begin{equation}
        \label{EqFmuk}
            \begin{array}{lll}
                F_{\mu,k} \left( t \right) & = & \log K_t \left( 2 \pi \left \vert k + \alpha \right \vert a \right) - \log K_{\sqrt{\frac{1}{4} + \mu}} \left( 2 \pi \left \vert k + \alpha \right \vert a \right) \\[1em]
                    
                && \qquad \qquad \qquad - \frac{t^2 - \left( 1/4 + \mu \right)}{\sqrt{4 \mu + 1}} \frac{\partial}{\partial t}_{\left \vert t = \sqrt{\frac{1}{4} + \mu} \right.} \log K_t \left( 2 \pi \left \vert k + \alpha \right \vert a \right).
            \end{array}
        \end{equation}
        
        \noindent for any integer $k \in \mathbb{Z}$ where, as always, the case $k = 0$ should only be considered if we have $\alpha \neq 0$. Let us state precisely the consequences of corollary \ref{CorAsymptoticExpansionModifiedBessel} we need.
        
        \begin{proposition}
        \label{PropLogK1}
        
            For every integer $k \neq 0$, and any real number $\mu \geqslant 0$, we have
            
            \[ \scalemath{0.97}{\begin{array}{lll}
                    \multicolumn{3}{l}{\log K_{2 \left \vert k \right \vert^{\delta} \sqrt{\frac{1}{4}+\mu}} \left( 2 \pi \left \vert k + \alpha \right \vert a \right)} \\[1em]
                    
                    & = & - \sqrt{\left( 2 \pi \left \vert k + \alpha \right \vert a \right)^2 + \left( 4 \mu + 1 \right) \left \vert k \right \vert^{2 \delta}} + \left \vert k \right \vert^{\delta} \sqrt{4 \mu + 1} \; \arcsinh \left( \frac{\left \vert k \right \vert^{\delta} \sqrt{4 \mu + 1}}{2 \pi \left \vert k + \alpha \right \vert a} \right) \\[1em]
                    
                    && \quad - \frac{1}{4} \log \left( \left( 2 \pi \left \vert k + \alpha \right \vert a \right)^2 + \left( 4 \mu + 1 \right) \left \vert k \right \vert^{2 \delta} \right) - \frac{1}{\left \vert k \right \vert^{\delta} \sqrt{4 \mu + 1}} U_1 \left( p \left( \frac{2 \pi \left \vert k + \alpha \right \vert a}{\left \vert k \right \vert^{\delta} \sqrt{4 \mu + 1}} \right) \right) \\[1em]
                    
                    && \quad \quad + \frac{1}{2} \log \left( \frac{\pi}{2} \right) + \widetilde{\eta_2} \left( \sqrt{4 \mu + 1} \left \vert k \right \vert^{\delta}, \; \; \frac{2 \pi \left \vert k + \alpha \right \vert a}{\left \vert k \right \vert^{\delta} \sqrt{4 \mu + 1}} \right),
                \end{array}} \]
            
            \noindent where the notations are made clear in appendix \ref{AppModifiedBessel}.
        
        \end{proposition}
        
        \begin{proof}
        
            This is a consequence of corollary \ref{CorAsymptoticExpansionModifiedBessel}.
        
        \end{proof}
        
        \begin{proposition}
        \label{PropLogK2}
        
            Assume we have $\alpha \neq 0$. For any real number $\mu \geqslant 0$, we have
            
            \[ \scalemath{0.97}{\begin{array}{lll}
                    \multicolumn{3}{l}{\log K_{2 \sqrt{\frac{1}{4}+\mu}} \left( 2 \pi \alpha a \right)} \\[0.5em]
                    
                    \qquad & = & \frac{1}{2} \log \left( \frac{\pi}{2} \right) - \sqrt{\left( 2 \pi \alpha a \right)^2 + \left( 4 \mu + 1 \right)} + \sqrt{4 \mu + 1} \; \arcsinh \left( \frac{\sqrt{4 \mu + 1}}{2 \pi \alpha a} \right) \\[0.5em]
                    
                    && \qquad - \frac{1}{4} \log \left( \left( 2 \pi \alpha a \right)^2 + \left( 4 \mu + 1 \right) \right) - \frac{1}{\sqrt{4 \mu + 1}} \; U_1 \left( p \left( \frac{2 \pi \alpha a}{\sqrt{4 \mu + 1}} \right) \right) \\[0.5em]
                    
                    && \qquad \qquad + \widetilde{\eta_2} \left( \sqrt{4 \mu + 1}, \frac{2 \pi \alpha a}{\sqrt{4 \mu + 1}} \right),
                \end{array}} \]
            
            \noindent where the notations are made clear in appendix \ref{AppModifiedBessel}.
        
        \end{proposition}
        
        \begin{proof}
        
            This is a consequence of corollary \ref{CorAsymptoticExpansionModifiedBessel}.
        
        \end{proof}
        
        \begin{proposition}
        \label{PropLogK3}
        
            For every integer $k$, with the exception of $k = 0$ should $\alpha$ vanish, and any real number $\mu \geqslant 0$, we have
            
            \[ \scalemath{0.97}{\begin{array}{lll}
                    \multicolumn{3}{l}{\log K_{\sqrt{\frac{1}{4}+\mu}} \left( 2 \pi \left \vert k + \alpha \right \vert a \right)} \\[1em]
                    
                    & = &  - \sqrt{\left( 2 \pi \left \vert k + \alpha \right \vert a \right)^2 + \frac{1}{4} + \mu} + \sqrt{\frac{1}{4}+\mu} \arcsinh \left( \frac{\sqrt{1/4+\mu}}{2 \pi \left \vert k + \alpha \right \vert a} \right) \\[1em]
                    
                    && \qquad - \frac{1}{4} \log \left( \left( 2 \pi \left \vert k + \alpha \right \vert a \right)^2 + \frac{1}{4} + \mu \right) - \frac{1}{\sqrt{1/4 + \mu}} U_1 \left( p \left( \frac{2 \pi \left \vert k + \alpha \right \vert a}{\sqrt{1/4+\mu}} \right) \right) \\[1em]
                    
                    && \qquad \qquad + \frac{1}{2} \log \left( \frac{\pi}{2} \right) + \widetilde{\eta_2} \left( \sqrt{\frac{1}{4}+\mu}, \frac{2 \pi \left \vert k + \alpha \right \vert a}{\sqrt{1/4+\mu}} \right),
                \end{array}} \]
            
            \noindent where the notations are made clear in appendix \ref{AppModifiedBessel}.
        
        \end{proposition}
        
        \begin{proof}
        
            This is a consequence of corollary \ref{CorAsymptoticExpansionModifiedBessel}.
        
        \end{proof}
        
        Having these expansions, we can study the series with general terms $A_{\mu,k}$. We will prove that the associated function has a holomorphic continuation to a region which contains the origin, and have a partial understanding of the asymptotic behavior of its derivative at $0$ as $\mu$ goes to infinity. The parts left uncomputed and unexplicit will be canceled in the overall study of the series with general terms $I_{\mu,k}$.

        \subsubsection*{\textbf{First part}}
        \label{AFirstPart}
        
            We begin the study by dealing with one of the remainder terms.
            
            \begin{proposition}
            \label{PropA1}
            
                The function
                
                \[ \begin{array}{lll}
                        s & \longmapsto & \frac{\sin \left( \pi s \right)}{\pi} \left( 4 \mu + 1 \right)^{-s} \sum\limits_{\left \vert k \right \vert \geqslant 1} \left( \left \vert k \right \vert^{2 \delta} - \frac{1}{4} \right)^{-s} \widetilde{\eta_2} \left( \left \vert k \right \vert^{\delta} \sqrt{4 \mu + 1}, \frac{2 \pi \left \vert k + \alpha \right \vert a}{\left \vert k \right \vert^{\delta} \sqrt{4 \mu + 1}} \right)
                    \end{array} \]
                
                \noindent is well-defined and holomorphic on the half-plane $\Re s > - 1 / \left( 2 \delta \right)$, and its derivative at $s=0$ satisfies, as $\mu$ goes to infinity,
                
                \[ \begin{array}{lll}
                        \scriptstyle \frac{\partial}{\partial s}_{\left \vert s=0 \right.} \left[ \frac{\sin \left( \pi s \right)}{\pi} \left( 4 \mu + 1 \right)^{-s} \sum\limits_{\left \vert k \right \vert \geqslant 1} \left( \left \vert k \right \vert^{2 \delta} - \frac{1}{4} \right)^{-s} \widetilde{\eta_2} \left( \left \vert k \right \vert^{\delta} \sqrt{4 \mu + 1}, \frac{2 \pi \left \vert k + \alpha \right \vert a}{\left \vert k \right \vert^{\delta} \sqrt{4 \mu + 1}} \right) \right] & \scriptstyle = & \scriptstyle o \left( 1 \right)
                    \end{array} . \]
            
            \end{proposition}
            
            \begin{proof}
            
                The key to this result is the bound given on $\widetilde{\eta_2}$ in corollary \ref{CorAsymptoticExpansionModifiedBessel}, whose notations will be used in this proof. For any non-zero integer $k$, we have
                
                \[ \begin{array}{lllclll}   
                        \nu & = & \left \vert k \right \vert^{\delta} \sqrt{4 \mu + 1} & \text{and} & x \nu & = & 2 \pi \left \vert k + \alpha \right \vert a
                    \end{array} . \]
                
                \noindent For $k$ with large enough absolute value, say with $\left \vert k \right \vert > K_0$ the hypotheses of corollary \ref{CorAsymptoticExpansionModifiedBessel} are satisfied, and for such integers, we have
                
                \[ \begin{array}{lll}
                        \left \vert \widetilde{\eta_2} \left( \left \vert k \right \vert^{\delta} \sqrt{4 \mu + 1}, \frac{2 \pi \left \vert k + \alpha \right \vert a}{\left \vert k \right \vert^{\delta} \sqrt{4 \mu + 1}} \right) \right \vert & \leqslant & \frac{C}{4 \pi^2 a^2} \cdot \frac{1}{\left( k + \alpha \right)^2}
                    \end{array} . \]
                
                \noindent The dominated convergence theorem proves that the function studied here is holomorphic on the half-plane $\Re s > - 1/\left( 2 \delta \right)$. Its derivative at $0$ is given by
                
                \[ \begin{array}{l}
                        \sum\limits_{\left \vert k \right \vert \geqslant 1} \widetilde{\eta_2} \left( \left \vert k \right \vert^{\delta} \sqrt{4 \mu + 1}, \frac{2 \pi \left \vert k + \alpha \right \vert a}{\left \vert k \right \vert^{\delta} \sqrt{4 \mu + 1}} \right)
                    \end{array} . \]
                
                \noindent To get the estimate on $\widetilde{\eta_2}$ for all non-zero integers, we will use corollary \ref{CorAsymptoticExpansionModifiedBessel} slightly differently. We have $\nu \geqslant \sqrt{4 \mu + 1}$, which means that the hypotheses of the corollary are satisfied for $\mu$ large enough, and all non-zero integers $k$. We then have
                
                \[ \begin{array}{lll}
                        \left \vert \sum\limits_{\left \vert k \right \vert \geqslant 1} \widetilde{\eta_2} \left( \left \vert k \right \vert^{\delta} \sqrt{4 \mu + 1}, \frac{2 \pi \left \vert k + \alpha \right \vert a}{\left \vert k \right \vert^{\delta} \sqrt{4 \mu + 1}} \right) \right \vert & \leqslant & \frac{C}{4 \pi^2 a^2} \sum\limits_{\left \vert k \right \vert \geqslant 1} \frac{1}{\left( k + \alpha \right)^2}
                    \end{array} . \]
                
                \noindent This allows us to use the dominated convergence theorem for the limit as $\mu$ goes to infinity. Using the second estimate provided by corollary \ref{CorAsymptoticExpansionModifiedBessel}, namely
                
                \[ \begin{array}{lll}
                        \left \vert \widetilde{\eta_2} \left( \left \vert k \right \vert^{\delta} \sqrt{4 \mu + 1}, \frac{2 \pi \left \vert k + \alpha \right \vert a}{\left \vert k \right \vert^{\delta} \sqrt{4 \mu + 1}} \right) \right \vert & \leqslant & C \left \vert k \right \vert^{- \delta} \cdot \frac{1}{\sqrt{4 \mu + 1}}
                    \end{array} , \]
                
                \noindent we see that the left-hand side of this inequality vanishes as $\mu$ goes to infinity, for all non-zero integers $k$. This concludes the proof.
            
            \end{proof}
            
            In order to complete this first part, let us take care of the case $k=0$.
            
            \begin{proposition}
            \label{PropA1-0}
            
                Assume we have $\alpha \neq 0$. The function
                
                \[ \begin{array}{lll}
                        s & \longmapsto & 3^{-s} \frac{\sin \left( \pi s \right)}{\pi} \left( \frac{1}{4} + \mu \right)^{-s} \widetilde{\eta_2} \left( \sqrt{4 \mu + 1}, \frac{2 \pi \alpha a}{\sqrt{4 \mu + 1}} \right)
                    \end{array} \]
                
                \noindent is entire, and its derivative at $0$ satisfies, as $\mu$ goes to infinity,
                
                \[ \begin{array}{lll}
                        \frac{\partial}{\partial s}_{\left \vert s = 0 \right.} \left[ 3^{-s} \frac{\sin \left( \pi s \right)}{\pi} \left( \frac{1}{4} + \mu \right)^{-s} \widetilde{\eta_2} \left( \sqrt{4 \mu + 1}, \frac{2 \pi \alpha a}{\sqrt{4 \mu + 1}} \right) \right] & = & o \left( 1 \right)
                    \end{array} . \]
            
            \end{proposition}
            
            \begin{proof}
            
                We begin by noting that the function studied in this proposition is entire, since there is no series involved. Then, we note that for $\mu$ large enough, the hypotheses of corollary \ref{CorAsymptoticExpansionModifiedBessel} are satisfied, and we conclude by noting we have
                
                \[ \begin{array}{lll}
                        \left \vert \widetilde{\eta_2} \left( \sqrt{4 \mu + 1}, \frac{2 \pi \alpha a}{\sqrt{4 \mu + 1}} \right) \right \vert & \leqslant & \frac{C}{4 \mu + 1}
                    \end{array} . \]
            
            \end{proof}

        \subsubsection*{\textbf{Second part}}
        \label{ASecondPart}
        
            We now move on to the other term involving a remainder $\widetilde{\eta_2}$.
            
            \begin{proposition}
            \label{PropA2}
            
                The function
                
                \[ \begin{array}{lll}
                        s & \longmapsto & - \frac{\sin \left( \pi s \right)}{\pi} \left( 4 \mu + 1 \right)^{-s} \sum\limits_{\left \vert k \right \vert \geqslant 1} \left( \left \vert k \right \vert^{2 \delta} - \frac{1}{4} \right)^{-s} \widetilde{\eta_2} \left( \sqrt{\frac{1}{4} + \mu}, \frac{2 \pi \left \vert k + \alpha \right \vert a}{\sqrt{1/4 + \mu}} \right)
                    \end{array} \]
                
                \noindent is well-defined and holomorphic on the half-plane $\Re s > - 1 / \left( 2 \delta \right)$, and its derivative at $s=0$ satisfies, as $\mu$ goes to infinity,
                
                \[ \begin{array}{lll}
                        \frac{\sin \left( \pi s \right)}{\pi} \left( 4 \mu + 1 \right)^{-s} \sum\limits_{\left \vert k \right \vert \geqslant 1} \left( \left \vert k \right \vert^{2 \delta} - \frac{1}{4} \right)^{-s} \widetilde{\eta_2} \left( \sqrt{\frac{1}{4} + \mu}, \frac{2 \pi \left \vert k + \alpha \right \vert a}{\sqrt{1/4 + \mu}} \right) & = & o \left( 1 \right)
                    \end{array} . \]
            
            \end{proposition}
            
            \begin{proof}
            
                This proof is similar to that of proposition \ref{PropA1}.
            
            \end{proof}
            
            Here again, we need to take care of the case $k = 0$, assuming $\alpha$ is not zero.
            
            \begin{proposition}
            \label{PropA2-0}
            
                Assume we have $\alpha \neq 0$. The function
                
                \[ \begin{array}{lll}
                        s & \longmapsto & - 3^{-s} \frac{\sin \left( \pi s \right)}{\pi} \left( \frac{1}{4} + \mu \right)^{-s} \widetilde{\eta_2} \left( \sqrt{\frac{1}{4} + \mu}, \frac{2 \pi \alpha a}{\sqrt{1/4 + \mu}} \right)
                    \end{array} \]
                
                \noindent is entire, and its derivative at $0$ satisfies, as $\mu$ goes to infinity,
                
                \[ \begin{array}{lll}
                        \frac{\partial}{\partial s}_{\left \vert s = 0 \right.} \left[ 3^{-s} \frac{\sin \left( \pi s \right)}{\pi} \left( \frac{1}{4} + \mu \right)^{-s} \widetilde{\eta_2} \left( \sqrt{\frac{1}{4} + \mu}, \frac{2 \pi \alpha a}{\sqrt{1/4 + \mu}} \right) \right] & = & o \left( 1 \right)
                    \end{array} . \]
            
            \end{proposition}
            
            \begin{proof}
            
                The proof is similar to that of proposition \ref{PropA1-0}.
            
            \end{proof}

        \subsubsection*{\textbf{Third part}}
        \label{AThirdPart}
        
            Having dealt with the ``remainder terms'', we come to a more complicated term, which cannot be fully studied. However, the part which will remain uncomputed will cancel another one later in this paper. In this third part, we will assume that we have $\delta < 1/2$, and that $1 / \left( 2 \delta \right)$ is not an integer.
            
            \begin{proposition}
            \label{PropA3}
            
                The function
                
                \[ \begin{array}{lll}
                        s & \longmapsto & - \frac{\sin \left( \pi s \right)}{\pi} \left( 4 \mu + 1 \right)^{-s} \sum\limits_{\left \vert k \right \vert \geqslant 1} \left( \left \vert k \right \vert^{2 \delta} - \frac{1}{4} \right)^{-s} \sqrt{\left( 2 \pi \left \vert k + \alpha \right \vert a \right)^2 + \left( 4 \mu + 1 \right) \left \vert k \right \vert^{2 \delta}}
                    \end{array} \]
                
                \noindent is well-defined and holomorphic on the half-plane $\Re s > 1 / \delta$. It has a holomorphic continuation to a region containing the origin, whose derivative there satisfies
                
                \[ \begin{array}{lll}
                        \multicolumn{3}{l}{\scriptstyle \frac{\partial}{\partial s}_{\left \vert s = 0 \right.} \left[ - \frac{\sin \left( \pi s \right)}{\pi} \left( 4 \mu + 1 \right)^{-s} \sum\limits_{\left \vert k \right \vert \geqslant 1} \left( \left \vert k \right \vert^{2 \delta} - \frac{1}{4} \right)^{-s} \sqrt{\left( 2 \pi \left \vert k + \alpha \right \vert a \right)^2 + \left( 4 \mu + 1 \right) \left \vert k \right \vert^{2 \delta}} \right]} \\[1em]
                        
                        \quad & \scriptstyle = & \scriptstyle - \frac{1}{4 \pi a \delta} \mu - \frac{1}{16 \pi a \delta} + \frac{\partial}{\partial s}_{\left \vert s = 0 \right.} \left[ - \frac{\sin \left( \pi s \right)}{\pi} \left( 4 \mu + 1 \right)^{-s} \sum\limits_{\left \vert k \right \vert \geqslant 1} \left \vert k \right \vert^{- 2 \delta s} \sqrt{\left( 2 \pi \left \vert k + \alpha \right \vert a \right)^2 + \left( 4 \mu + 1 \right) \left \vert k \right \vert^{2 \delta}} \right].
                    \end{array} \]
            
            \end{proposition}
            
            \begin{remark}
            
                Several points must be noted with regard to this last proposition.
                
                \begin{enumerate}
                    \item In the second part of the proposition, one should read ``the derivative at $0$ of the continuation of'' instead of ``derivative at $0$ of''. This kind of abuse of notation allows us to keep track of which part is being considered.
                    
                    \item The asymptotic expansion contains the derivative at $0$ of a holomorphic function which is not computed as $\mu$ goes to infinity. However, it will be canceled in the overall asymptotic study of $I_{\mu,k}$.
                \end{enumerate}
            
            \end{remark}
            
            \begin{proof}[Proof of proposition \ref{PropA3}]
            
                Note that the function we study is indeed holomorphic on the half-plane $\Re s > 1/\delta$. The binomial formula (see proposition \ref{PropBinomialFormula}) then gives
                
                \[ 
\begin{array}{lllll}
                        \left( \left \vert k \right \vert^{2 \delta} - \frac{1}{4} \right)^{-s} & = & \left \vert k \right \vert^{-2 \delta s} \left( 1 - \frac{1}{4} \left \vert k \right \vert^{-2\delta} \right)^{-s}  & = & \sum\limits_{j = 0}^{+ \infty} \frac{1}{j! 4^{j}} \left( s \right)_j \left \vert k \right \vert^{-2 \delta \left( s + j \right)}
                    \end{array} . \]
                
                \noindent We can plug this into the function we study, and interchange both sums, to get
                
                \[ \begin{array}{lll}
                        \multicolumn{3}{l}{\scriptstyle \sum\limits_{\left \vert k \right \vert \geqslant 1} \left( \left \vert k \right \vert^{2 \delta} - \frac{1}{4} \right)^{-s} \sqrt{\left( 2 \pi \left \vert k + \alpha \right \vert a \right)^2 + \left( 4 \mu + 1 \right) \left \vert k \right \vert^{2 \delta}}} \\[0.5em] 
                        
                        \qquad \qquad \qquad \qquad & \scriptstyle = & \scriptstyle \sum\limits_{j = 0}^{+ \infty} \frac{1}{j! 4^{j}} \left( s \right)_j \sum\limits_{\left \vert k \right \vert \geqslant 1} \left \vert k \right \vert^{-2 \delta \left( s + j \right)} \sqrt{\left( 2 \pi \left \vert k + \alpha \right \vert a \right)^2 + \left( 4 \mu + 1 \right) \left \vert k \right \vert^{2 \delta}}.
                    \end{array} \]
                
                \noindent We will deal with the terms $j=0$ and $j=1$ separately, they are the only ones who contribute to the proposition. For any integer $k > 0$, we have
                
                \[ \begin{array}{lllll}
                        \multicolumn{3}{l}{\scriptstyle \sqrt{\left( 2 \pi \left \vert k + \alpha \right \vert a \right)^2 + \left( 4 \mu + 1 \right) \left \vert k \right \vert^{2 \delta}}} & \scriptstyle = & \scriptstyle 2 \pi a k \sqrt{ \left( 1 + \frac{\alpha}{k} \right)^2 + \frac{4 \mu + 1}{4 \pi^2 a^2} \cdot \frac{1}{k^{2-2\delta}}} \\[1em]
                        
                        \quad & = & \multicolumn{3}{l}{\scriptstyle 2 \pi a k \left[ \sqrt{ \left( 1 + \frac{\alpha}{k} \right)^2 + \frac{4 \mu + 1}{4 \pi^2 a^2} \cdot \frac{1}{k^{2-2\delta}}} - \left( 1 + \frac{\alpha}{k} + \frac{4 \mu + 1}{8 \pi^2 a^2} \cdot \frac{1}{k^{2-2\delta}} \right) \right] + 2 \pi a k + 2 \pi a \alpha + \frac{4 \mu + 1}{4 \pi a} \cdot \frac{1}{k^{1-2\delta}}}.
                    \end{array} \]
                
                \noindent We further note that we have
                
                \[ \begin{array}{ll}
                        \multicolumn{2}{l}{\scriptstyle \left \vert \sqrt{ \left( 1 + \frac{\alpha}{k} \right)^2 + \frac{4 \mu + 1}{4 \pi^2 a^2} \cdot \frac{1}{k^{2-2\delta}}} - \left( 1 + \frac{\alpha}{k} + \frac{4 \mu + 1}{8 \pi^2 a^2} \cdot \frac{1}{k^{2-2\delta}} \right) \right \vert} \\[1em]
                        
                        \scriptstyle = & \scriptstyle \frac{1}{\sqrt{ \left( 1 + \frac{\alpha}{k} \right)^2 + \frac{4 \mu + 1}{4 \pi^2 a^2} \cdot \frac{1}{k^{2-2\delta}}} + 1 + \frac{\alpha}{k} + \frac{4 \mu + 1}{8 \pi^2 a^2} \cdot \frac{1}{k^{2-2\delta}}} \left \vert \left( 1 + \frac{\alpha}{k} \right)^2 + \frac{4 \mu + 1}{4 \pi^2 a^2} \cdot \frac{1}{k^{2-2\delta}} - \left( 1 + \frac{\alpha}{k} + \frac{4 \mu + 1}{8 \pi^2 a^2} \cdot \frac{1}{k^{2-2\delta}} \right)^2 \right \vert \\[1.5em]
                        
                        \scriptstyle \leqslant & \scriptstyle \frac{\left( 4 \mu + 1 \right)^2}{64 \pi^4 a^4} \cdot \frac{1}{k^{4-2 \delta}} + \frac{\left( 4 \mu + 1 \right) \alpha}{4 \pi^2 a2} \cdot \frac{1}{k^{3 - 2 \delta}},
                    \end{array} \]
                
                \noindent the last inequality being obtained after bounding the fraction by $1$, and computing the difference. This proves that the function
                
                \[ \begin{array}{lll}
                        s &\longmapsto & \sum\limits_{k \geqslant 1} k^{1 - 2 \delta s} \left[ \sqrt{ \left( 1 + \frac{\alpha}{k} \right)^2 + \frac{4 \mu + 1}{4 \pi^2 a^2} \cdot \frac{1}{k^{2-2\delta}}} - \left( 1 + \frac{\alpha}{k} + \frac{4 \mu + 1}{8 \pi^2 a^2} \cdot \frac{1}{k^{2-2\delta}} \right) \right]
                    \end{array} \]
                
                \noindent is holomorphic on the half-plane $\Re s > 1 - 1/ \left( 2 \delta \right)$, which contains the origin since we have $\delta < 1/2$. For any complex number $s$ with $\Re s > 1/\delta$, we have
                
                \[ \begin{array}{lll}
                        \multicolumn{3}{l}{\scriptstyle \sum\limits_{k \geqslant 1} k^{-2 \delta s} \sqrt{\left( 2 \pi \left \vert k + \alpha \right \vert a \right)^2 + \left( 4 \mu + 1 \right) \left \vert k \right \vert^{2 \delta}}} \\[0.8em]
                        
                        \qquad \qquad & \scriptstyle = & \scriptstyle 2 \pi a \sum\limits_{k \geqslant 1} k^{1 - 2 \delta s} \left[ \sqrt{ \left( 1 + \frac{\alpha}{k} \right)^2 + \frac{4 \mu + 1}{4 \pi^2 a^2} \cdot \frac{1}{k^{2-2\delta}}} - \left( 1 + \frac{\alpha}{k} + \frac{4 \mu + 1}{8 \pi^2 a^2} \cdot \frac{1}{k^{2-2\delta}} \right) \right] \\[0.8em]
                        
                        && \scriptstyle \qquad \qquad \qquad \qquad + 2 \pi a \zeta \left( 2 \delta s - 1 \right) + 2 \pi a \alpha \zeta \left( 2 \delta s \right) + \frac{4 \mu + 1}{4 \pi a} \zeta \left( 1 + 2 \delta \left( s - 1 \right) \right),
                    \end{array} \]
                
                \noindent and the associated function has a holomorphic continuation near the origin. The sum over integers $k \leqslant -1$ can be dealt with similarly. After multiplication by the relevant factor, the derivative at $s=0$ is the one left uncomputed in the statement of the proposition. Let us move on to the term corresponding to $j=1$, which is
                
                \[ \begin{array}{lll}
                        - s \frac{\sin \left(\pi s \right)}{4\pi} \left( 4 \mu + 1 \right)^{-s} \sum\limits_{\left \vert k \right \vert \geqslant 1} \left\vert k \right \vert^{-2 \delta \left( s+1 \right)} \sqrt{\left( 2 \pi \left \vert k + \alpha \right \vert a \right)^2 + \left( 4 \mu + 1 \right) \left \vert k \right \vert^{2 \delta}}
                    \end{array} . \]
                
                \noindent Using the same method as above, we have, on the half-plane $\Re s > 1/\delta - 1$,
                
                \[ \begin{array}{lll}
                        \multicolumn{3}{l}{ \scriptstyle - s \frac{\sin \left(\pi s \right)}{4\pi} \left( 4 \mu + 1 \right)^{-s} \sum\limits_{k \geqslant 1} k^{-2 \delta \left( s+1 \right)} \sqrt{\left( 2 \pi \left \vert k + \alpha \right \vert a \right)^2 + \left( 4 \mu + 1 \right) \left \vert k \right \vert^{2 \delta}}} \\[1em]
                        
                        & \scriptstyle = & \scriptstyle - \frac{1}{2} a s \sin \left( \pi s \right) \left( 4 \mu + 1 \right)^{-s} \sum\limits_{k \geqslant 1} k^{1 - 2 \delta \left( s + 1 \right)} \left( \sqrt{ \left( 1 + \frac{\alpha}{k} \right)^2 + \frac{4 \mu + 1}{4 \pi^2 a^2} \cdot \frac{1}{k^{2-2\delta}}} - \left( 1 + \frac{\alpha}{k} + \frac{4 \mu + 1}{8 \pi^2 a^2} \cdot \frac{1}{k^{2-2\delta}} \right) \right) \\[1em]
                        
                        && \qquad \scriptstyle - \frac{1}{2} a s \sin \left( \pi s \right) \left( 4 \mu + 1 \right)^{-s} \left[ \zeta \left( 1 + 2 \delta \left( s+1 \right) \right) + \alpha \zeta \left( 2 \delta \left( s+1 \right) \right) + \frac{4\mu + 1}{8 \pi^2 a^2} \zeta \left( 1 + 2 \delta s \right) \right].
                    \end{array} \]
                
                \noindent Previously found estimates show that the first part on the right-hand side induces a holomorphic function on the half-plane $\Re s > - 1 / \left( 2 \delta \right)$, whose derivative at $0$ vanishes, due to the factor $s \sin \left( \pi s \right)$. The Riemann zeta function having a unique pole, which is located at $s=1$ and simple, we note that the function
                
                \[ \begin{array}{lll}
                        s & \longmapsto & - \frac{1}{2} a s \sin \left( \pi s \right) \left( 4 \mu + 1 \right)^{-s} \left[ \zeta \left( 1 + 2 \delta \left( s+1 \right) \right) + \alpha \zeta \left( 2 \delta \left( s+1 \right) \right) \right]
                    \end{array} \]
                
                \noindent is entire, and that its derivative at the origin vanishes. Similarly, the function
                
                \[ \begin{array}{lll}
                        s & \longmapsto & - \frac{1}{2} a s \sin \left( \pi s \right) \left( 4 \mu + 1 \right)^{-s} \cdot \frac{4 \mu + 1}{8 \pi^2 a ^2} \zeta \left( 1 + 2 \delta s \right)
                    \end{array} \]
                
                \noindent is entire, but its derivative at $0$ does not vanish. We have
                
                \[ \begin{array}{lll}
                        \frac{\partial}{\partial s}_{\left \vert s = 0 \right.} \left[ - \frac{1}{2} a s \sin \left( \pi s \right) \left( 4 \mu + 1 \right)^{-s} \cdot \frac{4 \mu + 1}{8 \pi^2 a ^2} \zeta \left( 1 + 2 \delta s \right) \right] & = & - \frac{1}{8 \pi a \delta} \mu - \frac{1}{32 \pi a \delta}
                    \end{array} . \]
                
                \noindent One obtains similar results for the sum bearing over integers $k \leqslant -1$. It only remains to study the remaining sum over integers $j \geqslant 2$, given by
                
                \[ \begin{array}{l}
                        - \frac{\sin \left( \pi s \right)}{\pi} \left( 4 \mu + 1 \right)^{-s} \sum\limits_{j = 2}^{+ \infty} \frac{\left( s \right)_j}{j! 4^{j}} \sum\limits_{\left \vert k \right \vert \geqslant 1} \left \vert k \right \vert^{-2 \delta \left( s + j \right)} \sqrt{\left( 2 \pi \left \vert k + \alpha \right \vert a \right)^2 + \left( 4 \mu + 1 \right) \left \vert k \right \vert^{2 \delta}}
                    \end{array} . \]
                
                \noindent We only study the sum over integers $k > 0$, the other one being similar. We have
                
                \[ \begin{array}{lll}
                        \multicolumn{3}{l}{\scriptstyle \sum\limits_{j = 2}^{+ \infty} \frac{\left( s \right)_j}{j! 4^{j}} \sum\limits_{k \geqslant 1} k^{-2 \delta \left( s + j \right)} \sqrt{\left( 2 \pi \left \vert k + \alpha \right \vert a \right)^2 + \left( 4 \mu + 1 \right) k^{2 \delta}}} \\[1em]
                        
                        \qquad & = & \scriptstyle 2 \pi a \sum\limits_{j = 2}^{+ \infty} \frac{\left( s \right)_j}{j! 4^{j}} \sum\limits_{k \geqslant 1} k^{1 - 2 \delta \left( s + j \right)} \left[ \sqrt{ \left( 1 + \frac{\alpha}{k} \right)^2 + \frac{4 \mu + 1}{4 \pi^2 a^2} \cdot \frac{1}{k^{2-2\delta}}} - \left( 1 + \frac{\alpha}{k} + \frac{4 \mu + 1}{8 \pi^2 a^2} \cdot \frac{1}{k^{2-2\delta}} \right) \right] \\[1em]
                        
                        && \qquad \scriptstyle + 2 \pi a \sum\limits_{j = 2}^{+ \infty} \frac{\left( s \right)_j}{j! 4^{j}} \left[ \zeta \left( - 1 + 2 \delta \left( s+j \right) \right) + \alpha \zeta \left( 2 \delta \left( s + j \right) \right) + \frac{4 \mu + 1}{8 \pi^2 a^2} \zeta \left( 1 + 2 \delta \left( s + j - 1 \right) \right) \right].
                    \end{array} \]
                
                \noindent The first series on the right-hand side is holomorphic on $\Re s > -1 + 1/ \left( 2 \delta \right)$, and, after multiplication by the appropriate factor, the derivative at $s=0$ vanishes, because of the Pochhammer symbols $\left( s \right)_j$. For the second part, the function
                
                \[ \begin{array}{lll}
                        \scriptstyle s & \scriptstyle \longmapsto & \scriptstyle 2 \pi a \sum\limits_{j = 2}^{+ \infty} \frac{\left( s \right)_j}{j! 4^{j}} \left[ \zeta \left( - 1 + 2 \delta \left( s+j \right) \right) + \alpha \zeta \left( 2 \delta \left( s + j \right) \right) + \frac{4 \mu + 1}{8 \pi^2 a^2} \zeta \left( 1 + 2 \delta \left( s + j - 1 \right) \right) \right]
                    \end{array} , \]
                
                \noindent is holomorphic around $s=0$, since $1/\left( 2 \delta \right)$ is not an integer. After multiplication by the relevant factor, the derivative at $s=0$ vanishes, because of the Pochhammer symbols. The sum over $k \leqslant -1$ is taken care of similarly, thus completing the proof.
            
            \end{proof}
            
            As always, we need the corresponding result for $k=0$, when $\alpha$ does not vanish.
            
            \begin{proposition}
            \label{PropA3-0}
            
                Assume that $\alpha$ does not vanish. The function
                
                \[ \begin{array}{lll}
                        s & \longmapsto & - 3^{-s} \frac{\sin \left(\pi s \right)}{\pi} \left( \frac{1}{4} + \mu \right)^{-s} \sqrt{4 \pi^2 \alpha^2 a^2 + 4 \mu + 1}
                    \end{array} \]
                
                \noindent is entire, and its derivative at $s=0$ satisfies, as $\mu$ goes to infinity,
                
                \[ \begin{array}{lll}
                        \frac{\partial}{\partial s}_{\left \vert s = 0 \right.} \left[ - 3^{-s} \frac{\sin \left(\pi s \right)}{\pi} \left( \frac{1}{4} + \mu \right)^{-s} \sqrt{4 \pi^2 \alpha^2 a^2 + 4 \mu + 1} \right] & = & - 2 \sqrt{\mu} + o \left( 1 \right)
                    \end{array} . \]
            
            \end{proposition}
            
            \begin{proof}
            
                This is a direct computation.
            
            \end{proof}

        \subsubsection*{\textbf{Fourth part}}
        \label{AFourthPart}
            
            We continue our study with the part containing the inverse hyperbolic function $\arcsinh$. We assume in this fourth part that we have $\delta < 1/2$.
            
            \begin{proposition}
            \label{PropA4}
            
                The function
                
                \[ \begin{array}{lll}
                        s & \longmapsto & \frac{\sin \left(\pi s \right)}{\pi} \left( 4 \mu + 1 \right)^{-s+1/2} \sum\limits_{\left \vert k \right \vert \geqslant 1} \left( \left \vert k \right \vert^{2 \delta} - \frac{1}{4} \right)^{-s} \left \vert k \right \vert^{\delta} \arcsinh \left( \frac{\left \vert k \right \vert^{\delta} \sqrt{4 \mu + 1}}{2 \pi \left \vert k + \alpha \right \vert a} \right)
                    \end{array} , \]
                
                \noindent is well-defined and holomorphic on the half-plane $\Re s > 1$, has a holomorphic continuation to a region which contains $0$, and its derivative there satisfies
                
                \[ \begin{array}{lll}
                        \multicolumn{3}{l}{\scriptstyle \frac{\partial}{\partial s}_{\left \vert s = 0 \right.} \big[ \frac{\sin \left(\pi s \right)}{\pi} \left( 4 \mu + 1 \right)^{-s+1/2} \sum\limits_{\left \vert k \right \vert \geqslant 1} \left( \left \vert k \right \vert^{2 \delta} - \frac{1}{4} \right)^{-s} \left \vert k \right \vert^{\delta} \arcsinh \left( \frac{\left \vert k \right \vert^{\delta} \sqrt{4 \mu + 1}}{2 \pi \left \vert k + \alpha \right \vert a} \right) \big]} \\[1em]
                        
                        \qquad & \scriptstyle = & \scriptstyle \frac{1}{2 \pi a \delta} \mu + \frac{1}{8 \pi a \delta} + \frac{\partial}{\partial s}_{\left \vert s = 0 \right.} \big[ \frac{\sin \left(\pi s \right)}{\pi} \left( 4 \mu + 1 \right)^{-s+1/2} \sum\limits_{\left \vert k \right \vert \geqslant 1} \left \vert k \right \vert^{\delta - 2 \delta s} \left \vert k \right \vert^{\delta} \arcsinh \left( \frac{\left \vert k \right \vert^{\delta} \sqrt{4 \mu + 1}}{2 \pi \left \vert k + \alpha \right \vert a} \right) \big].
                    \end{array} \]
            
            \end{proposition}
            
            \begin{proof}
            
                Using Taylor's formula, we have, for every non-zero integer
                
                \[ \begin{array}{lll}
                        \scriptstyle \arcsinh \left( \frac{\left \vert k \right \vert^{\delta} \sqrt{4 \mu + 1}}{2 \pi \left \vert k + \alpha \right \vert a} \right) & \scriptsize = & \scriptsize \frac{\left \vert k \right \vert^{\delta} \sqrt{4 \mu + 1}}{2 \pi \left \vert k + \alpha \right \vert a} - \int_0^{\frac{\left \vert k \right \vert^{\delta} \sqrt{4 \mu + 1}}{2 \pi \left \vert k + \alpha \right \vert a}} \frac{x}{\left( 1 + x^2 \right)^{3/2}} \left( \frac{\left \vert k \right \vert^{\delta} \sqrt{4 \mu + 1}}{2 \pi \left \vert k + \alpha \right \vert a} - x \right) \mathrm{d}x
                    \end{array} . \]
                
                \noindent We further have the estimate
                
                \[ \begin{array}{lll}
                        \scriptstyle \left \vert \int_0^{\frac{\left \vert k \right \vert^{\delta} \sqrt{4 \mu + 1}}{2 \pi \left \vert k + \alpha \right \vert a}} \frac{x}{\left( 1 + x^2 \right)^{3/2}} \left( \frac{\left \vert k \right \vert^{\delta} \sqrt{4 \mu + 1}}{2 \pi \left \vert k + \alpha \right \vert a} - x \right) \mathrm{d}x \right \vert & \scriptstyle \leqslant & \scriptstyle \frac{1}{2} \left( \frac{\left \vert k \right \vert^{\delta} \sqrt{4 \mu + 1}}{2 \pi \left \vert k + \alpha \right \vert a} \right)^3
                    \end{array} . \]
                
                \noindent This proves that the function studied in this proposition is indeed holomorphic on the half-plane $\Re s > 1$. The binomial formula then gives
                
                \[ \begin{array}{lll}
                        \scriptstyle \sum\limits_{\left \vert k \right \vert \geqslant 1} \left( \left \vert k \right \vert^{2 \delta} - \frac{1}{4} \right)^{-s} \left \vert k \right \vert^{\delta} \arcsinh \left( \frac{\left \vert k \right \vert^{\delta} \sqrt{4 \mu + 1}}{2 \pi \left \vert k + \alpha \right \vert a} \right) & \hspace{-1pt} \scriptstyle = \hspace{-1pt} & \scriptstyle \sum\limits_{j = 0}^{+ \infty} \frac{\left( s \right)_j}{j! 4^{j}} \sum\limits_{\left \vert k \right \vert \geqslant 1}  \left \vert k \right \vert^{-2 \delta \left( s + j \right) + \delta} \arcsinh \left( \frac{\left \vert k \right \vert^{\delta} \sqrt{4 \mu + 1}}{2 \pi \left \vert k + \alpha \right \vert a} \right).
                    \end{array} \]
                
                \noindent We will deal with the terms corresponding to $j=0$ and $j=1$ separately. We have
                
                \[ \begin{array}{ll}
                        \multicolumn{2}{l}{\scriptstyle \sum\limits_{\left \vert k \right \vert \geqslant 1}  \left \vert k \right \vert^{-2 \delta s + \delta} \arcsinh \left( \frac{\left \vert k \right \vert^{\delta} \sqrt{4 \mu + 1}}{2 \pi \left \vert k + \alpha \right \vert a} \right)} \\[0.6em]
                        
                        \scriptstyle = & \scriptstyle \frac{1}{2 \pi a} \sqrt{4 \mu + 1} \sum\limits_{\left \vert k \right \vert \geqslant 1} \frac{1}{\left \vert k + \alpha \right \vert} \left \vert k \right \vert^{-2 \delta s + 2 \delta} - \sum\limits_{\left \vert k \right \vert \geqslant 1}  \left \vert k \right \vert^{-2 \delta s + \delta} \int_0^{\frac{\left \vert k \right \vert^{\delta} \sqrt{4 \mu + 1}}{2 \pi \left \vert k + \alpha \right \vert a}} \frac{x}{\left( 1 + x^2 \right)^{3/2}} \left( \frac{\left \vert k \right \vert^{\delta} \sqrt{4 \mu + 1}}{2 \pi \left \vert k + \alpha \right \vert a} - x \right) \mathrm{d}x.
                    \end{array} \]
                
                \noindent The estimate already obtained on the integral proves that the second term on the right-hand side induces a holomorphic function on the half-plane $\Re s > 2 - 1/\delta$, which contains $0$ since we have $\delta < 1/2$. Let us study the first term. We have
                
                \[ \begin{array}{lll}
                        \sum\limits_{\left \vert k \right \vert \geqslant 1} \frac{1}{\left \vert k + \alpha \right \vert} \left \vert k \right \vert^{-2 \delta s + 2 \delta} & = & 2 \zeta \left( 1 - 2 \delta + 2 \delta s \right) + 2 \sum\limits_{k \geqslant 1} \left[ \frac{1}{1 - \frac{\alpha^2}{k^2}} - 1 \right] k^{-1+2\delta - 2 \delta s}.
                    \end{array} \]
                
                \noindent Using the fact that we have
                
                \[ \begin{array}{lll}
                        \left \vert \left[ \frac{1}{1 - \frac{\alpha^2}{k^2}} - 1 \right] k^{- 2 \delta s + 2 \delta} \right \vert & \leqslant & \frac{\alpha^2}{1 - \alpha^2} k^{- 3 + 2 \delta - 2 \delta \Re s}
                    \end{array} , \]
                
                \noindent and the continuation of the Riemann zeta function, the term we study has a holomorphic continuation near the origin. Its derivative there, after multiplication by the relevant factor, is uncomputed in the proposition. Similarly, we can take care of the term $j=1$, and find its contribution. The term corresponding to $j \geqslant 2$ is
                
                \[ \begin{array}{l}
                        \frac{\sin \left(\pi s \right)}{\pi} \left( 4 \mu + 1 \right)^{-s+1/2} \sum\limits_{j = 2}^{+ \infty} \frac{\left( s \right)_j}{j! 4^{j}} \sum\limits_{\left \vert k \right \vert \geqslant 1}  \left \vert k \right \vert^{-2 \delta \left( s + j \right) + \delta} \arcsinh \left( \frac{\left \vert k \right \vert^{\delta} \sqrt{4 \mu + 1}}{2 \pi \left \vert k + \alpha \right \vert a} \right)
                    \end{array} , \]
                
                \noindent can be dealt with in a similar fashion . The derivative at $s=0$ of its continuation vanishes, because of the Pochhammer symbols $\left( s \right)_j$. This completes the proof.
            
            \end{proof}
            
            \begin{proposition}
            \label{PropA4-0}
            
                Assume we have $\alpha \neq 0$. The function
                
                \[ \begin{array}{lll}
                        s & \longmapsto & 2 \cdot 3^{-s} \frac{\sin \left( \pi s \right)}{\pi} \left( \frac{1}{4} + \mu \right)^{-s+1/2} \arcsinh \left( \frac{\sqrt{4 \mu + 1}}{2 \pi \alpha a} \right)
                    \end{array} \]
                
                \noindent is entire, and its derivative at $0$ satisfies, as $\mu$ goes to infinity,
                
                \[ \begin{array}{lll}
                        \scriptstyle \frac{\partial}{\partial s}_{\left \vert s = 0 \right.} \left[ 2 \cdot 3^{-s} \frac{\sin \left( \pi s \right)}{\pi} \left( \frac{1}{4} + \mu \right)^{-s+1/2} \arcsinh \left( \frac{\sqrt{4 \mu + 1}}{2 \pi \alpha a} \right) \right] & \scriptstyle = & \scriptstyle \sqrt{\mu} \log \mu - 2 \log \left( \frac{\pi}{2} \alpha a \right) \sqrt{\mu} + o \left( 1 \right)
                    \end{array} . \]
            
            \end{proposition}
            
            \begin{proof}
            
                This is a direct computation.
            
            \end{proof}

        \subsubsection*{\textbf{Fifth part}}
        \label{AFifthPart}
            
            We now turn our attention to the study of the logarithmic term from proposition \ref{PropLogK1}. Once again, we assume that we have $\delta < 1/2$.
            
            \begin{proposition}
            \label{PropA5}
            
                The function
                
                \[ \begin{array}{lll}
                        s & \hspace{-4pt} \mapsto & \hspace{-4pt} - \frac{\sin \left(\pi s \right)}{4 \pi} \left( 4 \mu + 1 \right)^{-s} \hspace{-1pt} \sum\limits_{\left \vert k \right \vert \geqslant 1} \left( \left \vert k \right \vert^{2 \delta} - \frac{1}{4} \right)^{-s} \hspace{-1pt} \log \left( \left( 2 \pi \left \vert k + \alpha \right \vert a \right)^2 \hspace{-1pt} + \left( 4 \mu + 1 \right) \left \vert k \right \vert^{2 \delta} \right)
                    \end{array} \]
                
                \noindent is holomorphic on the half-plane $\Re s > 1/\left( 2 \delta \right)$, has a holomorphic continuation to a region containing $0$, and its derivative there satisfies
                
                \[ \begin{array}{lll}
                        \multicolumn{3}{l}{\scriptstyle \frac{\partial}{\partial s}_{\left \vert s = 0 \right.}  \left[ - \frac{\sin \left(\pi s \right)}{4 \pi} \left( 4 \mu + 1 \right)^{-s} \sum\limits_{\left \vert k \right \vert \geqslant 1} \left( \left \vert k \right \vert^{2 \delta} - \frac{1}{4} \right)^{-s} \log \left( \left( 2 \pi \left \vert k + \alpha \right \vert a \right)^2 + \left( 4 \mu + 1 \right) \left \vert k \right \vert^{2 \delta} \right) \right]} \\[1em]
                        
                        \qquad \qquad \qquad & \scriptstyle = & \scriptstyle \frac{\partial}{\partial s}_{\left \vert s = 0 \right.}  \left[ - \frac{\sin \left(\pi s \right)}{4 \pi} \left( 4 \mu + 1 \right)^{-s} \sum\limits_{\left \vert k \right \vert \geqslant 1} \left \vert k \right \vert^{-2 \delta s} \log \left( \left( 2 \pi \left \vert k + \alpha \right \vert a \right)^2 + \left( 4 \mu + 1 \right) \left \vert k \right \vert^{2 \delta} \right) \right]
                    \end{array} \]
            
            \end{proposition}
            
            \begin{proof}
            
                The argument is essentially the same as in propositions \ref{PropA3} and \ref{PropA4}.
            
            \end{proof}
            
            The corresponding term for $k = 0$ is stated below.
            
            \begin{proposition}
            \label{PropA5-0}
            
                Assume we have $\alpha \neq 0$. The function
                
                \[ \begin{array}{lll}
                        s & \longmapsto & - \frac{1}{4} 3^{-s} \frac{\sin \left( \pi s \right)}{\pi} \left( \frac{1}{4} + \mu \right)^{-s} \log \left( \left( 2 \pi \alpha a \right)^2 + 4 \mu + 1 \right)
                    \end{array} \]
                
                \noindent is entire, and its derivative at $s=0$ satisfies, as $\mu$ goes to infinity,
                
                \[ \begin{array}{lll}
                        \scriptstyle \frac{\partial}{\partial}_{\left \vert s = 0 \right.} \left[ - \frac{1}{4} 3^{-s} \frac{\sin \left( \pi s \right)}{\pi} \left( \frac{1}{4} + \mu \right)^{-s} \log \left( \left( 2 \pi \alpha a \right)^2 + 4 \mu + 1 \right) \right] & \scriptstyle = & \scriptstyle - \frac{1}{4} \log \mu - \frac{1}{2} \log 2 + o \left( 1 \right)
                    \end{array} . \]
            
            \end{proposition}
            
            \begin{proof}
            
                This is a direct computation.            
            
            \end{proof}

        \subsubsection*{\textbf{Sixth part}}
        \label{ASixthPart}
        
            We now come to the last term from proposition \ref{PropLogK1}, which is the one involving the polynomial $U_1 \left( t \right) = \left( 3t - 5t^3 \right)/24$. For simplicity, we will split it into two monomial parts corresponding to $t$ and $t^3$. The coefficients will be taken into account in a subsequent proposition. We assume that we have $\delta < 1$.
            
            \begin{proposition}
            \label{PropA6Power1}
            
                The function
                
                \[ \begin{array}{lll}
                        s & \longmapsto & - \frac{\sin \left( \pi s \right)}{\pi} \left( 4 \mu + 1 \right)^{-s-1/2} \sum\limits_{\left \vert k \right \vert \geqslant 1} \left( \left \vert k \right \vert^{2 \delta} - \frac{1}{4} \right)^{-s} \left \vert k \right \vert^{- \delta} p \left( \frac{2 \pi \left \vert k + \alpha \right \vert a}{\left \vert k \right \vert^{\delta} \sqrt{4 \mu + 1}} \right)
                    \end{array} \]
                    
                \noindent which is well-defined on the half-plane $\Re s > 0$, has a holomorphic continuation to a region of the complex plane containing $0$.
            
            \end{proposition}
            
            \begin{proof}
            
                Recall that definition \ref{Defpxi} gives $p \left( x \right) = \left( 1+x^2 \right)^{-1/2}$ for all $x>0$. For any integer $k \geqslant 1$, we have
                
                \[ \begin{array}{llll}
                        \multicolumn{2}{l}{p \left( \frac{2 \pi \left \vert k + \alpha \right \vert a}{\left \vert k \right \vert^{\delta} \sqrt{4 \mu + 1}} \right)} = & \frac{k^{\delta} \sqrt{4 \mu + 1}}{2 \pi \left( k + \alpha \right) a} + \frac{k^{\delta} \sqrt{4 \mu + 1}}{2 \pi \left( k + \alpha \right) a} \left[ \left( 1 + \frac{\left( 4 \mu + 1 \right) k^{2 \delta}}{\left( 2 \pi \left( k + \alpha \right) a \right)^2} \right)^{-1/2} - 1 \right] \\[1em]
                        
                        = & \multicolumn{3}{l}{\frac{k^{\delta} \sqrt{4 \mu + 1}}{2 \pi \left( k + \alpha \right) a} - \frac{\left( 4 \mu + 1 \right) k^{2 \delta}}{\left( 2 \pi \left( k + \alpha \right) a \right)^2} \underbrace{\scriptstyle \left( 1 + \frac{\left( 4 \mu + 1 \right) k^{2 \delta}}{\left( 2 \pi \left( k + \alpha \right) a \right)^2} \right)^{-1/2} \left( 1 + \left( 1 + \frac{\left( 4 \mu + 1 \right) k^{2 \delta}}{\left( 2 \pi \left( k + \alpha \right) a \right)^2} \right)^{1/2} \right)^{-1}}_{\leqslant 1}}.
                    \end{array} \]
                
                \noindent In particular, the function
                
                \[ \begin{array}{lll}
                        s & \longmapsto & \sum\limits_{k \geqslant 1} \left( k^{2 \delta} - \frac{1}{4} \right)^{-s} k^{- \delta} \frac{k^{\delta} \sqrt{4 \mu + 1}}{2 \pi \left( k + \alpha \right) a} \left[ \left( 1 + \frac{\left( 4 \mu + 1 \right) k^{2 \delta}}{\left( 2 \pi \left( k + \alpha \right) a \right)^2} \right)^{-1/2} - 1 \right]
                    \end{array} \]
                
                \noindent is holomorphic on the half-plane $\Re s > \left( \delta - 1 \right) / \left( 2 \delta \right)$, which contains $0$, having assumed that we have $\delta < 1$. Let us now study the function
                
                \[ \begin{array}{lll}
                        s & \longmapsto & - \frac{\sin \left( \pi s \right)}{\pi} \left( 4 \mu + 1 \right)^{-s-1/2} \sum\limits_{k \geqslant 1} \left( k^{2 \delta} - \frac{1}{4} \right)^{-s} k^{-\delta} \cdot \frac{k^{\delta} \sqrt{4 \mu + 1}}{2 \pi \left( k + \alpha \right) a}
                    \end{array} , \]
                
                \noindent which is well-defined and holomorphic and the half-plane $\Re s > 0$.  Using the binomial formula, which is the content of proposition \ref{PropBinomialFormula}, we have
                
                \[ \begin{array}{lll}
                        \multicolumn{3}{l}{- \frac{\sin \left( \pi s \right)}{\pi} \left( 4 \mu + 1 \right)^{-s-1/2} \sum\limits_{k \geqslant 1} \left( k^{2 \delta} - \frac{1}{4} \right)^{-s} k^{-\delta} \cdot \frac{k^{\delta} \sqrt{4 \mu + 1}}{2 \pi \left( k + \alpha \right) a}} \\[1em]
                        
                        \qquad \qquad \qquad \qquad \qquad & = & - \frac{1}{2 \pi a} \frac{\sin \left( \pi s \right)}{\pi} \left( 4 \mu + 1 \right)^{-s} \sum\limits_{j = 0}^{+ \infty} \frac{\left( s \right)_j}{j! 4^j} \sum\limits_{k \geqslant 1} k^{-2 \delta \left( s + j \right)} \frac{1}{k + \alpha}.
                    \end{array} \]
                
                \noindent The sum over integers $j \geqslant 1$ inducing a holomorphic function around $0$, let us focus on the term corresponding to $j=0$, which is given by
                
                \[ \begin{array}{lll}
                        \multicolumn{3}{l}{- \frac{1}{2 \pi a} \frac{\sin \left( \pi s \right)}{\pi} \left( 4 \mu + 1 \right)^{-s} \sum\limits_{k \geqslant 1} k^{-2 \delta s} \frac{1}{k + \alpha}} \\[1em]
                        
                        \quad & = & - \frac{1}{2 \pi a} \frac{\sin \left( \pi s \right)}{\pi} \left( 4 \mu + 1 \right)^{-s} \Big( \zeta \left( 2 \delta s + 1 \right) + \sum\limits_{k \geqslant 1} k^{-2 \delta s + 1} \left( \left( 1 + \frac{\alpha}{k} \right)^{-1} - 1 \right) \Big).
                    \end{array} \]
                
                \noindent The second term above induces a holomorphic function around $0$, as we have
                
                \[ \begin{array}{lllll}
                        \left \vert \left( 1 + \frac{\alpha}{k} \right)^{-1} - 1 \right \vert & = & \frac{1}{k} \left( 1 + \frac{\alpha}{k} \right)^{-1} & \leqslant & \frac{1}{k}
                    \end{array} . \]
                
                 \noindent In the first term, the simple pole of $\zeta \left( 2 \delta s + 1 \right)$ is canceled by the factor $\sin \left( \pi s \right)$, and the result is holomorphic around the origin. Even though we only dealt with the sum over integers $k \geqslant 1$, only small modifications are required for the sum over integers $k \leqslant -1$, which is seen by switching the sign of $\alpha$. This concludes the proof.

            \end{proof}
            
            \begin{proposition}
            \label{PropA6Power1-0}
            
                Assume we have $\alpha \neq 0$. The function
                
                \[ \begin{array}{lll}
                        s & \longmapsto & - 3^{-s} \frac{\sin \left( \pi s \right)}{2 \pi} \left( \frac{1}{4} + \mu \right)^{-s-1/2} p \left( \frac{2 \pi \alpha a}{\sqrt{4 \mu + 1}} \right)
                    \end{array} \]
                
                \noindent is entire, and its derivative at $s=0$ satisfies, as $\mu$ goes to infinity,
                
                \[ \begin{array}{lll}
                        \scriptstyle \frac{\partial}{\partial s}_{\left \vert s = 0 \right.} \left[ - 3^{-s} \frac{\sin \left( \pi s \right)}{2 \pi} \left( \frac{1}{4} + \mu \right)^{-s-1/2} p \left( \frac{2 \pi \alpha a}{\sqrt{4 \mu + 1}} \right) \right] & \scriptstyle = & \scriptstyle o \left( 1 \right)
                    \end{array} . \]
            
            \end{proposition}
            
            \begin{proof}
            
                This is a direct computation.
            
            \end{proof}
            
            \begin{proposition}
            \label{PropA6Power3}
            
                The function
                
                \[ \begin{array}{lll}
                        s & \longmapsto & - \frac{\sin \left( \pi s \right)}{\pi} \left( 4 \mu + 1 \right)^{-s-1/2} \sum\limits_{\left \vert k \right \vert \geqslant 1} \left( \left \vert k \right \vert^{2 \delta} - \frac{1}{4} \right)^{-s} \left \vert k \right \vert^{- \delta} p \left( \frac{2 \pi \left \vert k + \alpha \right \vert a}{\left \vert k \right \vert^{\delta} \sqrt{4 \mu + 1}} \right)^3
                    \end{array} \]
                    
                \noindent is holomorphic on the half-plane $\Re s > \left( \delta - 1 \right) / \delta$, which contains $0$.
            
            \end{proposition}
            
            \begin{proof}
            
                For any non-zero integer $k$, we have
                
                \[ \begin{array}{lllll}
                        p \left( \frac{2 \pi \left \vert k + \alpha \right \vert a}{\left \vert k \right \vert^{\delta} \sqrt{4 \mu + 1}} \right)^3 & = & \left( 1 + \frac{4 \pi^2 \left( k + \alpha \right)^2 a^2}{\left( 4 \mu + 1 \right) \left \vert k \right \vert^{2 \delta}} \right)^{-3/2} & \leqslant & \frac{\left( 4 \mu + 1 \right)^3 \left \vert k \right \vert^{3 \delta}}{\left( 2 \pi \left \vert k + \alpha \right \vert a \right)^3}
                    \end{array} . \]
                
                \noindent The inequality $3 - 2 \delta + 2 \delta \Re s > 1$ being equivalent to $\Re s \left( \delta - 1 \right)/\delta$, we get the result, after having noted that the half-plane in question contains $0$, since we have $\delta < 1$.
            
            \end{proof}
            
            \begin{proposition}
            \label{PropA6Power3-0}
            
                Assume we have $\alpha \neq 0$. The function
                
                \[ \begin{array}{lll}
                        s & \longmapsto & - 3^{-s} \frac{\sin \left( \pi s \right)}{2 \pi} \left( \frac{1}{4} + \mu \right)^{-s-1/2} p \left( \frac{2 \pi \alpha a}{\sqrt{4 \mu + 1}} \right)^3
                    \end{array} \]
                
                \noindent is entire, and its derivative at $s=0$ satisfies, as $\mu$ goes to infinity,
                
                \[ \begin{array}{lll}
                        \scriptstyle \frac{\partial}{\partial s}_{\left \vert s = 0 \right.} \left[ - 3^{-s} \frac{\sin \left( \pi s \right)}{2 \pi} \left( \frac{1}{4} + \mu \right)^{-s-1/2} p \left( \frac{2 \pi \alpha a}{\sqrt{4 \mu + 1}} \right)^3 \right] & \scriptstyle = & \scriptstyle o \left( 1 \right)
                    \end{array} . \]
            
            \end{proposition}
            
            \begin{proof}
            
                This is a direct computation.
            
            \end{proof}
            
            Let us now put these contributions together.
            
            \begin{proposition}
            \label{PropA6}
            
                The function
                
                \[ \begin{array}{lll}
                        s & \longmapsto & - \frac{\sin \left( \pi s \right)}{\pi} \left( 4 \mu + 1 \right)^{-s-1/2} \sum\limits_{\left \vert k \right \vert \geqslant 1} \left( \left \vert k \right \vert^{2 \delta} - \frac{1}{4} \right)^{-s} \left \vert k \right \vert^{- \delta} U_1 \left( p \left( \frac{2 \pi \left \vert k + \alpha \right \vert a}{\left \vert k \right \vert^{\delta} \sqrt{4 \mu + 1}} \right) \right)
                    \end{array} \]
                    
                \noindent which is well-defined on the half-plane $\Re s > 0$, has a holomorphic continuation to a region of the complex plane containing $0$.
            
            \end{proposition}
            
            \begin{proposition}
            \label{PropA6-0}
            
                Assume we have $\alpha \neq 0$. The function
                
                \[ \begin{array}{lll}
                        s & \longmapsto & - 3^{-s} \frac{\sin \left( \pi s \right)}{2 \pi} \left( \frac{1}{4} + \mu \right)^{-s-1/2} U_1 \left( p \left( \frac{2 \pi \alpha a}{\sqrt{4 \mu + 1}} \right) \right)
                    \end{array} \]
                
                \noindent is entire, and its derivative at $s=0$ satisfies, as $\mu$ goes to infinity,
                
                \[ \begin{array}{lll}
                        \scriptstyle \frac{\partial}{\partial s}_{\left \vert s = 0 \right.} \left[ - 3^{-s} \frac{\sin \left( \pi s \right)}{2 \pi} \left( \frac{1}{4} + \mu \right)^{-s-1/2} U_1 \left( p \left( \frac{2 \pi \alpha a}{\sqrt{4 \mu + 1}} \right) \right) \right] & \scriptstyle = & \scriptstyle o \left( 1 \right)
                    \end{array} . \]
            
            \end{proposition}

        \subsubsection*{\textbf{Seventh part}}
        \label{ASeventhPart}
            
            We now come to a complicated term we need to deal with, where we will use the Ramanujan summation, presented in apprendix \ref{AppRamanujan}, following \cite{MR3677185}.
            
            \begin{proposition}
            \label{PropA7}
            
                The function
                
                \[ \begin{array}{lll}
                        s & \longmapsto & \frac{\sin \left( \pi s \right)}{\pi} \left( 4 \mu + 1 \right)^{-s} \sum\limits_{\left \vert k \right \vert \geqslant 1} \left( \left \vert k \right \vert^{2 \delta} - \frac{1}{4} \right)^{-s} \sqrt{\left( 2 \pi \left \vert k + \alpha \right \vert a \right)^2 + \frac{1}{4} + \mu}
                    \end{array} , \]
                
                \noindent which is holomorphic on the half-plane $\Re s > 1/\delta$, has a holomorphic continuation to a region containing $0$. Furthermore, the (continuation of) the function 
                
                \[ \begin{array}{lll}
                        \scriptstyle s & \scriptstyle \longmapsto & \scriptstyle \frac{\sin \left( \pi s \right)}{\pi} \left( 4 \mu + 1 \right)^{-s} \sum\limits_{\left \vert k \right \vert \geqslant 1} \left( \left \vert k \right \vert^{2 \delta} - \frac{1}{4} \right)^{-s} \sqrt{\left( 2 \pi \left \vert k + \alpha \right \vert a \right)^2 + \frac{1}{4} + \mu} \\[1em]
                        
                        && \qquad \scriptstyle - \frac{\sin \left( \pi s \right)}{\pi} \cdot \frac{\Gamma \left( \delta s \right) \Gamma \left( 3/2 - \delta s \right)}{\sqrt{\pi}} \left( 4 \pi a \right)^{2 \delta s - 1} \left( 4 \mu + 1 \right)^{1 - \left( 1 + \delta \right) s} \cdot \frac{1}{\left( 2 \delta s - 1 \right) \left( 2 \delta s - 2 \right)}
                    \end{array} \]
                
                \noindent vanishes at $s=0$, and its derivative at this point satisfies, as $\mu$ goes to infinity,
                
                \[ \scalemath{0.95}{\begin{array}{lll}
                        \scriptstyle \frac{\partial}{\partial s}_{\left \vert s = 0 \right.} \Big[ \frac{\sin \left( \pi s \right)}{\pi} \left( 4 \mu + 1 \right)^{-s} \sum\limits_{\left \vert k \right \vert \geqslant 1} \left( \left \vert k \right \vert^{2 \delta} - \frac{1}{4} \right)^{-s} \sqrt{\left( 2 \pi \left \vert k + \alpha \right \vert a \right)^2 + \frac{1}{4} + \mu} \\[1em]
                        
                        \qquad  \qquad \scriptstyle - \frac{\sin \left( \pi s \right)}{\pi} \cdot \frac{\Gamma \left( \delta s \right) \Gamma \left( 3/2 - \delta s \right)}{\sqrt{\pi}} \left( 4 \pi a \right)^{2 \delta s - 1} \left( 4 \mu + 1 \right)^{1 - \left( 1 + \delta \right) s} \cdot \frac{1}{\left( 2 \delta s - 1 \right) \left( 2 \delta s - 2 \right)} \Big] & \hspace{-3pt} \scriptstyle = \hspace{-3pt} & \scriptstyle - \sqrt{\mu} + 2 \pi \alpha^2 a + o \left( 1 \right).
                    \end{array}} \]
            
            \end{proposition}
            
            \begin{proof}[Proof of proposition \ref{PropA7}]
            
                We begin by noting that the function
                
                \[ \begin{array}{lll}
                        s & \longmapsto & \sum\limits_{\left \vert k \right \vert \geqslant 1} \left( \left \vert k \right \vert^{2 \delta} - \frac{1}{4} \right)^{-s} \sqrt{\left( 2 \pi \left \vert k + \alpha \right \vert a \right)^2 + \frac{1}{4} + \mu}
                    \end{array} \]
                
                \noindent is indeed holomorphic on the half-plane $\Re s > 1 / \delta$. Though the existence of a holomorphic continuation can be proved using Taylor's formula, this method does not yield the $\mu$-asymptotics of the derivative at $0$. We will instead use the Ramanujan summation, which provides both results. Using the binomial formula, we get
                
                \begin{equation}
                \label{Part7PosNegInteger}
                    \scalemath{0.98}{\begin{array}{lll}
                        \multicolumn{3}{l}{\scriptstyle \sum\limits_{\left \vert k \right \vert \geqslant 1} \left( \left \vert k \right \vert^{2 \delta} - \frac{1}{4} \right)^{-s} \sqrt{\left( 2 \pi \left \vert k + \alpha \right \vert a \right)^2 + \frac{1}{4} + \mu}} \\[1em]
                        
                        & \scriptstyle = & \scriptstyle \sum\limits_{j = 0}^{+ \infty} \frac{\left( s \right)_j}{j! 4^j} \Big[ \sum\limits_{k \geqslant 1} k^{-2 \delta \left( s + j \right)} \sqrt{\left( 2 \pi \left( k + \alpha \right) a \right)^2 + \frac{1}{4} + \mu} + \sum\limits_{k \geqslant 1} k^{-2 \delta \left( s + j \right)} \sqrt{\left( 2 \pi \left( k - \alpha \right) a \right)^2 + \frac{1}{4} + \mu} \Big].
                    \end{array}}
                \end{equation}
                
                \noindent The sum has been split into two parts because there will be a small variation in the argument. We begin by dealing with the first part of \eqref{Part7PosNegInteger}, induced by
                
                \[ \begin{array}{c}
                        \sum\limits_{k \geqslant 1} k^{-2 \delta \left( s + j \right)} \sqrt{\left( 2 \pi \left( k + \alpha \right) a \right)^2 + \frac{1}{4} + \mu}
                    \end{array} \]
                
                \noindent for every integer $j \geqslant 0$. The first step in using the Ramanujan sumamtion is to find a function which interpolates, at integers, the general terms of the series we consider. We consider a fixed complex number $s$ with $\Re s > 1 / \delta$, and set
                
                \[ \begin{array}{lllll}
                        f_{s,j} & : & z & \longmapsto & z^{-2 \delta \left( s + j \right)} \sqrt{\left( 2 \pi \left( z + \alpha \right) a \right)^2 + \frac{1}{4} + \mu}
                    \end{array} , \]
                
                \noindent which is holomorphic on the half-plane $\Re z > 0$, where we have denoted by $\sqrt{\cdot}$ the principal branch of the complex square root. The function $f_{s,j}$ is then of \textit{moderate growth}, in the sense of definition \ref{DefExpModGrowth}. We now need to check that both hypotheses from theorem \ref{ThmRamanujanSum} are satisfied. The first one, which is
                
                \[ \begin{array}{lll}
                        \lim\limits_{k \rightarrow + \infty} f_{s,j} \left( k \right) & = & 0
                    \end{array} \]
                
                \noindent holds, since we have $\Re s > 1 / \delta$. We now need to prove that we have
                
                \[ \begin{array}{lll}
                        \lim\limits_{k \rightarrow + \infty} \int_0^{+ \infty} \frac{f_{s,j} \left( k + i t \right) - f_{s,j} \left( k - i t \right)}{e^{2 \pi t} - 1} \mathrm{d}t & = & 0
                    \end{array} . \]
                
                \noindent The idea here is to use the dominated convergence theorem. In order to bound the function within the integral by an integrable one, we need to take some care with the denominator $e^{2 \pi t} - 1$ at $t = 0$. We have
                
                \[ \begin{array}{lll}
                        f_{s,j} \left( k + i t \right) & = & f_{s,j} \left( k \right) + i \int_0^t f_{s,j}' \left( k + ix \right) \mathrm{d}x, \\[0.5em]
                        
                        f_{s,j} \left( k - i t \right) & = & f_{s,j} \left( k \right) - i \int_0^t f_{s,j}' \left( k - ix \right) \mathrm{d}x,
                    \end{array} \]
                
                \noindent for any integer $k \geqslant 1$ and any $t \in \left] 0, + \infty \right[$. Taking the difference then yields
                
                \[ \begin{array}{lll}
                        f_{s,j} \left( k + i t \right) - f_{s,j} \left( k - i t \right) & = & i \int_0^t \left( f_{s,j}' \left( k + ix \right) + f_{s,j}' \left( k - ix \right) \right) \mathrm{d}x
                    \end{array} , \]
                
                \noindent giving an extra factor $t$ which cancels the singular behavior of $e^{2 \pi t} - 1$ at $t = 0$. Through an explicit evaluation of $f_{s,j}'$, one proves that we have
                
                \begin{equation}
                \label{EqHypothesisRamanujan}
                    %
                \end{equation}
                
                \noindent with implicit constants independent of $k$. The details are not written here, as they amount to cumbersome estimates which make the overall argument even more intricate. This allows us to use the dominated convergence theorem, which gives
                
                \[ \begin{array}{lllll}
                        \lim\limits_{k \rightarrow + \infty} \int_0^{+ \infty} \frac{f_{s,j} \left( k + i t \right) - f_{s,j} \left( k - i t \right)}{e^{2 \pi t} - 1} \mathrm{d}t & = & \int_0^{+ \infty} \lim\limits_{k \rightarrow + \infty} \frac{f_{s,j} \left( k + i t \right) - f_{s,j} \left( k - i t \right)}{e^{2 \pi t} - 1} \mathrm{d}t & = & 0
                    \end{array} . \]
                
                \noindent We can therefore apply theorem \ref{ThmRamanujanSum}, which states that we have
                
                \[ \begin{array}{lll}
                        \sum\limits_{k = 1}^{+ \infty} f_{s,j} \left( k \right) & = & \sum\limits_{k \geqslant 1}^{\left( \mathcal{R} \right)} f_{s,j} \left( k \right) + \int_1^{+ \infty} f_{s,j} \left( x \right) \mathrm{d}x
                    \end{array} . \]
                
                \noindent We will now study the first term on the right-hand side of the equality above, which is a Ramanujan sum. The main difference between classical and Ramanujan sums lies with the fact that theorem \ref{ThmAnalycityRamanujanSum} assure us that the function
                
                \[ \begin{array}{lll}
                        s & \longmapsto & \sum\limits_{k \geqslant 1}^{\left( \mathcal{R} \right)} f_{s,j} \left( k \right)
                    \end{array} \]
                
                \noindent is actually entire, as it is a Ramanujan sum of entire functions. Using estimates such as those hinted at to get \eqref{EqHypothesisRamanujan}, one proves that the function
                
                \[ \begin{array}{lll}
                        s & \longmapsto & \sum\limits_{j = 0}^{+ \infty} \frac{\left( s \right)_j}{j! 4^j} \sum\limits_{k \geqslant 1}^{\left( \mathcal{R} \right)} f_{s,j} \left( k \right)
                    \end{array} \]
                
                \noindent is entire, and, after multiplication by the relevant factor, its derivative at $0$ is
                
                \[ \begin{array}{lll}
                        \scriptstyle \frac{\partial}{\partial s}_{\left \vert s = 0 \right.} \left[ \frac{\sin \left( \pi s \right)}{\pi} \left( 4 \mu + 1 \right)^{-s} \sum\limits_{j = 0}^{+ \infty} \frac{\left( s \right)_j}{j! 4^j} \sum\limits_{k \geqslant 1}^{\left( \mathcal{R} \right)} f_{s,j} \left( k \right) \right] & \scriptstyle = & \scriptstyle \frac{1}{2} f_{0,0} \left( 1 \right) + i \int_0^{+ \infty} \frac{f_{0,0} \left( 1 + it \right) - f_{0,0} \left( 1 - it \right)}{e^{2 \pi t} - 1} \mathrm{d}t.
                    \end{array} \]
                
                \noindent It should now be remembered that each function $f_{s,j}$ depends implicitly on $\mu$, and that the behavior as $\mu$ goes to infinity of the Ramanujan sum above must be studied. We have
                
                \[ \begin{array}{lllll}
                        \frac{1}{2} f_{0,0} \left( 1 \right) & = & \frac{1}{2} \sqrt{\left( 2 \pi \left( 1 + \alpha \right) a \right)^2 + \frac{1}{4} + \mu} & = & \frac{1}{2} \sqrt{\mu} + o \left( 1 \right)
                    \end{array} \]
                
                \noindent as $\mu$ goes to infinity. We will now prove that we have
                
                \begin{equation}
                \label{EqA7RamanujanSumAsymptotic}
                    \begin{array}{lll}
                        \lim\limits_{\mu \rightarrow + \infty} \int_0^{+ \infty} \frac{f_{0,0} \left( 1 + it \right) - f_{0,0} \left( 1 - it \right)}{e^{2 \pi t} - 1} \mathrm{d}t & = & 0
                    \end{array} .
                \end{equation}
                
                \noindent To do this, we must once again call upon the dominated convergence theorem. The estimates we need for the domination hypothesis are the same as the ones alluded to earlier in this proof. Let us briefly expand on them. For any $t \geqslant 0$, we have
                
                \[ \begin{array}{lll}
                        f_{0,0} \left( 1 \pm i t \right) & = & \sqrt{\left( 2 \pi \left( 1 \pm i t + \alpha \right) a \right)^2 + \frac{1}{4} + \mu}
                    \end{array} , \]
                
                \noindent and we can add the two terms, to get
                
                \[ \begin{array}{lll}
                        f_{0,0} \left( 1 + i t \right) - f_{0,0} \left( 1 - i t \right) & = & \frac{16 i \pi^2 a^2 t \left( 1 + \alpha \right)}{\sqrt{\left( 2 \pi \left( 1 + \alpha + it \right) a \right)^2 + \frac{1}{4} + \mu} + \sqrt{\left( 2 \pi \left( 1 + \alpha - it \right) a \right)^2 + \frac{1}{4} + \mu}}
                    \end{array} . \]
                
                \noindent We can find a lower bound for the modulus of the denominator, by writing
                
                \[ \begin{array}{lll}
                        \multicolumn{3}{l}{\scriptstyle \left \vert \sqrt{\left( 2 \pi \left( 1 + \alpha + it \right) a \right)^2 + \frac{1}{4} + \mu} + \sqrt{\left( 2 \pi \left( 1 + \alpha - it \right) a \right)^2 + \frac{1}{4} + \mu} \right \vert} \\[1em]
                        
                        \qquad \qquad \qquad \qquad \qquad & \scriptstyle \geqslant & \scriptstyle \Re \sqrt{\left( 2 \pi \left( 1 + \alpha + it \right) a \right)^2 + \frac{1}{4} + \mu} + \Re \sqrt{\left( 2 \pi \left( 1 + \alpha - it \right) a \right)^2 + \frac{1}{4} + \mu},
                    \end{array} \]
                
                \noindent and each real part can be explicitely computed, yielding for instance
                
                \[ \begin{array}{lll}
                        \scriptstyle \Re \sqrt{\left( 2 \pi \left( 1 + \alpha + it \right) a \right)^2 + \frac{1}{4} + \mu} & \scriptstyle \geqslant & \scriptstyle \frac{1}{\sqrt{2}} \big[ \underbrace{\big \vert \scriptstyle \frac{1}{4} + \mu - 4 \pi^2 a^2 t^2 \big \vert + \frac{1}{4} + \mu - 4 \pi^2 a^2 t^2}_{\geqslant 0} + 4 \pi^2 a^2 \left( 1 + \alpha \right)^2 \big]^{1/2} \\[1em]
                        
                        & \scriptstyle \geqslant & \scriptstyle \sqrt{2} \pi a \left( 1 + \alpha \right).
                    \end{array} \]
                
                \noindent We therefore get the upper bound
                
                \[ \begin{array}{lll}
                        \left \vert f_{0,0} \left( 1 + i t \right) - f_{0,0} \left( 1 - i t \right) \right \vert & \leqslant & 8 \sqrt{2} \pi a t
                    \end{array} . \]
                
                \noindent This allows us to use the dominated convergence theorem, which yields \eqref{EqA7RamanujanSumAsymptotic}. We can now move to the core of this proof, \textit{i.e.} the study of the function
                
                \[ \begin{array}{lll}
                        s & \longmapsto & \frac{\sin \left( \pi s \right)}{\pi} \left( 4 \mu + 1 \right)^{-s} \sum\limits_{j = 0}^{+ \infty} \frac{\left( s \right)_j}{j! 4^j} \int_1^{+ \infty} \frac{1}{x^{2 \delta \left( s+j \right)}} \sqrt{\left( 2 \pi \left( x + \alpha \right) a \right)^2 + \frac{1}{4} + \mu} \; \mathrm{d}x
                    \end{array} , \]
                
                \noindent which is well-defined and holomorphic on the half-plane $\Re s > 1 / \delta$. The advantage of dealing with an integral rather than a series is that we can perform integration by parts and change of variables. For any integer $j \geqslant 0$, we have
                
                \[ \begin{array}{lll}
                        \multicolumn{3}{l}{\scriptstyle \int_1^{+ \infty} \frac{1}{x^{2 \delta \left( s+j \right)}} \sqrt{\left( 2 \pi \left( x + \alpha \right) a \right)^2 + \frac{1}{4} + \mu} \; \mathrm{d}x} \\[0.5em]
                        
                        \qquad & = & \scriptstyle \left[ - \frac{1}{2 \delta \left( s + j \right) - 1} \cdot \frac{1}{x^{2 \delta \left( s + j \right) - 1}} \sqrt{\left( 2 \pi \left( x + \alpha \right) a \right)^2 + \frac{1}{4} + \mu} \right]_{1}^{+ \infty} \\[0.em]
                        
                        && \qquad \qquad \qquad \scriptstyle + \frac{1}{2 \delta \left( s + j \right) - 1} \cdot \frac{1}{2} \int_1^{+ \infty} \frac{1}{x^{2 \delta \left( s + j \right) - 1}} \cdot \frac{8 \pi^2 a^2 \left( x + \alpha \right)}{\sqrt{\left( 2 \pi \left( x + \alpha \right) a \right)^2 + \frac{1}{4} + \mu}} \mathrm{d}x \\[0.5em]
                        
                        & = & \scriptstyle \frac{1}{2 \delta \left( s + j \right) - 1} \sqrt{\left( 2 \pi \left( 1 + \alpha \right) a \right)^2 + \frac{1}{4} + \mu} + \frac{2 \pi a}{2 \delta \left( s + j \right) - 1} \int_1^{+ \infty} \frac{1}{x^{2 \delta \left( s + j \right) - 1}} \left( 1 + \frac{1/4 + \mu}{\left( 2 \pi \left( x + \alpha \right) a \right)^2} \right)^{-1/2} \mathrm{d}x.
                    \end{array} \]
                
                \noindent We can now compute this last integral, using hypergeometric functions, and more precisely using corollary \ref{CorIntegralFormulaHypergeometric}. First, let us prepare the computation. We have
                
                \[ \begin{array}{lll}
                        \multicolumn{3}{l}{\scriptstyle \int_1^{+ \infty} \frac{1}{x^{2 \delta \left( s + j \right) - 1}} \left( 1 + \frac{1/4 + \mu}{\left( 2 \pi \left( x + \alpha \right) a \right)^2} \right)^{-1/2} \mathrm{d}x} \\[0.5em]
                        
                        & \scriptstyle = & \scriptstyle \frac{1}{2} \left( 2 \pi a \right)^{2 \delta \left( s + j \right) - 2} \left( \frac{1}{4} + \mu \right)^{- \delta \left( s + j \right) + 1} \int_0^{\frac{1/4 + \mu}{4 \pi^2 \left( 1 + \alpha \right)^2 a^2}} \frac{t^{\delta \left( s + j \right) - 2}}{\sqrt{1 + t}} \left( 1 - \frac{2 \pi \alpha a}{\sqrt{1/4 + \mu}} t^{1/2} \right)^{- 2 \delta \left( s + j \right) + 1} \mathrm{d}t.
                    \end{array} \]
                
                \noindent On the interval of integration, we now have
                
                \[ \begin{array}{lllllllll}
                        0 & \leqslant & \frac{2 \pi \alpha a}{\sqrt{1/4 + \mu}} t^{1/2} & \leqslant & \frac{2 \pi \alpha a}{\sqrt{1/4 + \mu}} \cdot \frac{\sqrt{1/4 + \mu}}{2 \pi \left( 1 + \alpha \right) a} & = & \frac{\alpha}{1 + \alpha} & < & 1
                    \end{array} . \]
                
                \noindent We now use the binomial formula, and interchange sum and integral, to get
                
                \begin{equation}
                \label{EqA7Ramanujan-01}
                    \begin{array}{lll}
                        \multicolumn{3}{l}{\scriptstyle \scriptstyle \int_1^{+ \infty} \frac{1}{x^{2 \delta \left( s+j \right)}} \sqrt{\left( 2 \pi \left( x + \alpha \right) a \right)^2 + \frac{1}{4} + \mu} \; \mathrm{d}x} \\[0.5em]
                        
                        & \scriptstyle = & \scriptstyle \frac{1}{2 \delta \left( s + j \right) - 1} \sqrt{\left( 2 \pi \left( 1 + \alpha \right) a \right)^2 + \frac{1}{4} + \mu} \; \; + \Big[ \frac{1}{2 \left( 2 \delta \left( s + j \right) - 1 \right)} \left( 2 \pi a \right)^{2 \delta \left( s + j \right) - 1} \left( \frac{1}{4} + \mu \right)^{- \delta \left( s + j \right) + 1} \\[0.5em]
                        
                        && \qquad \qquad \scriptstyle \cdot \sum\limits_{n=0}^{+ \infty} \frac{\Gamma \left( 2 \delta \left( s + j \right) + n - 1 \right)}{n! \Gamma \left( 2 \delta \left( s + j \right) - 1 \right)} \cdot \frac{\left( 2 \pi \alpha a \right)^n}{\left( 1/4 + \mu \right)^{n/2}} \int_0^{\frac{1/4 + \mu}{4 \pi^2 \left( 1 + \alpha \right)^2 a^2}} \frac{t^{\left( \delta \left( s + j \right) + n/2 - 1 \right) - 1}}{\sqrt{1+t}} \mathrm{d}t \Big].
                    \end{array}
                \end{equation}
                
                \noindent Each integral appearing in the last sum of \eqref{EqA7Ramanujan-01} is now seen to be of the form presented in corollary \ref{CorIntegralFormulaHypergeometric}. For any integer $n \geqslant 0$, we have indeed
                
                \begin{equation}
                \label{EqA7Ramanujan-02}
                    \begin{array}{lll}
                        \multicolumn{3}{l}{\scriptstyle \int_0^{\frac{1/4 + \mu}{4 \pi^2 \left( 1 + \alpha \right)^2 a^2}} \frac{t^{\left( \delta \left( s + j \right) + n/2 - 1 \right) - 1}}{\sqrt{1+t}} \mathrm{d}t} \\[0.5em]
                        
                        \qquad \qquad & \scriptstyle = & \scriptstyle \frac{2}{2 \delta \left( s + j \right) + n - 2} \left( \frac{1/4 + \mu}{4 \pi^2 \left( 1 + \alpha \right)^2 a^2} \right)^{\delta \left( s + j \right) + n/2 - 1} \left( 1 + \frac{1/4 + \mu}{4 \pi^2 \left( 1 + \alpha \right)^2 a^2} \right)^{-1/2} \\[0.5em]
                        
                        && \qquad \qquad \qquad \qquad \qquad \scriptstyle \cdot F \left( \frac{1}{2}, 1; \delta \left( s + j \right) + \frac{n}{2}; \frac{1/4 + \mu}{4 \pi^2 \left( 1 + \alpha \right)^2 a^2 + 1/4 + \mu} \right),
                    \end{array}
                \end{equation}
                
                \noindent and plugging \eqref{EqA7Ramanujan-02} into \eqref{EqA7Ramanujan-01} yields
                
                \[ \begin{array}{lll}
                        \multicolumn{3}{l}{\scriptstyle \scriptstyle \int_1^{+ \infty} \frac{1}{x^{2 \delta \left( s+j \right)}} \sqrt{\left( 2 \pi \left( x + \alpha \right) a \right)^2 + \frac{1}{4} + \mu} \; \mathrm{d}x} \\[0.5em]
                        
                        \qquad & = & \scriptstyle \frac{1}{2 \delta \left( s + j \right) - 1} \sqrt{\left( 2 \pi \left( 1 + \alpha \right) a \right)^2 + \frac{1}{4} + \mu} \; \; + \big[ 4 \pi^2 a^2 \left( 1 + \alpha \right)^{3 - 2 \delta \left( s + j \right)} \sum\limits_{n = 0}^{+ \infty} \frac{\Gamma \left( 2 \delta \left( s + j \right) + n - 2 \right)}{n! \Gamma \left( 2 \delta \left( s + j \right) \right)} \\[0.5em]
                        
                        && \qquad \scriptstyle \cdot \left( \frac{\alpha}{1 + \alpha} \right)^n \frac{1}{\sqrt{4 \pi^2 \left( 1 + \alpha \right)^2 a^2 + 1/4 + \mu}} F \left( \frac{1}{2}, 1; \delta \left( s + j \right) + \frac{n}{2}; \frac{1/4 + \mu}{4 \pi^2 \left( 1 + \alpha \right)^2 a^2 + 1/4 + \mu} \right) \big].
                    \end{array} \]
                
                \noindent We will deal with these two terms separately. Prior to beginning this study, one should recall that $1/\left( 2 \delta \right)$ is assumed not to be an integer, which means that there is no integer $j > 0$ such that we have $2 \delta j = 1$. Regarding the first term, the function
                
                \[ \begin{array}{lll}
                        s & \longmapsto & \frac{\sin \left( \pi s \right)}{\pi} \left( 4 \mu + 1 \right)^{-s} \sum\limits_{j = 0}^{+ \infty} \frac{\left( s \right)_j}{j! 4^j} \frac{1}{2 \delta \left( s + j \right) - 1} \sqrt{\left( 2 \pi \left( 1 + \alpha \right) a \right)^2 + \frac{1}{4} + \mu}
                    \end{array} \]
                
                \noindent is holomorphic around $0$, and its derivative there satisfies
                
                \[ \begin{array}{lll}
                        \scriptstyle \frac{\partial}{\partial s}_{\left \vert s = 0 \right.} \left[ \frac{\sin \left( \pi s \right)}{\pi} \left( 4 \mu + 1 \right)^{-s} \sum\limits_{j = 0}^{+ \infty} \frac{\left( s \right)_j}{j! 4^j} \frac{1}{2 \delta \left( s + j \right) - 1} \sqrt{\left( 2 \pi \left( 1 + \alpha \right) a \right)^2 + \frac{1}{4} + \mu} \right] & \scriptstyle = & \scriptstyle - \sqrt{\mu} \; + \; o \left( 1 \right)
                    \end{array} \]
                
                \noindent as $\mu$ goes to infinity. We can now study the second term, given by
                
                \[ \begin{array}{l}
                        \scriptstyle 4 \pi^2 a^2 \sum\limits_{j = 0}^{+ \infty} \frac{\left( s \right)_j}{j! 4^j} \left( 1 + \alpha \right)^{3 - 2 \delta \left( s + j \right)} \sum\limits_{n = 0}^{+ \infty} \frac{\Gamma \left( 2 \delta \left( s + j \right) + n - 2 \right)}{n! \Gamma \left( 2 \delta \left( s + j \right) \right)} \left( \frac{\alpha}{1 + \alpha} \right)^n \frac{1}{\sqrt{4 \pi^2 \left( 1 + \alpha \right)^2 a^2 + 1/4 + \mu}} \\[0.7em] 
                        
                        \scriptstyle \qquad \qquad \qquad \qquad \qquad \qquad \qquad \qquad \qquad \qquad \cdot F \left( \frac{1}{2}, 1; \delta \left( s + j \right) + \frac{n}{2}; \frac{1/4 + \mu}{4 \pi^2 \left( 1 + \alpha \right)^2 a^2 + 1/4 + \mu} \right). 
                    \end{array} \]
                
                \noindent We will need to break apart this term according to the value of the integer $n$.
                
                \medskip
                
                \hspace{10pt} $\bullet$ For $n \geqslant 4$, we have  $\delta \left( \Re s + j \right) + \frac{n}{2} - \frac{1}{2} - 1 > \delta \Re s + \frac{1}{2}$, which means that
                
                    \[ \begin{array}{lll}
                            t & \longmapsto & F \left( \frac{1}{2}, 1; \delta \left( s + j \right) + \frac{n}{2}; t \right)
                        \end{array} \]
                
                    \noindent is bounded on $\left[ 0, 1 \right[$, uniformly in $j$ and in $s$ near $0$, by proposition \ref{PropHypergeometricZ1}. The sum over $n \geqslant 4$ induces a holomorphic function around $0$. After having multiplied by the relevant factor, its derivative at $0$ vanishes, because of the Pochhammer symbol for non-zero integers $j$, and because of $\Gamma \left( 2 \delta \left( s + j \right) \right)$ for $j = 0$.
                
                    \medskip
                
                \hspace{10pt} $\bullet$ For $n = 3$, we consider, after simplifications,
                
                    \[ \begin{array}{l}
                            \scriptstyle \frac{4}{3} \cdot \frac{\pi^2 a^2 \alpha^3}{\sqrt{4 \pi^2 \left( 1 + \alpha \right)^2 a^2 + 1/4 + \mu}} \sum\limits_{j=0}^{+ \infty} \frac{\left( s \right)_j}{j! 4^j} \left( 1 + \alpha \right)^{-2 \delta \left( s + j \right)} \cdot \delta \left( s + j \right) F \left( \frac{1}{2}, 1; \delta \left( s + j \right) + \frac{3}{2}; \frac{1/4 + \mu}{4 \pi^2 \left( 1 + \alpha \right)^2 a^2 + 1/4 + \mu} \right)
                        \end{array} . \]
                
                    \noindent For $j \geqslant 1$, we can use proposition \ref{PropHypergeometricZ1} to bound the hypergeometric function uniformly, for $s$ near $0$. The term corresponding to $j = 0$ inducing an entire function, the whole sum is holomorphic around $0$, and, after multiplication by the relevant factor, its derivative there vanishes, because of the Pochhammer symbol for non-zero integers $j$, and because of $s + j$ for $j = 0$.
                
                    \medskip
                
                \hspace{10pt} $\bullet$ For $n = 2$, we consider, after simplifications,
                
                    \[ \begin{array}{l}
                            \scriptstyle \frac{2 \pi^2 a^2 \alpha^2}{\sqrt{4 \pi^2 \left( 1 + \alpha \right)^2 a^2 + 1/4 + \mu}} \sum\limits_{j = 0}^{+ \infty} \frac{\left( s \right)_j}{j! 4^j} \left( 1 + \alpha \right)^{1 - 2 \delta \left( s + j \right)} F \left( \frac{1}{2}, 1; \delta \left( s + j \right); \frac{1/4 + \mu}{4 \pi^2 \left( 1 + \alpha \right)^2 a^2 + 1/4 + \mu} \right)
                        \end{array} . \]
                        
                    \noindent For $j > 1/\left( 2 \delta \right)$, we use proposition \ref{PropHypergeometricZ1} to bound the hypergeometric function independently of the parameters, for $s$ around $0$. The remaining terms, in finite number, induce entire functions, so the sum over $j \geqslant 0$ is holomorphic around $0$. The Pochhammer symbol vanishing at $0$ for all $j > 0$, we have, as $\mu$ goes to infinity,
                
                    \[ \begin{array}{lllll}
                            \multicolumn{5}{l}{\scriptstyle \frac{\partial}{\partial s}_{\left \vert s = 0 \right.} \left[ \frac{\sin \left( \pi s \right)}{\pi} \frac{2 \pi^2 a^2 \alpha^2 \left( 4 \mu + 1 \right)^{-s}}{\sqrt{4 \pi^2 \left( 1 + \alpha \right)^2 a^2 + 1/4 + \mu}} \sum\limits_{j = 0}^{+ \infty} \frac{\left( s \right)_j}{j! 4^j} \left( 1 + \alpha \right)^{1 - 2 \delta \left( s + j \right)} F \left( \frac{1}{2}, 1; \delta \left( s + j \right); \frac{1/4 + \mu}{4 \pi^2 \left( 1 + \alpha \right)^2 a^2 + 1/4 + \mu} \right) \right]} \\[1em]
                        
                            \qquad \qquad \qquad & = & \scriptstyle \frac{2 \pi^2 a^2 \alpha^2 \left( 4 \mu + 1 \right)^{-s}}{\sqrt{4 \pi^2 \left( 1 + \alpha \right)^2 a^2 + 1/4 + \mu}} \left( 1 + \alpha \right) F \left( \frac{1}{2}, 1; 1; \frac{1/4 + \mu}{4 \pi^2 \left( 1 + \alpha \right)^2 a^2 + 1/4 + \mu} \right) & = & \scriptstyle \pi \alpha^2 a,
                        \end{array} \]
                
                    \noindent the hypergeometric function being given by remark \ref{RmkHypergeometricInterchangeAB} and proposition \ref{PropBinomialFormula}.
                
                    \medskip
                
                \hspace{10pt} $\bullet$ For $n = 1$, we consider, after simplifications,
                
                    \[ \begin{array}{l}
                            \scriptstyle \frac{2 \pi^2 a^2 \alpha}{\sqrt{4 \pi^2 \left( 1 + \alpha \right)^2 a^2 + 1/4 + \mu}} \sum\limits_{j=0}^{+\infty} \frac{\left( s \right)_j}{j! 4^j} \cdot \frac{\left( 1 + \alpha \right)^{2 - 2 \delta \left( s + j \right)}}{2 \delta \left( s + j \right) - 1} F \left( \frac{1}{2}, 1; \delta \left( s + j \right) + \frac{1}{2}; \frac{1/4 + \mu}{4 \pi^2 \left( 1 + \alpha \right)^2 a^2 + 1/4 + \mu} \right)
                        \end{array} . \]
                
                    \noindent For $j > 1/\delta$, we use proposition \ref{PropHypergeometricZ1} to bound the hypergeometric function independently of the parameters, for $s$ around $0$. The remaining terms, in finite number, induce entire functions, so the whole sum over $j \geqslant 0$ is holomorphic around $0$. The Pochhammer symbol vanishing at $0$ for all $j > 0$, we have, as $\mu$ goes to infinity,
                
                    \[ \begin{array}{lllll}
                            \multicolumn{5}{l}{\scriptstyle \frac{\partial}{\partial s}_{\left \vert s = 0 \right.} \hspace{-2pt} \left[ \frac{\sin \left( \pi s \right)}{\pi} \frac{2 \pi^2 a^2 \alpha \left( 4 \mu + 1 \right)^{-s}}{\sqrt{4 \pi^2 \left( 1 + \alpha \right)^2 a^2 + 1/4 + \mu}} \hspace{-1pt} \sum\limits_{j=0}^{+\infty} \frac{\left( s \right)_j \left( 1 + \alpha \right)^{2 - 2 \delta \left( s + j \right)}}{j! 4^j \left( 2 \delta \left( s + j \right) - 1 \right)} F \left( \frac{1}{2}, 1; \delta \left( s + j \right) + \frac{1}{2}; \frac{1/4 + \mu}{4 \pi^2 \left( 1 + \alpha \right)^2 a^2 + 1/4 + \mu} \right) \hspace{-1pt} \right]} \\[1em]
                        
                            \qquad \qquad \qquad \qquad & \scriptstyle = & \multicolumn{3}{l}{\scriptstyle - \frac{2 \pi^2 a^2 \alpha}{\sqrt{4 \pi^2 \left( 1 + \alpha \right)^2 a^2 + 1/4 + \mu}} \left( 1 + \alpha \right)^2 F \left( \frac{1}{2}, 1; \frac{1}{2}; \frac{1/4 + \mu}{4 \pi^2 \left( 1 + \alpha \right)^2 a^2 + 1/4 + \mu} \right)} \\[1em]
                        
                            & \scriptstyle = & \scriptstyle - \frac{1}{2} \alpha \sqrt{4 \pi^2 \left( 1 + \alpha \right)^2 a^2 + 1/4 + \mu} & \scriptstyle = & \scriptstyle - \frac{1}{2} \alpha \sqrt{\mu} + o \left( 1 \right),
                        \end{array} \]
                    
                    \noindent where proposition \ref{PropBinomialFormula} is used to compute the hypergeometric function.
                
                    \medskip
                
                \hspace{10pt} $\bullet$ For $n = 0$, we consider, after simplifications,
                
                    \[ \begin{array}{l}
                            \scriptstyle \frac{4 \pi^2 a^2}{\sqrt{4 \pi^2 \left( 1 + \alpha \right)^2 a^2 + 1/4 + \mu}} \sum\limits_{j = 0}^{+ \infty} \frac{\left( s \right)_j}{j! 4^j} \frac{\left( 1 + \alpha \right)^{3 - 2 \delta \left( s + j \right)}}{\left( 2 \delta \left( s + j \right) - 2 \right) \left( 2 \delta \left( s + j \right) - 1 \right)} F \left( \frac{1}{2}, 1; \delta \left( s + j \right); \frac{1/4 + \mu}{4 \pi^2 \left( 1 + \alpha \right)^2 a^2 + 1/4 + \mu} \right)
                        \end{array} . \]
                    
                    \noindent For $j > 3/\left( 2 \delta \right)$, we use proposition \ref{PropHypergeometricZ1} to bound the hypergeometric function independently of the parameters, for $s$ around $0$. The sum over $j \geqslant 1$, given by
                    
                    \[ \begin{array}{l}
                            \scriptstyle \frac{\sin \left( \pi s \right)}{\pi} \frac{4 \pi^2 a^2 \left( 4 \mu + 1 \right)^{-s}}{\sqrt{4 \pi^2 \left( 1 + \alpha \right)^2 a^2 + 1/4 + \mu}} \sum\limits_{j = 1}^{+ \infty} \frac{\left( s \right)_j}{j! 4^j} \frac{\left( 1 + \alpha \right)^{3 - 2 \delta \left( s + j \right)}}{\left( 2 \delta \left( s + j \right) - 2 \right) \left( 2 \delta \left( s + j \right) - 1 \right)} F \left( \frac{1}{2}, 1; \delta \left( s + j \right); \frac{1/4 + \mu}{4 \pi^2 \left( 1 + \alpha \right)^2 a^2 + 1/4 + \mu} \right)
                        \end{array} \]
                    
                    \noindent is holomorphic around $0$, and its derivative there vanishes, because of the Pochhammer symbol. The term corresponding to $j = 0$ requires more care. It is given by
                    
                    \[ \begin{array}{lll}
                            \scriptstyle \frac{\sin \left( \pi s \right)}{\pi} \frac{4 \pi^2 a^2 \left( 4 \mu + 1 \right)^{-s}}{\sqrt{4 \pi^2 \left( 1 + \alpha \right)^2 a^2 + 1/4 + \mu}} \frac{\left( 1 + \alpha \right)^{3 - 2 \delta s}}{\left( 2 \delta s - 2 \right) \left( 2 \delta s - 1 \right)} F \left( \frac{1}{2}, 1; \delta s; \frac{1/4 + \mu}{4 \pi^2 \left( 1 + \alpha \right)^2 a^2 + 1/4 + \mu} \right)
                        \end{array} . \]
                    
                    \noindent We need to work a bit more on this term, since the hypergeometric function is going to induce a simple pole at $0$, which will be compensated by the sine. We have
                    
                    \[ \begin{array}{lll}
                            \multicolumn{3}{l}{\scriptstyle F \left( \frac{1}{2}, 1; \delta s; \frac{1/4 + \mu}{4 \pi^2 \left( 1 + \alpha \right)^2 a^2 + 1/4 + \mu} \right)} \\[1em]
                            
                            \quad & \scriptstyle = & \scriptstyle \frac{\Gamma \left( \delta s \right) \Gamma \left( 3/2 - \delta s \right)}{\Gamma \left( 1/2 \right)} \left( \frac{4 \pi^2 \left( 1 + \alpha \right)^2 a^2}{4 \pi^2 \left( 1 + \alpha \right)^2 a^2 + 1/4 + \mu} \right)^{\delta s - 3/2} F \left( \delta s - \frac{1}{2}, \delta s - 1; \delta s - \frac{1}{2}; \frac{4 \pi^2 \left( 1 + \alpha \right)^2 a^2}{4 \pi^2 \left( 1 + \alpha \right)^2 a^2 + 1/4 + \mu} \right) \\[1em]
                            
                            && \qquad \scriptstyle + \frac{\Gamma \left( \delta s \right) \Gamma \left( \delta s - 3/2 \right)}{\Gamma \left( \delta s - 1/2 \right) \Gamma \left( \delta s - 1 \right)} F \left( \frac{1}{2}, 1; \frac{5}{2} - \delta s; \frac{4 \pi^2 \left( 1 + \alpha \right)^2 a^2}{4 \pi^2 \left( 1 + \alpha \right)^2 a^2 + 1/4 + \mu} \right)
                        \end{array} \]
                    
                    \noindent using proposition \ref{PropTransformationLinear}. The first of these last two hypergeometric functions is given by proposition \ref{PropBinomialFormula}, and we can simplify some of the Gamma functions, yielding
                    
                    \[ %
 \]
                    
                    \noindent After multiplication by the appropriate factor, we get
                    
                    \[ %
 \]
                    
                    \noindent We now deal with each of these terms separately, beginning with
                    
                    \[ \begin{array}{l}
                            \scriptstyle \frac{\sin \left( \pi s \right)}{\pi} \left( 4 \mu + 1 \right)^{-s} \frac{4 \pi^2 a^2}{\sqrt{4 \pi^2 \left( 1 + \alpha \right)^2 a^2 + 1/4 + \mu}} \cdot \frac{\left( 1 + \alpha \right)^{3 - 2 \delta s}}{\left( 2 \delta s - 1 \right) \left( 2 \delta s - 3 \right)} F \left( \frac{1}{2}, 1; \frac{5}{2} - \delta s; \frac{4 \pi^2 \left( 1 + \alpha \right)^2 a^2}{4 \pi^2 \left( 1 + \alpha \right)^2 a^2 + 1/4 + \mu} \right)
                        \end{array} . \]
                    
                    \noindent This term induces a holomorphic function around $0$, whose derivative here satisfies
                    
                    \[ %
 \]
                    
                    \noindent as $\mu$ goes to infinity. The second term, given by
                    
                    \[ \begin{array}{l}
                            \scriptstyle \frac{\sin \left( \pi s \right)}{\pi \sqrt{\pi}} 4^{\delta s - 1} \left( 4 \mu + 1 \right)^{1 - \left( 1 + \delta \right) s} \frac{\left( 2 \pi a \right)^{2 \delta s - 1}}{\left( 2 \delta s - 2 \right) \left( 2 \delta s - 1 \right)}  \Gamma \left( \delta s \right) \Gamma \left( \frac{3}{2} - \delta s \right)
                        \end{array} \]
                    
                    \noindent induces a holomorphic function near $0$, since the simple pole of the Gamma function is compensated by the simple zero of the sine function. Its derivative at $s=0$ can be computed, but it is not necessary to do that to get this proposition. We now turn our attention to the other part from \eqref{Part7PosNegInteger}, induced by
                    
                    \[ \begin{array}{c}
                            \sum\limits_{k \geqslant 1} k^{-2 \delta \left( s + j \right)} \sqrt{\left( 2 \pi \left( k - \alpha \right) a \right)^2 + \frac{1}{4} + \mu}
                        \end{array} . \]
                    
                    \noindent The method used here is similar, and the final result can be obtained formally by switching the sign of $\alpha$. The only difference is that we need to deal with the integral
                    
                    \[ %
 \]
                    
                    \noindent and we have, on the interval of integration,
                    
                    \[ %
 . \]
                    
                    \noindent Unlike the case we have dealt with in detail, the last quantity above cannot be bounded, since there is no limit to how close $\alpha$ can get to $1$. This prevents us from using the binomial formula. To avoid that problem, we should instead study
                    
                    \[ %
 \]
                    
                    \noindent On this new interval of integration, we have
                    
                    \[ %
 . \]
                    
                    \noindent The rest is similar to what we have already seen, and the remaining term, given by
                    
                    \[ %
 , \]
                    
                    \noindent induces a holomorphic function around $0$, whose derivative there equals $\sqrt{\mu} + o \left( 1 \right)$, as $\mu$ goes to infinity. Putting all these results together yield the proposition.
            
            \end{proof}
            
            We will now take care of the case $k = 0$, assuming that we have $\alpha \neq 0$.
            
            \begin{proposition}
            \label{PropA7-0}
            
                Assume we have $\alpha \neq 0$. The function
                
                \[ %
 \]
                
                \noindent is entire, and its derivative at $0$ satisfies, as $\mu$ goes to infinity,
                
                \[ %
 . \]
            
            \end{proposition}
            
            \begin{proof}
            
                This is a direct computation.
            
            \end{proof}

        \subsubsection*{\textbf{Eighth part}}
        \label{AEighthPart}
            
            The term in question here can be studied using the Ramanujan summation, along the lines of what is done in proposition \ref{PropA7}.
            
            \begin{proposition}
            \label{PropA8}
            
                The function
                
                \[ %
 \]
                
                \noindent is holomorphic on the half-plane $\Re s > 1/\left( 2 \delta \right)$, has a holomorphic continuation to an open neighborhood of $0$. Furthermore, the (continuation of) the function 
                
                \[ %
 \]
                
                \noindent vanishes at $s=0$, and its derivative at this point satisfies, as $\mu$ goes to infinity,
                
                \[ %
 \]
            
            \end{proposition}
            
            \begin{proof}
            
                The proof is entirely similar to that of proposition \ref{PropA7}. Let us simply mention that the integral remaining in the derivative at $0$ is found when asymptotically studying the Ramanujan sum as $\mu$ goes to infinity.
            
            \end{proof}
            
            Let us now state the result for the case $k = 0$, supposing we have $\alpha \neq 0$.
            
            \begin{proposition}
            \label{PropA8-0}
            
                Assume we have $\alpha \neq 0$. The function
                
                \[ %
 \]
                
                \noindent is entire, and its derivative at $s = 0$ is given by, as $\mu$ goes to infinity,
                
                \[ %
 . \]
            
            \end{proposition}
            
            \begin{proof}
            
                This is a direct computation.
            
            \end{proof}

        \subsubsection*{\textbf{Ninth part}}
        \label{ANinthPart}
        
            The next step in this section is to study the logarithmic term coming from proposition \ref{PropLogK3}. The arguments are the same as in proposition \ref{PropA7}.
            
            \begin{proposition}
            \label{PropA9}
            
                The function
                
                \[ %
 \]
                
                \noindent is holomorphic on the half-plane $\Re s > 1/\left( 2 \delta \right)$, has a holomorphic continuation to an open neighborhood of $0$. Its derivative there satisfies, as $\mu$ goes to infinity,
                
                \[ %
 . \]
            
            \end{proposition}
            
            \begin{proof}
            
                The proof is similar to that of proposition \ref{PropA7}.
            
            \end{proof}
            
            For this part too, we need to account for the case $k = 0$, should $\alpha$ not vanish.
            
            \begin{proposition}
            \label{PropA9-0}
            
                Assume we have $\alpha \neq 0$. The function
                
                \[ %
 \]
                
                \noindent is entire, and its derivative at $s = 0$ is given by, as $\mu$ goes to infinity,
                
                \[ %
 . \]
            
            \end{proposition}
            
            \begin{proof}
            
                This is a direct computation.
            
            \end{proof}

        \subsubsection*{\textbf{Tenth part}}
        \label{ATenthPart}
        
            Finally, we must take care of the polynomial term $U_1$ appearing in proposition \ref{PropLogK3}. The arguments used here are the same as in proposition \ref{PropA7}.
            
            \begin{proposition}
            \label{PropA10}
            
                The function
                
                \[ %
 \]
                
                \noindent is holomorphic on the half-plane $\Re s > 0$, has a holomorphic continuation to a neighborhood of $0$. Furthermore, the (continuation of) the function 
                
                \[ %
 \]
                
                \noindent vanishes at $s=0$, and its derivative at this point satisfies, as $\mu$ goes to infinity,
                
                \[ %
 \]
            
            \end{proposition}
            
            \begin{proof}
            
                The proof is similar to that of proposition \ref{PropA7}.
            
            \end{proof}
            
            Let us combine propositions \ref{PropA6} and \ref{PropA10} into one, in order to make some simplifications.
            
            \begin{proposition}
            \label{PropA6and10}
            
                The function
                
                \[ %
 \]
                
                \noindent is holomorphic on the half-plane $\Re s > 0$, has a holomorphic continuation to a neighborhood of $0$, whose value at $s=0$ vanishes, and whose derivative at $s=0$ satisfies, as $\mu$ goes to infinity,
                
                \[ %
 \]
                
            \end{proposition}
            
            \begin{proof}
            
                This result is a direct consequence of propositions \ref{PropA6} and \ref{PropA10}, after having noted that the function
                
                \[ %
 \]
                
                \noindent extends holomorphically near $0$, vanishes at $s=0$, and that we have
                
                \[ %
 \]
                
                \noindent as $\mu$ goes to infinity. These statements are obtained by using Laurent series expansions, as in proposition \ref{PropR}.
            
            \end{proof}
            
            Once more, let us take care of the case $k = 0$, should we have $\alpha \neq 0$.
            
            \begin{proposition}
            \label{PropA10-0}
            
                Assume we have $\alpha \neq 0$. The function
                
                \[ %
 \]
                
                \noindent is entire, and its derivative at $0$ satisfies, as $\mu$ goes to infinity,
                
                \[ %
 . \]
            
            \end{proposition}
            
            \begin{proof}
            
                This is a direct computation.
            
            \end{proof}

        \subsubsection*{\textbf{Eleventh part}}
        \label{AEleventhPart}
            
            Going back to equation \eqref{EqFmuk} and definition \ref{DefAmuk}, we have so far dealt with every term coming from the difference
            
            \[ %
 , \]
            
            \noindent with terms in $\log \left( \pi / 2 \right)$ being canceled. The remaining term we need to take care comes from the last part of equation \eqref{EqFmuk}, which is given by
            
            \[ %
 . \]
            
            \noindent There is actually no new content here, but stating a proposition will still help when putting all the pieces together to get the results proved in this paper.
            
            \begin{proposition}
            \label{PropA11}
            
                The function
                
                \[ %
 \]
                
                \noindent is holomorphic on the half-plane $\Re s > 1/\delta$, and has a holomorphic continuation to an open neighborhood of $0$.
            
            \end{proposition}
            
            \begin{proof}
            
                This is a direct consequence of propositions \ref{PropAmukAsymptA}, \ref{PropA1}, \ref{PropA2}, \ref{PropA3}, \ref{PropA4}, \ref{PropA5}, \ref{PropA7}, \ref{PropA8}, \ref{PropA9}, \ref{PropA6and10}.
            
            \end{proof}
            
            \begin{remark}
            
                The point of having removed the two terms from the logarithmic derivative of the Bessel function is that we get a holomorphic function around $0$ without having to multiply by the factor $\sin \left( \pi s \right)$.
            
            \end{remark}
            
            For this final part related to the study of the terms $A_{\mu,k}$, we have to see what happens in the case $k = 0$, which plays a role when $\alpha$ does not vanish.
            
            \begin{proposition}
            \label{PropA11-0}
            
                Assume we have $\alpha \neq 0$. The following function is entire
                
                \[ \begin{array}{lll}
                        s & \longmapsto & - \frac{3^{-s+1}}{2} \frac{\sin \left( \pi s \right)}{\pi} \left( \frac{1}{4} + \mu \right)^{- s + \frac{1}{2}} \frac{\partial}{\partial t}_{\left \vert t = \sqrt{\frac{1}{4} + \mu} \right.} \log K_t \left( 2 \pi \alpha a \right)
                    \end{array} . \]
            
            \end{proposition}
            
            \begin{proof}
            
                The result is direct.
            
            \end{proof}

    \subsection{\texorpdfstring{Study of the integrals $M_{\mu,k}$}{Study of the integrals Mmu,k}}
    \label{SubSecStudyMmuk}
    
        Having studied the integrals~$L_{\mu,k}$ coming from 
        
        \[ \begin{array}{lll}
                I_{\mu,k} \left( s \right) & = & L_{\mu,k} + M_{\mu,k} \left( s \right),
            \end{array} \]
        
        \noindent we turn our attention to $M_{\mu,k}$. Recall that, according to definition \ref{DefMmuk}, we have
        
        \[ \begin{array}{lll}
                M_{\mu,k} \left( s \right) & = & \displaystyle \frac{\sin \left( \pi s \right)}{\pi} \int_{2 \left \vert k \right \vert^{\delta} \sqrt{\frac{1}{4} + \mu}}^{+ \infty} \; \left( t^2 - \left( \frac{1}{4} + \mu \right) \right)^{-s} \; f_{\mu,k} \left( t \right) \; \text{d}t
            \end{array} \]
            
        \noindent for non-zero integers $k$, and that, should $\alpha$ not vanish, we also have
            
        \[ \begin{array}{lll}
                M_{\mu,0} \left( s \right) & = & \displaystyle \frac{\sin \left( \pi s \right)}{\pi} \int_{2 \sqrt{\frac{1}{4} + \mu}}^{+ \infty} \; \left( t^2 - \left( \frac{1}{4} + \mu \right) \right)^{-s} \; f_{\mu,0} \left( t \right) \; \text{d}t
            \end{array} . \]
            
        \noindent As indicated in definition \ref{Deffmuk}, we have
        
        \[ \begin{array}{lll}
                f_{\mu,k} \left( t \right) & = & \frac{\partial}{\partial t} \log K_t \left( 2\pi \left \vert k + \alpha \right \vert a \right) - \frac{2t}{\sqrt{4 \mu + 1}} \frac{\partial}{\partial t}_{\left \vert t = \sqrt{\frac{1}{4} + \mu} \right.} \log K_t \left( 2 \pi \left \vert k + \alpha \right \vert a \right).
            \end{array} \]
        
        \noindent We will now split $M_{\mu,k}$, for every integer $k$, according to the expression of $f_{\mu,k}$.
        
        \begin{definition}
        \label{DefMtildemuk}
                
            On the strip $1 < \Re s < 2$, and for any real number $\mu \geqslant 0$, we set
            
            \[ \begin{array}{lll}
                    \widetilde{M}_{\mu,k} \left( s \right) & = & \frac{\sin \left( \pi s \right)}{\pi} \int_{2 \left \vert k \right \vert^{\delta} \sqrt{\frac{1}{4} + \mu}}^{+ \infty} \left( t^2 - \left( \frac{1}{4} + \mu \right) \right)^{-s} \frac{\partial}{\partial t} \log K_t \left( 2\pi \left \vert k + \alpha \right \vert a \right) \mathrm{d}t
                \end{array} \]
            
            \noindent for any non-zero integer $k$. Assuming we have $\alpha \neq 0$, we also set, on the same strip
            
            \[ \begin{array}{lll}
                    \widetilde{M}_{\mu,0} \left( s \right) & = & \frac{\sin \left( \pi s \right)}{\pi} \int_{2 \sqrt{\frac{1}{4} + \mu}}^{+ \infty} \left( t^2 - \left( \frac{1}{4} + \mu \right) \right)^{-s} \frac{\partial}{\partial t} \log K_t \left( 2 \pi \alpha a \right) \mathrm{d}t
                \end{array} . \]
                
        \end{definition}
        
        \begin{definition}
        \label{DefRmuk}
                
            On the strip $1 < \Re s < 2$, and for any real number $\mu \geqslant 0$, we set
            
            \[ \begin{array}{lll}
                    R_{\mu,k} \left( s \right) & = & - \frac{\sin \left( \pi s \right)}{\pi} \cdot \frac{2}{\sqrt{4 \mu + 1}} \int_{2 \left \vert k \right \vert^{\delta} \sqrt{\frac{1}{4} + \mu}}^{+ \infty} t \left( t^2 - \left( \frac{1}{4} + \mu \right) \right)^{-s} \\[0.5em]
                    
                    && \qquad \qquad \qquad \qquad \qquad \qquad \quad \cdot \frac{\partial}{\partial t}_{\left \vert t = \sqrt{\frac{1}{4} + \mu} \right.} \log K_t \left( 2 \pi \left \vert k + \alpha \right \vert a \right) \mathrm{d}t
                \end{array} \]
            
            \noindent for any non-zero integer $k$. We also set, on the same strip
            
            \[ \begin{array}{lll}
                    R_{\mu,0} \left( s \right) & = & - \frac{\sin \left( \pi s \right)}{\pi} \cdot \frac{2}{\sqrt{4 \mu + 1}} \int_{2 \sqrt{\frac{1}{4} + \mu}}^{+ \infty} t \left( t^2 - \left( \frac{1}{4} + \mu \right) \right)^{-s} \\[0.5em]
                    
                    && \qquad \qquad \qquad \qquad \qquad \qquad \quad \cdot \frac{\partial}{\partial t}_{\left \vert t = \sqrt{\frac{1}{4} + \mu} \right.} \log K_t \left( 2 \pi \alpha a \right) \mathrm{d}t
                \end{array} \]
            
            \noindent assuming we have $\alpha \neq 0$.
                
        \end{definition}
        
        \begin{remark}
        
            We have $M_{\mu,k} \left( s \right) = \widetilde{M}_{\mu,k} \left( s \right) + R_{\mu,k} \left( s \right)$ on the strip $1 < \Re s < 2$.
        
        \end{remark}

        \subsubsection{Study of the integrals $R_{\mu,k}$}
        \label{SubSubSecStudyRmuk}
        
            We begin this section by taking care of the remainder terms $R_{\mu,k}$. The relevant derivatives at $s = 0$ will be studied together with those from propositions \ref{PropA11} and \ref{PropA11-0}. We begin by a couple of lemmas.
            
            \begin{lemma}
            \label{LemRmuk}
            
                On the strip $1 < \Re s < 2$, we have
                
                \[ \begin{array}{lll}
                        R_{\mu,k} \left( s \right) & = & \frac{\left( 4 \mu + 1 \right)^{-s+ \frac{1}{2}}}{1-s} \frac{\sin \left( \pi s \right)}{\pi} \left( \left \vert k \right \vert^{2 \delta} - \frac{1}{4} \right)^{-s+1} \frac{\partial}{\partial t}_{\left \vert t = \sqrt{\frac{1}{4} + \mu} \right.} \log K_t \left( 2 \pi \left \vert k + \alpha \right \vert a \right)
                    \end{array} . \]
            
            \end{lemma}
            
            \begin{proof}
            
                On the strip $1 < \Re s < 2$, we have
                
                \[ \begin{array}{lll}
                        \scriptstyle R_{\mu,k} \left( s \right) & \scriptstyle = & \scriptstyle - \frac{\sin \left( \pi s \right)}{\pi} \cdot \frac{2}{\sqrt{4 \mu + 1}} \left( \int_{2 \left \vert k \right \vert^{\delta} \sqrt{\frac{1}{4} + \mu}}^{+ \infty} t \left( t^2 - \left( \frac{1}{4} + \mu \right) \right)^{-s} \mathrm{d}t \right) \frac{\partial}{\partial t}_{\left \vert t = \sqrt{1/4 + \mu} \right.} \log K_t \left( 2 \pi \left \vert k + \alpha \right \vert a \right) \\[1em]
                        
                        & \scriptstyle = & \scriptstyle - \frac{\sin \left( \pi s \right)}{\pi} \cdot \frac{1}{\sqrt{4 \mu + 1}} \left[ \frac{1}{1 - s} \left( t^2 - \left( \frac{1}{4} + \mu \right) \right)^{-s+1} \right]_{2 \left \vert k \right \vert^{\delta} \sqrt{1/4 + \mu}}^{+ \infty} \frac{\partial}{\partial t}_{\left \vert t = \sqrt{1/4 + \mu} \right.} \log K_t \left( 2 \pi \left \vert k + \alpha \right \vert a \right) \\[1em]
                        
                        & \scriptstyle = & \scriptstyle \frac{1}{1-s} \cdot \frac{\sin \left( \pi s \right)}{\pi} \left( 4 \mu + 1 \right)^{-s+ 1/2} \left( \left \vert k \right \vert^{2 \delta} - \frac{1}{4} \right)^{-s+1} \frac{\partial}{\partial t}_{\left \vert t = \sqrt{1/4 + \mu} \right.} \log K_t \left( 2 \pi \left \vert k + \alpha \right \vert a \right).
                    \end{array} \]
                
                \noindent The proof of the proposition is thus complete.
            
            \end{proof}
            
            We further have the version of this lemma corresponding to the term $k = 0$.
            
            \begin{lemma}
            \label{LemRmuk-0}
            
                Assume we have $\alpha \neq 0$. On the strip $1 < \Re s < 2$, we have
                
                \[ \begin{array}{lll}
                        R_{\mu, 0} \left( s \right) & = & \frac{1}{2} \cdot \frac{3^{-s + 1}}{1 - s} \cdot \frac{\sin \left( \pi s \right)}{\pi} \left( \frac{1}{4} + \mu \right)^{-s+\frac{1}{2}} \frac{\partial}{\partial t}_{\left \vert t = \sqrt{1/4 + \mu} \right.} \log K_t \left( 2 \pi \left \vert k + \alpha \right \vert a \right)
                    \end{array} \]
            
            \end{lemma}
            
            \begin{proof}
            
                This is a direct computation, similar to the proof of lemma \ref{LemRmuk}.
            
            \end{proof}
            
            Having these two computations, we can study the terms $R_{\mu,k}$.
            
            \begin{proposition}
            \label{PropR}
            
                The function
                
                \[ \begin{array}{lll}
                        s & \longmapsto & \sum\limits_{\left \vert k \right \vert \geqslant 1} R_{\mu,k} \left( s \right)
                    \end{array} \]
                
                \noindent is holomorphic on the half-plane $\Re s > 1/\delta$, and has a holomorphic continuation to an open neighborhood of zero. Its derivative at $s = 0$ satisfies, as $\mu$ goes to infinity,
                
                \[ \begin{array}{lll}
                        \multicolumn{3}{l}{\scriptstyle \frac{\partial}{\partial s}_{\left \vert s = 0 \right.} \left[ \sum\limits_{\left \vert k \right \vert \geqslant 1} R_{\mu,k} \left( s \right) \right]} \\[1em]
                        
                        & \scriptstyle = & \hspace{-3pt} \scriptstyle - \frac{\partial}{\partial s}_{\left \vert s = 0 \right.}  \left[ \frac{\sin \left( \pi s \right)}{\pi} \left( 4 \mu + 1 \right)^{-s} \left( - \sqrt{4 \mu + 1} \sum\limits_{\left \vert k \right \vert \geqslant 1} \left( \left \vert k \right \vert^{2 \delta} - \frac{1}{4} \right)^{- s + 1} \frac{\partial}{\partial t}_{\left \vert t = \sqrt{\frac{1}{4} + \mu} \right.} \log K_t \left( 2 \pi \left \vert k + \alpha \right \vert a \right) \right. \right. \\[1em]
                        
                        && \qquad \qquad \qquad \qquad \qquad \; \scriptstyle + \frac{1}{\sqrt{\pi}} \Gamma \left( \delta s \right) \Gamma \left( \frac{3}{2} - \delta s \right) \left( 4 \pi a \right)^{2 \delta s - 1} \left( 4 \mu + 1 \right)^{1 - \delta s} \cdot \frac{1}{\left( 2 \delta s - 1 \right) \left( 2 \delta s - 2 \right)} \\[1em]
                        
                        && \qquad \qquad \qquad \qquad \qquad \; \; \left. \left. \scriptstyle + \frac{1}{\sqrt{\pi}} \Gamma \left( \delta s \right) \Gamma \left( \frac{1}{2} - \delta s \right) 4^{\delta s - 1} \left( 2 \pi a \right)^{2 \delta s - 1} \left( 4 \mu + 1 \right)^{1 - \delta s} \cdot \frac{1}{2 \delta s - 1} \right) \right] \\[1em]
                        
                        && \scriptstyle + \frac{1}{4 \pi a} \left( 1 + \frac{1}{\delta} \right) \mu \log \mu - \frac{1}{4 \pi a \delta} \left( 1 + 2 \delta \log \left( 4 \pi a \right) + 3 \delta - 2 \log 2 \right) \mu  + \frac{1}{16 \pi a} \left( 1 + \frac{1}{\delta} \right) \log \mu \\[1em]
                        
                        && \qquad \scriptstyle - \frac{1}{8 \pi a} \log \left( 4 \pi a \right) + \frac{1}{8 \pi a \delta} \log 2 - \frac{1}{8 \pi a} + o \left( 1 \right).
                    \end{array} \]
                
                \noindent Furthermore, the same derivative, this time with $\mu = 0$, has the following asymptotic expansion, as $a$ goes to infinity,
                
                \[ \begin{array}{lll}
                        \frac{\partial}{\partial s}_{\left \vert s = 0 \right.} \left[ \sum\limits_{\left \vert k \right \vert \geqslant 1} R_{0,k} \left( s \right) \right] & = & o \left( 1 \right)
                    \end{array} . \]
                
            \end{proposition}
            
            \begin{proof}
            
                This result relies on proposition \ref{PropA11}, which tells us that the function
                
                \[ \begin{array}{lll}
                        \scriptstyle s & \scriptstyle \longmapsto & \scriptstyle \left( s - 1 \right) \sum\limits_{\left \vert k \right \vert \geqslant 1} R_{\mu,k} \left( s \right) \\[1em]
                        
                        && \scriptstyle \quad + \frac{\sin \left( \pi s \right)}{\pi} \left( 4 \mu + 1 \right)^{-s} \left[ \frac{1}{\sqrt{\pi}} \Gamma \left( \delta s \right) \Gamma \left( \frac{3}{2} - \delta s \right) \left( 4 \pi a \right)^{2 \delta s - 1} \left( 4 \mu + 1 \right)^{1 - \delta s} \cdot \frac{1}{\left( 2 \delta s - 1 \right) \left( 2 \delta s - 2 \right)} \right. \\[1em]
                        
                        && \scriptstyle \qquad \qquad \qquad \qquad \quad \left. + \frac{1}{\sqrt{\pi}} \Gamma \left( \delta s \right) \Gamma \left( \frac{1}{2} - \delta s \right) 4^{\delta s - 1} \left( 2 \pi a \right)^{2 \delta s - 1} \left( 4 \mu + 1 \right)^{1 - \delta s} \cdot \frac{1}{2 \delta s - 1} \right] \\[1em]
                    \end{array} \]
                
                \noindent is holomorphic on the half-plane $\Re s > 1/\delta$, and has a holomorphic continuation around $0$, whose derivative there equals that of the (continuation of) the function 
                
                \[ \begin{array}{lll}
                        \scriptstyle s & \scriptstyle \longmapsto & \scriptstyle - \sum\limits_{\left \vert k \right \vert \geqslant 1} R_{\mu,k} \left( s \right) \\[1em]
                        
                        && \scriptstyle \quad - \frac{1}{s-1} \cdot \frac{\sin \left( \pi s \right)}{\pi} \left( 4 \mu + 1 \right)^{-s} \left[ \frac{1}{\sqrt{\pi}} \Gamma \left( \delta s \right) \Gamma \left( \frac{3}{2} - \delta s \right) \left( 4 \pi a \right)^{2 \delta s - 1} \left( 4 \mu + 1 \right)^{1 - \delta s} \cdot \frac{1}{\left( 2 \delta s - 1 \right) \left( 2 \delta s - 2 \right)} \right. \\[1em]
                        
                        && \scriptstyle \qquad \qquad \qquad \qquad \qquad \left. + \frac{1}{\sqrt{\pi}} \Gamma \left( \delta s \right) \Gamma \left( \frac{1}{2} - \delta s \right) 4^{\delta s - 1} \left( 2 \pi a \right)^{2 \delta s - 1} \left( 4 \mu + 1 \right)^{1 - \delta s} \cdot \frac{1}{2 \delta s - 1} \right] \\[1em]
                    \end{array} \]
                
                \noindent Furthermore, the common value of these derivatives equals the one left uncomputed on the right-hand side of the equality stated in the current proposition. In order to prove the part of the proposition related to the $\mu$-asymptotic study, it only remains to evaluate three derivatives at $s=0$ as $\mu$ goes to infinity. The first one is
                
                \[ \begin{array}{c}
                        \frac{\partial}{\partial s}_{\left \vert s = 0 \right.} \left[ - \frac{1}{s-1} \frac{\sin \left( \pi s \right)}{\pi \sqrt{\pi}} \Gamma \left( \delta s \right) \Gamma \left( \frac{3}{2} - \delta s \right) \left( 4 \pi a \right)^{2 \delta s - 1} \left( 4 \mu + 1 \right)^{1 - \left( 1+\delta \right) s} \frac{1}{\left( 2 \delta s - 1 \right) \left( 2 \delta s - 2 \right)} \right].
                    \end{array} \]
                
                \noindent This is done by using Laurent series expansions. We have
                
                \[ \begin{array}{lll}
                        \multicolumn{3}{l}{\scriptstyle - \frac{1}{s-1} \frac{\sin \left( \pi s \right)}{\pi \sqrt{\pi}} \Gamma \left( \delta s \right) \Gamma \left( \frac{3}{2} - \delta s \right) \left( 4 \pi a \right)^{2 \delta s - 1} \left( 4 \mu + 1 \right)^{1 - \left( 1+\delta \right) s} \frac{1}{\left( 2 \delta s - 1 \right) \left( 2 \delta s - 2 \right)}} \\[1em]
                        
                        & \scriptstyle = & \scriptstyle \frac{4 \mu + 1}{16 \pi a \delta} \left( 1 + s + O \left( s^2 \right) \right) \left( 1 + O \left( s^2 \right) \right) \left( 1 - \delta \gamma s + O \left( s^2 \right) \right) \left( 1 + \delta \left( 2 \log 2 + \gamma - 2 \right) s + O \left( s^2 \right) \right) \\[0.5em]
                        
                        && \scriptstyle \quad \cdot \left( 1 + 2 \delta \log \left( 4 \pi a \right) s + O \left( s^2 \right) \right) \left( 1 - \left( 1 + \delta \right) \log \left( 4 \mu + 1 \right) s + O \left( s^2 \right) \right) \left( 1 + 2 \delta s + O \left( s^2 \right) \right) \left( 1 + \delta s + O \left( s^2 \right) \right),
                    \end{array} \]
                
                \noindent and the required derivative is given, as $\mu$ goes to infinity, by
                
                \[ \begin{array}{lll}
                        \multicolumn{3}{l}{\scriptstyle \frac{\partial}{\partial s}_{\left \vert s = 0 \right.} \left[ - \frac{1}{s-1} \frac{\sin \left( \pi s \right)}{\pi \sqrt{\pi}} \Gamma \left( \delta s \right) \Gamma \left( \frac{3}{2} - \delta s \right) \left( 4 \pi a \right)^{2 \delta s - 1} \left( 4 \mu + 1 \right)^{1 - \left( 1+\delta \right) s} \frac{1}{\left( 2 \delta s - 1 \right) \left( 2 \delta s - 2 \right)} \right]} \\[1em]
                        
                        \qquad & \scriptstyle = & \scriptstyle - \frac{1}{4 \pi a} \left( 1 + \frac{1}{\delta} \right) \mu \log \mu + \frac{1}{4 \pi a \delta} \left[ 1 - 2 \log 2 + 2 \delta \log \left( 4 \pi a \right) + \delta \right] \mu - \frac{1}{16 \pi a} \left( 1 + \frac{1}{\delta} \right) \log \mu \\[0.5em]
                        
                        && \scriptstyle \qquad + \frac{1}{8 \pi a \delta} \left[ \delta \log \left( 4 \pi a \right) - \log 2 \right] + o \left( 1 \right).
                    \end{array} \]
                
                \noindent The last derivative we need to deal with is
                
                \[ \begin{array}{c}
                        \frac{\partial}{\partial s}_{\left \vert s = 0 \right.} \left[ - \frac{1}{s-1} \frac{\sin \left( \pi s \right)}{\pi \sqrt{\pi}} \Gamma \left( \delta s \right) \Gamma \left( \frac{1}{2} - \delta s \right) 4^{\delta s - 1} \left( 2 \pi a \right)^{2 \delta s - 1} \left( 4 \mu + 1 \right)^{1 - \left( 1+\delta \right) s} \frac{1}{2 \delta s - 1} \right].
                    \end{array} \]
                
                \noindent We have the following Laurent series expansion
                
                \[ \begin{array}{lll}
                        \multicolumn{3}{l}{\scriptstyle - \frac{1}{s-1} \frac{\sin \left( \pi s \right)}{\pi \sqrt{\pi}} \Gamma \left( \delta s \right) \Gamma \left( \frac{1}{2} - \delta s \right) 4^{\delta s - 1} \left( 2 \pi a \right)^{2 \delta s - 1} \left( 4 \mu + 1 \right)^{1 - \left( 1+\delta \right) s} \frac{1}{2 \delta s - 1}} \\[1em]
                        
                        & \scriptstyle = & \scriptstyle - \frac{4 \mu + 1}{8 \pi a \delta} \left( 1 + s + O \left( s^2 \right) \right) \left( 1 + O \left( s^2 \right) \right) \left( 1 - \delta \gamma s + O \left( s^2 \right) \right) \left( 1 + \delta \left( 2 \log 2 + \gamma \right) s + O \left( s^2 \right) \right) \left( 1 + 2 \delta s + O \left( s^2 \right) \right) \\[0.5em]
                        
                        && \quad \scriptstyle \cdot \left( 1 + 2 \delta \log \left( 2 \right) s + O \left( s^2 \right) \right) \left( 1 + 2 \delta \log \left( 2 \pi a \right) s + O \left( s^2 \right) \right) \left( 1 - \left( 1+\delta \right) \log \left( 4 \mu + 1 \right) s + O \left( s^2 \right) \right), 
                    \end{array} \]
                
                \noindent and the required derivative is given, as $\mu$ goes to infinity, by
                
                \[ \begin{array}{lll}
                        \multicolumn{3}{l}{\scriptstyle \frac{\partial}{\partial s}_{\left \vert s = 0 \right.} \left[ - \frac{1}{s-1} \frac{\sin \left( \pi s \right)}{\pi \sqrt{\pi}} \Gamma \left( \delta s \right) \Gamma \left( \frac{1}{2} - \delta s \right) 4^{\delta s - 1} \left( 2 \pi a \right)^{2 \delta s - 1} \left( 4 \mu + 1 \right)^{1 - \left( 1+\delta \right) s} \frac{1}{2 \delta s - 1} \right]} \\[1em]
                        
                        \qquad & \scriptstyle = & \scriptstyle \frac{1}{2 \pi a} \left( 1 + \frac{1}{\delta} \right) \mu \log \mu - \frac{1}{2 \pi a \delta} \left[ 1 + 2 \delta \log \left( 4 \pi a \right) - 2 \log 2 + 2 \delta \right] \mu  + \frac{1}{8 \pi a} \left( 1 + \frac{1}{\delta} \right) \log \mu \\[0.5em]
                        
                        && \scriptstyle \qquad - \frac{1}{8 \pi a \delta} \left[ 2 \delta \log \left( 4 \pi a \right) - 2 \log 2 + \delta \right] + o \left( 1 \right).
                    \end{array} \]
                
                \noindent Let us now study the behavior when $\mu = 0$, as $a$ goes to infinity. For any relative integer $k$, except $0$ should $\alpha$ vanish, we have
                
                \[ \begin{array}{lll}
                        R_{0,k} \left( s \right) & = & \frac{1}{1-s} \cdot \frac{\sin \left( \pi s \right)}{\pi} \left( \left \vert k \right \vert^{2 \delta} - \frac{1}{4} \right)^{-s + 1} \frac{\partial}{\partial t}_{\left \vert t = 1/2 \right.} \log K_t \left( 2 \pi \left \vert k + \alpha \right \vert a \right) \\[1em]
                        
                        & = & \frac{1}{1-s} \cdot \frac{\sin \left( \pi s \right)}{\pi} \left( \left \vert k \right \vert^{2 \delta} - \frac{1}{4} \right)^{-s + 1} \mathbb{E}_1 \left( 4 \pi \left \vert k + \alpha \right \vert a \right) e^{4 \pi \left \vert k + \alpha \right \vert a},
                    \end{array} \]
                
                \noindent using proposition \ref{PropDerivativeOneHalfModifiedBessel}, where $\mathbb{E}_1$ stands for the \textit{exponential integral} function. The asymptotic expansion given in proposition \ref{PropDerivativeOneHalfModifiedBessel} allows us to conclude.
            
            \end{proof}

            As always, we must take care of the term $k = 0$, whenever it makes sense.
            
            \begin{proposition}
            \label{PropR-0}
            
                Assume we have $\alpha \neq 0$. The function $s \longmapsto R_{\mu,0} \left( s \right)$ is holomorphic on $\mathbb{C} \setminus \left \{ 1 \right \}$, with a simple pole at $1$. Its derivative at $s=0$ is given by
                
                \[ \begin{array}{c}
                        \frac{\partial}{\partial s}_{\left \vert s = 0 \right.} \left[ \frac{3^{-s+1}}{2} \frac{\sin \left( \pi s \right)}{\pi} \left( \frac{1}{4} + \mu \right)^{- s + \frac{1}{2}} \frac{\partial}{\partial t}_{\left \vert t = \sqrt{\frac{1}{4} + \mu} \right.} \log K_t \left( 2 \pi \alpha a \right) \right].
                    \end{array} \]
                
                \noindent For $\mu = 0$, this derivative vanishes as $a$ goes to infinity.
            
            \end{proposition}
            
            \begin{proof}
            
                This is a direct computation, which is done using the expression of $R_{\mu,0}$. The last part is proved using proposition \ref{PropDerivativeOneHalfModifiedBessel}.
            
            \end{proof}

        \subsubsection{Study of the integrals $\widetilde{M}_{\mu,k}$}
        \label{SubSubSecStudyMTildemuk}
        
            Now that we are done with the study of the remainder terms $R_{\mu,k}$, we turn our attention to the core of this subsection, which is comprised of the terms presented in definition \ref{DefMtildemuk}.
            
            \begin{lemma}
            \label{LemLogDerivativeModifiedBessel}
            
                For any real number $t > 0$, and any relative integer $k$, with the exception of $0$ should $\alpha$ vanish, we have
                
                \[ \begin{array}{lll}
                        \multicolumn{3}{l}{\scriptstyle \frac{\partial}{\partial t} \log K_t \left( 2 \pi \left \vert k + \alpha \right \vert a \right)} \\[0.5em] 
                        
                        \qquad & \scriptstyle = & \scriptstyle \arcsinh \left( \frac{t}{2 \pi \left \vert k + \alpha \right \vert a} \right) - \frac{1}{2} \cdot \frac{t}{t^2 + 4 \pi^2 \left( k + \alpha \right)^2 a^2} - \frac{1}{8} \frac{\partial}{\partial t} \left( \frac{1}{t} \left( 1 + \frac{1}{t^2} \cdot 4 \pi^2 \left( k + \alpha \right)^2 a^2 \right)^{-1/2} \right) \\[1em]
                        
                        && \qquad \scriptstyle + \frac{5}{24} \frac{\partial}{\partial t} \left( \frac{1}{t} \left( 1 + \frac{1}{t^2} \cdot 4 \pi^2 \left( k + \alpha \right)^2 a^2 \right)^{-3/2} \right) + \frac{\partial}{\partial t} \left( \frac{1}{t^2} \widetilde{\eta_2} \left( t, \frac{1}{t} \cdot 2 \pi \left \vert k + \alpha \right \vert a \right) \right),
                    \end{array} \]
                
                \noindent the remainder term $\widetilde{\eta_2}$ being introduced in corollary \ref{CorAsymptoticExpansionModifiedBessel}.
            
            \end{lemma}
            
            \begin{proof}
            
                Let us begin by recalling that we have
                
                \[ \begin{array}{lll}
                        \multicolumn{3}{l}{\scriptstyle \log K_t \left( 2 \pi \left \vert k + \alpha \right \vert a \right)} \\[0.5em]
                        
                        \quad & \scriptstyle = & \scriptstyle \frac{1}{2} \log \frac{\pi}{2} + t \arcsinh \left( \frac{t}{2 \pi \left \vert k + \alpha \right \vert a} \right) - \sqrt{t^2 + 4 \pi^2 \left( k + \alpha \right)^2 a^2} - \frac{1}{4} \log \left( t^2 + 4 \pi^2 \left( k + \alpha \right)^2 a^2 \right) \\[0.5em]
                        
                        && \quad \scriptstyle - \frac{1}{8t} \left( 1 + \frac{1}{t^2} \cdot 4 \pi^2 \left( k + \alpha \right)^2 a^2 \right)^{-1/2} + \frac{5}{24t} \left( 1 + \frac{1}{t^2} \cdot 4 \pi^2 \left( k + \alpha \right)^2 a^2 \right)^{-3/2} + \frac{1}{t^2} \widetilde{\eta_2} \left( t, \frac{1}{t} \cdot 2 \pi \left \vert k + \alpha \right \vert a \right),
                    \end{array} \]
                
                \noindent under the hypotheses presented in this lemma, this equality being a direct consequence of corollary \ref{CorAsymptoticExpansionModifiedBessel}. After having differentiated with respect to $t$, and making some simplifications, we get the required formula.
            
            \end{proof}
            
            The strategy is now to take the integral defining $\widetilde{M}_{\mu,k}$, and to substitute the logarithmic derivative of the Bessel function by the expression above. This results in considering four terms separately.

            \subsubsection*{\textbf{First part}}
            \label{MTildeFirstPart}
            
                The first term we study is associated to the remainder $\widetilde{\eta_2}$.
                
                \begin{proposition}
                \label{PropMTilde1}
                
                    The function
                    
                    \[ 
 \]
                    
                    \noindent is holomorphic on the half-plane $\Re s > -1/4$, and its derivative at $0$ satisfies
                    
                    \[ %
 \]
                    
                    \noindent as $\mu$ goes to infinity. Furthermore, the same derivative, this time for $\mu = 0$, has the following asymptotic expansion, as $a$ goes to infinity,
                    
                    \[ %
 . \]
                
                \end{proposition}
                
                \begin{proof}
                
                    We begin by performing an integration by parts, which gives
                    
                    \[ %
 \]
                    
                    \noindent The first term on the right-hand side has already been studied in proposition \ref{PropA1}, where the associated $\mu$-asymptotic behavior is proved, which uses estimates that also give the $a$-asymptotic behavior. Hence, we need only look at the second term. On the interval of integration, using remark \ref{RmkEstimateEta2Tilde}, we have
                    
                    \[ %
 \]
                    
                    \noindent Still on the interval of integration, we have
                    
                    \[ %
 \]
                    
                    \noindent on the half-plane $\Re s > -1/4$.  We can use this estimate to get
                    
                    \[ %
 \]
                    
                    \noindent This proves that the function
                    
                    \[ %
 \]
                    
                    \noindent is holomorphic around $0$, and that its derivative there vanishes. This concludes the proof of this proposition.
                
                \end{proof}
                
                We can now deal with the associated $k = 0$ case, assuming $\alpha$ does not vanish.
                
                \begin{proposition}
                \label{PropMTilde1-0}
                
                    Suppose we have $\alpha \neq 0$. The function
                    
                    \[ %
 \]
                    
                    \noindent is holomorphic on the half-plane $\Re s > -3/2$, and we have, as $\mu$ goes to infinity,
                    
                    \[ %
 . \]
                    
                    \noindent Furthermore, the same derivative, taken with $\mu = 0$, satisfies, as $a$ goes to infinity,
                    
                    \[ %
 . \]
                
                \end{proposition}
                
                \begin{proof}
                
                    The proof is similar to that of proposition \ref{PropMTilde1}, though no series is involved.
                
                \end{proof}

            \subsubsection*{\textbf{Second part}}
            \label{MTildeSecondPart}
                
                We will now deal with the $\arcsinh$ term from lemma \ref{LemLogDerivativeModifiedBessel}.
                
                \begin{proposition}
                \label{PropMTilde2}
                
                    The function
                    
                    \[ %
 , \]
                    
                    \noindent is holomorphic on a half-plane of complex numbers with large enough real part, has a holomorphic continuation near $0$, and we have, as $\mu$ goes to infinity,
                    
                    \[ %
 \]
                    
                    \noindent Furthermore, the same derivative, taken with $\mu = 0$, satisfies, as $a$ goes to infinity,
                    
                    \[ %
 . \]
                
                \end{proposition}
                
                \begin{remark}
                
                    In the $\mu$-asymptotic study above, two derivatives were left uncomputed. They correspond, up to sign, to derivatives left aside in propositions \ref{PropA3} and \ref{PropA4}. In the final result, these uncomputable derivatives will cancel one another.
                
                \end{remark}
                
                \begin{proof}[Proof of proposition \ref{PropMTilde2}]
                
                    We begin by using the binomial formula (the reader is referred to proposition \ref{PropBinomialFormula}), which holds on the interval of integration. We get
                    
                    \[ %
 \]
                    
                    \noindent since we can interchange the sums and the integral. Ultimately, we want to compute each integral using hypergeometric functions. First, an integration by parts yields
                    
                    \[ %
 \]
                    
                    \noindent After summation over $k$ and $j$, we get
                    
                    \begin{equation}
                    \label{EqArcsinhBreakingApart}
                        %
                    \end{equation}
                    
                    \noindent We will now study these three terms separately, beginning with
                    
                    \[ %
 \]
                    
                    \noindent The first step is to split the sum over $k$ into one bearing on positive integers, and another on negative integers. After change of sign, we get
                    
                    \begin{equation}
                    \label{EqArcsinhTwoParts}
                        \scalemath{0.98}{%
}
                    \end{equation}
                    
                    \noindent Using Taylor's formula, the first of these two series yields
                    
                    \[ %
 \]
                    
                    \noindent After multiplication by the appropriate factor, the sum over $j \geqslant 2$ induces a holomorphic function around $0$, whose derivative there vanishes because of the Pochhammer symbol. The term $j=1$ also induces a holomorphic function around~$0$, though its derivative at this point does not vanish entirely, and is given by
                    
                    \[ %
 \]
                    
                    \noindent The term corresponding to $j=0$, \textit{i.e.}
                    
                    \[ %
 \]
                    
                    \noindent has a holomorphic continuation near $0$, as can be seen by using Taylor's formula. The computation as $\mu$ goes to infinity of its derivative at $s=0$ is not necessary, as it is left aside in the statement of this proposition. Thus, we only need to study it as $a$ goes to infinity, with $\mu = 0$. We have
                    
                    \[ %
 \]
                    
                    \noindent from which we get the following asymptotic estimate, as $a$ goes to infinity,
                    
                    \[ %
 . \]
                    
                    \noindent This concludes the study of the first term on the right-hand side of \eqref{EqArcsinhTwoParts}. The second term can be studied similarly, and the results can be obtained by switching the sign of $\alpha$. To sum up what we have seen so far, the function
                    
                    \[ %
 \]
                    
                    \noindent has a holomorphic continuation near $0$, whose derivative there is given by
                    
                    \[ %
 \]
                        
                    \noindent and vanishes as $a$ goes to infinity, when $\mu$ equals zero. We move on to the next term from \eqref{EqArcsinhBreakingApart}, which is, after multiplication by $\sin \left( \pi s \right)/\pi$,
                    
                    \[ %
 \]
                        
                    \noindent First, for any $t$ in the interval of integration, and any integers $j \geqslant 2$, $k \geqslant 1$, we have
                    
                    \[ %
 \]
                    
                    \noindent For any integer $j \geqslant 2$, the function
                    
                    \[ %
 \]
                    
                    \noindent is thus holomorphic near $0$, and for these integers $j$, we further have
                    
                    \[ %
 \]
                    
                    \noindent This proves that the function
                    
                    \[ %
 \]
                    
                    \noindent is holomorphic near $0$. Its derivative there vanishes, because of the Pochhamer symbol. Only the terms $j = 0$ and $j = 1$ remain. We compute the integrals using corollary \ref{CorIntegralFormula} and the change of variable $x = 4 \pi^2 \left( k + \alpha \right)^2 a^2/t^2$. We have
                    
                    \[ %
 \]
                    
                    \noindent Having this formula, we can now take care of both integers $j$ yet to be studied.
                    
                    \medskip
                    
                    \hspace{10pt} $\bullet$ We begin with the case $j = 1$. We consider
                    
                        \begin{equation}   
                        \label{EqArcsinhTermJ1}
                            %
                        \end{equation}
                        
                        \noindent The second of these two parts is simpler, as it does not involve any hypergeometric function, and is thus considered first. It is given by
                        
                        \[ %
 \]
                        
                        \noindent where $\zeta_H$ stands for the Hurwitz zeta function, which is meromorphic on the complex plane, and has a single pole, of order $1$, located at $1$. The term we are studying thus induces a holomorphic function around $0$, and we have
                        
                        \[ %
 \]
                        
                        \noindent where $\psi$ denotes the \textit{Digamma function}. We thus have
                        
                        \[ %
 \]
                        
                        \noindent We move on to the first term of \eqref{EqArcsinhTermJ1}, \textit{i.e.} we consider
                        
                        \[ %
 \]
                        
                        \noindent The hypergeometric function needs to be modified before we can work with it. Using proposition \ref{PropTransformationFraction}, we have, for every integer $k \geqslant 1$,
                        
                        \[ \scalemath{0.94}{%
} \]
                        
                        \noindent the advantage being that the last parameter of the hypergeometric function is now strictly between $0$ and $1$. To further simplify this factor, we need to extract the first term of the hypergeometric series. We use proposition \ref{PropExtractionTermsGeneralizedHypergeometric}, which gives
                        
                        \[ %
 \]
                        
                        \noindent Having $\Re \left( \left( -s+1 \right) - \left( -s \right) - 1/2 \right) = 1/2 > 0$, proposition \ref{PropExtractionTermsGeneralizedHypergeometric} allows us to bound this last generalized hypergeometric function, uniformly in every parameter, for $s$ in some neighborhood of $0$. After multiplying by the appropriate factor from \eqref{EqArcsinhTermJ1}, the associated term induces a holomorphic function around $0$, whose derivative there vanishes, because of the factor $s \sin \left( \pi s \right)$. Hence, we need only deal with
                        
                        \[ %
 \]
                        
                        \noindent We can further simplify the complex power, by writing
                        \[ %
 \]
                        
                        \noindent The term containing the integral behaves nicely around $0$, as we have
                        \[ %
 \]
                        
                        \noindent on the half-plane $\Re s < 1$. The term
                        \[ %
 \]
                        
                        \noindent thus induces a holomorphic function near $0$, and its derivative there vanishes, because of the factor $s \sin \left( \pi s \right)$. The term which remains to be studied is therefore
                        \begin{equation}
                        \label{EqArcsinhRemainingTermJ1}
                            \textstyle \frac{\sin \left( \pi s \right)}{\pi} \cdot \frac{1}{2s+1} \cdot \frac{1}{16 \pi a} \left( 4 \mu + 1 \right)^{-s+1} \sum\limits_{k \geqslant 1} \frac{k^{-2 \delta s}}{k + \alpha}.
                        \end{equation}
                        
                        \noindent We will break apart this series, to relate it to the Riemann zeta function. We have
                        
                        \[ %
 \]
                        
                        \noindent and the second part of the right-hand side induces a holomorphic function around~$0$, whose value at $0$ can be computed by relating it to the constant terms in the Laurant series expansions at $1$ of the Hurwitz and Riemann zeta functions. We have
                        
                        \[ %
 \]
                        
                        \noindent where $\psi$ stands as always for the Digamma function, \textit{i.e.} the logarithmic derivative of the Gamma function, and $\gamma$ is the Euler-Mascheroni constant. After multiplication by the appropriate factor from \eqref{EqArcsinhRemainingTermJ1}, we thus have
                        
                        \[ %
 \]
                        
                        \noindent The part from \eqref{EqArcsinhRemainingTermJ1} involving the Riemann zeta function also induces a holomorphic function around $0$, since the pole of $\zeta$ is canceled by the factor $\sin \left( \pi s \right)$. The derivative at $s=0$ is found by considering the Laurent series expansions. We have
                        
                        \[ %
 \]
                        
                        \noindent The derivative at $s=0$ is thus given, as $\mu$ goes to infinity, by
                        \[ %
 \]
                        
                        \noindent and the beginning of this last computation also shows that the derivative vanishes as $a$ goes to infinity. This concludes the study of the case $j = 1$.

                        \medskip

                    \hspace{10pt} $\bullet$ We now turn to the remaining case, which is $j = 0$. We consider
                    
                        \begin{equation}
                        \label{EqArcsinhTermJ0}
                            %
                        \end{equation}
                        
                        \noindent The second term of \eqref{EqArcsinhTermJ0} is simpler to deal with, and we thus begin by considering
                        
                        \[ %
 \]
                        
                        \noindent This term is closely related to the Hurwitz zeta function, as we have
                        
                        \[ %
 \]
                        
                        \noindent This term is holomorphic function around $0$, and its derivative there is given by
                        
                        \[ %
 \]
                        
                        \noindent This computation being exact, it can used in the $\mu$ and $a$ asymptotic studies. Let us now move on to the first term of \eqref{EqArcsinhTermJ0}, given by
                        
                        \begin{equation}
                        \label{EqArcsinhRemainingTermJ0}
                            %
                        \end{equation}
                        
                        \noindent Unlike the previous case $j=1$, there is no $s$ in the first two parameters of the hypergeometric function which we could extract using proposition \ref{PropExtractionTerms} or \ref{PropExtractionTermsGeneralizedHypergeometric}. To do something like that, we need to lower the second parameter $-s+1$ to $-s$, using proposition \ref{PropContiguousFunctions} and remembering that the first two parameters in a hypergeometric function can be interchanged. We have
                        
                        \[ %
 \]
                        
                        \noindent Note that the first hypergeometric on the right-hand side above is the one we want to study, that the second one can be computed using the binomial formula, here presented as proposition \ref{PropBinomialFormula}, since its first and third parameters are equal, and that the third one contains $-s$ as its first parameter, allowing its simplification using proposition \ref{PropExtractionTermsGeneralizedHypergeometric}. The aim being to make factors $s$ appear, we thus write
                        
                        \[ %
 \]
                        
                        \noindent The hypergeometric function to be studied appears on both sides of this last equality, but the occurence on the right-hand side is much simpler, because of the additional factor $s$. Plugging this into \eqref{EqArcsinhRemainingTermJ0} yields
                        
                        \begin{equation}
                        \label{EqMTildePart2J0Contiguous}
                            %
                        \end{equation}
                        
                        \noindent We begin by studying the third term of \eqref{EqMTildePart2J0Contiguous}, which is
                        
                        \[ %
 \]
                        
                        \noindent Using a Taylor expansion, similarly to what was done in proposition \ref{PropA3}, the function associated to this term has a holomorphic continuation to a neighborhood of $0$, and its derivative there is given by
                        
                        \[ %
 \]
                        
                        \noindent This derivative needs not be computed as $\mu$ goes to infinity, as it is left aside in the statement of the proposition. However, we still need to find its $a$-asymptotic behavior after having set $\mu = 0$. We have
                        
                        \[ %
 \]
                        
                        \noindent as $a$ goes to infinity, using the arguments and computations done in the proof of proposition \ref{PropA3}. We move on to the next term from \eqref{EqMTildePart2J0Contiguous}, which is
                        
                        \[ \scalemath{0.95}{%
} \]
                        
                        \noindent Using proposition \ref{PropTransformationFraction}, we work on the hypergeometric function to turn its last parameter into a real number strictly between $0$ and $1$. We have
                        
                        \[ %
 \]
                        
                        \noindent The term we wish to study therefore becomes
                        \[ %
 \]
                        
                        \noindent This last hypergeometric function not only has its last parameter strictly between~$0$ and $1$, its second parameter equals $1$. This allows us to extract the first few terms of the hypergeometric series using proposition \ref{PropExtractionTerms}. We have
                        
                        \begin{equation}
                        \label{EqMTildePart2HypergeometricBreaking}
                            %
                        \end{equation}
                        
                        \noindent For $s$ in a neighborhood of $0$, we have $\Re \left( \left( -s+4 \right) - 1 - 5/2 \right) = \Re \left( -s + 1/2 \right) > 0$, which means we can bound from above both hypergeometric functions in this last equality, uniformly in every parameter. Thus, the function associated to the term
                        
                        \[ %
 \]
                        
                        \noindent is holomorphic for $s$ around $0$, and its derivative there vanishes, due to the presence of the factor $s \sin \left( \pi s \right)$. The next term from \eqref{EqMTildePart2HypergeometricBreaking} is given by
                        
                        \[ %
 \]
                        
                        \noindent For any integer $k \geqslant 1$, we have
                        
                        \[ %
 \]
                        
                        \noindent This proves the holomorphy around $0$ of the function associated to
                        \[ %
 \]
                        
                        \noindent Furthermore, its derivative at $s=0$ vanishes. Thus, we need only consider
                        
                        \[ %
 \]
                        
                        \noindent For any integer $k \geqslant 1$, we have
                        
                        \[ %
 \]
                        
                        \noindent Thus, the function associated to the term
                        
                        \[ %
 \]
                        
                        \noindent is holomorphic around $0$, and its derivative at $s=0$ vanishes. The next term is
                        
                        \[ %
 \]
                        
                        \noindent which has a holomorphic continuation to an open neighborhood of $0$. Its derivative at this point vanishes. The last term induced by decomposition \eqref{EqMTildePart2HypergeometricBreaking} is given by
                        
                        \[ %
 \]
                        
                        \noindent Using the same computations as those performed in the proof of proposition \ref{PropA3}, this term has a holomorphic continuation to a neighborhood of $0$, and its derivative at $s=0$ vanishes, due to the presence of the extra factor $s$. Going back to \eqref{EqMTildePart2J0Contiguous}, the only remaining part we have yet to take care of is given by
                        
                        \[ %
 \]
                        
                        \noindent Once again, we need to modify the hypergeometric function, so its last parameter becomes a real number strictly between $0$ and $1$. Using proposition \ref{PropTransformationFraction}, we have
                        
                        \[ \scalemath{0.94}{%
} \]
                        
                        \noindent and proposition \ref{PropExtractionTermsGeneralizedHypergeometric} allows us to get the equality
                        
                        \begin{equation}
                        \label{EqMTildePart2GeneralizedHypergeometricBreaking}
                            \scalemath{0.91}{%
}
                        \end{equation}
                        
                        \noindent Having $\Re \left( \left( -s+2 \right) - \left( -s+3/2 \right) - \left( -s \right) \right) = \Re \left( s + 1/2 \right) > 0$ for $s$ in a neighborhood of $0$, the generalized hypergeometric function can be bounded, uniformly in every parameter. This means that the function associated to the term
                        
                        \[ \scalemath{0.92}{%
} \]
                        
                        \noindent is holomorphic on a neighborhood of $0$, and that its derivative at $s=0$ vanishes, because of the extra factor $s$. The next term coming from \eqref{EqMTildePart2GeneralizedHypergeometricBreaking} is given by
                        
                        \[ %
 \]
                        
                        \noindent We will break apart the last two factor above. For any integer $k \geqslant 1$, we have
                        
                        \[ %
 \]
                        \noindent thereby proving that the function associated to
                        
                        \[ \scalemath{0.93}{%
} \]
                        
                        \noindent is holomorphic on a neighborhood of $0$, and its derivative at $s=0$ vanishes, because of the factor $s \sin \left( \pi s \right)$. We are left with studying
                        
                        \[ %
 \]
                        
                        \noindent We will now simplify the complex power, by writing
                        \begin{equation}
                        \label{EqMTildePart2ComplexPower}
                            %
                        \end{equation}
                        
                        \noindent and the integral remainder satisfies the estimate
                        
                        \[ %
 \]
                        
                        \noindent for $s$ in a neighborhood of $0$. This means that the function associated to
                        
                        \[ %
 \]
                        
                        \noindent is holomorphic around $0$, and its derivative at $s=0$ vanishes. The last term induced by the decomposition \eqref{EqMTildePart2ComplexPower} of the complex power is given by
                        
                        \[ %
 \]
                        
                        \noindent For any integer $k \geqslant 1$, we have
                        
                        \[ %
 \]
                        
                        \noindent which proves that the holomorphic function around $s=0$ associated to
                        
                        \[ %
 , \]
                        
                        \noindent has a vanishing derivative at $s=0$. Finally, note that the function associated to
                        
                        \[ \scalemath{0.98}{%
} \]
                        
                        \noindent has a holomorphic continuation near $0$, and that its derivative at $s=0$ vanishes, due to the presence of the extra factor $s$. The last term induced by \eqref{EqMTildePart2GeneralizedHypergeometricBreaking} is
                        
                        \[ %
 \]
                        
                        \noindent Once again, we need to expand the complex power using Taylor's formula. We have
                        
                        \[ %
 \]
                        
                        \noindent and the integral remainder satisfies
                        
                        \[ %
 \]
                        
                        \noindent for $s$ in a neighborhood of $0$. Thus, the function associated to the term
                        
                        \[ %
 \]
                        
                        \noindent is holomorphic around $0$, and its derivative at $s=0$ vanishes, because of the extra factor $s$. Next, we note that
                        
                        \[ %
 \]
                        
                        \noindent has a holomorphic continuation around $0$, and that its derivative at $s=0$ vanishes, using computations previously explained related to the expansion of $1/\left( k + \alpha \right)$. The last term we will study in this proof is thus given
                        
                        \[ %
 \]
                        
                        \noindent Using the known behavior of the Riemann zeta function, this term has a holomorphic continuation around $0$, and its derivative at $s=0$ satisfies
                        
                        \[ %
 \]
                        
                        \noindent To sum up what we have proved, the function associated to
                        
                        \[ %
 \]
                        
                        \noindent has a holomorphic continuation near $0$, and its derivative at $s=0$ satisfies
                        
                        \[ \scalemath{0.98}{%
} \]
                        
                        \noindent as $\mu$ goes to infinity, and, after having set $\mu = 0$, satisfies, as $a$ goes to infinity,
                        
                        \[%
 \]
                        
                        \noindent The same methods can be used to prove that the last part induced by \eqref{EqArcsinhBreakingApart}, which is the function associated to
                        
                        \[ %
 \]
                        
                        \noindent has a holomorphic continuation near $0$, and to compute its derivative at $s=0$, asymptotically as $\mu$ goes to infinity for all $a > 0$, and as $a$ goes to infinity for $\mu = 0$. This contribution can be obtained by switching the sign of $\alpha$ in the results above. The study of all the terms from \eqref{EqArcsinhBreakingApart}, once put together, yield the full proposition.
                    
                \end{proof}
                
                We will take care of the term corresponding to $k=0$, which, as always, is only considered should $\alpha$ not vanish.
                
                \begin{proposition}
                \label{PropMTilde2-0}
                
                    Assume we have $\alpha \neq 0$. The function
                    
                    \[ 
 \]
                    
                    \noindent which is holomorphic on a half-plane consisting of complex numbers with large enough real part, has a holomorphic continuation to a neighborhood of $0$, whose derivative at $s=0$ satisfies
                    
                    \[ %
 \]
                    
                    \noindent as $\mu$ goes to infinity. The same derivative further satisfies, after having set $\mu = 0$,
                    
                    \[ %
 \]
                    
                    \noindent as $a$ goes to infinity.
                
                \end{proposition}
                
                \begin{proof}
                
                    This result can be shown following the reasoning performed in the proof of proposition \ref{PropMTilde2}, though it is simpler here, as there are no series involved.
                
                \end{proof}

            \subsubsection*{\textbf{Third part}}
            \label{MTildeThirdPart}
            
                We are now ready to deal with the third term induced by lemma \ref{LemLogDerivativeModifiedBessel}, which involves a rational fraction in $t$.
                
                \begin{proposition}
                \label{PropMTilde3}
                
                    The function
                    
                    \[ %
 \]
                    
                    \noindent which is holomorphic on a half-plane consisting of complex numbers with large enough real part, has a holomorphic continuation to a neighborhood of $0$, and its derivative at $s=0$ satisfies
                    
                    \[ %
 \]
                    
                    \noindent as $\mu$ goes to infinity. The same derivative further satisfies, after having set $\mu = 0$,
                    
                    \[ %
 \]
                    
                    \noindent as $a$ goes to infinity.
                
                \end{proposition}
                
                \begin{proof}
                
                    Using proposition \ref{PropBinomialFormula}, which is the binomial formula, we have
                    \begin{equation}
                    \label{EqMTildePart3BreakingPositiveNegativeIntegers}
                        %
                    \end{equation}
                    
                    \noindent We begin by studying the first term induced by \ref{EqMTildePart3BreakingPositiveNegativeIntegers}, which is given by
                    
                    \[ %
 \]
                    
                    \noindent For any integer $j \geqslant 2$, we have
                    
                    \[ %
 \]
                    
                    \noindent Therefore, the function
                    \[ %
 \]
                    
                    \noindent is holomorphic around $0$, and its derivative at $s=0$ vanishes, due to the presence of the Pochhammer symbol. Thus, only deal the cases $j = 0$ and $j = 1$ need to be dealt with. First, let us compute the integrals more precisely, using corollary \ref{CorIntegralFormula}. For any $j \in \left \{ 0, 1 \right \}$ and any integer $k \geqslant 1$, we have
                    
                    \[ %
 \]
                    
                    \bigskip

                    \hspace{0.5cm} $\bullet$ We first deal with the case $j = 1$, and thus consider
                    
                        \begin{equation}
                        \label{EqMTildePart3TermJ1}
                            %
                        \end{equation}
                        
                        \noindent The second term above is simpler than the other, since we have
                        
                        \[ %
 \]
                        
                        \noindent where $\zeta_H$ denotes, as always, the Hurwitz zeta function. This term induces a holomorphic function around $0$, and its derivative at $s=0$ is given by
                        
                        \[ %
 \]
                        
                        \noindent This computation, being exact, is used in both the $\mu$-asymptotic and $a$-asymptotic studies. It needs not be more explicit, as it will be canceled by another contribution shortly. We now turn our attention to the first term from \eqref{EqMTildePart3TermJ1}, given by
                        
                        \[ %
 \]
                        
                        \noindent For any integer $k \geqslant 1$, we can use proposition \ref{PropTransformationFraction}, which yields
                        
                        \[ %
 \]
                        
                        \noindent Using proposition \ref{PropHypergeometricZ1}, the hypergeometric function on the right-hand side above can be bounded, uniformly in every parameter, for $s$ in a neighborhood of $0$. This proves that the function associated to the term
                        
                        \[ \scalemath{0.98}{%
} \]
                        
                        \noindent is holomorphic around $0$, and that its derivative at $s=0$ is given by
                        
                        \[ %
 \]
                        
                        \noindent thus cancelling the term found earlier. This concludes the study of the case $j = 1$.
                        
                        \medskip

                    \hspace{10pt} $\bullet$ We now move on to the case $j = 0$, and thus consider
                    
                        \begin{equation}
                        \label{EqMTildePart3TermJ0}
                            \scalemath{0.97}{%
}
                        \end{equation}
                        
                        \noindent We first need to take care of the second term from \eqref{EqMTildePart3TermJ0}, given by
                        
                        \[ %
 \]
                        
                        \noindent This term has a holomorphic continuation near $0$, and its derivative at $s=0$ is
                        
                        \[ %
 \]
                        
                        \noindent This computation is exact, and is thus used in both asymptotic studies, as $\mu$ goes to infinity for all $a > 0$, and as $a$ goes to infinity for $\mu = 0$. We then move on to the first term of \eqref{EqMTildePart3TermJ0}, given by
                        
                        \begin{equation}
                        \label{EqMTildePart3Term1J0}
                            \scalemath{0.98}{%
}
                        \end{equation}
                        
                        \noindent Using proposition \ref{PropTransformationFraction}, we work on the hypergeometric function, to turn its last parameter into a real number strictly between $0$ and $1$. We have
                        
                        \[ \scalemath{0.96}{%
} \]
                        
                        \noindent and we will now relate this last hypergeometric function to its value at $s=0$, using proposition \ref{PropEulerIntegralFormula}, which is known as Euler's integral formula. We have
                        
                        \[ %
 \]
                        
                        \noindent Taking the difference between the left-hand side and its value at $s=0$ then yields
                        
                        \[ %
 \]
                        
                        \noindent this last equality being obtained by using Taylor's formula and Fubini's theorem. We then have the following estimate
                        
                        \[ %
 \]
                        
                        \noindent on the half-plane $\Re s < 1$. Therefore, the function associated to the term
                        
                        \[ %
 \]
                        
                        \noindent is holomorphic around $0$, and its derivative at $s=0$ vanishes, since a factor $s$ appears when expliciting the difference of hypergeometric functions. The study of the term \eqref{EqMTildePart3Term1J0} is then reduced to that of
                        
                        \[ \scalemath{0.94}{%
} \]
                        
                        \noindent Now that the parameters of the hypergeometric function have been simplified, it can be computed explicitely, using proposition \ref{PropExpleLog}. We have
                        
                        \[ %
 \]
                        
                        \noindent and the term we need to study becomes
                        
                        \[ %
 \]
                        
                        \noindent We will now study the second term above, using Taylor's formula. This is similar to the proofs of propositions \ref{PropA3} and \ref{PropA4}. To begin with, let us study the behavior as $a$ goes to infinity, having set $\mu = 0$. We have
                        
                        \[ %
 \]
                        
                        \noindent and each of the integrals above satisfies
                        
                        \[ %
 \]
                        
                        \noindent This means that the term associated to
                        
                        \[ %
 \]
                        
                        \noindent is holomorphic around $0$. Its derivative at $s=0$ further vanishes as $a$ goes to infinity. We now move on to the study as $\mu$ goes to infinity. We have
                        
                        \[ %
 \]
                        
                        \noindent After noting that we have
                        \[ %
 \]
                        
                        \noindent we see that the function associated to
                        
                        \[ %
 \]
                        
                        \noindent is holomorphic around $0$. Furthermore, its derivative at $s=0$ is not to be evaluated as $\mu$ goes to infinity. The first few terms of the Taylor expansion above have a holomorphic continuation near $0$, as we can see by writing
                        
                        \[ %
 \]
                        
                        \noindent It is not necessary to evaluate the derivative at $s=0$ as $\mu$ goes to infinity. Thus, we need only take care of
                        
                        \begin{equation}
                        \label{EqMTildePart3TermJ0LogPositive}
                            %
                        \end{equation}
                        
                        \noindent In order to avoid any unnecessary computation, let us note that the second series induced by \eqref{EqMTildePart3BreakingPositiveNegativeIntegers} can be dealt with similarly, and the contributions to the asymptotics behaviors are obtained by switching the sign of $\alpha$, up to the term
                        
                        \begin{equation}
                        \label{EqMTildePart3TermJ0LogNegative}
                            %
                        \end{equation}
                        
                        \noindent Let us study the terms \eqref{EqMTildePart3TermJ0LogPositive} and \eqref{EqMTildePart3TermJ0LogNegative} added together, and thus consider
                        
                        \[ %
 \]
                        
                        \noindent For any integer $k \geqslant 1$, we have
                        
                        \[ \scalemath{0.9}{%
} \]
                        
                        \noindent and this term therefore induces a holomorphic function around $0$, whose derivative at $s=0$, after multiplication by the relevant factor, is given by
                        
                        \[ %
 \]
                        
                        \noindent It only remains to compute the sum above. We have
                        
                        \[ %
 \]
                        
                        \noindent Let us mention that this last equality is a consequence of the fact that we have
                        
                        \[ %
 \]
                        
                        \noindent and that this formula stems from the power series expansion of the cotangent. Putting these results together yields the proposition, having noted that we have
                        
                        \[ %
 \]
                        
                        \noindent which is a direct consequence of the reflection formula for the Gamma function.
                        
                \end{proof}
                
                When $\alpha$ does not vanish, the case $k = 0$ must also be considered.
                
                \begin{proposition}
                \label{PropMTilde3-0}
                
                    Assume we have $\alpha \neq 0$. The function
                    
                    \[ %
 \]
                    
                    \noindent which is holomorphic on a half-plane consisting of complex numbers with large enough real part, has a holomorphic continuation to a neighborhood of $0$, whose derivative at $s=0$ satisfies
                    
                    \[ %
 \]
                    
                    \noindent as $\mu$ goes to infinity. The same derivative further satisfies, after having set $\mu = 0$,
                    
                    \[ %
 \]
                    
                    \noindent as $a$ goes to infinity.
                
                \end{proposition}
                
                \begin{proof}
                
                    This result can be proved along the line of proposition \ref{PropMTilde3}, the arguments being simpler as there are no series involved.
                
                \end{proof}

            \subsubsection*{\textbf{Fourth Part}}
            \label{MTildeFourthPart}
            
                Let us deal with the last two terms from \ref{LemLogDerivativeModifiedBessel}, which are given by
                
                \[ %
 \]
                    
                    \noindent which is holomorphic on a half-plane consisting of complex numbers with large enough real part, has a holomorphic continuation to a neighborhood of $0$, and its derivative at $s=0$ satisfies, as $\mu$ goes to infinity,
                    
                    \[ %
 \]
                    
                    \noindent After having set $\mu = 0$, the same derivative satisfies, as $a$ goes to infinity,
                    
                    \[ %
 \]
                
                \end{proposition}
                
                \begin{proof}
                
                    The proof is similar to that of propositions \ref{PropMTilde2} and \ref{PropMTilde3}.
                
                \end{proof}
                
                In this final case as well, we must deal with $k = 0$, should $\alpha$ not vanish.
                
                \begin{proposition}
                \label{PropMTilde4-0}
                
                    Assume we have $\alpha \neq 0$. The function
                    
                    \[ %
 \]
                    
                    \noindent which is holomorphic on a half-plane consisting of complex numbers with large enough real part, has a holomorphic continuation to a neighborhood of $0$, whose derivative at $s=0$ satisfies, as $\mu$ goes to infinity,
                    
                    \[ %
 \]
                    
                    \noindent After having set $\mu = 0$, the same derivative satisfies, as $a$ goes to infinity,
                    
                    \[ %
 \]
                
                \end{proposition}
                
                \begin{proof}
                
                    The proof is similar to that of propositions \ref{PropMTilde2} and \ref{PropMTilde3}.
                
                \end{proof}

    \subsection{Asymptotic study of the determinant of the pseudo-Laplacian}
    
        Having performed all the necessary computations, we can now state the asymptotic behavior of the determinant of the pseudo-Laplacian on a cuspidal end with a flat unitary line bundle, for the Dirichlet boundary conditions. This is done in four theorems, having two asymptotic studies, and both cases $\alpha = 0$ and $\alpha \neq 0$, which correspond respectively to a trivial line bundle $L$, and a non trivial one.

        \subsubsection{As $\mu$ goes to infinity, for all $a>0$}
        
            We begin by putting together all the results from section \ref{SecDetPseudoLaplacian} regarding the $\mu$-asymptotic expansion. The proofs only refer to the relevant propositions, whose contributions should all be summed directly, any sign or coefficient having already been taken into account.
        
            \begin{theorem}
            \label{ThmAlphaNeq0Mu}
        
                Assume we have $\alpha \neq 0$. We have, as $\mu$ goes to infinity,
            
                \[ \begin{array}{lll}
                        \multicolumn{3}{l}{\log \det \left( \Delta_{L,0} + \mu \right)} \\[0.5em]
                    
                        & = & - \frac{1}{4 \pi a} \mu \log \mu + \frac{1}{4 \pi a} \mu + \sqrt{\mu} \log \mu \\[0.5em]
                    
                        && \quad - \big[ 2 \int_0^{+ \infty} \frac{1}{e^{2 \pi t} - 1} \left( \arctan \left( \frac{t}{1 + \alpha} \right) + \arctan \left( \frac{t}{1 - \alpha} \right) \right) \mathrm{d}t - \log 2 \\[0.5em]
                    
                        && \qquad \quad + \alpha \log \left( \frac{1 + \alpha}{1 - \alpha} \right) + \frac{1}{4a} + \frac{1}{2} \log \left( 4 \pi^2 \left( 1 - \alpha^2 \right) a^2 \right) + \log \left( \pi \alpha a \right) \big] \sqrt{\mu} \\[0.5em]
                    
                        && \qquad \qquad \quad - \frac{3}{4} \log \mu + o \left( 1 \right).
                    \end{array} \]
        
            \end{theorem}
        
            \begin{proof}
        
                This result is a combination of propositions \ref{PropA1}, \ref{PropA1-0}, \ref{PropA2}, \ref{PropA2-0}, \ref{PropA3}, \ref{PropA3-0}, \ref{PropA4}, \ref{PropA4-0}, \ref{PropA5}, \ref{PropA5-0}, \ref{PropA7}, \ref{PropA7-0}, \ref{PropA8}, \ref{PropA8-0}, \ref{PropA9}, \ref{PropA9-0}, \ref{PropA6and10}, \ref{PropA10-0}, \ref{PropA11}, \ref{PropA11-0}, \ref{PropR}, \ref{PropR-0}, \ref{PropMTilde1}, \ref{PropMTilde1-0}, \ref{PropMTilde2}, \ref{PropMTilde2-0}, \ref{PropMTilde3}, \ref{PropMTilde3-0}, \ref{PropMTilde4}, \ref{PropMTilde4-0}.
        
            \end{proof}
        
            \begin{theorem}
            \label{ThmAlpha0Mu}
        
                Assume we have $\alpha = 0$. We have, as $\mu$ goes to infinity,
            
                \[ \begin{array}{lll}
                        \multicolumn{3}{l}{\log \det \left( \Delta_{L,0} + \mu \right)} \\[0.5em]
                    
                        & = & - \frac{1}{4 \pi a} \mu \log \mu + \frac{1}{4 \pi a} \mu + \frac{1}{2} \sqrt{\mu} \log \mu \\[0.5em]
                    
                        && \qquad - \big[ 4 \int_0^{+ \infty} \frac{1}{e^{2 \pi t} - 1} \arctan \left( t \right) \mathrm{d}t - \log 2 + 1 + \frac{1}{4a} + \log \left( 2 \pi a \right) \big] \sqrt{\mu} \\[0.5em]
                    
                        && \qquad \qquad - \frac{1}{2} \log \mu + o \left( 1 \right).
                    \end{array} \]
        
            \end{theorem}
        
            \begin{proof}
        
                This result is a combination of propositions \ref{PropA1}, \ref{PropA2}, \ref{PropA3}, \ref{PropA4}, \ref{PropA5}, \ref{PropA7}, \ref{PropA8}, \ref{PropA9}, \ref{PropA6and10}, \ref{PropA11}, \ref{PropR}, \ref{PropMTilde1}, \ref{PropMTilde2}, \ref{PropMTilde3}, \ref{PropMTilde4}.
        
            \end{proof}

        \subsubsection{As $a$ goes to infinity, with $\mu = 0$}
        
            We now deal with the $a$-asymptotic study, for which we only consider the case $a = 0$. Once again, only the relevant propositions and results are given in the proofs, and all the relevant contributions are to be added directly, with every sign and coefficient having been taken into account.
            
            \begin{theorem}
            \label{ThmAlphaNeq0A}
            
                Assume we have $\alpha \neq 0$. We have, as $a$ goes to infinity,
            
                \[ \scalemath{0.96}{\begin{array}{lll}
                        \log \det \Delta_{L,0} & = & 2 \pi \alpha^2 a - 2 \pi \alpha a + \frac{\pi}{3} a - \frac{1}{2} \log \frac{\sin \left( \pi \alpha \right)}{\pi \alpha} - \frac{1}{2} \log \left( 2 \pi \alpha \right) + o \left( 1 \right).
                    \end{array}} \]
            
            \end{theorem}
            
            \begin{proof}
            
                This is a direct consequence of propositions \ref{PropAmukAsymptA}, \ref{PropAmukAsymptA-0}, \ref{PropR}, \ref{PropR-0}, \ref{PropMTilde1}, \ref{PropMTilde1-0}, \ref{PropMTilde2}, \ref{PropMTilde2-0}, \ref{PropMTilde3}, \ref{PropMTilde3-0}, \ref{PropMTilde4}, \ref{PropMTilde4-0}.
            
            \end{proof}
            
            \begin{theorem}
            \label{ThmAlpha0A}
            
                Assume we have $\alpha = 0$. We have, as $a$ goes to infinity,
            
                \[ \begin{array}{lll}
                        \log \det \Delta_{L,0} & = & \frac{\pi}{3} a + \frac{1}{2} \log a + o \left( 1 \right).
                    \end{array} \]
            
            \end{theorem}
            
            \begin{proof}
            
                This is a direct consequence of propositions \ref{PropAmukAsymptA}, \ref{PropR}, \ref{PropMTilde1}, \ref{PropMTilde2}, \ref{PropMTilde3}, \ref{PropMTilde4}.
            
            \end{proof}

\appendix

\section{Self-adjoint operators}
\label{AppSelfAdjointOp}

    The aim of this appendix is to quickly present the notion of self-adjoint operators between Hilbert spaces, as well as some useful results. For more details, the reader is referred to \cite[Appendix C]{MR2339952} and in \cite{MR1335452}, which delves deeply in the theory of linear operators. We denote by $H$ a separable Hilbert space, by $\left< \cdot, \cdot \right>$ its inner product, and by $\left \Vert \cdot \right \Vert$ the associated norm.
    
    \subsection{Quadratic forms}
    
        Before we can move to the study of self-adjoint extensions of symmetric operators, we need to review the notion of quadratic forms and their relationship to self-adjoint operators.
    
        \begin{definition}
    
            A \textit{quadratic form} on $H$ is the datum of a vector subspace $V$ of~$H$, and of a sesquilinear map $Q : V \times V \longrightarrow \mathbb{C}$. The subspace $V$ is called the \textit{domain} of $Q$, and denoted by $\mathrm{Dom} \, Q$.
    
        \end{definition}
        
        For the remainder of this paragraph, we consider such a quadratic form $Q$.
        
        \begin{definition}
        \label{DefQuadFormNorm}
        
            The quadratic form $Q$ is said to be \textit{positive} if, for all $x \in \mathrm{Dom} \, Q$, we have $Q \left( x, x \right) \geqslant 0$. In this case, it is further said to be \textit{closed} if its domain is complete for the norm
            
            \[ \begin{array}{lll}
                    \left \Vert x \right \Vert_Q & = & \sqrt{Q \left( x,x \right) + \left \Vert x \right \Vert^2}
                \end{array} . \]
        
        \end{definition}
        
        \begin{definition}
        
            Let $T$ be a positive symmetric operator. The \textit{associated quadratic form} $Q_T$ is defined on $\mathrm{Dom} \, Q_t = \mathrm{Dom} \, T$ by $Q_T \left( x,y \right) = \left< Tx, y \right>$.
        
        \end{definition}
        
        \begin{proposition}
        
            Let $T$ be a positive symmetric operator. The set of positive closed quadratic forms extending $Q_T$ has a smallest element in terms of inclusion of domains, denoted by $\overline{Q_T}$.
        
        \end{proposition}
        
        \begin{proof}
        
            This is stated as proposition $C.1.6$ in \cite[Appendix C]{MR2339952}. Let us simply note that $\mathrm{Dom} \, \overline{Q_T}$ is the completion of $\mathrm{Dom} \, Q_T$ with respect to $\left \Vert \cdot \right \Vert_Q$.
        
        \end{proof}
        
        \begin{proposition}
        
            Let $Q$ be a closed positive quadratic form. Then there exists a positive self-adjoint operator $T$ such that we have $Q = \overline{Q_T}$.
        
        \end{proposition}
        
        \begin{proof}
        
            This is proposition $C.1.5$ from \cite[Appendix C]{MR2339952}.
        
        \end{proof}

    \subsection{Spectrum of self-adjoint operators}
    
        In this section, we will study the \textit{spectrum} of a densely defined positive self-adjoint operator $T$ on the Hilbert space~$H$. For more general considerations, the reader is referred to \cite{MR1335452}.
        
        \begin{definition}
        
            The \textit{resolvent set} of $T$ is defined as
            
            \[ \begin{array}{lll}
                    \rho \left( T \right) & = & \left \{ \lambda \in \mathbb{C}, \; T - \lambda : \mathrm{Dom} \, T \longrightarrow H \text{ is bijective} \right \}
                \end{array} , \]
            
            \noindent and its spectrum as the complement set $\Sigma \left( T \right) = \mathbb{C} \setminus \rho \left( T \right)$.
        
        \end{definition}
        
        \begin{definition}
        
            A complex number $\lambda$ is said to be an \textit{eigenvalue} of $T$ if the operator $T - \lambda : \mathrm{Dom} \, T \longrightarrow H$ is not injective, in which case its kernel is said to be an \textit{eigenspace} of $T$, whose dimension is called the \textit{multiplicity} of the eiganvalue $\lambda$.
        
        \end{definition}
        
        \begin{remark}
        
            Unlike the finite dimensional case, there can be complex numbers $\lambda$ for which $T - \lambda : \mathrm{Dom} \, T \longrightarrow H$ is injective, but not surjective. Thus, the spectrum is not limited to the eigenvalues, even including those of infinite multiplicity.
        
        \end{remark}
        
        \begin{proposition}
        
            The spectrum of $T$ is included in $\mathbb{R}_+$.
        
        \end{proposition}
        
        \begin{proof}
        
            The spectrum of a self-adjoint operator is real, as in \cite[Sec. 5.3.4]{MR1335452}. Let us show that its elements are positive. For any real number $\lambda < 0$, we have
            
            \[ \begin{array}{ccccc}
                    \left< \left( T - \lambda \right) u, u \right> & \geqslant & - \lambda \left \Vert u \right \Vert^2 & = & \left \vert \lambda \right \vert \left \Vert u \right \Vert^2
                \end{array} \]
            
            \noindent for any $u \in \mathrm{Dom} \, T$, and the Cauchy-Schwarz inequality yields
            
            \[ \begin{array}{ccc}
                    \left \vert \lambda \right \vert \left \Vert u \right \Vert & \leqslant & \left \Vert \left( T - \lambda \right) u \right \Vert
                \end{array} . \]
            
            \noindent The kernel of $T - \lambda$ is thus reduced to $0$, and we further have
            
            \[ \begin{array}{ccccc}
                    \mathrm{Im} \left( T - \lambda \right) & = & \ker \left( T - \lambda \right)^{\perp} & = & H
                \end{array} , \]
            
            \noindent so $\lambda$ cannot be in the spectrum of $T$, thereby completing the proof of the proposition.
        
        \end{proof}
        
        \begin{definition}
        
            The \textit{discrete sprectrum} of $T$, denoted by $\Sigma_{\mathrm{dis}} \left( T \right)$, is defined as the set of eigenvalues $\lambda$ of $T$ which have finite multiplicity and are isolated, in the sense that there is a neighborhood of $\lambda$ in $\mathbb{C}$ disjoint from $\Sigma \left( T \right)$.
        
        \end{definition}
        
        \begin{definition}
        
            The \textit{essential spectrum} of $T$, denoted by $\Sigma_{\mathrm{ess}} \left( T \right)$ is defined as the complement set of $\Sigma_{\mathrm{dis}} \left( T \right)$ in the spectrum, \textit{i.e.} as $\Sigma_{\mathrm{ess}} \left( T \right) = \Sigma \left( T \right) \setminus \Sigma_{\mathrm{dis}} \left( T \right)$.
        
        \end{definition}
        
        \begin{definition}
        
            The lower bound of $\Sigma_{\mathrm{ess}} \left( T \right)$ is denoted by $\sigma_{\mathrm{ess}} \left( T \right)$. If the essential spectrum is empty, this lower bound is set at $+ \infty$.
        
        \end{definition}
        
        \begin{theorem}[Inf-sup theorem]
        \label{ThmInfSup}
        
            For any positive integer $n \in \mathbb{N}^{\ast}$, the quantity
            
            \[ \begin{array}{lll}
                    \mu_n \left( T \right) & = & \inf\limits_{\psi_1, \dots, \psi_n \in \mathrm{Dom} \, \overline{Q_T}} \; \sup \left \{ \frac{\overline{Q_T} \left( \psi, \psi \right)}{\left< \psi, \psi \right>}, \; \psi \in \mathrm{span} \left( \psi_1, \dots, \psi_n \right), \; \psi \neq 0 \right \}
                \end{array} , \]
            
            \noindent is well-defined, the infimum being taken on linearly independent elements. This sequence is increasing, covers the eigenvalues of $T$ lying below $\sigma_{\mathrm{ess}} \left( T \right)$ with multiplicity, and becomes constant equal to this lower bound if $\Sigma_{\mathrm{ess}} \left( T \right)$ is non-empty.
        
        \end{theorem}
        
        \begin{proof}
        
            This is theorem $4.5.2$ from \cite{MR1349825}, with a small and direct adjustment, since we consider $\overline{Q_T} \left( \psi, \psi \right)$ and not $\left< T \psi, \psi \right>$ in the definition of $\mu_n \left( T \right)$.
        
        \end{proof}
        
        \begin{definition}
        \label{DefSpectralCountingFunction}
        
            The \textit{spectral counting function} is defined for any $\lambda > 0$ by
            
            \[ \begin{array}{lll}
                    N \left( T, \lambda \right) & = & \# \left \{ n \in \mathbb{N}^{\ast}, \; \mu_n \left( T \right) \leqslant \lambda \right \}
                \end{array} . \]
        
        \end{definition}
        
        \begin{remark}
        \label{RmkSpectralCountingEssentialSpectrum}
        
            The spectral counting function becomes infinite if and only if the essential spectrum is non-empty.
        
        \end{remark}

    \subsection{Self-adjoint extensions}
    
        It should be noted that a symmetric operator may have several self-adjoint extensions, but comparing them cannot be done by looking at their domains. Instead, we will define an order on the set of self-ajoint operators, which also for instance be found in \cite[Def I.5.4]{MR1102675}.
        
        \begin{definition}
        \label{DefOrderSelfAdjoint}
        
            Let $T_1$ and $T_2$ be two positive self-adjoint operators on $H$. We write $T_1 \preccurlyeq T_2$ if we have the inclusion $\mathrm{Dom} \, \overline{Q_{T_2}} \subseteq \mathrm{Dom} \, \overline{Q_{T_1}}$ and the inequality
            
            \[ \begin{array}{lll}
                    \overline{Q_{T_1}} \left( x, x \right) & \leqslant & \overline{Q_{T_2}} \left( x, x \right)
                \end{array} \]
            
            \noindent for every $x \in \mathrm{Dom} \, \overline{Q_{T_2}}$.
        
        \end{definition}
        
        \begin{definition}
        
            Let $T$ be a symmetric positive operator. The \textit{Friedrichs extension} of $T$ is the only self-adjoint extension $T_F$ of $T$ such that we have $\overline{Q_T} = \overline{Q_{T_F}}$.
        
        \end{definition}
        
        \begin{remark}
        \label{RmkDomainFriedrichs}
        
            The domain of the Friedrichs extension can be expressed as
            
            \[ \begin{array}{lll}
                    \mathrm{Dom} \, T_F & = & \mathrm{Dom} \, \overline{Q_T} \cap \mathrm{Dom} \, T^{\ast}
                \end{array} . \]
        
        \end{remark}
        
        \begin{proposition}
        
            The Friedrichs extension of a symmetric positive operator is its maximal self-adjoint extension with respect to $\preccurlyeq$.
        
        \end{proposition}
        
        \begin{proof}
        
            Let $T$ be a symmetric positive operator, denote by $T_F$ its Friedrichs extension, and by $T_{SA}$ any self-adjoint extension of $T$. We have
            
            \[ \begin{array}{ccccccccc}
                    \mathrm{Dom} \, Q_T & = & \mathrm{Dom} \, T & \subseteq & \mathrm{Dom} \, T_{SA} & = & \mathrm{Dom} \, Q_{T_{SA}} & \subseteq & \mathrm{Dom} \overline{Q_{T_{SA}}}.
                \end{array} \]
            
            \noindent The closed quadratic form $\overline{Q_{T_{SA}}}$ being an extension of $Q_T$, we have
            
            \[ \begin{array}{ccccc}
                    \mathrm{Dom} \overline{Q_{T_F}} & = & \mathrm{Dom} \overline{Q_T} & \subseteq & \mathrm{Dom} \overline{Q_{T_{SA}}},
                \end{array} \]
            
            \noindent and we further have $\overline{Q_{T_{SA}}} \left( x,x \right) = \overline{Q_{T_F}} \left( x,x \right)$ on the domain of $\overline{Q_{T_F}}$.
        
        \end{proof}
        
        \begin{proposition}
        \label{PropSpectralCountingOrderRelation}
        
            Let $T_1$ and $T_2$ be two self-adjoint operators, with $T_1 \preccurlyeq T_2$. We have $\mu_n \left( T_1 \right) \leqslant \mu_n \left( T_2 \right)$ for any integer $n > 0$. As a consequence, for any real number $\lambda > 0$, we have $N \left( T_2, \lambda \right) \leqslant N \left( T_1, \lambda \right)$.
        
        \end{proposition}
        
        \begin{proof}
        
            Let $n > 0$ be a positive integer and $\psi_1$, \dots, $\psi_n$ be elements of $\mathrm{Dom} \, \overline{Q_{T_2}}$. For any non-zero element $\psi$ in the vector space spanned by every $\psi_i$, we have
            
            \[ \begin{array}{lll}
                    \frac{\overline{Q_{T_1}} \left( \psi, \psi \right)}{\left< \psi, \psi \right>} & \leqslant & \frac{\overline{Q_{T_2}} \left( \psi, \psi \right)}{\left< \psi, \psi \right>}
                \end{array} \]
            
            \noindent since we assumed $T_1 \preccurlyeq T_2$. After taking the upper bound, we get
            
            \[ \begin{array}{lll}
                    \multicolumn{3}{l}{\sup \left \{ \frac{\overline{Q_{T_1}} \left( \psi, \psi \right)}{\left< \psi, \psi \right>}, \; \psi \in \mathrm{span} \left( \psi_1, \dots, \psi_n \right), \; \psi \neq 0 \right \}} \\[1em]
                    
                    \qquad \qquad \qquad \qquad \qquad & \leqslant & \sup \left \{ \frac{\overline{Q_{T_2}} \left( \psi, \psi \right)}{\left< \psi, \psi \right>}, \; \psi \in \mathrm{span} \left( \psi_1, \dots, \psi_n \right), \; \psi \neq 0 \right \}.
                \end{array} \]
            
            \noindent Taking the lower bound on elements $\psi_1, \dots, \psi_n \in \mathrm{Dom} \, \overline{Q_{T_2}}$, together with the inclusion of domains $\mathrm{Dom} \, \overline{Q_{T_2}} \subseteq \mathrm{Dom} \, \overline{Q_{T_1}}$, gives the result. The second part, related to the spectral counting function, is a direct consequence.
        
        \end{proof}
        
        For the remainder of this appendix, we consider two Hilbert spaces $H_1$ and $H_2$, with respective inner products $\left< \cdot, \cdot \right>_1$ and $\left< \cdot, \cdot \right>_2$. We also consider two operators $T_1$ and $T_2$ on $H_1$ and $H_2$.
        
        \begin{definition}
        
            The direct sum operator $T_1 \oplus T_2$ is defined as
            
            \[ \begin{array}{lllll}
                    T_1 \oplus T_2 & : & \mathrm{Dom} \, T_1 \oplus \mathrm{Dom} \, T_2 & \longrightarrow & H_1 \oplus H_2 \\[0.5em]
                    && \left( x_1, x_2 \right) & \longmapsto & \left( T_1 x_1, T_2 x_2 \right)
                \end{array} . \]
        
        \end{definition}
        
        \begin{proposition}
        \label{PropDomainDirectSum}
        
            Assume that $T_1$ and $T_2$ are self-adjoint. Then $T_1 \oplus T_2$ is also self-adjoint. Furthermore, we have
            
            \[ \begin{array}{lll}
                    \mathrm{Dom} \, \overline{Q_{T_1 \oplus T_2}} & = & \mathrm{Dom} \, \overline{Q_{T_1}} \oplus \mathrm{Dom} \, \overline{Q_{T_2}}
                \end{array} . \]
        
        \end{proposition}
        
        \begin{proof}
        
            This is a classical result.
        
        \end{proof}
        
        \begin{proposition}
        \label{PropSpectralCountingDirectSum}
        
            Let $\lambda > 0$ be a strictly positive real number. We have
            
            \[ \begin{array}{lll}
                    N \left( T_1 \oplus T_2, \lambda \right) & = & N \left( T_1, \lambda \right) + N \left( T_2, \lambda \right)
                \end{array} . \]
        
        \end{proposition}
        
        \begin{proof}
        
            Consider a real number $\lambda > 0$. We have
            
            \[ \begin{array}{lll}
                    \ker \left( T_1 \oplus T_2 - \lambda \id \right) & = & \ker \left( T_1 - \lambda \id \right) \oplus \ker \left( T_2 - \lambda \id \right)
                \end{array} \]
            
            \noindent Furthermore, the operator $T_1 \oplus T_2 - \lambda \id$ is surjective if and only if both $T_1 - \lambda \id$ and $T_2 - \lambda \id$ are. Put together, these two facts yield the equality
            
            \[ \begin{array}{lll}
                    \Sigma_{\mathrm{ess}} \left( T_1 \oplus T_2 \right) & = & \Sigma_{\mathrm{ess}} \left( T_1 \right) \cup \Sigma_{\mathrm{ess}} \left( T_2 \right)
                \end{array} . \]
            
            \noindent The proof of the proposition is then completed by combining these results.
        
        \end{proof}

\section{The Ramanujan summation}
\label{AppRamanujan}

    The main reference here is \cite{MR3677185}, where Candelpergher presents a technique, known as the \textit{Ramanujan summation}, to study series of holomorphic functions. Its purpose is to prove the existence of meromorphic continuations of such functions, and affords a greater control on these extensions than Taylor's formula.
    
    \subsection{Introduction to the method}
    
        At the core of the Ramanujan summation lies the idea of comparing a sum of values taken by a function to the integral of said function. This notion is embodied in two well-known formulae.
        
        \begin{theorem}[Euler--Maclaurin formula]
        
            Let $f : \left[ a,b \right] \longrightarrow \mathbb{C}$ be a $\mathcal{C}^p$ function defined on a segment whose endpoints are integers. We have
            
            \[ \scalemath{0.97}{\begin{array}{lll}
                    \sum\limits_{k=a}^b \; f \left( k \right) & = & \int_a^b \; f \left( x \right) \; \mathrm{d}x + \frac{1}{2} \left( f \left( a \right) + f \left( b \right) \right) + \sum\limits_{k=2}^p \; \frac{1}{k!} B_k \left( f^{\left( k-1 \right)} \left( b \right) - f^{\left( k-1 \right)} \left( a \right) \right) \\[1em]
                    
                    &&  \qquad + \left( -1 \right)^{p+1} \int_a^b \frac{1}{p!} b_p \left( x - \left[ x \right] \right) f^{\left( p \right)} \left( x \right) \mathrm{d}x
                \end{array}} \]
            
            \noindent where $b_p \left( x \right)$ is the $p$-th Bernoulli polynomial, and $B_k$ denotes the $k$-th Bernoulli number, with the convention $B_1 = 1/2$.
        
        \end{theorem}
        
        \begin{proof}
        
            This formula is proved in \cite[Sec. 1.1]{MR3677185}, though it should be noted that the factor $\left( -1 \right)^k$ in the sum is not necessary, as Bernoulli numbers with odd indices all vanish except $B_1$.
        
        \end{proof}
        
        While this formula only deals with integrals over segments, the second one, presented below, tackles the problem of integrals over unbounded intervals.
        
        \begin{theorem}[Abel--Plana formula]
        
            Let $f$ be a holomorphic function on the half-plane $\Re z > 1$, continuous up to the boundary, such that we have
            
            \[ \begin{array}{ccc}
                    \left \vert z \right \vert^{1+\varepsilon} \left \vert f \left( z \right) \right \vert & \leqslant & C
                \end{array} \]
            
            \noindent for some $\varepsilon > 0$ and $C > 0$. We then have
            
            \[ \begin{array}{ccc}
                    \sum\limits_{k=1}^{+ \infty} \; f \left( k \right) & = & \int_1^{+\infty} \; f \left( x \right) \; \mathrm{d}x + \frac{1}{2} f \left( 1 \right) + i \int_0^{+ \infty} \; \frac{f \left( 1+it \right) - f \left( 1-it \right)}{e^{2 \pi t} - 1} \; \mathrm{d}t
                \end{array} . \]
        
        \end{theorem}
        
        \begin{proof}
        
            As we will soon see, this result is close to the Ramanujan summation, and is presented here for historical purposes. Nevertheless, it is proved in \cite[Sec. 1.4.2]{MR3677185} and also in \cite[Sec.8.3.1]{MR1429619}, under weaker assumptions than those stated above.
        
        \end{proof}

    \subsection{Functions with exponential growth}
    
        We will now define a growth condition on functions under which the last integral on the right-hand side of the Abel--Plana formula makes sense.
        
        \begin{definition}
        \label{DefExpModGrowth}
        
            A holomorphic function $f$ defined on a half-plane $\Re z > a$ for some real number $0 < a < 1$ is said to be \textit{of exponential type} at most $\alpha > 0$ if there is a constant $0 < \beta < \alpha$ such that $f$ satisfies a bound
            
            \[ \begin{array}{ccc}
                    \left \vert f \left( z \right) \right \vert & \leq & C e^{\beta \left \vert z \right \vert}
                \end{array} \]
            
            \noindent for some constant $C > 0$. This space is denoted by $\mathcal{O}^{\alpha}$. Such a function is further said to be of \textit{moderate growth} if it is exponential of type $\varepsilon$ for all $\varepsilon > 0$.
        
        \end{definition}
        
        \begin{remark}
        
            For an element $f \in \mathcal{O}^{2 \pi}$, the function
            
            \[ \begin{array}{ccc}
                    t & \longmapsto & \frac{f \left( 1+it \right) - f \left( 1-it \right)}{e^{2 \pi t} - 1}
                \end{array} \]
            
            \noindent is integrable on the unbounded interval $\left] 0, + \infty \right[$.
        
        \end{remark}
    
        \begin{theorem}[Carlson]
        
            Let $f \in \mathcal{O}^{\pi}$. Assume that we have $f \left( k \right) = 0$ for every integer $k > 0$. Then the function $f$ vanishes identically.
        
        \end{theorem}
        
        \begin{proof}
        
            This theorem is proved in \cite[App. B]{MR3677185}.
        
        \end{proof}

    \subsection{The Ramanujan summation}
    
        \begin{definition}
        
            Let $f \in \mathcal{O}^{2 \pi}$. The \textit{Ramanujan sum} of $f$ is defined as
            
            \[ \begin{array}{lll}
                    \sum\limits_{k \geqslant 1}^{\left( \mathcal{R} \right)} \; f \left( k \right) & = & \frac{1}{2} f \left( 1 \right) + i \int_0^{+ \infty} \; \frac{f \left( 1 + it \right) - f \left( 1 - it \right)}{e^{2 \pi t} - 1} \; \mathrm{d}t
                \end{array} . \]
        
        \end{definition}
        
        \begin{remark}
        
            Unlike the traditionnal sum, this Ramanujan sum does not depend only on the value of $f$ at integers, but on the whole function $f$. However, using Carlon's theorem, the function $f$ is entirely determined by its value at integers if it is assumed to be of exponential type at most $\pi$.
        
        \end{remark}
    
        \begin{theorem}
        \label{ThmRamanujanSum}
        
            Let $f$ be an element of $\mathcal{O}^{\pi}$ such that we have
            
            \[ \begin{array}{lll}
                    \lim\limits_{k \rightarrow + \infty} \; f \left( k \right) & = & 0,
                \end{array} \]
            
            \noindent and further assume that we have
            
            \[ \begin{array}{lll}
                    \lim\limits_{k \rightarrow + \infty} \; \int_0^{+ \infty} \; \frac{f \left( k + it \right) + f \left( k - it \right)}{e^{2 \pi t} - 1} \; \mathrm{d}t & = & 0.
                \end{array} \]
            
            \noindent Then $f$ is integrable over $\left] 1, + \infty \right[$ if and only if the series of general term $f \left( k \right)$ is absolutely convergent, and in this case, we have
            
            \[ \begin{array}{lll}
                    \sum\limits_{k=1}^{+ \infty} \; f \left( k \right) & = & \int_1^{+\infty} f \left( x \right) \; \mathrm{d}x + \sum\limits_{k \geqslant 1}^{\left( \mathcal{R} \right)} \; f \left( k \right)
                \end{array} . \]
        
        \end{theorem}
        
        \begin{proof}
        
            This is the content of \cite[Sec. 1.4.3]{MR3677185}.
        
        \end{proof}
        
        \begin{remark}
        
            As explained by Candelpherger in \cite[Sec. 1.4.3]{MR3677185}, this result is obtained by proving we can write
            
            \[ \begin{array}{lll}
                    \sum\limits_{k=1}^n f \left( k \right) & = & R_f \left( n \right) - R_f \left( 1 \right)
                \end{array} , \]
            
            \noindent where $R_f$ is the function defined by
            
            \[ \begin{array}{lllll}
                    R_f & : & x & \longmapsto & - \int_1^x f \left( t \right) \; \mathrm{d}t + \frac{1}{2} f \left( x \right) + i \int_0^{+ \infty} \; \frac{f \left( x + it \right) + f \left( x - it \right)}{e^{2 \pi t} - 1} \; \mathrm{d}t.
                \end{array} \]
            
            \noindent Writing the partial sum above can be done by using the residue formula, and constitutes theorem $1$ from \cite[Sec. 1.3.2]{MR3677185}.
        
        \end{remark}

    \subsection{Properties of the Ramanujan summation}
        
        In order to state the difference between the Ramanujan sum and the classical sum, and to get a regularity result, we need to introduce the notion of \textit{functions locally uniformly in} $\mathcal{O}^{\pi}$.
        
        \begin{definition}
        
            Let $U$ be an open subset of $\mathbb{C}$ and $a \in \left] 0,1 \right[$. A function
            
            \[ \begin{array}{lllll}
                    f & : & \left \{ z \in \mathbb{C}, \; \Re z > a\right \} \times U & \longrightarrow & \mathbb{C}
                \end{array} \]
            
            \noindent is said to be \textit{locally uniformly in} $\mathcal{O}^{\pi}$ if $f$ is holomorphic in $z$, and if, for any compact subset $K$ of $U$, there exists a real number $\beta$ with $0 < \beta < \pi$ and a constant $C > 0$ such that we have
            
            \[ \begin{array}{ccc}
                    \left \vert f \left( z, s \right) \right \vert & \leqslant & C e^{\beta \left \vert z \right \vert}
                \end{array} . \]
            
            \noindent This space is denoted by $\mathcal{O}^{\pi}_{\mathrm{loc}} \left( U \right)$.
        
        \end{definition}
        
        \begin{theorem}
        \label{ThmAnalycityRamanujanSum}
        
            Let $U$ be an open subset of the complex plane, and $f \in \mathcal{O}^{\pi}_{\mathrm{loc}} \left( U \right)$. Furthermore, assume that $f$ is holomorphic in $s$ on $U$. The function
            
            \[ \begin{array}{lll}
                    s & \longmapsto & \sum\limits_{k \geqslant 1}^{\left( \mathcal{R} \right)} \; f \left( k, s \right)
                \end{array} \]
            
            \noindent is then holomorphic on $U$, and one may differentiate term by term.
        
        \end{theorem}
        
        \begin{proof}
        
            This result constitutes theorem $9$ of \cite[Sec. 3.1.1]{MR3677185}, and relies on the dominated convergence theorem.
        
        \end{proof}

\section{Special functions}

    This appendix is devoted to compiling information on two particular types of functions: the modified Bessel functions, and the hypergeometric functions. We will begin by reviewing the notion of total variation of a function, which is used in the study of Bessel functions.
    
    \subsection{Total variation of a function}
    
        In this section, we will follow Olver's presentation from \cite[Sec. 1.11]{MR1429619}.
        
        \begin{definition}
        
            Let $f$ be a real function defined on a segment $\left[ a,b \right]$. The \textit{total variation} of $f$ between $a$ and $b$ is defined as
            
            \[ \begin{array}{ccc}
                    V_{a,b} \left( f \right) & = & \sup\limits_{a=a_0 < a_1 < \dots < a_{n-1} < a_n = b} \; \sum\limits_{k=0}^{n-1} \; \left \vert f \left( a_{k+1} \right) - f \left( a_k \right) \right \vert
                \end{array} , \]
            
            \noindent the supremum being taken on all subdivisions of the segment $\left[ a,b \right]$. If this quantity is finite, then the function $f$ is said to be of \textit{bounded variation}.
        
        \end{definition}
        
        \begin{proposition}
        \label{PropTotalVariationReal}
        
            Assume $f$ is a $\mathcal{C}^1$-function on $\left[ a,b \right]$. Then we have
            
            \[ \begin{array}{ccc}
                    V_{a,b} \left( f \right) & = & \int_a^b \; \left \vert f' \left( x \right) \right \vert \; \mathrm{d}x
                \end{array} . \]
        
        \end{proposition}
        
        \begin{proof}
        
            This follows from the Taylor-Lagrange theorem (\textit{i.e.} with mean-value remainder), and the definition of the Riemann integral.
        
        \end{proof}
        
        \begin{remark}
        
            If $f$ is a complex-valued differentiable function, its variation on $\left[ a,b \right]$ is defined by the integral in the proposition above. Once again, the function is said to be of \textit{bounded variation} if the integral converges.
        
        \end{remark}
        
        \begin{definition}
        
            Let $f$ be a holomorphic function on an open domain $U$ of the complex plane, and $z$, $w$ two points in $U$. Consider a piecewise smooth path $\gamma \left( z,w \right)$ joining $z$ to $w$, parametrized by $t \in \left[ a,b \right] \mapsto z \left( t \right)$. The \textit{variation} of $f$ along $\gamma \left( z,w \right)$ is defined as
            
            \[ \begin{array}{ccc}
                    V_{\gamma \left( z,w \right)} \left( f \right) & = & \sum\limits_{k=0}^{n-1} \; \int_{t_k}^{t_{k+1}} \; \left \vert f' \left( z \left( t \right) \right) z' \left( t \right) \right \vert \; \mathrm{d}t
                \end{array} , \]
            
            \noindent where $t_0 = a < t_1 < \dots < t_{n-1} < t_n = b$ are the points of $\left[ a,b \right]$ corresponding to the non-smooth points of $\gamma \left( z,w \right)$.
        
        \end{definition}
        
        \begin{proposition}
        
            As in the definition above, let $\gamma \left( z,w \right)$ be a piecewise smooth curve joining two points $z$ and $w$. We have
            
            \[ \begin{array}{ccc}
                    V_{\gamma \left( z,w \right)} \left( f \right) & = & \sup \sum\limits_{k=0}^{n-1} \; \left \vert f \left( z_{k+1} \right) - f \left( z_k \right) \right \vert
                \end{array} , \]
            
            \noindent where the supremum is taken on points $z_0 = z$, $z_1$, \dots, $z_{n-1}$, $z_n = w$ arranged in the order in which they are reached by $\gamma \left( z, w\right)$.
        
        \end{proposition}
        
        \begin{proof}
        
            This form of the total variation can be obtained by applying the same results used in the proof of proposition \ref{PropTotalVariationReal} on each $\left[ t_k, t_{k+1} \right]$.
        
        \end{proof}

    \subsection{Modified Bessel functions}
    \label{AppModifiedBessel}
    
        We begin with the first category of functions we need in this paper: modified Bessel functions. Consider the differential equation
        
        \begin{equation}
        \label{EquationModifiedBessel}
            \begin{array}{lll}
                z^2 \frac{\mathrm{d}^2 w}{\mathrm{d}z^2} + z \frac{\mathrm{d} w}{\mathrm{d} z} - \left( z^2 + \nu^2 \right) w & = & 0,
            \end{array}
        \end{equation}
        
        \noindent where $\nu$ is a complex number. The \textit{modified Bessel functions}, defined as particular solutions of this equation, are studied in \cite[Sec. 2.10 \& 7.8]{MR1429619}, with useful asymptotic properties in \cite[Sec. 10.7]{MR1429619}. These objects are also dealt with in \cite[Sec. 10.25]{MR2723248}.
        
        \begin{definition}
        
            The \textit{modified Bessel function of the second kind} $K_{\nu}$ is defined as the only solution of \ref{EquationModifiedBessel} with, as $z$ goes to infinity in the sector $\left \vert \arg z \right \vert < \pi / 2$,
            
            \[ \begin{array}{lll}
                    K_{\nu} \left( z \right) & \sim & \sqrt{\frac{\pi}{2z}} e^{-z}
                \end{array} . \]
            
            \noindent The complex number $\nu$ is called the \textit{order} of the modified Bessel function.
        
        \end{definition}
        
        \begin{proposition}
        \label{PropIntegralRepModifiedBessel}
        
            The modified Bessel function of the second kind admits the following integral representation
            
            \[ \begin{array}{lll}
                    K_{\nu} \left( x \right) & = & \frac{\sqrt{\pi}}{\Gamma \left( 1/2 + \nu \right)} \left( \frac{1}{2} x \right)^{\nu} \int_1^{+ \infty} e^{-xt} \left( t^2 - 1 \right)^{\nu - 1/2} \mathrm{d}t
                \end{array} \]
            
            \noindent for any $\nu \geqslant 1/2$ and any $x > 0$. As a consequence, we have, for every $x > 0$,
            
            \[ \begin{array}{lll}
                    K_{1/2} \left( x \right) & = & \sqrt{\frac{\pi}{2x}} e^{-x}
                \end{array} . \]
        
        \end{proposition}
        
        \begin{proof}
        
            This is given as exercice $8.4$ in \cite[Sec. 7.8]{MR1429619}, where it is hinted that one should use the previously proved integral representation of the Hankel functions. It should be noted that this is result is stated there for complex order and parameters, and that the extension to $\nu = 1/2$ can be made by continuity since only real parameters are involved here. The special value at $1/2$ is a direct consequence.
        
        \end{proof}
        
        \begin{proposition}
        \label{PropDerivativeOneHalfModifiedBessel}
        
            For any $x > 0$, the logarithmic derivative at $1/2$ with respect to the order of the modified Bessel function of the second kind $K_{\nu} \left( x \right)$ is given by
            
            \[ \begin{array}{lll}
                    \frac{\partial}{\partial \nu}_{\left \vert \nu = 1/2 \right.} \log K_{\nu} \left( x \right) & = & \mathbb{E}_1 \left( 2x \right) e^{2x}
                \end{array} , \]
            
            \noindent where $\mathbb{E}_1$ is the exponential integral. Hence we have, for every integer $N \geqslant 0$,
            
            \[ \begin{array}{lll}
                    \frac{\partial}{\partial \nu}_{\left \vert \nu = 1/2 \right.} \log K_{\nu} \left( x \right) & = & \frac{1}{2x} \sum\limits_{n = 0}^{N} \left( -1 \right)^n n! \, x^{-n} + O \left( x^{-N-1} \right)
                \end{array} . \]
        
        \end{proposition}
        
        \begin{proof}
        
            The computation of the order-derivative at $1/2$ of $K_{\nu} \left( x \right)$ is detailed in \cite{MR92871}, based of the integral representation from proposition \ref{PropIntegralRepModifiedBessel}. Together with the special value $K_{1/2} \left( x \right)$ for $x > 0$ given above, we get the first part of the proposition. Let us now prove the rest. A change of variable gives, for every $x > 0$,
            
            \[ \begin{array}{lllll}
                    \frac{\partial}{\partial \nu}_{\left \vert \nu = 1/2 \right.} \log K_{\nu} \left( x \right) & = & \mathbb{E}_1 \left( 2x \right) e^{2x} & = & \int_0^{+ \infty} \frac{e^{-t}}{t + 2x} \mathrm{d}t
                \end{array} , \]
            
            \noindent which we can we use to get the required asymptotic expansion by induction.
        
        \end{proof}
        
        \begin{proposition}
        \label{PropAppZerosModifiedBessel}
        
            Let $u > 0$ be a strictly positive real number. The function
            
            \[ \begin{array}{lll}
                    \nu & \longmapsto & K_{i \nu} \left( u \right)
                \end{array} \]
            
            \noindent is entire, has only simple zeros, all of whom are located on the real line.
        
        \end{proposition}
        
        \begin{proof}
        
            The argument that follows is adapted from \cite[Appendix A]{MR2439244}, where Saharian deals with  Legendre functions. Using the differential equation \ref{EquationModifiedBessel} satisfied by modified Bessel functions of the second kind, we have
            
            \[ \begin{array}{lll}
                    K_{\beta} \left( z \right) & = & K_{- \beta} \left( z \right)
                \end{array} . \]
            
            \noindent Furthermore, we see that $K_{\beta} \left( t \right)$ is real whenever $\beta$ and $t > 0$ are real. The Schwarz reflection principle then gives the following identities
            
            \[ \begin{array}{lllll}
                    \overline{K_{i \nu} \left( u \right)} & = & K_{- i \overline{\nu}} \left( u \right) & = & K_{i \overline{\nu}} \left( u \right)
                \end{array} . \]
            
            \noindent Using once more equation \ref{EquationModifiedBessel}, we have
            
            \[ \begin{array}{lll}
                    u K_{i \overline{\nu}}' \left( u \right) & = & \left( u^2 - \overline{\nu}^2 \right) K_{i \overline{\nu}} \left( u \right) - u^2 K_{i \overline{\nu}}'' \left( u \right), \\[0.5em]
                    
                    u K_{i \nu}' \left( u \right) & = & \left( u^2 - \nu^2 \right) K_{i \nu} \left( u \right) - u^2 K_{i \nu}'' \left( u \right),
                \end{array} \]
            
            \noindent the derivatives being taken with respect to the parameter. We get
            
            \[ \begin{array}{lll}
                    \multicolumn{3}{l}{u \left[ K_{i \nu} \left( u \right) K_{i \overline{\nu}}' \left( u \right) -  K_{i \overline{\nu}} \left( u \right) K_{i \nu}' \left( u \right) \right]} \\[0.5em] 
                    
                    \qquad & = & \left( \nu^2 - \overline{\nu}^2 \right) K_{i \nu} \left( u \right) K_{i \overline{\nu}} \left( u \right) + u^2 \left[ K_{i \overline{\nu}} \left( u \right) K_{i \nu}'' \left( u \right) - K_{i \nu} \left( u \right) K_{i \overline{\nu}}'' \left( u \right) \right].
                \end{array} \]
            
            \noindent If $\nu$ is neither real nor purely imaginary, we have $\nu^2 - \overline{\nu}^2 \neq 0$, and get
            
            \[ \begin{array}{lll}
                    \frac{1}{v} K_{i \nu} \left( v \right) K_{i \overline{\nu}} \left( v \right) & = & - \frac{1}{\overline{\nu}^2 - \nu^2} \left[ K_{i \nu} \left( v \right) K_{i \overline{\nu}}' \left( v \right) -  K_{i \overline{\nu}} \left( v \right) K_{i \nu}' \left( v \right) \right] \\[0.5em]
                    
                    && \qquad + \frac{v}{\overline{\nu}^2 - \nu^2} \left[ K_{i \overline{\nu}} \left( v \right) K_{i \nu}'' \left( v \right) - K_{i \nu} \left( v \right) K_{i \overline{\nu}}'' \left( v \right) \right]
                \end{array} \]
            
            \noindent for every real number $v$ between $0$ and $u$. After integrating for $v$ real between $0$ and $u$, as well as integrating by parts, we get
            
            \begin{equation}
            \label{FormulaModifiedBesselAwayFromAxes}
                \begin{array}{lll}
                    \int_0^u \; \frac{1}{v} \left \vert K_{i \nu} \left( v \right) \right \vert^2 \; \mathrm{d}v & = & \int_0^u \; \frac{1}{v} K_{i \nu} \left( v \right) K_{i \overline{\nu}} \left( v \right) \; \mathrm{d}v \\[0.5em]
                    
                    & = & - \frac{u}{\overline{\nu}^2 - \nu^2} \left[ K_{i \nu} \left( u \right) K_{i \overline{\nu}}' \left( u \right) -  K_{i \overline{\nu}} \left( u \right) K_{i \nu}' \left( u \right) \right].
                \end{array}
            \end{equation}
            
            \noindent If $\nu \in\mathbb{C}$ is such that we have $K_{i \nu} \left( u \right) = 0$, then we also have $K_{i \overline{\nu}} \left( u \right) = 0$ by the conjugation properties stated above. This gives
            
            \[ \begin{array}{lll}
                    \int_0^u \; \frac{1}{v} \left \vert K_{i \nu} \left( v \right) \right \vert^2 \; \mathrm{d}v & = & 0
                \end{array} , \]
            
            \noindent which implies that we have $K_{i \nu} \left( v \right) = 0$ for every $v \in \left[ 0, u \right]$. This is absurd, since the function $z \longmapsto K_{i \nu} \left( z \right)$ is holomorphic and does not vanish identically. Thus, the zeros of the function $\nu \longmapsto K_{i \nu} \left( u \right)$ are either real or purely imaginary. However, modified Bessel functions of the second kind are strictly positive when both the order and the parameter are real. Hence, the zeros of $\nu \longmapsto K_{i \nu} \left( u \right)$ can only be real. We will now prove that such zeros, if they exist, can only be simple. To that effect, we note that formula \ref{FormulaModifiedBesselAwayFromAxes} gives, by a difference quotient argument,
            
            \begin{equation}
            \label{FormulaModifiedBesselRealOrder}
                \begin{array}{lll}
                    \int_0^u \; \frac{1}{v} \left \vert K_{i \nu} \left( v \right) \right \vert^2 \; \mathrm{d}v & = & - i \frac{u}{2 \nu} \left[ K_{i \nu} \left( u \right) \frac{\partial}{\partial \beta}_{\vert \beta = i \nu} K_{\beta}' \left( u \right) - K_{i \nu}' \left( u \right) \frac{\partial}{\partial \beta}_{\vert \beta = i \nu} K_{\beta} \left( u \right) \right]
                \end{array}
            \end{equation}
            
            \noindent for any non-zero real number $\nu$. If such $\nu \in \mathbb{R}^{\ast}$ was a zero of order at least two, then we would have
            
            \[ \begin{array}{lllll}
                    K_{i \nu} \left( u \right) & = & \frac{\partial}{\partial \beta}_{\vert \beta = i \nu} K_{\beta} \left( u \right) & = & 0
                \end{array} , \]
            
            \noindent and the integral on the left-hand side of formula \ref{FormulaModifiedBesselRealOrder} vanishes. This is absurd, since the function $z \longmapsto K_{i \nu} \left( z \right)$ is holomorphic and does not vanish identically. The proof of the proposition is now complete.
        
        \end{proof}
        
        In this paper, we need two types of asymptotics for modified Bessel functions of the second kind: one for large orders, and one for large parameters, both with some control of the remainder. Let us now see these results, which can be found in \cite{MR1429619}.
        
        \begin{definition}
        \label{Defpxi}
        
            Let $\delta > 0$. We define the cone $C_{\delta}$ by
            
            \[ \begin{array}{lll}
                    C_{\delta} & = & \left \{ z \in \mathbb{C}, \; \; \left \vert \arg z \right \vert \leqslant \frac{\pi}{2} - \delta \right \}
                \end{array} . \]
            
            \noindent We now define two functions $p$ and $\xi$ by
            
            \[ \begin{array}{cccclcccccl}
                    p & : & C_{\delta} & \longmapsto & \mathbb{C} & \quad \text{and} \quad & \xi & : & C_{\delta} & \longmapsto & \mathbb{C} \\
                    && z & \longmapsto & \frac{1}{\sqrt{1 + z^2}} & & && z & \longmapsto & \sqrt{1 + z^2} + \log \frac{z}{1 + \sqrt{1 + z^2}}
                \end{array} , \]
            
            \noindent where the functions $\sqrt{\cdot}$ and $\log$ denote the principal branches of the square root and the logarithm, in keeping with equations $\left( 7.07 \right)$ and $\left( 7.09 \right)$ of \cite[Sec. 10.7.3]{MR1429619}.
        
        \end{definition}
        
        \begin{definition}
        
            The polynomials $U_k$ are defined inductively by $U_0 = 1$ and
            
            \[ \begin{array}{lll}
                    U_{k+1} \left( t \right) & = & \frac{1}{2} t^2 \left( 1 - t^2 \right) U_k' \left( t \right) + \frac{1}{8} \int_0^t \; U_k \left( x \right) \; \mathrm{d}x
                \end{array} . \]
        
        \end{definition}
        
        \begin{remark}
        
            This definition is presented in equation $\left( 7.10 \right)$ of \cite[Sec. 10.7.3]{MR1429619}. In this paper, we need the first three terms of this sequence, which are explicitely given in $\left( 7.11 \right)$ of \cite[Sec. 10.7.3]{MR1429619}. We have
            
            \[ \begin{array}{lll}
                    U_1 \left( t \right) & = & \frac{1}{24} \left( 3t - 5t^3 \right), \\[0.5em]
                    U_2 \left( t \right) & = & \frac{1}{1152} \left( 81 t^2 - 462 t^4 + 385 t^6 \right).
                \end{array} \]
        
        \end{remark}
        
        \begin{proposition}
        \label{PropAsymptoticExpansionModifiedBessel}
        
            For any integer $n \in \mathbb{N}^{\ast}$ and fixed parameter $z \in C_{\delta}$, we have
            
            \[ \begin{array}{lll}
                    K_ {\nu} \left( \nu z \right) & = & \sqrt{\frac{\pi}{2 \nu}} \cdot \frac{e^{- \nu \xi \left( z \right)}}{\left( 1 + z^2 \right)^{1/4}} \cdot \left[ \; \sum\limits_{k=0}^{n-1} \left( -1 \right)^k \frac{1}{\nu^k} U_k \left( p \left( z \right) \right) + \eta_n \left( \nu, z \right) \right]
                \end{array} , \]
            
            \noindent where the remainder $\eta_n \left( \nu, z \right)$ satisfies the bound
            
            \[ \begin{array}{lll}
                    \left \vert \eta_n \left( \nu, z \right) \right \vert & \leqslant & \frac{2}{\nu^n} \exp \left[ \, \frac{2}{\nu} V_{\gamma \left( 0, p \left( z \right) \right)} \left( U_1 \right) \, \right]  V_{\gamma \left( 0, p \left( z \right) \right)} \left( U_n \right)
                \end{array} \]
            
            \noindent for any $\xi$-progressive path $\gamma$ joining $0$ to $p \left( z \right)$, that is any path between these points for which $\Re \xi \left( z \right)$ is increasing.
        
        \end{proposition}
        
        \begin{proof}
        
            This proposition is proved in \cite[Sec. 10.7]{MR1429619}.
        
        \end{proof}
        
        \begin{remark}
        
            This proposition will be used for real parameters $z > 0$, for which the considered path is the segment from $0$ to $p \left( z \right)$. The definition of the variation along such paths coincides with the total variation of a function of a real variable.
        
        \end{remark}
        
        This paper actually calls for a more restrictive version of this asymptotic expansion, which we state here for clarity.
        
        \begin{corollary}
        \label{CorAsymptoticExpansionModifiedBessel}
        
            For every $x, \nu > 0$, we can write
            
            \[ \begin{array}{lll}
                    \log K_{\nu} \left( \nu x \right) & = & \frac{1}{2} \log \frac{\pi}{2 \nu} - \nu \xi \left( x \right) - \frac{1}{4} \log \left( 1 + x^2 \right) - \frac{1}{\nu} U_1 \left( p \left( x \right) \right) + \widetilde{\eta_2} \left( \nu, x \right)
                \end{array} . \]
                
            \noindent If we have $\nu \geqslant A$ or $\nu x \geqslant B$ large enough, then the remainder $\widetilde{\eta_2}$ satisfies
            
            \[ \begin{array}{lll}
                    \left \vert \widetilde{\eta_2} \left( \nu, x \right) \right \vert & \leqslant & C \min \left( \frac{1}{\nu^2 x^2}, \frac{1}{\nu^2} \right)
                \end{array} , \]
            
            \noindent where $C > 0$ does not depend on $x$ or $\nu$, but depends on $A$ or $B$.
        
        \end{corollary}
        
        \begin{proof}
        
            Taking the logarithm of the expansion from proposition \ref{PropAsymptoticExpansionModifiedBessel}, we get
            
            \[ \begin{array}{lll}
                    \log K_{\nu} \left( \nu x \right) & = & \frac{1}{2} \log \frac{\pi}{2 \nu} - \nu \xi \left( x \right) - \frac{1}{4} \log \left( 1 + x^2 \right) \\[0.5em]
                    
                    && \qquad \qquad \qquad \qquad \qquad + \log \left( 1 - \frac{1}{\nu} U_1 \left( p \left( x \right) \right) + \eta_2 \left( \nu, x \right) \right).
                \end{array} \]
            
            \noindent We can then set
            
            \[ \begin{array}{lll}
                    \widetilde{\eta_2} \left( \nu, x \right) & = & \frac{1}{\nu} U_1 \left( p \left( x \right) \right) + \log \left( 1 - \frac{1}{\nu} U_1 \left( p \left( x \right) \right) + \eta_2 \left( \nu, x \right) \right)
                \end{array} , \]
                
            \noindent which gives the required formula. We now need to prove the important part of the result, which is the bound on $\widetilde{\eta_2}$. First, we note that we have
            
            \[ \begin{array}{lllllll}
                    \left \vert \frac{1}{\nu} U_1 \left( p \left( x \right) \right) \right \vert & = & \left \vert \frac{1}{24 \nu} \left( 3 p \left( x \right) - 5 p \left( x \right)^3 \right) \right \vert & \leqslant & \frac{1}{\nu} p \left( x \right) & \leqslant & \min \left( \frac{1}{\nu x}, \frac{1}{\nu} \right)
                \end{array} , \]
            
            \noindent since we have $p \left( x \right) = \left( 1 + x^2 \right)^{-1/2}$. Furthermore, we have
            
            \[ \begin{array}{lll}
                    \left \vert \eta_2 \left( \nu, x \right) \right \vert & \leqslant & \frac{2}{\nu^2} \exp \left[ \frac{2}{\nu} V_{0,1} \left( U_1 \right) \right] V_{0,p\left( x \right)} \left( U_2 \right) \\[0.5em]
                    
                    & \leqslant & \frac{2}{\nu^2} \exp \left[ \frac{2}{A} V_{0,1} \left( U_1 \right) \right] \cdot \frac{1}{1152} \left \vert 81 p \left( x \right)^2 - 462 p \left( x \right)^4 + 385 p \left( x \right)^6 \right \vert \\[0.5em]
                    
                    & \leqslant & \underbrace{2 \exp \left[ 2 V_{0,1} \left( U_1 \right) \right]}_{A'} \cdot \min \left( \frac{1}{\nu^2 x^2}, \frac{1}{\nu^2} \right), 
                \end{array} \]
            
            \noindent assuming we have $A \geqslant 1$, and where the total variation is understood in the real sense of proposition \ref{PropTotalVariationReal}. If $A$ is large enough so as to have $A'/\nu^2 + 1/\nu \leqslant 1/2$ for every $\nu \geqslant A$, then we have in particular
            
            \[ \begin{array}{lll}
                    \left \vert - \frac{1}{\nu} U_1 \left( p \left( x \right) \right) + \eta_2 \left( \nu, x \right) \right \vert & < & 1
                \end{array} , \]
            
            \noindent and one can then use the power series expansion of the logarithm to get
            
            \[ \begin{array}{lll}
                    \widetilde{\eta_2} \left( \nu, x \right) & = & \eta_2 \left( \nu, x \right) + \sum\limits_{n=2}^{+ \infty} \frac{\left( -1 \right)^{n+1}}{n} \left[ \eta_2 \left( \nu, x \right) - \frac{1}{\nu} U_1 \left( p \left( x \right) \right) \right]^n
                \end{array} . \]
            
            \noindent This allows us to properly bound $\widetilde{\eta_2}$, as we have
            
            \[ \begin{array}{lllll}
                    \multicolumn{3}{l}{\left \vert \widetilde{\eta_2} \left( \nu, x \right) \right \vert} & \leqslant & \left \vert \eta_2 \left( \nu, x \right) \right \vert + \sum\limits_{n=2}^{+ \infty} \frac{1}{n} \left[ \left \vert \eta_2 \left( \nu, x \right) \right \vert + \left \vert \frac{1}{\nu} U_1 \left( p \left( x \right) \right) \right \vert \right]^n \\[0.5em]
                    
                    \quad & \leqslant & \multicolumn{3}{l}{A' \min \left( \frac{1}{\nu^2 x^2}, \frac{1}{\nu^2} \right) + \frac{1}{2} \left[ A' \min \left( \frac{1}{\nu^2 x^2}, \frac{1}{\nu^2} \right) + \min \left( \frac{1}{\nu x}, \frac{1}{\nu} \right) \right]^2 \sum\limits_{n=0}^{+ \infty} \left[ \frac{A'}{\nu^2} + \frac{1}{\nu} \right]^n} \\[0.5em]
                    
                    & \leqslant & \multicolumn{3}{l}{\min \left( \frac{1}{\nu^2 x^2}, \frac{1}{\nu^2} \right) \left[ A' + \frac{1}{2} \left( A' \min \left( \frac{1}{\nu x}, \frac{1}{\nu} \right) + 1 \right)^2 \right] \sum\limits_{n=0}^{+ \infty} \left[ \frac{A'}{\nu^2} + \frac{1}{\nu} \right]^n} \\[0.5em]
                    
                    & \leqslant & \multicolumn{3}{l}{\min \left( \frac{1}{\nu^2 x^2}, \frac{1}{\nu^2} \right) \left[ A' + \frac{1}{2} \left( \frac{A'}{\nu} + 1 \right)^2 \right] \sum\limits_{n=0}^{+ \infty} 2^{-n}} \\[0.5em]
                    
                    & \leqslant & \multicolumn{3}{l}{\underbrace{2 \left[ A' + \frac{1}{2} \left( A'A + 1 \right)^2 \right]}_{C} \min \left( \frac{1}{\nu^2 x^2}, \frac{1}{\nu^2} \right).}
                \end{array} \]
            
            \noindent This completes the proof of the corollary, the case $\nu x \geqslant B$ large enough being entirely similar.
        
        \end{proof}
        
        \begin{remark}
        \label{RmkEstimateEta2Tilde}
        
            In the proof above, bounding $p \left( x \right)^2$ by $1/\left( 1 + x^2 \right)$ instead of $1/x^2$, actually gives the estimate
            
            \[ \begin{array}{lll}
                    \left \vert \widetilde{\eta_2} \left( \nu, x \right) \right \vert & \leqslant & \frac{C}{\nu^2 \left( 1 + x^2 \right)}
                \end{array} . \]
            
            \noindent It is not stated as such in the corollary, as this variant is used in a much smaller portion of this paper.
        
        \end{remark}
        
        \begin{proposition}
        \label{PropAsymptoticsModifiedBesselParameter}
        
            For any integer $n \in \mathbb{N}^{\ast}$, and any order $\nu$, we have
            
            \[ \begin{array}{lll}
                    K_{\nu} \left( z \right) & = & \sqrt{\frac{\pi}{2 z}} e^{-z} \left[ \; \sum\limits_{k = 0}^{n-1} \frac{1}{z^k} A_k \left( \nu \right) + \gamma_n \left( \nu, z \right) \; \right]
                \end{array} , \]
            
            \noindent for $z \in C_{\delta}$, where each polynomial $A_k$ is defined by
            
            \[ \begin{array}{lll}
                    A_k \left( \nu \right) & = & \frac{1}{8^k k!} \prod\limits_{j=1}^k \; \left( 4 \nu^2 - \left( 2j-1 \right)^2 \right)
                \end{array} , \]
            
            \noindent and the remainder $\gamma_n$ satisfies the bound
            
            \[ \begin{array}{lll}
                    \left \vert \gamma_n \left( \nu, z \right) \right \vert & \leqslant & 2 \left \vert \frac{A_n \left( \nu \right)}{z^n} \right \vert \exp \left( \left \vert \frac{1}{z} \left( \nu^2 - \frac{1}{4} \right) \right \vert \right)
                \end{array} . \]
            
        \end{proposition}
        
        \begin{proof}
        
            This result constitutes exercise $13.2$ from \cite[Sec. 7.13]{MR1429619}, and follows from the expansion of Hankel's functions which are given there.
        
        \end{proof}

    \subsection{Hypergeometric functions}
    \label{AppHypergeometric}
    
        The last part of this appendix is devoted to the presentation of the \textit{hypergeometric functions}, which are used throughout this paper. We will follow \cite[Sec. 5.9]{MR1429619}, though the required content can be found, without proofs, in \cite[Sec. 15]{MR2723248}.
        
        \subsubsection{Hypergeometric series}
        
            One of the easiest introduction to hypergeometric functions is through the \textit{hypergeometric series}.
        
            \begin{definition}
        
                Let $s$ be a complex number, and $k$ be a positive integer. Assume neither $s$ nor $s+k$ is a negative integer. The \textit{Pochhammer symbol} is defined as
            
                \[ \begin{array}{lllll}
                        \left( s \right)_k & = & \frac{\Gamma \left( s+k \right)}{\Gamma \left( s \right)} & = & s \left( s+1 \right) \dots \left( s+k \right)
                    \end{array} \]
        
            \end{definition}
        
            \begin{remark}
        
                The notation adopted here is the one used in \cite{MR1429619, MR2723248}. However, these symbols are sometimes denoted by $s^{\left( k \right)}$ and referred to as \textit{rising factorials}, with the notation $\left( s \right)_k$ being reserved for the so-called \textit{falling factorials}.
        
            \end{remark}
        
            \begin{remark}
        
                Using the fact the Gamma function has a single pole at every negative integer, and is holomorphic elsewhere, the definition above can be extended to the case where $s$ is a negative integer, still with $k$ positive.
        
            \end{remark}
        
            \begin{proposition-definition}
        
                Let $a$, $b$, and $c$ be complex numbers. Assume $c$ is not a negative integer. The hypergeometric series
            
                \[ \begin{array}{lll}
                        F \left( a, b; c; z \right) & = & \sum\limits_{n=0}^{+\infty} \; \frac{\left( a \right)_n \left( b \right)_n}{\left( c \right)_n} \cdot \frac{z^n}{n!}
                    \end{array} \]
            
                \noindent converges on the disk $\left \vert z \right \vert < 1$, where it induces a hypergeometric function.
        
            \end{proposition-definition}
            
            \begin{proof}
            
                This result directly follows from d'Alembert's ratio test for series, and from the fact that we have $\left( a \right)_{n+1} = \left( a+n \right) \left( a \right)_n$.
            
            \end{proof}
        
            \begin{remark}
            \label{RmkHypergeometricInterchangeAB}
        
                In this definition, the parameters $a$ and $b$ can be interchanged. Furthermore, the following notation is sometimes used in Olver's work
                
                \[ \begin{array}{lll}
                        \bm{F} \left( a, b; c; z \right) & = & \frac{1}{\Gamma \left( c \right)} F \left( a, b; c; z \right)
                    \end{array} ,\]
                
                \noindent and will never be used here. The results proved in \cite{MR1429619} will be adapted when needed to remain with the more classical notation $F$.
        
            \end{remark}

        \subsubsection{The hypergeometric differential equation}
        
            The functions defined above naturally satisfy an ordinary differential equation.
        
            \begin{proposition}
        
                Let $a$, $b$, and $c$ be complex numbers. Assume $c$ is not a negative integer. The hypergeometric function $z \longmapsto F \left( a, b; c; z \right)$ satisfies
                
                \begin{equation}
                \label{EqHypergeometric}
                    \begin{array}{lll}
                        z \left( 1 - z \right) \frac{\mathrm{d}^2 w}{\mathrm{d}z^2} + \left[ c - \left( a + b + 1 \right) z \right] \frac{\mathrm{d}w}{\mathrm{d}z} - ab \, w & = & 0
                    \end{array} .
                \end{equation}
        
            \end{proposition}
        
            \begin{proof}
        
                This can be checked directly on the unit disk $\left \vert z \right \vert < 1$, using the definition of the Pochhammer symbol as a quotient of Gamma functions, and the known property $\Gamma \left( s+1 \right) = s \Gamma \left( s \right)$, which translates into $\left( s \right)_{n+1} = \left( s+n \right) \left( s \right)_n$.
        
            \end{proof}
            
            \begin{corollary}
        
                The hypergeometric function $z \longmapsto F \left( a, b; c; z \right)$ has a holomorphic continuation to $\mathbb{C} \setminus \left[ 1, + \infty \right[$.
        
            \end{corollary}
            
            \begin{proof}
            
                The regularity with regard to the parameter $z$ of solutions to \eqref{EqHypergeometric} is studied in \cite[Sec. 5.3.1]{MR1429619}.
            
            \end{proof}
            
            \begin{corollary}
        
                For any $z \in \mathbb{C} \setminus \left[ 1, +\infty \right[$, the hypergeometric function $F \left( a, b; c; z \right)$ is entire in $a$ and $b$, is holomorphic in $c$ on $\mathbb{C} \setminus \mathbb{Z}_{\leqslant 0}$, with a simple pole at every negative integer $k \in \mathbb{Z}_{\leqslant 0}$.
        
            \end{corollary}
            
            \begin{proof}
            
                The regularity with regard to the parameters $a$, $b$, and $c$ of solutions to \eqref{EqHypergeometric} is studied in \cite[Sec. 5.3.3]{MR1429619}.
            
            \end{proof}

        \subsubsection{Examples of hypergeometric functions}
        
            There are two particular instances of hypergeometric functions which are of importance for this paper.
        
            \begin{proposition}
            \label{PropBinomialFormula}
        
                For any complex numbers $a$ and $b$, with $a \not \in \mathbb{Z}_{\leqslant 0}$, we have
            
                \[ \begin{array}{lllll}
                        \sum\limits_{k=0}^{+\infty} \; \left( b \right)_k \cdot \frac{z^k}{k!} & = & F \left( a, b; a; z \right) & = & \left( 1 - z \right)^{-b}
                    \end{array} \]
            
                \noindent on the open unit disk. The value of $a$ actually plays no role.
        
            \end{proposition}
            
            \begin{proof}
            
                Proving this proposition amounts to noting that both sides verify the differential equation $ \left( 1-z \right) u' = b u$ with initial condition $u \left(0 \right) = 1$.
            
            \end{proof}
        
            \begin{proposition}
            \label{PropExpleLog}
        
                For any complex number $z$ with $\left \vert z \right \vert < 1$, we have
            
                \[ \begin{array}{lll}
                        z F \left( 1, 1; 2; z \right) & = & - \log \left( 1 - z \right)
                    \end{array} . \]
        
            \end{proposition}
            
            \begin{proof}
            
                Integrating the power series expansion of $\left( 1-z \right)^{-1}$ gives this result.
            
            \end{proof}

        \subsubsection{Transformation of the variable}
        
            This paragraph is devoted to the presentation of some transformations on the last variables of hypergeometric functions.
            
            \begin{proposition}
            \label{PropTransformationFraction}
            
                Let $a, b, c$ be complex numbers. Assume we have $c \not \in \mathbb{Z}_{\leqslant 0}$. For any $z \in \mathbb{C} \setminus \left[ 1, + \infty \right[$, we have
                
                \[ \begin{array}{lllll}
                        F \left( a, b; c; z \right) & \hspace{-1.5pt} = \hspace{-1.5pt} & \left( 1 - z \right)^{-a} F \left( a, c - b; c; \frac{z}{z - 1} \right) & \hspace{-1.5pt} = \hspace{-1.5pt} & \left( 1 - z \right)^{-b} F \left( b, c - a; c; \frac{z}{z - 1} \right)
                    \end{array} . \]
            
            \end{proposition}
            
            \begin{proof}
            
                This can be found in \cite[Sec. 5.10.3]{MR1429619}.
            
            \end{proof}
            
            \begin{proposition}
            \label{PropTransformationLinear}
            
                Let $a, b, c$ be complex numbers. Assume we have $c \not \in \mathbb{Z}_{\leqslant 0}$, and that $a + b - c$ is not an integer. For any $z \in \mathbb{C} \setminus \left[ 1, + \infty \right[$, we have
                
                \[ \begin{array}{lll}
                        F \left( a, b; c; z \right) & = & \frac{\Gamma \left( c \right) \Gamma \left( a + b - c \right)}{\Gamma \left( a \right) \Gamma \left( b \right)} \left( 1 - z \right)^{c - a - b} F \left( c - a, c - b; 1 + c - a - b; 1 - z \right) \\[1em]
                        
                        && \qquad + \frac{\Gamma \left( c \right) \Gamma \left( c - a - b \right)}{\Gamma \left( c - a \right) \Gamma \left( c - b \right)} F \left( a, b; 1 + a + b - c; 1 - z \right).
                    \end{array} \]
            
            \end{proposition}
            
            \begin{proof}
            
                This can be found in \cite[Sec. 5.10.4]{MR1429619}.
            
            \end{proof}

        \subsubsection{Integral representations}
        
            We will now see how to compute certain integrals using hypergeometric functions. The first step is to prove the \textit{Euler integral formula}.
        
            \begin{proposition}[Euler integral formula]
            \label{PropEulerIntegralFormula}
        
                Assume we have $\Re c > \Re b > 0$. The hypergeometric function $F \left( a, b; c; z \right)$ then has the following integral representation
            
                \[ \begin{array}{lll}
                        F \left( a, b; c; z \right) & = & \frac{\Gamma \left( c \right)}{\Gamma \left( b \right) \Gamma \left( c - b \right)} \int_0^1 \; t^{b-1} \left( 1 - t \right)^{c - b - 1} \left( 1 - zt \right)^{-a} \; \mathrm{d}t
                    \end{array} . \]
        
            \end{proposition}
            
            \begin{proof}
            
                This integral representation is based on one that holds for the beta function
                
                \[ \begin{array}{lllll}
                        \frac{\Gamma \left( s \right) \Gamma \left( w \right)}{\Gamma \left( s+w \right)} & = & B \left( s, w \right) & = & \int_0^1 t^{s-1} \left( 1-t \right)^{w-1} \mathrm{d}t,
                    \end{array} \]
                
                \noindent which is proved in \cite[Sec. 2.1.6]{MR1429619}, and holds for any complex numbers $s$ and $w$ with strictly positive real parts. The rest of the proof is found in \cite[Sec. 5.9.4]{MR1429619}.
            
            \end{proof}
        
            \begin{corollary}
            \label{CorIntegralFormula}
        
                For any complex numbers $\mu$ and $\nu$, with $\Re \mu > 0$ and $\nu - \mu \not \in \mathbb{Z}_{\leqslant 0}$, and any real number $u > 0$, we have
            
                \[ \begin{array}{lll}
                        \int_0^u \frac{y^{\mu-1}}{\left( 1 + y \right)^{\nu}} \mathrm{d}y & = & \frac{u^{\mu - \nu}}{\mu - \nu} F \left( \nu, \nu - \mu; \nu - \mu + 1; - \frac{1}{u} \right) + \frac{\Gamma \left( \mu \right) \Gamma \left( \nu - \mu \right)}{\Gamma \left( \nu \right)}
                    \end{array} \]
        
            \end{corollary}
        
            \begin{proof}
        
                Let $\mu, \nu$ be complex numbers with $\Re \mu > 0$, and $u > 0$ be a strictly positive real number. First, using the definition of the \textit{beta function}, we have
            
                \begin{equation}
                \label{EqBetaFunction}
                    \begin{array}{lllll}
                        \frac{\Gamma \left( \mu \right) \Gamma \left( \nu - \mu \right)}{\Gamma \left( \nu \right)} & = & \int_0^1 t^{\mu - 1} \left( 1 - t \right)^{\nu - \mu - 1} \mathrm{d}t & = & \int_0^{+ \infty} \frac{y^{\mu - 1}}{\left( 1 + y \right)^{\nu}} \mathrm{d}t
                    \end{array}
                \end{equation}
            
                \noindent using the change of variable $t = y/\left( y+1 \right)$. Using proposition \ref{PropEulerIntegralFormula}, we have
            
                \[ \begin{array}{lll}
                        F \left( \nu, \nu - \mu; \nu - \mu + 1; z \right) & = & \frac{\Gamma \left( \nu - \mu + 1 \right)}{\Gamma \left( \nu - \mu \right)} \int_0^1 \; t^{\nu - \mu -1} \left( 1 + \frac{t}{u} \right)^{- \nu} \; \mathrm{d}t \\[1em]
                    
                        & = & \left( \nu - \mu \right) u^{\nu - \mu} \int_0^{1/u} x^{\nu - \mu - 1} \left( 1 + x \right)^{-\nu} \mathrm{d}x
                    \end{array} \]
            
                \noindent the last equality being obtained by setting $x=t/u$. This yields
            
                \begin{equation}
                \label{EqHyperGeoIntegral}
                    \begin{array}{lll}
                        \frac{u^{\mu - \nu}}{\mu - \nu} F \left( \nu, \nu - \mu; \nu - \mu + 1; - \frac{1}{u} \right) & = & - \int_0^{1/u} x^{\nu - \mu - 1} \left( 1 + x \right)^{-\ nu} \mathrm{d}x \\[1em]
                    
                        & = & - \int_u^{+ \infty} \frac{y^{\mu - 1}}{\left( 1 + y \right)^{\nu}} \mathrm{d}y
                    \end{array}
                \end{equation}
                
                \noindent after having performed the change of variable $y = 1/x$. The result then stems from equations \eqref{EqBetaFunction} and \eqref{EqHyperGeoIntegral}. The proposition can be obtained in the more general case, \textit{i.e.} without assuming that we have $\Re \nu > \Re \mu$, by analytic continuation.
        
            \end{proof}
            
            \begin{corollary}
            \label{CorIntegralFormulaHypergeometric}
            
                For any complex numbers $\mu$ and $\nu$, with $\Re \mu > 0$ and $\nu - \mu \not \in \mathbb{Z}_{\leqslant 0}$, and any real number $u > 0$, we have
                
                \[ \begin{array}{lll}
                        \int_0^u \frac{y^{\mu-1}}{\left( 1 + y \right)^{\nu}} \mathrm{d}y & = & \frac{1}{\mu} u^{\mu} \left( 1 + u \right)^{-\nu} F \left( \nu, 1; \mu + 1; \frac{u}{1 + u} \right)
                    \end{array} . \]
            
            \end{corollary}
            
            \begin{proof}
            
                To prove this result, one begins by applying corollary \ref{CorIntegralFormula}, then uses proposition \ref{PropTransformationFraction} on the hypergeometric function, and finally proposition \ref{PropTransformationLinear}.
            
            \end{proof}

        \subsubsection{Hypergeometric functions at $z = 1$}
        
            As we have already seen, the hypergeometric function $F \left( a, b; c; z \right)$ does not, in general, make sense at $z = 1$. We will now see when it actually does.
        
            \begin{proposition}
            \label{PropHypergeometricZ1}
        
                Assume $a,b,c$ are complex numbers with $c \not \in \mathbb{Z}_{\leqslant 0}$ and that we have $\Re \left( c - a - b \right) > 0$. The hypergeometric series converges at $z = 1$, and we have
            
                \[ \begin{array}{lllll}
                        \sum\limits_{k=0}^{+ \infty} \frac{\left( a \right)_k \left( b \right)_k}{\left( c \right)_k} \cdot \frac{1}{k!} & = & \lim\limits_{z \rightarrow 1} F \left( a, b; c; z \right) & = & \frac{\Gamma \left( c \right) \Gamma \left( c - a - b \right)}{\Gamma \left( c - a \right) \Gamma \left( c - b \right)}
                    \end{array} . \]
            
                \noindent In particular, the function $x \longmapsto F \left( a, b; c; x \right)$ is continuous on $\left[ 0,1 \right[$, and can be bounded, locally uniformly in the parameters $a,b,c$.
        
            \end{proposition}
            
            \begin{proof}
            
                This is the content of \cite[Sec. 5.9.5]{MR1429619}.
            
            \end{proof}

        \subsubsection{Contiguous functions}
        
            It is sometimes required in this paper to apply some linear transformation on the first three parameters of hypergeometric functions. There are multiple such relations, though we will only present one.
        
            \begin{proposition}
            \label{PropContiguousFunctions}
        
                Let $a,b,c$ be complex numbers, with $c \in \mathbb{Z}_{\leqslant 1}$. For any complex number $z \in \mathbb{C} \setminus \left[ 1, + \infty \right[$, we have
                
                \[ \begin{array}{lll}
                        \left( a - 1 + \left( b + 1 - c \right) z \right) F \left( a, b; c; z \right) + \left( c - a \right) F \left( a - 1, b; c; z \right) \\[0.5em]
                        
                        \qquad \qquad \qquad \qquad \qquad \qquad \qquad \qquad - \left( c - 1 \right) \left( 1 - z \right) F \left( a, b; c-1; z \right) & = & 0.
                    \end{array} \]
        
            \end{proposition}
            
            \begin{proof}
            
                This equality can be verified directly using the hypergeometric series and the properties of the Gamma function (seen as properties on the Pochhammer symbols) on the region $\left \vert z \right \vert < 1$. The unicity of analytic continuation then completes the proof.
            
            \end{proof}

        \subsubsection{Extraction of terms and generalized hypergeometric functions}
        
            Going back to the definition of the hypergeometric series, it is sometimes necessary in this paper to set aside the first few terms of the series, and still recognize some hypergeometric function in the remainder.
            
            \begin{proposition}
            \label{PropExtractionTerms}
                    
                Let $a, c$ be complex numbers, with $c \not \in \mathbb{Z}_{\leqslant 0}$. For any complex number $z \in \mathbb{C} \setminus \left[ 1, + \infty \right[$, we have
                
                \[ \begin{array}{lll}
                        F \left( a, 1; c; z \right) & = & 1 + \frac{a}{c} z + \frac{a \left( a + 1 \right)}{c \left( c + 1 \right)} z^2 F \left( a + 2, 1; c + 2; z \right)
                    \end{array} . \]
                    
            \end{proposition}
            
            \begin{proof}
            
                It is enough to prove the proposition for $z$ in the open unit disk, by unicity of meromorphic extensions. We have
                
                \[ \begin{array}{lllll}
                        F \left( a, 1; c; z \right) & = & 1 + \frac{a}{c} z + \sum\limits_{n=2}^{+ \infty} \frac{\left( a \right)_n \left( 1 \right)_n}{\left( c \right)_n} \cdot \frac{z^n}{n!} & = & 1 + \frac{a}{c} z + \sum\limits_{n=2}^{+ \infty} \frac{\left( a \right)_n}{\left( c \right)_n} z^n \\[1em]
                        
                        &&& = & 1 + \frac{a}{c} z + z^2 \sum\limits_{n=0}^{+ \infty} \frac{\left( a \right)_{n+2}}{\left( c \right)_{n+2}} z^n.
                    \end{array} \]
                
                \noindent We can then conclude, using the equality $\left( a \right)_{n+2} = a \left( a+1 \right) \left( a+2 \right)_n$, and its version for $c$, which can be proved using the expression of the Pochhammer symbol as a quotient of Gamma functions.
            
            \end{proof}
            
            There is a similar result which holds without assuming that one of the first two parameters equals $1$. This will require the introduction of \textit{generalized hypergeometric functions}.
            
            \begin{proposition-definition}
            \label{DefGeneralizedHypergeometric}
                    
                Let $a_1, a_2, a_3, b_1, b_2$ be complex numbers. Assume that we have $b_1, b_2 \not \in \mathbb{Z}_{\leqslant 0}$. The generalized hypergeometric series
                
                \[ \begin{array}{lll}
                        F \left( a_1, a_2, a_3; b_1, b_2; z \right) & = & \sum\limits_{n = 0}^{+ \infty} \frac{\left( a_1 \right)_n \left( a_2 \right)_n \left( a_3 \right)_n}{\left( b_1 \right)_n \left( b_2 \right)_n} \cdot \frac{z^n}{n!}
                    \end{array} \]
                
                \noindent is absolutely convergent on the open unit disk.
                    
            \end{proposition-definition}
            
            \begin{proof}
            
                This is a consequence of d'Alembert's ratio test for series.
            
            \end{proof}
            
            \begin{proposition}
            \label{PropExtractionTermsGeneralizedHypergeometric}
                    
                Let $a, b, c$ be complex numbers, with $c \not \in \mathbb{Z}_{\leqslant 0}$. For any complex number $z$ in the open unit disk, we have
                
                \[ \begin{array}{lll}
                        F \left( a, b; c; z \right) & = & 1 + \frac{ab}{c} z F \left( a+1, b+1, 1; c+1, 2; z \right) \\[1em]
                        
                        & = & 1 + \frac{a b}{c} z + \frac{a \left( a + 1 \right) b \left( b + 1 \right)}{c \left( c + 1 \right)} \cdot \frac{z^2}{2} F \left( a + 2, b + 2, 1; c + 2, 3; z \right).
                    \end{array} \]
                
                \noindent In particular, assuming we have $\Re \left( c - a - b \right) > 0$, the generalized hypergeometric function $x \longmapsto F \left( a + 1, b + 1, 1; c + 1, 2; x \right)$ is bounded on $\left[ 0, 1 \right[$, locally uniformly in the parameters $a, b, c$. A similar result holds for the generalized hypergeometric function $x \longmapsto F \left( a + 2, b + 2, 1; c + 2, 3; x \right)$.
                    
            \end{proposition}
            
            \begin{proof}
            
                This is a direct computation, similar to the one performed in the proof of proposition \ref{PropExtractionTerms}, using the series defining the hypergeometric function and its generalized versions.
            
            \end{proof}

\bibliographystyle{amsplain}

\bibliography{bibliography}

\providecommand{\bysame}{\leavevmode\hbox to3em{\hrulefill}\thinspace}
\providecommand{\MR}{\relax\ifhmode\unskip\space\fi MR }
\providecommand{\MRhref}[2]{%
  \href{http://www.ams.org/mathscinet-getitem?mr=#1}{#2}
}
\providecommand{\href}[2]{#2}
\begin{thebibliography}{10}

\bibitem{MR1165865}
D.~Burghelea, L.~Friedlander, and T.~Kappeler, \emph{Meyer-{V}ietoris type
  formula for determinants of elliptic differential operators}, J. Funct. Anal.
  \textbf{107} (1992), no.~1, 34--65. \MR{1165865}

\bibitem{MR3677185}
B.~Candelpergher, \emph{Ramanujan summation of divergent series}, Lecture Notes
  in Mathematics, vol. 2185, Springer, Cham, 2017. \MR{3677185}

\bibitem{MR1102675}
R.~Carmona and J.~Lacroix, \emph{Spectral theory of random {S}chr\"{o}dinger
  operators}, Probability and its Applications, Birkh\"{a}user Boston, Inc.,
  Boston, MA, 1990. \MR{1102675}

\bibitem{MR1890995}
G.~Carron, \emph{D\'{e}terminant relatif et la fonction {X}i}, Amer. J. Math.
  \textbf{124} (2002), no.~2, 307--352. \MR{1890995}

\bibitem{MR688031}
Y.~Colin~de Verdi\`ere, \emph{Pseudo-laplaciens. {I}}, Ann. Inst. Fourier
  (Grenoble) \textbf{32} (1982), no.~3, xiii, 275--286. \MR{688031}

\bibitem{MR699488}
\bysame, \emph{Pseudo-laplaciens. {II}}, Ann. Inst. Fourier (Grenoble)
  \textbf{33} (1983), no.~2, 87--113. \MR{699488}

\bibitem{MR1070713}
J.~B. Conway, \emph{A course in functional analysis}, second ed., Graduate
  Texts in Mathematics, vol.~96, Springer-Verlag, New York, 1990. \MR{1070713}

\bibitem{MR1349825}
E.~B. Davies, \emph{Spectral theory and differential operators}, Cambridge
  Studies in Advanced Mathematics, vol.~42, Cambridge University Press,
  Cambridge, 1995. \MR{1349825}

\bibitem{MR0417174}
P.~Deligne, \emph{\'{E}quations diff\'{e}rentielles \`a points singuliers
  r\'{e}guliers}, Lecture Notes in Mathematics, Vol. 163, Springer-Verlag,
  Berlin-New York, 1970. \MR{0417174}

\bibitem{MR2343536}
J.~Eichhorn, \emph{Global analysis on open manifolds}, Nova Science Publishers,
  Inc., New York, 2007. \MR{2343536}

\bibitem{MR4167014}
G.~Freixas~i Montplet and A.-M. von Pippich, \emph{Riemann--{R}och isometries
  in the non-compact orbifold setting}, J. Eur. Math. Soc. (JEMS) \textbf{22}
  (2020), no.~11, 3491--3564. \MR{4167014}

\bibitem{MR1335452}
T.~Kato, \emph{Perturbation theory for linear operators}, Classics in
  Mathematics, Springer-Verlag, Berlin, 1995, Reprint of the 1980 edition.
  \MR{1335452}

\bibitem{MR2339952}
X.~Ma and G.~Marinescu, \emph{Holomorphic {M}orse inequalities and {B}ergman
  kernels}, Progress in Mathematics, vol. 254, Birkh\"{a}user Verlag, Basel,
  2007. \MR{2339952}

\bibitem{MR1617554}
W.~M\"uller, \emph{Relative zeta functions, relative determinants and
  scattering theory}, Comm. Math. Phys. \textbf{192} (1998), no.~2, 309--347.
  \MR{1617554}

\bibitem{MR92871}
F.~Oberhettinger, \emph{On the derivative of {B}essel functions with respect to
  the order}, J. Math. and Phys. \textbf{37} (1958), 75--78. \MR{92871}

\bibitem{MR1429619}
F.~W.~J. Olver, \emph{Asymptotics and special functions}, AKP Classics, A K
  Peters, Ltd., Wellesley, MA, 1997, Reprint of the 1974 original [Academic
  Press, New York; MR0435697 (55 \#8655)]. \MR{1429619}

\bibitem{MR2723248}
F.~W.~J. Olver, D.~W. Lozier, R.~F. Boisvert, and C.~W. Clark (eds.),
  \emph{N{IST} handbook of mathematical functions}, U.S. Department of
  Commerce, National Institute of Standards and Technology, Washington, DC;
  Cambridge University Press, Cambridge, 2010. \MR{2723248}

\bibitem{MR2393625}
C.~A.~M. Peters and J.~H.~M. Steenbrink, \emph{Mixed {H}odge structures},
  Ergebnisse der Mathematik und ihrer Grenzgebiete. 3. Folge. A Series of
  Modern Surveys in Mathematics [Results in Mathematics and Related Areas. 3rd
  Series. A Series of Modern Surveys in Mathematics], vol.~52, Springer-Verlag,
  Berlin, 2008. \MR{2393625}

\bibitem{MR1157815}
W.~Rudin, \emph{Functional analysis}, second ed., International Series in Pure
  and Applied Mathematics, McGraw-Hill, Inc., New York, 1991. \MR{1157815}

\bibitem{MR2439244}
A.~A. Saharian, \emph{A summation formula over the zeros of the associated
  {L}egendre function with a physical application}, J. Phys. A \textbf{41}
  (2008), no.~41, 415203, 17. \MR{2439244}

\bibitem{MR1976398}
E.~M. Stein and R.~Shakarchi, \emph{Complex analysis}, Princeton Lectures in
  Analysis, vol.~2, Princeton University Press, Princeton, NJ, 2003.
  \MR{1976398}

\bibitem{MR2451566}
C.~Voisin, \emph{Hodge theory and complex algebraic geometry. {I}}, english
  ed., Cambridge Studies in Advanced Mathematics, vol.~76, Cambridge University
  Press, Cambridge, 2007, Translated from the French by Leila Schneps.
  \MR{2451566}

\end{thebibliography}

\end{document}